\documentclass[11pt]{article}

\usepackage{amssymb,amsbsy,amsthm,amsmath,graphicx,mathrsfs,epsfig,graphics,times}

\input xy

\xyoption{all}

\begin{document}

\title{\textbf{Pleasant extensions retaining algebraic structure, II}}
\author{Tim Austin}
\date{}

\maketitle


\newenvironment{nmath}{\begin{center}\begin{math}}{\end{math}\end{center}}

\newtheorem{thm}{Theorem}[section]
\newtheorem*{thm*}{Theorem}
\newtheorem{lem}[thm]{Lemma}
\newtheorem{prop}[thm]{Proposition}
\newtheorem{cor}[thm]{Corollary}
\newtheorem{conj}[thm]{Conjecture}
\newtheorem{dfn}[thm]{Definition}
\newtheorem{prob}[thm]{Problem}
\newtheorem{ques}[thm]{Question}
\theoremstyle{remark}
\newtheorem*{ex}{Example}


\newcommand{\s}{\sigma}
\renewcommand{\O}{\Omega}
\renewcommand{\S}{\Sigma}
\newcommand{\co}{\mathrm{co}}
\newcommand{\e}{\mathrm{e}}
\newcommand{\eps}{\varepsilon}
\renewcommand{\d}{\mathrm{d}}
\newcommand{\im}{\mathrm{i}}
\renewcommand{\k}{\kappa}
\renewcommand{\l}{\lambda}
\newcommand{\G}{\Gamma}
\newcommand{\g}{\gamma}
\renewcommand{\L}{\Lambda}
\renewcommand{\a}{\alpha}
\renewcommand{\b}{\beta}
\newcommand{\Sone}{\mathrm{S}^1}

\newcommand{\Aut}{\mathrm{Aut}}
\renewcommand{\Pr}{\mathrm{Pr}}
\newcommand{\Hom}{\mathrm{Hom}}
\newcommand{\id}{\mathrm{id}}
\newcommand{\BSp}{\mathsf{BSp}}
\newcommand{\WBSp}{\mathsf{WBSp}}
\newcommand{\Sys}{\underline{\mathsf{Sys}}}
\newcommand{\Rep}{\underline{\mathsf{Rep}}}
\newcommand{\Lat}{\mathrm{Lat}}
\newcommand{\FLat}{\mathrm{FLat}}
\newcommand{\SG}{\mathrm{SG}}
\newcommand{\CSG}{\mathrm{CSG}}
\newcommand{\Clos}{\mathrm{Clos}}
\newcommand{\pro}{\mathrm{pro}}
\newcommand{\nil}{\mathrm{nil}}
\newcommand{\rat}{\mathrm{rat}}
\newcommand{\Ab}{\mathrm{Ab}}
\newcommand{\CIS}{\mathrm{CIS}}
\newcommand{\DDelta}{{\Delta\!\!\!\!\bullet\ }}

\newcommand{\bbN}{\mathbb{N}}
\newcommand{\bbR}{\mathbb{R}}
\newcommand{\bbZ}{\mathbb{Z}}
\newcommand{\bbQ}{\mathbb{Q}}
\newcommand{\bbT}{\mathbb{T}}
\newcommand{\bbD}{\mathbb{D}}
\newcommand{\bbC}{\mathbb{C}}
\newcommand{\bbP}{\mathbb{P}}

\newcommand{\A}{\mathcal{A}}
\newcommand{\B}{\mathcal{B}}
\newcommand{\C}{\mathcal{C}}
\newcommand{\E}{\mathcal{E}}
\newcommand{\F}{\mathcal{F}}
\newcommand{\I}{\mathcal{I}}
\newcommand{\calL}{\mathcal{L}}
\renewcommand{\P}{\mathcal{P}}
\newcommand{\U}{\mathcal{U}}
\newcommand{\W}{\mathcal{W}}
\newcommand{\Y}{\mathcal{Y}}
\newcommand{\Z}{\mathcal{Z}}

\newcommand{\frH}{\mathfrak{H}}
\newcommand{\frM}{\mathfrak{M}}

\newcommand{\bfX}{\mathbf{X}}
\newcommand{\bfY}{\mathbf{Y}}
\newcommand{\bfZ}{\mathbf{Z}}
\newcommand{\bfW}{\mathbf{W}}
\newcommand{\bfV}{\mathbf{V}}

\newcommand{\rmC}{\mathrm{C}}
\newcommand{\rmE}{\mathrm{E}}
\newcommand{\rmH}{\mathrm{H}}
\newcommand{\T}{\mathrm{T}}

\newcommand{\sfC}{\mathsf{C}}
\newcommand{\sfD}{\mathsf{D}}
\newcommand{\sfE}{\mathsf{E}}
\newcommand{\sfV}{\mathsf{V}}
\newcommand{\sfW}{\mathsf{W}}
\newcommand{\sfY}{\mathsf{Y}}
\newcommand{\sfZ}{\mathsf{Z}}

\newcommand{\uhr}{\!\!\upharpoonright}
\newcommand{\into}{\hookrightarrow}
\newcommand{\onto}{\twoheadrightarrow}
\newcommand{\actson}{\curvearrowright}

\newcommand{\bb}[1]{\mathbb{#1}}
\newcommand{\bs}[1]{\boldsymbol{#1}}
\newcommand{\fr}[1]{\mathfrak{#1}}
\renewcommand{\bf}[1]{\mathbf{#1}}
\renewcommand{\sf}[1]{\mathsf{#1}}
\renewcommand{\rm}[1]{\mathrm{#1}}
\renewcommand{\cal}[1]{\mathcal{#1}}

\renewcommand{\t}[1]{\tilde{#1}}
\newcommand{\ol}[1]{\overline{#1}}
\renewcommand{\hat}[1]{\widehat{#1}}

\newcommand{\mvee}{\hbox{$\bigvee$}}

\newcommand{\lfr}{\lfloor}
\newcommand{\rfr}{\rfloor}

\newcommand{\fin}{\nolinebreak\hspace{\stretch{1}}$\lhd$}
\newcommand{\tick}[1]{\nolinebreak\hspace{\stretch{1}}$\surd_{\mathrm{#1}}$}

\begin{abstract}
In this paper we combine the general tools developed in~\cite{Aus--lindeppleasant1} with several ideas taken from
earlier work on one-dimensional nonconventional ergodic averages by
Furstenberg and Weiss~\cite{FurWei96}, Host and Kra~\cite{HosKra05}
and Ziegler~\cite{Zie07} to study the averages
\[\frac{1}{N}\sum_{n=1}^N(f_1\circ T^{n\bf{p}_1})(f_2\circ T^{n\bf{p}_2})(f_3\circ T^{n\bf{p}_3})\quad\quad f_1,f_2,f_3 \in L^\infty(\mu)\]
associated to a triple of directions $\bf{p}_1,\bf{p}_2,\bf{p}_3 \in
\bbZ^2$ that lie in general position along with $\bs{0}\in \bbZ^2$.
We will show how to construct a `pleasant' extension of an
initially-given $\bbZ^2$-system for which these averages admit
characteristic factors with a very concrete description, involving
the same structure as for those in~\cite{Aus--nonconv} together with two-step pro-nilsystems (reminiscent of~\cite{HosKra05} and its predecessors).

We will also use this
analysis to construct pleasant extensions and then prove norm
convergence for the polynomial nonconventional ergodic averages
\[\frac{1}{N}\sum_{n=1}^N(f_1\circ T_1^{n^2})(f_2\circ T_1^{n^2}T_2^n)\]
associated to two commuting transformations $T_1$, $T_2$.
\end{abstract}

\parskip 0pt

\tableofcontents

\parskip 7pt

\section{Introduction}

This paper continues the work of~\cite{Aus--lindeppleasant1}, and we
will freely refer to that paper for a detailed background discussion
and several necessary results.

We consider probability-preserving actions $T:\bbZ^2\curvearrowright
(X,\mu)$ on standard Borel spaces, and study the associated `nonconventional'
ergodic averages of the form
\[\frac{1}{N}\sum_{n=1}^N (f_1\circ T^{n\bf{p}_1})(f_2\circ T^{n\bf{p}_2})(f_3\circ T^{n\bf{p}_3})\quad\quad \hbox{for}\ f_1,f_2,f_3 \in L^\infty(\mu)\]
where $\bf{p}_1$, $\bf{p}_2$, $\bf{p}_3 \in \bbZ^2\setminus
\{\bs{0}\}$ are distinct and are such that together with $\bs{0}$
they lie in general position: that is, such that no three of the
points $\bs{0}$, $\bf{p}_1$, $\bf{p}_2$, $\bf{p}_3$ lie on a line.
Following the general terminology recalled
in~\cite{Aus--lindeppleasant1}, a triple of factors
$\xi_i:(X,\mu,T)\to (Y_i,\nu_i,S_i)$ is \textbf{characteristic} for
these averages if
\begin{multline*}
\frac{1}{N}\sum_{n=1}^N(f_1\circ T^{n\bf{p}_1})(f_2\circ
T^{n\bf{p}_2})(f_3\circ T^{n\bf{p}_3})\\ \sim
\frac{1}{N}\sum_{n=1}^N(\sfE_\mu(f_1\,|\,\xi_1)\circ
T^{n\bf{p}_1})(\sfE_\mu(f_2\,|\,\xi_2)\circ
T^{n\bf{p}_2})(\sfE_\mu(f_3\,|\,\xi_3)\circ T^{n\bf{p}_3}),
\end{multline*}
for any $f_1$, $f_2$, $f_3 \in L^\infty(\mu)$, where we write $f_N
\sim g_N$ to denote that $\|f_N - g_N\|_2 \to 0$ as $N\to\infty$.

Motivated by the approach to such averages developed
in~\cite{Aus--nonconv,Aus--newmultiSzem}
(building on several earlier contributions, discussed properly
in~\cite{Aus--lindeppleasant1}), we here seek an extension
$\pi:(\t{X},\t{\mu},\t{T})\to (X,\mu,T)$ of an arbitrary
initially-given system in which a characteristic triple of factors
can be found whose form is as simple as possible.  The new feature of
the present paper is that we insist on retaining the linear dependence among the directions $\bf{p}_i$
in the extended system, which creates difficulties that did not
arise in those earlier works owing to an implicit assumption of linear independence.

More precisely, the `best' extensions among $\bbZ^2$-systems for the study of these averages generally require
characteristic factors that are not as simple as the pure joins of
isotropy factors that emerge in the linearly independent case (see
Theorem 1.1 in~\cite{Aus--lindeppleasant1}).  The extra ingredients we need to construct these characteristic factors are two-step pro-nilsystems, which re-appear here after having taken centre stage in the study of the nonconventional averages of actions of $\bbZ$ through the works of Conze and Lesigne~\cite{ConLes84}, Zhang~\cite{Zha96}, Furstenberg and Weiss~\cite{FurWei96}, Host and Kra~\cite{HosKra05} and Ziegler~\cite{Zie07}. 

\begin{thm}[Pleasant extensions for linearly dependent triple linear averages]\label{thm:char-three-lines-in-2D}
Any system $T:\bbZ^2\curvearrowright (X,\mu)$ has an extension
$\pi:(\t{X},\t{\mu},\t{T}) \to (X,\mu,T)$ such that for any
$\bf{p}_1$, $\bf{p}_2$, $\bf{p}_3 \in \bbZ^2$ that are in general
position with the origin the averages
\[\frac{1}{N}\sum_{n=1}^N(f_1\circ \t{T}^{n\bf{p}_1})(f_2\circ \t{T}^{n\bf{p}_2})(f_3\circ \t{T}^{n\bf{p}_3}),\quad\quad f_1,f_2,f_3 \in L^\infty(\t{\mu}),\]
admit the characteristic factors
\[\xi_i = \zeta_0^{\t{T}^{\bf{p}_i}}\vee \zeta_0^{\t{T}^{\bf{p}_i}=\t{T}^{\bf{p}_j}}\vee\zeta_0^{\t{T}^{\bf{p}_i}=\t{T}^{\bf{p}_k}}\vee \zeta_{\nil,2}^{\t{T}}\quad\quad i=1,2,3\]
where the target of $\zeta_{\nil,2}^{\t{T}}$ is an inverse limit of direct integrals of two-step $\bbZ^2$-pro-nilsystems.
\end{thm}

Note that this result promises a single extension that
simultaneously enjoys simplified characteristic factors for every
triple of directions in general position with the origin. Motivated
by~\cite{Aus--nonconv,Aus--lindeppleasant1}, we will refer to such
an extension as a \textbf{pleasant extension for linearly dependent
triple linear nonconventional averages}.

In addition to its technical interest,
Theorem~\ref{thm:char-three-lines-in-2D} can be applied to prove a new case of $L^2$-convergence for
Bergelson and Leibman's polynomial nonconventional ergodic
averages~(\cite{BerLei02}):

\begin{thm}\label{thm:polyconv}
If $T_1,T_2:\bbZ\curvearrowright (X,\mu)$ commute then the averages
\[\frac{1}{N}\sum_{n=1}^N (f_1\circ T_1^{n^2})(f_2\circ T_1^{n^2}T_2^n)\]
converge in $L^2(\mu)$ as $N\to\infty$ for any $f_1,f_2 \in
L^\infty(\mu)$.
\end{thm}

We prove this in Section~\ref{sec:convproof}, making use of an extension of $(X,\mu,T_1,T_2)$ in which the
above quadratic averages admit quite concrete characteristic
factors, related to those we obtain in
Theorem~\ref{thm:char-three-lines-in-2D}.

Although the new convergence result of Theorem~\ref{thm:polyconv} is
modest in itself, the methods we develop in pursuit of
Theorem~\ref{thm:char-three-lines-in-2D} seem to indicate a much more far-reaching structure that may emerge in connexion with Bergelson and Leibman's conjecture of polynomial
nonconventional average convergence, and potentially in other
questions on the structure of joinings between different classes of
system in the ergodic theory of $\bbZ^d$-actions.

\textbf{Acknowledgements}\quad My thanks go to Vitaly Bergelson,
Bernard Host, Bryna Kra, Mariusz Lema\'nczyk, Emmanuel Lesigne, Terence Tao, Dave
Witte Morris and Tamar Ziegler for several helpful discussions and
to the Mathematical Sciences Research Institute (Berkeley) for its
hospitality during the 2008 program on Ergodic Theory and Additive
Combinatorics.

\section{Some preliminary results on isometric extensions}

In this paper we will make free use of the background results
recalled in Section 2 of~\cite{Aus--lindeppleasant1}, and of the
formalism of idempotent classes of system and satedness developed in
Section 3 of~\cite{Aus--lindeppleasant1}. However, in addition to
those we will now need to make quite extensive use of the theory of
isometric extensions of not-necessarily-ergodic
probability-preserving systems, as developed
in~\cite{Aus--ergdirint} building on classical works of Mackey,
Furstenberg and Zimmer (see that paper for more complete
references).  We recall some of the necessary statements here, and
also introduce the new property of `fibre-normality' (adapted from a
definition of Furstenberg and Weiss~\cite{FurWei96}) that will be
useful later.

Before all else, let us remind the reader that we work throughout in
the category of probability-preserving actions on standard Borel
spaces, and consequently that whenever an isometric extension of
such systems is coordinatized using a measurably-varying family of
compact homogeneous spaces, it will be implicit that these
homogeneous spaces are constructed from some measurable-varying
compact \emph{metrizable} groups, themselves drawn from within some
metrizable compact fibre repository group.  This may always be
assumed, even though for brevity we sometimes omit to mention
metrizability explicitly.  The definition of these extensions (along
with this convention concerning metrizability) can be found in
Section 3 of~\cite{Aus--ergdirint}.

\subsection{Mackey Theory and the Furstenberg-Zimmer Structure Theorem over a non-ergodic
base}\label{subs:MandFZ}

The classical Mackey Theory describing the ergodic decomposition of
a skew-product extension of an ergodic system by rotations on a
compact homogeneous space is extended to the case of a non-ergodic
base by allowing families of compact homogeneous space fibres over
the base that are invariant for the action but otherwise can vary
measurably (in a suitable sense made formal in Section 3
of~\cite{Aus--ergdirint}).  These results apply to jointly
measurable, probability-preserving actions of an arbitrary locally
compact second countable group $\G$. Referring to such families, the
main results of the extended Mackey Theory are the following:

\begin{thm}\label{thm:homo-nonergMackey-1}
Suppose that $(X,\mu,T) = (Y,\nu,S)\ltimes
(G_\bullet/H_\bullet,m_{G_\bullet/H_\bullet},\rho)$ is a
$\G$-system, $\zeta_0^S:Y\to Z_0^S$ a coordinatization of the base
isotropy factor and $P:Z_0^S \stackrel{\rm{p}}{\longrightarrow} Y$ a
version of the disintegration of $\nu$ over $\zeta_0^S$. Then there
are subgroup data $K_\bullet \leq G_\bullet$ and a cocycle-section
$b:Y\to G_\bullet$ such that the factor map
\[\phi:X \to Z_0^S\ltimes (K_\bullet\backslash G_\bullet/H_\bullet): (y,gH_{\zeta_0^S(y)})\mapsto (\zeta_0^S(y),K_{\zeta_0^S(y)}b(y)gH_{\zeta_0^S(y)})\]
is a coordinatization of the isotropy factor $\zeta_0^T:X\to Z_0^T$,
and the probability kernel
\[(s,K_sg'H_s)\stackrel{\rm{p}}{\mapsto} P(s,\,\cdot\,)\ltimes m_{b(\bullet)^{-1}K_sg'H_s/H_s}\]
is a version of the ergodic decomposition of $\mu$ over $\zeta_0^T$,
where for any subset $S \subseteq G_s$ we write $S/H_s := \{gH_s:\
g\in S\}$.

Moreover, the Mackey group data $K_\bullet$ is \textbf{conjugate-minimal}: if $K'_\bullet \leq G_\bullet$ is another measurable
assignment of compact subgroup data on $Z_0^S$ and $b':Y\to
G_\bullet$ another section such that the cocycle-section
$(\g,y)\mapsto b'(S^\g y)\rho(\g,y)b'(y)^{-1}$ takes a value in
$K'_{\zeta_0^S(y)}$ for $\nu$-almost every $y$ for every $\g$, then
there is a section $c:Z_0^S\to G_\bullet$ such that
\[c(s) \cdot K'_s \cdot c(s)^{-1} \geq K_s\] for
$(\zeta_0^S)_\#\nu$-almost every $s$\qed
\end{thm}

\begin{thm}\label{thm:homo-nonergMackey-2}
Suppose that $S:\G\curvearrowright (Y,\nu)$, $H_\bullet \leq
G_\bullet$ are $S$-invariant measurable compact group data and
$\rho:\G\times Y\to G_\bullet$ is a cocycle-section over $S$ and $X$
is the space $Y\ltimes G_\bullet/H_\bullet$ but equipped with some
unknown $(S\ltimes \rho)$-invariant and relatively ergodic lift
$\mu$ of $\nu$. Then there are subgroup data $K_\bullet \leq
G_\bullet$ and a section $b:Y\to G_\bullet$ such that $\mu =
\nu\ltimes m_{b(\bullet)^{-1}K_\bullet H_\bullet/H_\bullet}$. \qed
\end{thm}

As in the classical case of an ergodic base system, replacing some
given group data $G_\bullet$ with the Mackey group data $K_\bullet$
and recoordinatizing (see Corollary 3.27 in Glasner~\cite{Gla03})
gives the following corollary.

\begin{cor}\label{cor:implicit-covering-nonergMackey}
Given a $\G$-system $\bfY = (Y,\nu,S)$, measurable $S$-invariant
homogeneous space data $G_\bullet/K_\bullet$ over $Y$ and a
cocycle-section $\rho:\G\times Y\to G_\bullet$ over $S$, and
defining $X:= Y\ltimes G_\bullet/K_\bullet$ and $T:= S\ltimes\rho$,
any $(S\ltimes \rho)$-relatively ergodic lift $\mu$ of $\nu$ admits
a re-coordinatization of the canonical extension $(X,\mu,T)\to \bfY$
to $\bfY\ltimes
(G'_\bullet/H'_\bullet,m_{G'_\bullet/H'_\bullet},\rho')\to \bfY$
leaving the base system fixed (so the new lifted measure is just the
direct integral measure), and such that the implicit covering group
extension $\bfY\ltimes (G'_\bullet,m_{G'_\bullet},\rho')\to\bfY$ is
also relatively ergodic. \qed
\end{cor}

Extensions by measurable homogeneous space data acquire greater
significance through the non-ergodic version of the structure
theorems of Furstenberg~\cite{Fur77} and Zimmer~\cite{Zim76.1},
which identifies them as all the possible isometric extensions and
accounts for the overall isometric subextension of a relatively
independent join of extensions in terms of these.

\begin{thm}\label{thm:rel-ind-joinings}
Suppose that $\pi_i:\bfX_i\to\bfY_i$ are relatively ergodic
extensions for $i=1,2,\ldots,n$, that $\nu$ is a joining of
$\bfY_1$, $\bfY_2$, \ldots, $\bfY_n$ forming the system $\bfY =
(Y,\nu,S) := (Y_1\times Y_2\times\cdots\times Y_n,\nu,S_1\times
S_2\times \cdots\times S_n)$.  Suppose further that $\bfX =
(X,\mu,T)$ is similarly a joining of $\bfX_1$, $\bfX_2$, \ldots,
$\bfX_n$ that extends $\nu$ through the coordinatewise factor map
$\pi:\bfX\to\bfY$ assembled from the $\pi_i$, and such that under
$\mu$ the coordinate projections $\a_i:\bfX\to\bfX_i$ are relatively
independent over the tuple of further factors $\pi_i\circ\a_i$. Then
there are intermediate isometric extensions
\[\bfX_i\stackrel{\zeta_1}{\longrightarrow}\bfZ_i\stackrel{\pi_i|_{\zeta_i}}{\longrightarrow}\bfY_i\]
such that the intermediate factor map
\[\zeta_1\vee\zeta_2\vee\ldots\vee\zeta_3:\bfX\to\bfZ\]
whose target $\bfZ$ is the resulting joining of the systems $\bfZ_i$
is precisely a coordinatization of the maximal factor between $\bfX$
and $\bfY$ that defines an isometric extension of $\bfY$ (which
contains the relatively invariant extension $\bfZ_0^T\vee \bfY \to
\bfY$, which may be nontrivial). \qed
\end{thm}

As in~\cite{Aus--ergdirint}, and following well-known practice in
the ergodic case, we can define the maximal isometric and maximal
distal subextensions of an extension $\pi:\bfX\to \bfY$; we
generally denote the maximal $n$-step distal subextension by
\[\bfX\stackrel{\zeta_{n/\pi}^T}{\longrightarrow} \bfZ_n^T(\bfX/\pi)\stackrel{\pi|_{\zeta_{n/\pi}^T}}{\longrightarrow} \bfY\]
for some coordinatizing intermediate target system
$\bfZ_n^T(\bfX/\pi)$.

Given an ergodic system $\bfX$, the maximal isometric subextension
of the trivial factor is just the \textbf{Kronecker factor} of
$\bfX$, and as a simpler special case of
Theorem~\ref{thm:rel-ind-joinings} this factor $\pi:\bfX\to\bfZ_1^T$
can be coordinatized as an ergodic rotation action on a compact
group: that is, there are a compact core-free homogeneous space
$G/H$ and a homomorphism $\phi:\G\longrightarrow G$ with dense image such that
\[T|_\pi^\g(gH) = \phi(\g)gH\quad\quad \g\in\G,\ g\in G.\]
In this case we will sometimes write $(G/H,m_{G/H},\phi)$ in place
of $(G/H,m_{G/H},T|_\pi)$.

\subsection{Factors and automorphisms of isometric extensions}

Two of the main results of~\cite{Aus--ergdirint} are structure
theorems for factors and automorphisms of relatively ergodic
extensions by compact homogeneous space data.  These can be deduced
quite simply after an appropriate change of viewpoint: considering
instead the graphical joining associated to a factor map or
automorphism, we obtain an extension from a smaller to a larger
joining that is also coordinatized by compact homogeneous space
data, and so to which the non-ergodic Mackey Theory can be applied.
The structure of the factor map or automorphism can then be
recovered from the Mackey data for this joining.  We refer the
reader to Section 6 of~\cite{Aus--ergdirint} for details, only
recalling here some notation and the two particular results that we
need.

First, if $(X_i,\mu_i) = (Y,\nu)\ltimes
(G_{i,\bullet}/H_{i,\bullet},m_{G_{i,\bullet}/H_{i,\bullet}})$ for
$i=1,2$ are two different extensions of a standard Borel probability
space by homogeneous space data, $R$ is a probability-preserving
transformation of $(Y,\nu)$ and if in addition we are given a
section of homomorphisms $\Phi_\bullet:G_{1,\bullet}\longrightarrow
G_{2,R(\bullet)}$ such that $\Phi_\bullet(H_{1,\bullet}) =
H_{2,R(\bullet)}$ and another section $b:Y\to G_{2,R(\bullet)}$,
then we write
\[\a = R\ltimes (L_{b(\bullet)}\circ
\Phi_\bullet)|^{H_{1,\bullet}}_{H_{2,\bullet}}:(X_1,\mu_1)\to
(X_2,\mu_2)\] for the map defined as an extension of $R:(Y,\nu)\to
(Y,\nu)$ by the fibrewise action of the affine endomorphisms
associated to $\Phi_\bullet$ and left-multiplication by $b$.  Note
that the condition $\Phi_\bullet(H_{1,\bullet}) = H_{2,R(\bullet)}$
is needed for this formula for $\a$ to make sense at all.

Once again let $\G$ be an arbitrary locally compact second countable
group. Our main results here are that relative factors and
automorphisms can all be described in terms of such data.

\begin{thm}[Relative Factor Structure Theorem]\label{thm:RFST}
Suppose that $\bfY = (Y,\nu,S)$ is a $\G$-system, that
$G_{i,\bullet}/H_{i,\bullet}$ are $S$-invariant core-free
homogeneous space data on $Y$ and that $\s_i:\G\times Y\to
G_{i,\bullet}$ are ergodic cocycle-sections for the action $S$, and
let $\bfX_i = (X_i,\mu_i,T_i) := \bfY\ltimes
(G_{i,\bullet}/H_{i,\bullet},\s_i)$. Suppose further that $\bfX_2\to
\bfY$ admits insertion as a subextension of $\bfX_1\to\bfY$:
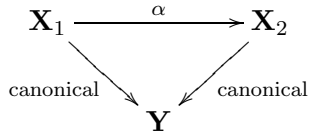
\begin{center}
$\phantom{i}$\xymatrix{ \bfX_1\ar[rr]^\a\ar[dr]_{\rm{canonical}} & &
\bfX_2\ar[dl]^{\rm{canonical}}\\ & \bfY }
\end{center}
Then there are an $S$-invariant measurable family of epimorphisms
$\Phi_\bullet:G_{1,\bullet} \to G_{2,\bullet}$ such that
$\Phi_\bullet(H_{1,\bullet}) = H_{2,\bullet}$ almost surely and a
section $b:Y\to G_{2,\bullet}$ such that $\a = \id_Y\ltimes
(L_{b(\bullet)}\circ
\Phi_\bullet)|^{H_{1,\bullet}}_{H_{2,\bullet}}$, $\mu_1$-almost
surely. \qed
\end{thm}

The conclusion for automorphisms is very similar.

\begin{thm}[Relative Automorphism Structure Theorem]\label{thm:RAST}
Suppose that $\bfY = (Y,\nu,S)$ is a $\G$-system, that
$G_\bullet/H_\bullet$ are $S$-invariant core-free homogeneous space
data on $Y$ and that $\s:\G\times Y\to G_\bullet$ is an ergodic
cocycle-section for the action $S$, and let $\bfX = (X,\mu,T) :=
\bfY\ltimes (G_\bullet/H_\bullet,\s)$. Suppose further that
$R:\L\curvearrowright (X,\mu)$ is an action of a discrete group $\L$
that commutes with $T$ and respects the canonical factor map
$\pi:\bfX\to\bfY$, and so defines an automorphism of this extension
of $\G$-actions. Then for each $h\in\L$ there are an $S$-invariant
measurable family of isomorphisms $\Phi_{h,\bullet}:G_\bullet \to
G_{R|^h_\pi(\bullet)}$ such that $\Phi_{h,\bullet}(H_\bullet) =
H_{{R|^h_{\pi}(\bullet)}}$ almost surely and a section $\rho_h:Y\to
G_{R|^h_{\pi}(\bullet)}$ such that \[R^h = R|^h_\pi\ltimes
(L_{\rho_h(\bullet)}\circ
\Phi_{h,\bullet})|^{H_\bullet}_{H_{R|^h_{\pi}(\bullet)}}\] for each
$h\in \L$, and then
\begin{itemize}
\item we have
\[\s(\g,R|_\pi^h(y)) = \rho_h(S^\g y)\cdot \Phi_{h,y}(\s(\g,y))\cdot \rho_h(y)^{-1}\]
for $\nu$-almost all $y$ for all $\g \in \G$ and $h \in \L$, and
\item we have
\[\Phi_{h_1h_2,y} = \Phi_{h_1,R|_{\pi}^{h_2}(y)}\circ \Phi_{h_2,y}\]
and
\[\rho_{h_1h_2}(y) = \rho_{h_1}(R|_\pi^{h_2}(y))\cdot \Phi_{h_1,R|_\pi^{h_2}(y)}(\rho_{h_2}(y))\]
for $\nu$-almost all $y$ for all $h_1,h_2 \in \L$. \qed
\end{itemize}
\end{thm}

\subsection{Some auxiliary notation for Abelian cocycles}

It will help us to collect here some convenient notation for the
more detailed study of Abelian cocycles over special kinds of
system. This is partly motivated by the recent paper of Bergelson,
Tao and Ziegler~\cite{BerTaoZie09}.

First, for any system $T:\G\curvearrowright (X,\mu)$, Polish Abelian group $A$ and measurable function $\s:X\to A$ we denote by $\DDelta_T\s:\G\times X\to A$ the
resulting coboundary:
\[\DDelta_T\s(\g,x) := \s(T^\g x)\cdot \s(x)^{-1}.\]
If $X$ has the structure of a compact Abelian group and $T = R_\phi$ is the rotation action corresponding to a homomorphism
$\phi:\G\longrightarrow X$, then we will generally abbreviate
$\DDelta_{R_\phi}$ to $\DDelta_\phi$.  The slightly unconventional notation `$\DDelta$' will be a reminder that we write the group operation of $A$ multiplicatively, even though it is Abelian.

In addition, we write $\C(X;A)$ for the group of all Borel maps
$X\to A$ under pointwise multiplication.  Given an action $T$, we write $\Z^1(T;A)$ for the collection of all its Borel cocycles $\G\times X\to A$ and $\B^1(T;A)$ for the subcollection of its $A$-valued coboundaries. As is standard, since $A$ is Abelian, $\B^1(T;A) \leq \Z^1(T;A)$ are groups under pointwise multiplication. If $\pi:(X,\mu)\to (Y,\nu)$ then we write $\Z^1(T|_\pi;A)\circ\pi$ for the subgroup of all $\phi \in \Z^1(T;A)$ for which $\phi = \psi\circ\pi$ for some $\psi \in \Z^1(T|_\pi;A)$.

\subsection{Fibre-normality}\label{subs:fibre-normal}

Alongside the notion of sated extensions that we have brought
from~\cite{Aus--lindeppleasant1}, we will now introduce another
general property enjoyed by some systems and show that we may always
pass to extensions where this property obtains.

Importantly, henceforth we will assume that $\G = \bbZ^d$, and will
consider also an arbitrary subgroup $\L \leq \bbZ^d$.  Many of the
results below could be extended unchanged to the setting of a
discrete group $\G$ and a \emph{central} subgroup $\L \leq \G$, but
even the case of $\L \unlhd \G$ with the conjugation action
$\G\curvearrowright \L$ nontrivial introduces new subtleties that we
do not wish to address here.

\begin{dfn}\label{dfn:fib-nor}
A relatively ergodic extension of systems $\a:\bfX\to \bfY$ is
\textbf{fibre-normal} if the maximal isometric subextension
$\a|_{\zeta_{1/\a}^T}:\bfZ_1^T(\bfX/\a) \to \bfY$ can be
coordinatized as an extension by measurable group data:
\begin{center}
$\phantom{i}$\xymatrix{
\bfZ_1^T(\bfX/\a)\ar[dr]_{\a|_{\zeta_{1/\a}^T}}\ar@{<->}[rr]^-\cong
& & \bfY\ltimes
(G_\bullet,m_{G_\bullet},\s)\ar[dl]^{\rm{canonical}}\\ & \bfY }
\end{center}
(rather than just homogeneous space data, as is always possible by
the results of Section 5 in~\cite{Aus--ergdirint}).

Equivalently, we will write that $\bfX$ is \textbf{fibre-normal over
the factor} $\a$ or that $(X,\mu)$ is $T$-\textbf{fibre-normal over
the factor} $\a$.  If $(\bfY,\a) = (\sfC \bfX,\zeta^\bfX_\sfC)$ for
some idempotent class $\sfC$ then we will write that $\bfX$ is
\textbf{fibre-normal over} $\sfC$.
\end{dfn}

This definition --- and the use to which we will put it --- is
strongly motivated by that of Furstenberg and Weiss' `normal'
systems in Section 8 of~\cite{FurWei96}. We will see its value very
concretely in the proof of Lemma~\ref{lem:2D-proj-full}, at which
point an analogous proof involving extensions by arbitrary
homogeneous space data would be considerably more grueling. The
ability to pass to fibre-normal extensions may also be of some
independent interest.

The main goals of this subsection are to show that for an order
continuous idempotent class $\sfC$ any system admits an extension
that is fibre-normal over $\sfC$ for any given subaction $T^{\
\uhr\L}$, and that fibre-normality over order continuous idempotent
classes is preserved under inverse limits.  We will eventually apply
these results to the idempotent class $\sfZ_0^{\bf{p}_1}\vee
\bigvee_{i=2}^k\sfZ_0^{\bf{p}_1- \bf{p}_i}$ and its relatives.  Our
proof partly follows that of Furstenberg and Weiss for their
instance of fibre-normality in~\cite{FurWei96}, although in other
ways we take a slightly different route (avoiding, in particular,
their use of the abstract characterization of extensions by compact
group data in terms of graphical self-joinings given in their Lemma
8.5, and originally traceable to work of Veech~\cite{Vee82}).

\begin{prop}\label{prop:fibre-normals-exist}
If $\sfC$ is an order continuous idempotent class and $\L\leq \G$ is
a fixed subgroup then every $\G$-system $\bfX_0$ admits an extension
$\pi:\bfX\to \bfX_0$ that is both $(\sfZ_0^\L\vee \sfC)$-sated and
$T^{\ \uhr\L}$-fibre-normal over $\sfC$.
\end{prop}

\begin{prop}\label{prop:fibre-normals-stable}
If $\sfC$ is an order continuous idempotent class, and a given
inverse sequence consists of systems all of which are both
$(\sfZ_0^\L\vee \sfC)$-sated and have $\L$-subaction fibre-normal
over $\sfC$, then this is also true of its inverse limit.
\end{prop}

\textbf{Example}\quad The assumption of satedness alongside
fibre-normality in Proposition~\ref{prop:fibre-normals-stable} is
essential.  We will give an example to show this with $\G = \L =
\bbZ^2$ and $\sfC := \sfZ_0^{\bf{e}_1}\vee \sfZ_0^{\bf{e}_2}$.

First let $\pi_m:\bbT^\bbN\to \bbT^m$ be the initial coordinate
projection, and
\[\bfX_{(0)} = (X_{(0)},\mu_{(0)},T_{(0)}):= (\bbT^\bbN,m_{\bbT^\bbN},\phi)\ltimes (G/H,m_{G/H},\s)\]
where
\begin{itemize}
\item $\phi$ is a dense homomorphic embedding
\[\phi:\bbZ^2\to \bbT^\bbN:(m,n)\mapsto (m+n)\cdot w\]
where $w = (w_1,w_2,\ldots)$ is a sequence of irrational and
rationally independent $w_i\in\bbT$ and so has a dense orbit in
$\bbT^\bbN$, \item $G/H$ is a core-free compact metrizable
homogeneous space with $H\neq \{1_G\}$, \item and $\s:\bbZ^2\times
\bbT^\bbN \to G$ is any ergodic cocycle over $T_{(0)}$ such that
$\s\cdot N$ is not cohomologous to a cocycle measurable with respect
to $\pi_m$ for any finite $m$ and proper $N \lhd G$ (it is easy to
see that a generic $\s$ has this property for many choices of $G$,
such as $G = \rm{O}(3)$).
\end{itemize}
Note that $\bfX_{(0)}$ is a $\sfZ_0^{\bf{e}_1 - \bf{e}_2}$-system,
but that $T_{(0)}^{\bf{e}_i}$ is ergodic for $i=1,2$. Let
$\a:\bbT^\bbN \times G/H\to \bbT^\bbN$ be the canonical factor map,
and note that $\pi_m\circ\a$ is a smaller factor map onto the
finite-dimensional Abelian group rotation
$(\bbT^m,m_{\bbT^m},\pi_m\circ\phi)$.

Now for each $m=1,2,\ldots,\infty$ let $\bfY_{(m)} := (\bbT^m\times
\bbT^m,m_{\bbT^m\times \bbT^m},\rho_m)$ and $\xi_m:\bfY_{(m)}\to
(\bbT^m,m_{\bbT^m},\pi_m\circ\phi)$ be the factor map $\xi_m(s,t) =
s+t$ with $\rho_m(\bf{e}_1) = (w_1,w_2,\ldots,w_m,0,0,\ldots,0)$ and
$\rho_m(\bf{e}_2) = (0,0,\ldots,0,w_1,w_2,\ldots,w_m)$ (with the
obvious interpretation when $m=\infty$). It is clear that this
defines an ergodic $\bbZ^2$-action and that the factor map $\xi_m$
does indeed map $\rho_m$ onto $\pi_m\circ\phi$ (and hence intertwine
the two corresponding rotation actions). Finally, let
\[\bfX_{(m)}:= \bfY_{(m)}\times_{\{\xi_m = \pi_m\circ\a\}}\bfX_{(0)}\]
and for $m < \infty$ let $\psi^{(m+1)}_{(m)}:\bfX_{(m+1)}\to
\bfX_{(m)}$ be the obvious factor map defined by lifting the map
$Y_{(m+1)}\to Y_{(m)}:(s,t)\mapsto (\pi_m(s),\pi_m(t))$.

Now it is easy to check, firstly, that $\bfX_{(m)}\to
\sfC\bfX_{(m)}$ is simply equivalent to the coordinate projection
factor map $\bfX_{(m)}\to \bfY_{(m)}$; and secondly that the maximal
isometric subextension of $\bfX_{(m)}\to \bfY_{(m)}$ is equivalent
to $\id_{Y_{(m)}}\times\a$: that is, that the fibre copies of $G/H$
in $\bfX_{(0)}$ are not retained in this maximal isometric
subextension, because this would require that for some proper $N\lhd
G$ the cocycle $\s\cdot N$ be measurable with respect to $\pi_m$. As
a result, each $\bfX_{(m)}$ is fibre-normal over $\sfC$.  On the
other hand, $\bfX_{(\infty)}$ can be identified with the inverse
limit of the inverse sequence $(\bfX_{(m)})_{m\geq 0}$,
$(\psi^{(m)}_{(k)})_{m\geq k\geq 0}$, but now $\sfC\bfX_{(\infty)}$
is the whole of the underlying group rotation
$(\bbT^\bbN\times\bbT^\bbN,m_{\bbT^\bbN\times
\bbT^\bbN},\psi_\infty)$, with respect to which the cocycle $\s$
\emph{is} measurable, and so now the maximal isometric subextension
of $\bfX_{(\infty)}\to\sfC\bfX_{(\infty)}$ is simply the whole of
$\bfX_{(\infty)}$, which involves the non-normal homogeneous space
fibres $G/H$ and so is not fibre-normal.

Intuitively, the phenomenon observed above is possible because, as
we ascend through the systems $\bfX_{(m)}$ for increasing $m$, their
maximal $\sfC$-factors determine increasingly large factors of the
original base system $\bfX_{(0)}$, until at precisely the point of
taking the inverse limit the $\sfC$-factor determines a large enough
factor of $\bfX_{(0)}$ that the maximal isometric extension can
capture some new, larger fibres that are not normal.  It is this
possibility, and some more complicated variations, that the
additional assumption of satedness prevents. In this connexion we
remark that this subtlety did not arise in Furstenberg and Weiss'
original use of fibre-normality (just `normality', in their
terminology) in~\cite{FurWei96}, because they were concerned only
with fibre-normality over Kronecker factors: that is, over the
idempotent class $\sfZ_1$, which is hereditary and hence
always-sating. By contrast, in this work we will often be concerned
with fibre-normality over joins of several isotropy factors, and we
have seen that in general such joins are \emph{not} always sating
and so genuinely require greater care, as shown by the above
example. \fin

The proofs of both of both
Propositions~\ref{prop:fibre-normals-exist}
and~\ref{prop:fibre-normals-stable} will involve heavy use of
inverse limits.

\begin{lem}\label{lem:order-cty-of-isom-exts}
Suppose that $\sfC$ is an order continuous idempotent class and that
$(\bfX_{(m)})_{m\geq 0}$, $(\psi^{(m)}_{(k)})_{m\geq k\geq 0}$ is a
$(\sfZ_0^\L\vee \sfC)$-sated inverse sequence with inverse limit
$\bfX_{(\infty)}$, $(\psi_{(m)})_{m\geq 0}$. Then
\[\zeta_{1/\sfC}^{T_{(\infty)}^{\ \uhr\L}}
\simeq \bigvee_{m\geq 0}\zeta_{1/\sfC}^{T_{(m)}^{\
\uhr\L}}\circ\psi_{(m)}.\]
\end{lem}

\textbf{Proof}\quad The relation
\[\zeta_{1/\sfC}^{T_{(\infty)}^{\ \uhr\L}} \succsim \bigvee_{m\geq 0}\zeta_{1/\sfC}^{T_{(m)}^{\ \uhr\L}}\circ\psi_{(m)}.\]
is clear by monotonicity, so it remains only to prove its reverse
$\precsim$.

By Lemma 3.6 of~\cite{Aus--lindeppleasant1} the class $\sfD :=
\sfZ_0^\L\vee \sfC$ is still order continuous. Also we have by
definition (see 5.11 in~\cite{Aus--ergdirint}) that $\zeta^{T^{\
\uhr\L}}_{1/\sfC} = \zeta^{T^{\ \uhr\L}}_{1/\sfD}$ for any $\bfX =
(X,\mu,T)$, and know from the non-ergodic Furstenberg-Zimmer Theory
that this is precisely the maximal factor of $\bfX$ generated by all
the finite-rank $T^{\ \uhr\L}$-invariant
$\zeta_\sfD^\bfX$-submodules of $L^2(\mu)$.  It will therefore
suffice to show that any $T_{(\infty)}^{\ \uhr\L}$-invariant
finite-rank $\zeta_\sfD^{\bfX_{(\infty)}}$-submodule $\frM \leq
L^2(\mu_{(\infty)})$ can be approximated by $\psi_{(m)}$-lifts of
$T_{(m)}^{\ \uhr\L}$-invariant finite-rank
$\zeta_\sfD^{\bfX_{(m)}}$-submodules of $L^2(\mu_{(m)})$ by taking
$m$ sufficiently large.

Since $\bfX_{(\infty)}$, $(\psi_{(m)})_{m\geq 0}$ is the inverse
limit we have
\[\id_{X_{(\infty)}} \succsim \bigvee_{m\geq 0}(\zeta_\sfD^{\bfX_{(\infty)}}\vee\psi_{(m)}) \succsim \bigvee_{m\geq 0}\psi_{(m)} \simeq \id_{X_{(\infty)}},\]
so in fact all these factor maps are equivalent. Let $\phi_1$,
\ldots, $\phi_d$ be an orthonormal basis for a $T_{(\infty)}^{\
\uhr\L}$-invariant finite-rank
$\zeta_\sfD^{\bfX_{(\infty)}}$-submodule $\frM \leq
L^2(\mu_{(\infty)})$. On the one hand, each $\phi_i$ can be
$L^2$-approximated by the
$(\zeta_\sfD^{\bfX_{(\infty)}}\vee\psi_{(m)})$-measurable functions
$\sfE_{\mu_{(\infty)}}(\phi_i\,|\,\zeta_\sfD^{\bfX_{(\infty)}}\vee\psi_{(m)})$
by taking $m$ sufficiently large. On the other, by definition there
is a $d\times d$ matrix of measurable functions $U_{i,j}:\G\times
\sfD X_{(\infty)}\to \bbC$ such that
\[\phi_i(T_{(\infty)}^\g(x)) = \sum_{j=1}^dU_{i,j}(\g,\zeta_\sfD^{\bfX_{(\infty)}}(x))\cdot\phi_j(x)\]
for $\g \in \L$ and $\mu_{(\infty)}$-a.e. $x\in X_{(\infty)}$.
Taking conditional expectation with respect to
$\zeta_\sfD^{\bfX_{(\infty)}}\vee\psi_{(m)}$ and bearing in mind
that $U_{i,j}$ is already $\zeta_\sfD^{\bfX_{(\infty)}}$-measurable
we obtain
\[\sfE_{\mu_{(\infty)}}(\phi_i\,|\,\zeta_\sfD^{\bfX_{(\infty)}}\vee\psi_{(m)})\circ T_{(\infty)}^\g = \sum_{j=1}^dU_{i,j}(\g,\zeta_\sfD^{\bfX_{(\infty)}}(\,\cdot\,))\cdot\sfE_{\mu_{(\infty)}}(\phi_j\,|\,\zeta_\sfD^{\bfX_{(\infty)}}\vee\psi_{(m)}),\]
so the conditional expectations
$\sfE_{\mu_{(\infty)}}(\phi_i\,|\,\zeta_\sfD^{\bfX_{(\infty)}}\vee\psi_{(m)})$
are lifted from a finite-rank
$(\zeta_\sfD^{\bfX_{(\infty)}}|_{\zeta_\sfD^{\bfX_{(\infty)}}\vee\psi_{(m)}})$-submodule
of
$L^2((\zeta_\sfD^{\bfX_{(\infty)}}\vee\psi_{(m)})_\#\mu_{(\infty)})$
that is invariant under the restriction to this factor of
$T_{(\infty)}^{\ \uhr\L}$, and as $m\to\infty$ these submodules
approximate $\frM$ in $L^2$.

So far we have not used the satedness of our inverse sequence; we
will need this to obtain a further approximation by finite-rank
submodules of $L^2(\mu_{(m)})$. This follows because by satedness
the joining of $\sfD\bfX_{(\infty)}$ and $\bfX_{(m)}$ as factors of
$\bfX_{(\infty)}$ must be relatively independent over
$\sfD\bfX_{(m)}$, and therefore by the Furstenberg-Zimmer Structure
Theorem~\ref{thm:rel-ind-joinings} any finite-rank
$(\zeta_\sfD^{\bfX_{(\infty)}}|_{\zeta_\sfD^{\bfX_{(\infty)}}\vee\psi_{(m)}})$-submodule
\[\frak{N} \leq
L^2((\zeta_\sfD^{\bfX_{(\infty)}}\vee\psi_{(m)})_\#\mu_{(\infty)})\]
that is invariant under the restricted $\L$-subaction must be
measurable with respect to $\zeta_\sfD^{\bfX_{(\infty)}}\vee
\zeta_{1/\sfD}^{T_{(m)}}$. Hence any $f \in \frak{N}$ can be
approximated arbitrarily well by finite sums of products of the form
$\sum_p g_p\cdot h_p$ with each $g_p$ being
$\zeta_\sfD^{\bfX_{(\infty)}}$-measurable and each $h_p$ being
$\zeta_{1/\sfD}^{T_{(m)}}$-measurable for some finite $m$. Now the
order continuity of $\sfD$ implies that by taking $m$ sufficiently
large we can further approximate each $g_p$ in this finite sum by
some $\zeta_\sfD^{\bfX_{(m)}}$-measurable function $g'_p$, and now
$\sum_p g'_p\cdot h_p$ is an approximation to $f$ that is a
$\psi_{(m)}$-measurable function obtained from a $T_{(m)}^{\
\uhr\L}$-invariant finite-rank $\zeta_\sfD^{\bfX_{(m)}}$-submodule
of $L^2(\mu_{(m)})$, as required. This completes the proof. \qed

\textbf{Remark}\quad An example similar to that given previously
shows that the hypothesis that each $\bfX_{(m)}$ is $(\sfZ_0^\L\vee
\sfC)$-sated (or at least that the factors $\zeta_{\sfZ_0^\L\vee
\sfC}^{\bfX_{(m+1)}}$ and $\psi^{(m+1)}_{(m)}$ of $\bfX_{(m+1)}$ be
relatively independent over $\zeta_{\sfZ_0^\L\vee
\sfC}^{\bfX_{(m)}}\circ \psi^{(m+1)}_{(m)}$ for each $m\geq 0$) is
not superfluous here. \fin

\begin{lem}\label{lem:still-f-n-after-joining}
If $\bfX = \bfY\ltimes (G_\bullet,m_{G_\bullet},\rho)$ is a
relatively ergodic extension by compact group data with canonical
factor $\pi:\bfX \to \bfY$, $\pi':\bfY'\to \bfY$ is any other
extension and $\l$ is any $(\mu,\nu')$-joining supported on $X' :=
X\times_{\{\pi = \pi'\}} Y'$ and relatively ergodic over the
canonical factor map onto $\bfY$, then the natural extension
$(X',\l,T\times S') \to \bfY'$ is also coordinatizable as an
extension by compact group data.
\end{lem}

\textbf{Proof}\quad This follows from the non-ergodic Mackey
Theorem~\ref{thm:homo-nonergMackey-2}.  As a standard Borel system
we have by definition that
\[(X',T') = (Y'\ltimes G_{\pi'(\bullet)},S'\ltimes (\rho\circ\pi')),\]
and so that theory gives us Mackey group data $M_\bullet \leq
G_{\pi'(\bullet)}$ and a section $b:Y' \to G_\bullet$ and an
$S'$-invariant section $g:Y' \to G_\bullet$ such that $\l =
\nu'\ltimes m_{b(\bullet)^{-1}M_\bullet g(\bullet)}$.  Now
re-coordinatizing by the fibre-wise isomorphism
\[(y',g') \mapsto (y',b(y')g'g(y'))\]
this gives a coordinatization of $(X',\mu',T')\to (Y',\nu',S')$ by
the compact group data $M_\bullet$ with the relatively ergodic
cocycle $(\g,y')\mapsto b((S')^\g y')\rho(\g,\pi'(y))b(y')$, which
is of the required form. \qed

\begin{lem}\label{lem:still-f-n-after-joining-2}
If $\pi:\bfX\to\bfY$ is a factor and
$\bfX\stackrel{\a_m}{\longrightarrow}\bfZ_m\stackrel{\pi|_{\a_m}}{\longrightarrow}\bfY$,
$m\geq 1$, is a family of intermediate factors each of which can be
coordinatized by compact group data, then so can their join
\[\bfX \stackrel{\a_1 \vee\a_2\vee\cdots}{\longrightarrow} \bfZ\stackrel{\pi|_{\a_1\vee \a_2\vee\cdots}}{\longrightarrow} \bfY.\]
\end{lem}

\textbf{Proof}\quad Having chosen coordinatizations by compact group
data
\begin{center}
$\phantom{i}$\xymatrix{
\bfZ_m\ar[dr]_{\pi|_{\a_m}}\ar@{<->}[rr]^-\cong & & \bfY\ltimes
(G_{m,\bullet},m_{G_{m,\bullet}},\s_m)\ar[dl]^{\rm{canonical}}\\ &
\bfY }
\end{center}
we can glue these together to coordinatize $\bfZ\to \bfY$ using the
compact group data $G_\bullet := \prod_{m\geq 1}G_{m,\bullet}$ and
cocycle-section $(\s_m)_{m\geq 1}$ and some invariant measure on
$Y\ltimes G_\bullet$ obtained from the joining.  Now the non-ergodic
Mackey Theory allows us to find some Mackey subgroup data
$M_\bullet$ for this extension and convert its coordinatization into
a coordinatization by a relatively ergodic cocycle-section using
that compact group data just as for the previous lemma, completing
the proof. \qed

\textbf{Proof of Proposition~\ref{prop:fibre-normals-exist}}\quad
Once again let $\sfD := \sfZ_0^\L\vee \sfC$. We specify recursively
an inverse sequence of extensions, similar to that in the proof of
Theorem 8.8 of Furstenberg and Weiss in~\cite{FurWei96}, as follows.
First set $\bfX_{(0)} := \bfX_0$, and now proceed as follows.
\begin{itemize}
\item When $m$ is even let
$\psi^{(m+1)}_{(m)}:\bfX_{(m+1)}\to\bfX_{(m)}$ be a $\sfC$-sated
extension.
\item When $m$ is odd, let
\begin{center}
$\phantom{i}$\xymatrix{ \big(\bfZ_1^{T^{\
\uhr\L}_{(m)}}(\bfX_{(m)}/\sfD)\big)^{\
\uhr\L}\ar[dr]_{\zeta_\sfD}\ar@{<->}[rr]^-\cong & &\sfD\bfX_{(m)}^{\
\uhr\L}\ltimes
(G_{m,\bullet}/H_{m,\bullet},m_{G_{m,\bullet}/H_{m,\bullet}},\s_m)\ar[dl]^{\rm{canonical}}\\
& (\sfD\bfX_{(m)})^{\ \uhr\L} }
\end{center}
be a coordinatization of the $T_{(m)}^{\ \uhr\L}$-isometric
extension of the $\L$-subactions using core-free homogeneous space
data and an ergodic cocycle-section $\s_m$. Implicitly this
coordinatization specifies a covering group extension
\[\pi':\sfD\bfX_{(m)}^{\ \uhr\L}\ltimes
(G_{m,\bullet},m_{G_{m,\bullet}},\s_m) \to \big(\bfZ_1^{T_{(m)}^{\
\uhr\L}}(\bfX_{(m)}/\sfD)\big)^{\ \uhr\L},\] and we now recall from
the Relative Factor Structure Theorem that the whole $\G$-action on
the target of this factor map can be lifted to give an action of the
whole group $\G$ upstairs, so that we may express
\[\sfD\bfX_{(m)}^{\ \uhr\L}\ltimes
(G_{m,\bullet},m_{G_{m,\bullet}},\s_m) =\bfY_{(m)}^{\ \uhr\L}\] for
some $\G$-system $\bfY_{(m)}$. Finally let
\[\bfX_{(m+1)}:= \bfY_{(m)} \otimes_{\{\pi' =
\zeta^{T_{(m)}}_{1/\sfD}\}} \bfX_{(m)}\] and
$\psi^{(m+1)}_{(m)}:\bfX_{(m+1)}\to\bfX_{(m)}$ be the second
coordinate factor map back onto $\bfX_{(m)}$.  In addition, let us
introduce the auxiliary notation
\[\eta_{(m+1)}:\bfX_{(m+1)}\to \bfY_{(m)}\] for the first coordinate
projection.  The important feature here is that by construction the
factor $\bfZ_1^{T^{\ \uhr\L}_{(m)}}(\bfX_{(m)}/\sfD)$ of
$\bfX_{(m)}$ which is $\L$-isometric over
$\zeta_{\sfD}^{\bfX_{(m)}}$ has now been swallowed by the factor
$\bfY_{(m)}$ which is $\L$-isometric and fibre-normal over the copy
$\zeta_{\sfD}^{\bfX_{(m)}}\circ\psi^{(m+1)}_{(m)}$.
\end{itemize}

The main difference between this construction and that of
Furstenberg and Weiss in~\cite{FurWei96} is that we must interleave
extensions that enlarge homogeneous space fibres to their covering
group fibres with extensions that recover full isotropy satedness.
Nevertheless, the proof we will offer that the final inverse limit
extension has the desired fibre normality essentially follows
theirs.

Let $\bfX_{(\infty)}$, $(\psi_{(m)})_{m\geq 0}$ be the inverse limit
of the above inverse sequence; we will show that it has the desired
satedness and fibre-normality.

On the one hand, the cofinal inverse subsequence
$(\bfX_{(m)})_{m\geq 0\ \rm{even}}$, $(\psi^{(m)}_{(k)})_{m\geq
k\geq 0\ \rm{even}}$ is $\sfD$-sated by construction.  It follows by
Lemma 3.12 of~\cite{Aus--lindeppleasant1} that $\bfX_{(\infty)}$ is
also $\sfD$-sated, and also by
Lemma~\ref{lem:order-cty-of-isom-exts} that
\[\zeta_{1/\sfD}^{T^{\ \uhr\L}_{(\infty)}}\simeq \bigvee_{m\geq 0\ \rm{even}}\zeta_{1/\sfD}^{T^{\ \uhr\L}_{(m)}}\circ\psi_{(m)}\]
(recall that this required the satedness assumption). Since
\[\zeta_{1/\sfD}^{T^{\ \uhr\L}_{(m)}}\circ\psi_{(m)}
\precsim \zeta_\sfD^{\bfX_{(\infty)}}\vee\big(\zeta_{1/\sfD}^{T^{\
\uhr\L}_{(m)}}\circ\psi_{(m)}\big) \precsim \zeta_{1/\sfD}^{T^{\
\uhr\L}_{(\infty)}}\] this implies by sandwiching that
\[\zeta_{1/\sfD}^{T^{\ \uhr\L}_{(\infty)}}\simeq \bigvee_{m\geq 0\ \rm{even}}\big(\zeta_\sfD^{\bfX_{(\infty)}}\vee\big(\zeta_{1/\sfD}^{T^{\ \uhr\L}_{(m)}}\circ\psi_{(m)}\big)\big),\]
and so since also
\begin{multline*}
\zeta_\sfD^{\bfX_{(\infty)}}\vee\big(\zeta_{1/\sfD}^{T^{\
\uhr\L}_{(m)}}\circ\psi_{(m)}\big) \precsim
\zeta_\sfD^{\bfX_{(\infty)}}\vee(\eta_{(m+1)}\circ\psi_{(m+1)})\\
\precsim \zeta_\sfD^{\bfX_{(\infty)}}\vee\big(\zeta_{1/\sfD}^{T^{\
\uhr\L}_{(m+2)}}\circ\psi_{(m+2)}\big)
\end{multline*}
for even $m$ we obtain that
\[\zeta_{1/\sfD}^{T^{\ \uhr\L}_{(\infty)}}\simeq \bigvee_{m\geq 0\ \rm{even}}\big(\zeta_\sfD^{\bfX_{(\infty)}}\vee(\eta_{(m+1)}\circ\psi_{(m+1)})\big).\]

On the other hand, the extension $\eta_{(m+1)}\circ\psi_{(m+1)}
\succsim \zeta_\sfD^{\bfX_{(m)}}\circ\psi_{(m)}$ is isometric and
fibre normal (we constructed it as a relatively ergodic covering
group data extension), and so by
Lemma~\ref{lem:still-f-n-after-joining} the extension
\[\zeta_\sfD^{\bfX_{(\infty)}}\vee(\eta_{(m+1)}\circ\psi_{(m+1)})
\to \zeta_\sfD^{\bfX_{(\infty)}}\] is also isometric and
fibre-normal. Therefore $\zeta_{1/\sfD}^{T_{(\infty)}^{\ \uhr\L}}$
can be expressed as a join of extensions of
$\zeta_\sfD^{\bfX_{(\infty)}}$ by compact group data, and so by
Lemma~\ref{lem:still-f-n-after-joining-2} it can itself be
coordinatized in that form.  This gives the desired fibre-normality.
\qed

\textbf{Remark}\quad In general, it can happen that the maximal
isometric extension $\zeta_{1/\sfD}^{T_{(m+1)}}\to
\zeta_{1/\sfD}^{T_{(m)}}\circ\psi^{(m+1)}_{(m)}$ is properly larger
than the extension $\eta_{(m+1)}\to
\zeta_{1/\sfD}^{T_{(m)}}\circ\psi^{(m+1)}_{(m)}$, hence the care we
had to exercise in obtaining the joining expression for
$\zeta_{1/\sfD}^{T^{\ \uhr\L}_{(m)}}$ that we eventually used in the
above proof.  This follows easily from constructions similar to the
example that follows the statement of
Proposition~\ref{prop:fibre-normals-stable}. As a result, the larger
maximal isometric extension can again require nontrivial homogeneous
space data (that is, it can fail to be fibre-normal), so we could
not use it directly in setting up the above appeal to
Lemma~\ref{lem:still-f-n-after-joining-2}. \fin

\textbf{Proof of Proposition~\ref{prop:fibre-normals-stable}}\quad
This essentially follows from the argument above: if
$(\bfX_{(m)})_{m\geq 0}$, $(\psi^{(m)}_{(k)})_{m\geq k \geq 0}$ is a
$(\sfZ_0^\L\vee\sfC)$-sated inverse sequence with all members
$T_{(m)}^{\ \uhr\L}$-fibre-normal over $\sfC$ and with inverse limit
$\bfX_{(\infty)}$, $(\psi_{(m)})_{m\geq 0}$, then the extension
$\zeta_{1/\sfC}^{\bfX_{(\infty)}} \to \zeta_{\sfZ_0^\L\vee
\sfC}^{\bfX_{(\infty)}}$ can be identified with the join of the
extensions
\[\zeta_{\sfZ_0^\L\vee \sfC}^{\bfX_{(\infty)}}\vee(\zeta_{1/\sfC}^{\bfX_{(m)}}\circ\psi_{(m)}) \to \zeta_{\sfZ_0^\L\vee \sfC}^{\bfX_{(\infty)}},\]
and each of these can be coordinatized by compact group data by
Lemma~\ref{lem:still-f-n-after-joining} and hence so can their join
by Lemma~\ref{lem:still-f-n-after-joining-2}. \qed

Now a final simple inverse-limit argument (very similar to that for
the existence of multiply sated extensions in Theorem 3.11
of~\cite{Aus--lindeppleasant1}) immediately gives the following.

\begin{cor}\label{cor:multiple-fibre-normality}
If $(\sfC_i)_{i\in I}$ is a countable family of order continuous
idempotent classes and $(\L_i)_{i\in I}$ is a countable family of
subgroups of $\bbZ^d$ then any $\bbZ^d$-system $(X_0,\mu_0,T_0)$
admits an extension $(X,\mu,T)\to (X_0,\mu_0,T_0)$ that is
$\sfC_i$-sated, $(\sfZ_0^{\L_i}\vee\sfC_i)$-sated and such that
$T^{\ \uhr \L_i}$ is fibre-normal over $\sfC_i$ for each $i\in I$.
\qed
\end{cor}

\begin{dfn}[FIS$^+$]
A $\bbZ^d$-system $(X,\mu,T)$ is \textbf{fully isotropy-sated with
fibre-normality} or \textbf{FIS$^{\bs{+}}$} if it is both
$(\sfZ_0^{p_1}\vee \sfZ_0^{p_2}\vee\cdots\vee\sfZ_0^{p_k})$-sated
and $T^{\ \uhr q}$-fibre-normal over $\sfZ_0^{p_1}\vee
\sfZ_0^{p_2}\vee\cdots\vee\sfZ_0^{p_k}$ for every choice of
homomorphisms $p_i:\bbZ^{r_i}\into \bbZ^d$ and $q:\bbZ^s \into
\bbZ^d$.
\end{dfn}

By the properties of isotropy factors established earlier we can
deduce the following strengthening of the existence of the fully
isotropy-sated (FIS) extensions of Definition 3.13
in~\cite{Aus--lindeppleasant1}.

\begin{cor}\label{cor:FIS+}
Any $\bbZ^d$-system admits an FIS$^+$ extension. \qed
\end{cor}

\section{Direct integrals of nilsystems and their inverse limits}\label{sec:nil}

Nilsystems have been an object of
study for ergodic theorists for some time: see, for instance, the
monograph of Auslander, Green and Hahn~\cite{AusGreHah63}, the
foundational papers of Parry~\cite{Par69,Par70} and the more recent
book of Starkov~\cite{Sta00}. In recent years they have come to
occupy a central place in the study of nonconventional averages
associated to powers of a single transformation, where they and
their higher-step analogs are now known to describe precisely the
characteristic factors for linear nonconventional averages (see the
papers of Host and Kra~\cite{HosKra05} and of Ziegler~\cite{Zie07}
and the references listed there). Moreover, pro-nilsystem factors of $\bbZ^d$-actions retain their
r\^ole as precise characteristic factors for nonconventional
averages associated to several commuting transformations, subject to
some additional ergodicity assumptions on various combinations of
these transformations (see Zhang~\cite{Zha96} and Frantzikinakis and
Kra~\cite{FraKra05}).

In view of these results it is not surprising that they re-appear in our Theorem~\ref{thm:char-three-lines-in-2D}.  However, we do now need a simple non-ergodic generalization of the pro-nilsystems studied in those earlier papers. Building on the machinery of extensions by measurably varying compact homogeneous spaces from~\cite{Aus--ergdirint}, in this section we introduce this generalization and establish some elementary properties that will be needed later.

\subsection{Nil-systems, cocycles and nil-selectors}

\textbf{Notation}\quad Given a compact Abelian group $Z$, we will routinely identify the semidirect product group $Z\ltimes \C(Z)$ with the group of all transformations of $Z\times \Sone$ that act as skew-product transformations over some rotation of $Z$ (so $(z,\s) \in Z\ltimes\C(Z)$ is identified with $R_z\ltimes \s$), and equip it with the restriction of the coarse topology on transformations (equivalent to the SOT on bounded linear operators on $L^2(m_{Z\times \Sone})$), under which it is clearly Polish.  In addition, we let $\rm{res}:Z\ltimes \C(Z) \to Z$ be the quotient map, which can be interpreted as restricting transformations to the factor $Z\times \Sone\to Z$. \fin

\begin{dfn}[Nil-cocycle]
If $Z$ is a compact Abelian Lie group (so having finitely many connected components, but possibly more than one), $\G$ a discrete Abelian group, $\phi:\G\to Z$ a homomorphism and $\s:\G\times Z\to \Sone$ a cocycle over $R_\phi$, then $\s$ is a \textbf{nil-cocycle} over $R_\phi$ if there is some transitive two-step nilpotent Lie subgroup $G\subseteq Z\ltimes \C(Z)$ such that
\[G \supseteq \{R_{\phi(\g)}\ltimes \s(\g,\cdot):\ \g \in \G\}.\]

If $Z$ is an arbitrary compact Abelian group and $\phi:\G\to Z$ a homomorphism then a \textbf{nil-cocycle} over $R_\phi$ is the lift to $Z$ of a nil-cocycle over some Lie quotient of $(Z,\phi)$.
\end{dfn}

\begin{dfn}[Nil-selectors]
If $Z$ is Lie as above then a \textbf{nil-selector over $Z$} is a Borel selection $z \mapsto b_z \in \C(Z)$ such that there is some transitive two-step nilpotent Lie group $G \subseteq Z\ltimes \C(Z)$ for which this selection is a cross-section of the quotient map $\rm{res}|_G:G\onto Z$ (which is still surjective because $G$ acts transitively).

If $Z$ is arbitrary then a \textbf{nil-selector over $Z$} is a Borel selection of the form $z\mapsto b_{q(z)}\circ q$ for some Lie quotient $q:Z\to Z_1$ and nil-selector $b_\bullet$ over $Z_1$; borrowing from the terminology of group cohomology we will sometimes refer to this as the \textbf{inflation} of $b$ to $Z$.
\end{dfn}

Note that if $Z$ is a compact Abelian Lie group and $G \subseteq Z\ltimes \C(Z)$ is as above then the space $Z\times\Sone$ is identified with some two-step nilmanifold $G/\G$ such that $[G,G]/([G,G]\cap \G) \cong \Sone$ and $G/[G,G]\G \cong Z$ (see, for instance, Green and Tao~\cite{GreTao--nildist} for a nice introduction to this kind of calculation).  If $R_\phi\ltimes \s$ acts on $Z\ltimes \Sone$ with $\s$ a nil-cocycle, then it is isomorphic to a $\bbZ^2$-action by rotations on such a nilmanifold.  This is the traditional definition of a `nilsystem', but for us it will prove more convenient to proceed via the above definition of a nil-cocycle.

\begin{dfn}[Nil-systems and pro-nilsystems]
For $\G$ a discrete Abelian group, an ergodic $\G$-system $\bfX$ is a \textbf{two-step nilsystem} if it is isomorphic to a two-step Abelian system $(Z\times A,m_{Z\times A},R_\phi\ltimes \s)$ with $Z$ and $A$ compact Abelian Lie groups and such that $\chi(\s)$ is a nil-cocycle over $R_\phi$ for every $\chi \in \hat{A}$.

More generally, $\bfX$ is an \textbf{ergodic pro-nilsystem} if it is an inverse limit of nilsystems.
\end{dfn}

\textbf{Remark}\quad Extending the previous observation, it easily seen that $\bfX$ is a $\bbZ^2$-nilsystem if and only if it is isomorphic to a $\bbZ^2$-action by rotations on a two-step nilmanifold $G/\G$ such that $[G,G]/([G,G]\cap \G) \cong A$ and $G/[G,G]\G \cong Z$. \fin

We will make one crucial appeal (in Subsection~\ref{subs:reduce-to-erg}) to the ability to make a measurable selection from isomorphism classes of nil-cocycles.

\begin{lem}[A Borel selection of canonical nil-cocycles]\label{lem:canon-nil-cocyc-measble}
For a fixed compact Abelian group $Z$ and discrete Abelian group $\G$, there is a Borel subset
\[\cal{A}(Z) \subseteq \Hom(\G,Z)\times \C(\G\times Z)\]
such that
\begin{itemize}
\item $\cal{A}(Z)$ intersects every fibre $\{\phi\}\times \C(\G\times Z)$, and
\item if $\phi \in \Hom(\G,Z)$ then a cocycle $\s:\G\times Z\to\Sone$ over $R_\phi$ is a nil-cocycle if and only if it is cohomologous over $R_\phi$ to some $\s_1$ for which $(\phi,\s_1) \in \cal{A}(Z)$.
\end{itemize}
\end{lem}

\textbf{Remark}\quad Note that we do not assume that $\phi$ has dense image. \fin

\textbf{Proof}\quad Suppose first that $Z$ is a Lie group.  Then there are only countably many possible transitive two-step nilpotent subgroups of $Z\ltimes \C(Z)$ up to fibrewise rotations of the extension $Z\times \Sone\to Z$ (see, for instance, Rudolph~\cite{Rud93} or Host and Kra~\cite{HosKra08} for a classification).  Picking a sequence of representatives $G_1,G_2,\ldots, \subseteq Z\ltimes \C(Z)$ of these isomorphism classes, we see easily that
\begin{itemize}
\item the set $I_\phi \subseteq \bbN$ of $i\geq 1$ such that the homomorphism $\phi:\G\to Z$ admits a lift $\G\to G_i$ varies measurably with $\phi \in \Hom(\G,Z)$;
\item for each $\phi$ we have $|I_\phi|\geq 1$, since one of the groups $G_i$ is simply isomorphic to the product $Z\times \Sone$ for which all lifts are possible.
\end{itemize}

Let $\A_i(Z) := \{(\phi,\s):\ I_\phi \ni i,\,R_\phi \ltimes \s \in G_i\}$ and $\A(Z) := \bigcup_{i\geq 1}\A_i(Z)$.  Each of these sets is clearly measurable and $\A(Z)$ intersects every fibre $\{\phi\}\times \C(\G\times Z)$.  Moreover, if $\s$ is a nil-cocycle over $\phi$ then there is some transitive two-step nilpotent Lie $G' \subseteq Z\ltimes\C(Z)$ containing $R_\phi\ltimes \s$, and now picking a fibrewise rotation of $Z\times \Sone$ that identifies $G'$ with some $G_i$, this correspondingly identifies $R_\phi\ltimes \s$ with $R_\phi\ltimes \s_1$ for some $\s_1 \in \A_i(Z)$, as required.

Finally, if $Z$ is not necessarily a Lie group, then we let $\A(Z)$ be the union of the inflations of the collections $\A(Z_1)$ corresponding to all Lie group quotients $Z\onto Z_1$. \qed

\textbf{Remark}\quad Some analog of the above result should hold in higher ranks, but it is made more complicated because the countability of isomorphism classes of acting nilpotent Lie groups that we have used can fail.  I have not examined this question carefully. \fin

We will later make central use of the following intrinsic characterization of nil-cocycles, which originates in the works~\cite{ConLes84,ConLes88.1,ConLes88.2} of Conze and Lesigne and is examined in depth in Rudolph's paper~\cite{Rud93} (see in particular his Theorem 3.8).

\begin{prop}\label{prop:CL}
Suppose that $\phi:\G\to Z$ is a dense homomorphism and that $\s:\G\times Z\to \Sone$ is a cocycle over $R_\phi$.  Then $\s$ is a nil-cocycle if and only if for every $z \in Z$ the Conze-Lesigne equation
\[\DDelta_z\s(\g,\cdot) = \DDelta_{\phi(\g)}b_z\cdot c_z(\g)\quad\quad \forall \g \in \G\]
has a solution in $b_z \in \C(Z)$ and $c_z \in \Hom(\G,\Sone)$.
\end{prop}

\textbf{Proof}\quad This is essentially as in the case of $\bbZ$-actions, which are well-treated in the above references, so we only sketch the proof here.  The forward implication follows at once by letting $b_z$ be a nil-selector from some transitive two-step nilpotent Lie group $G$ containing $R_\phi\ltimes \s$, for which the Conze-Lesigne equations simply become the assertion that $[G,G]$ consists of constant vertical rotations.

For the backwards implication, first observe that if $(b_z,c_z)$ and $(b'_z,c'_z)$ are two competing solutions of the above equation for some $z$, then by comparing these equations and using the ergodicity of $R_\phi$ we find that $b_z\cdot\ol{b'_z}$ must be an affine function on $Z$ (this argument will re-appear in Lemma~\ref{lem:cocycle-basics} below).

Making a measurable selection $z\mapsto b_z$, it follows that the map $\k:(z,z')\mapsto (b_z\circ R_{z'})\cdot\ol{b_{zz'}}\cdot b_{z'}$ is a $2$-cocycle $Z\times Z\to\E(Z)$.  Since $\E(Z) \cong \Sone\times \hat{Z}$ with a twisted action of $Z$, the continuity results of Theorems~\ref{thm:cohomA} and~\ref{thm:cohomB} show that this cocycle is inflated from some Lie group quotient $Z\onto Z_1$ up to cohomology.  Adjusting the maps $b_z$ themselves by the affine-map-valued cochain that gives this cohomology, and then adjusting $c_z$ accordingly, we find that both $b_z$ and $c_z$ may in fact be taken to have been inflated from $Z_1$, and hence we have reduced to the case in which $Z$ is a Lie group.

Given this assumption, we define the group
\[G := \{R_z\ltimes b_z:\ b_z\ \hbox{solves the Conze-Lesigne eq. at $z$ with some}\ c_z\},\]
so the cocycle equation for $\s$ implies that $R_{\phi(\g)}\ltimes \s(\g,\cdot) \in G$ for every $\g \in \G$.  The assumption that the Conze-Lesigne equations all have solutions gives that $G$ is transitive, and the continuity of those equations in $(z,b_z,c_z)$ together with the discreteness of $\Hom(\G,\Sone)$ imply that $G$ is a closed subgroup of $Z\ltimes \C(Z)$.  It is therefore Polish, and so the fact that it is an extension of $Z$ by $\E(Z)$ implies that it is a Lie group (for example by using the work of Gleason, Montgomery and Zippin on Hilbert's Fifth Problem, because $G$ has no small subgroups~\cite{MonZip55}, although the methods of~\cite{Rud93} yield a more elementary proof).  It remains only to show that it is two-step nilpotent.  This holds because if $b_z$ and $b_{z'}$ solve the Conze-Lesigne equations at $z$ and $z'$ respectively, then differencing the first of these equations by $z'$, the second by $z$ and comparing the results gives
\[\DDelta_{\phi(\g)}(\DDelta_{z'}b_z\cdot\DDelta_z\ol{b_{z'}}) = 1\quad\quad\forall \g \in \G,\]
so by the density of $\phi$ this implies that $\DDelta_{z'}b_z\cdot \DDelta_z\ol{b_{z'}}$ is a constant.  Re-writing this conclusion, we have shown that $[G,G]$ is the set of constant vertical rotations $\id_Z\ltimes (\rm{const.})$, and hence is central in $G$, as required. \qed

\subsection{Nilsystems from local nilsystems}

In the analysis of this paper we will also need a general ability to pass from a system for which the ergodic components of a finite-index subaction are nilsystems to a single, global nilsystem.

\begin{dfn}[Local nil-cocycles and nil-selectors]\label{dfn:local-nil-select}
If $Z$ is a compact Abelian group, $\G$ a discrete Abelian group, $\phi:\G\to Z$ a homomorphism and $Z_0 \leq Z$ a finite-index subgroup, then a cocycle $\s:\G\times Z\to \Sone$ over $R_\phi$ is a \textbf{$Z_0$-local nil-cocycle} if its restriction $\s|_{\L\times zZ_0}:\L\times zZ_0\to \Sone$ is a nil-cocycle for each coset $zZ_0 \in Z/Z_0$, where $\L := \phi^{-1}(Z_0)$.

Similarly, a \textbf{$Z_0$-local nil-selector} is a Borel map $\b:Z_0\to \C(Z)$ such that for each coset $z_0Z_0$ the translated restrictions $z\mapsto \b(z)|_{z_0Z_0}\circ R_{z_0}$ define a nil-selector on $Z_0$.

We will sometimes refer to ordinary nil-cocycles and selectors as \textbf{global} to emphasize that they are not merely local.
\end{dfn}

\begin{prop}\label{prop:loc-nil-to-nil-1}
Suppose that $\G$ is a discrete Abelian group, $\phi:\G\to Z$ is a
dense homomorphism, $Z_0 \leq Z$ and $\L \leq \G$ are finite-index
subgroups for which $\phi(\L)\subseteq Z_0$ and $\tau:\G\times
Z\to\Sone$ is a cocycle over $R_\phi$. If $\tau|_{\L\times
Z_0}:\L\times Z_0\to \Sone$ is a nil-cocycle on $Z_0$, then in fact
$\tau$ is a nil-cocycle on $Z$.
\end{prop}

\textbf{Remark}\quad This amounts to a higher-rank variant of the
result for $\bbZ$-actions that any root of a nil-cocycle is still a
nil-cocycle: see, for instance, Proposition 3.18 in
Meiri~\cite{Mei90}. \fin

\textbf{Proof}\quad This proof breaks naturally into two special
cases.  In both cases we will show that for any $z \in Z$ there are
$b_z \in \C(Z)$ and $c_z \in \Hom(\G,\Sone)$ such that
\[\DDelta_z\tau(\g,\cdot) = \DDelta_{\phi(\g)}b_z\cdot c_z(\g)\quad\quad \forall \g \in \G,\]
from which point Proposition~\ref{prop:CL} completes the proof.

\quad\textbf{Step 1}\quad Suppose first that $\ol{\phi(\L)} = Z_0 =
Z$. In this case we know that for any $z \in Z$ there are $b_z$ and
$c^\circ_z$ for which the above equation holds for every $\g \in \L$, and
will show that for some extension $c_z$ of $c^\circ_z$ to $\G$ it must in fact hold for every $\g \in \G$.

To see this, let $z \in Z$, $\g\in \G$ and $\l \in \L$ and consider
the translation of the Conze-Lesigne equation at $z$ and $\l$ by
$\phi(\g)$:
\[\DDelta_z(\tau(\l,\cdot)\circ R_{\phi(\g)}) = \DDelta_{\phi(\l)}(b_z\circ R_{\phi(\g)})\cdot c^\circ_z(\l).\]
By the cocycle equations for $\tau$, the left-hand side here equals
$\DDelta_z\big(\tau(\l,\cdot)\cdot
\DDelta_{\phi(\l)}\tau(\g,\cdot)\big)$, and substituting this
expression and re-arranging gives
\[\DDelta_z\tau(\l,\cdot) = \DDelta_{\phi(\l)}\big((b_z\circ R_{\phi(\g)})\cdot \DDelta_z\ol{\tau(\g,\cdot)}\big)\cdot c^\circ_z(\l).\]
It follows that the function $b'_z := (b_z\circ R_{\phi(\g)})\cdot
\DDelta_z\ol{\tau(\g,\cdot)}$ is a competing solution of the
Conze-Lesigne equation at $z$ and $\l$ with the same $c^\circ_z$, and so dividing these two versions of this equation gives
\[\DDelta_{\phi(\l)}\big(\DDelta_{\phi(\g)}b_z\cdot \DDelta_z\ol{\tau(\g,\cdot)}\big) = 1.\]
Sine this holds for all $\l\in\L$ and we have assumed that
$\ol{\phi(\L)} = Z$, it follows that $\DDelta_{\phi(\g)}b_z\cdot
\DDelta_z\ol{\tau(\g,\cdot)}$ is equal to a constant, say
$\ol{c_z(\g)}$, and now we can check directly that this must define
a homomorphism $c_z:\G\to \Sone$ extending $c_z^\circ$.

\quad\textbf{Step 2}\quad Now suppose that $Z_0$ is an arbitrary
finite-index subgroup of $Z$.  By initially shrinking $Z_0$ further
if necessary we may assume that $Z_0 = \ol{\phi(\L)}$. Now applying
Step 1 to the cocycle $\tau|_{\phi^{-1}(Z_0)\times Z_0}$ and the
inclusion $\L\leq \phi^{-1}(Z_0)$, we already know that
$\tau|_{\phi^{-1}(Z_0)\times Z_0}$ is a nil-cocycle, and hence we
may now assume that $\L = \phi^{-1}(Z_0)$.

We will next construct the global Conze-Lesigne solutions $b_z$ for
$z \in Z_0$. Let $\O \subset \G$ be a fundamental domain for $\L$
and $\lfloor\cdot\rfloor$ and $\{\cdot\}$ the corresponding integer-
and fractional-part maps. By assumption, we already have some
$b^\circ_z\in \C(Z_0)$ and $c^\circ_z \in \Hom(\L,\Sone)$ such that
\begin{eqnarray}\label{eq:CLagain}
\DDelta_z\tau(\g,w) = \DDelta_{\phi(\g)}b^\circ_z(w)\cdot c^\circ_z(\g)\quad\quad \forall \g \in \L,\ w \in Z_0.
\end{eqnarray}
Moreover, the cocycle equation for $\tau$ gives that
\[\tau(\g,\phi(\k)w) = \tau(\g,w)\cdot\DDelta_{\phi(\g)}\tau(\k,w)\]
for any $\k \in \O$, $z \in Z_0$, and so extending the definition of $b_z^\circ$ by setting $b^\circ_z(\phi(\k)w) := b^\circ_z(w)\cdot \DDelta_z\tau(\k,w)$ shows that we may assume that~(\ref{eq:CLagain}) holds for all $w \in Z$.

Let us also suppose that we have arbitrarily extended each $c^\circ_z$ to a homomorphism $c_z:\G\to \Sone$. Having done this, if $\g = \lfloor \g\rfloor + \{\g\}$ is an arbitrary element of $\G$ then the cocycle equations for $\tau$ and the extended definition of $b^\circ_z$ give
\begin{eqnarray*}
&&\DDelta_z\tau(\g,\phi(\k)w)\\
&&= \DDelta_z \tau(\lfloor \g + \k\rfloor,\phi(\k)w)\cdot \DDelta_z\tau(\{\g + \k\} - \k,\phi(\lfloor \g + \k\rfloor + \k)w)\\
&&= \DDelta_{\phi(\lfloor\g + \k\rfloor)}b^\circ_z(\phi(k)w)\cdot c_z(\lfloor\g+\k\rfloor)\cdot \DDelta_z\tau(\{\g + \k\} - \k,\phi(\lfloor \g + \k\rfloor + \k)w)\\
&&= \big(b^\circ_z(\phi(k)\phi(\lfloor\g + \k\rfloor)w)\cdot \DDelta_z\tau(\{\g + \k\} - \k,\phi(\k)\phi(\lfloor \g + \k\rfloor)w)\big)\\
&&\quad\quad\quad\quad\quad\quad\quad\quad\quad\quad\quad\quad\cdot\ol{b^\circ_z(\phi(\k)w)}\cdot c_z(\lfloor\g+\k\rfloor)\\
&&= \big(b^\circ_z(\phi(\lfloor\g + \k\rfloor)w)\\
&& \quad\quad\quad\quad\quad\quad\cdot \DDelta_z(\tau(\k,\phi(\lfloor\g + \k\rfloor)w)\cdot\tau(\{\g + \k\} - \k,\phi(\k)\phi(\lfloor \g + \k\rfloor)w))\big)\\
&&\quad\quad\quad\quad\quad\quad\quad\quad\quad\quad\quad\quad\cdot\ol{b^\circ_z(\phi(\k)w)}\cdot c_z(\lfloor\g+\k\rfloor)\\
&&= \big(b^\circ_z(\phi(\lfloor\g + \k\rfloor)w)\cdot \DDelta_z(\tau(\{\g +\k\},\phi(\lfloor\g + \k\rfloor)w)\big)\\
&&\quad\quad\quad\quad\quad\quad\quad\quad\quad\quad\quad\quad\cdot\ol{b^\circ_z(\phi(\k)w)}\cdot c_z(\lfloor\g+\k\rfloor)\\
&&= \big(b^\circ_z(\phi(\g)\phi(\k)w)\cdot\ol{b^\circ_z(\phi(\k)w)}\cdot c_z(\lfloor\g+\k\rfloor)\\
&&= \DDelta_{\phi(\g)}b^\circ_z(\phi(\k)w)\cdot c_z(\g)\cdot c_z(\lfloor \g + \k\rfloor - \g).
\end{eqnarray*}
Since $c_z(\lfloor \g + \k\rfloor - \g) = c_z(\{\g + \k\} - \k) = c_z(\{\g + \k\})\ol{c_z(\k)}$, if we now define $b_z(\phi(\k)u) := b^\circ_z(\phi(\k)u)\cdot c_z(\k)$ for $k \in \O$ and $u \in Z_0$ then the above calculation asserts that
\[\DDelta_z\tau(\g,w) = \DDelta_{\phi(\g)}b_z(w)\cdot c_z(\g)\quad\quad\forall \g\in \G,\,w\in Z,\]
as required.

To finish the proof we need only extend the definition of $b_\bullet$ to the whole of $Z$, and this can be done using one last appeal to the cocycle equations: if $\l\in \O$ and $z \in Z_0$ then we set $b_{\phi(\l)z}(w) := \tau(\l,zw)\cdot b_z(w)$ for $w\in Z$ and compute that
\begin{eqnarray*}
\DDelta_{\phi(\l)z}\tau(\g,w) &=& \DDelta_{\phi(\l)}\tau(\g,zw)\cdot \DDelta_z\tau(\g,w)\\
&=& \DDelta_{\phi(\g)}\tau(\l,zw)\cdot \DDelta_{\phi(\g)}b_z(w)\cdot c_z(\g)\\
&=& \DDelta_{\phi(\g)}b_{\phi(\l)z}(w)\cdot c_z(\g).
\end{eqnarray*}
Thus we have obtained global solutions to all the desired Conze-Lesigne equations, and hence verified the conditions of Proposition~\ref{prop:CL}. \qed

Specializing to $\bbZ^2$-actions, we next show how the above proposition can be used in conjunction with the ability to create new nil-cocycles upon extending a base group rotation over which we already have some nil-structure.

\begin{lem}\label{lem:making-nilrtns-commute}
If $\phi:\bbZ^2 \to Z$ is a dense homomorphism, $G\subseteq Z\ltimes \C(Z)$ is a transitive two-step nilpotent Lie group and $\s:Z\to \Sone$ is a Borel map such that $R_{\phi(\bf{e}_1)}\ltimes \s \in G$, then there are some product system $(\t{Z},\t{\phi}) := (Z\times Z',(\phi,\phi'))$ with first coordinate projection $q$ onto $(Z,\phi)$, a transformation $R_{\t{\phi}(\bf{e}_2)}\ltimes \s_2$ that commutes with $R_{\t{\phi}(\bf{e}_1)}\ltimes (\s\circ q)$ and a transitive two-step nilpotent Lie group $\t{G} \subseteq \t{Z}\ltimes \C(\t{Z})$ containing both of these transformations.
\end{lem}

\textbf{Proof}\quad Start by choosing any $R_{\phi(\bf{e}_2)}\ltimes \tau \in G$ that extends the rotation $R_{\phi(\bf{e}_2)}$ (this is possible by the transitivity of $G$).  Now the commutator
\[[R_{\phi(\bf{e}_1)}\ltimes \s,R_{\phi(\bf{e}_2)}\ltimes \tau] = \id_Z\ltimes (\DDelta_{\phi(\bf{e}_1)}\tau\cdot\ol{\DDelta_{\phi(\bf{e}_2)}\s})\]
is simply a constant vertical rotation, say by $\theta\in \Sone$.  Multiplying $(Z,\phi)$ by a rotation $(Z',\phi')$ on a subgroup of $\Sone$ for which $\theta$ is an $R_{\phi'(\bf{e}_2)}$-eigenvalue, we obtain a product extension of $(Z,\phi)$ through the coordinate projection $q$ so that $\theta$ is an $R_{\t{\phi}(\bf{e}_2)}$-eigenvalue, say with eigenvector $\chi \in \E(\t{Z})$. Thence by setting $\s_2 := \chi\cdot(\tau\circ q)$ we obtain a cocycle over $R_{\t{\phi}(\bf{e}_2)}$ that commutes with $R_{\t{\phi}(\bf{e}_1)}\ltimes(\s\circ q)$.

Finally, setting
\[\t{G} := \big\langle\{R_{(z,z')}\ltimes (\s'\circ q):\ z' \in Z',\ R_z\ltimes \s' \in G\}\cup \{R_{\t{\phi}(\bf{e}_2)}\ltimes \s_2\}\big\rangle,\]
this is clearly transitive and contains both $R_{\t{\phi}(\bf{e}_1)}\ltimes \s$ and $R_{\t{\phi}(\bf{e}_2)}\ltimes \s_2$, and a quick calculation shows that any two of its elements have commutator of the form $\id_{\t{Z}}\ltimes (\rm{const}.)$, hence central in $\t{G}$. \qed

\begin{prop}\label{prop:loc-nil-to-nil-2}
Suppose that $\phi:\bbZ^2\to Z$ is a homomorphism such that $\ol{\phi(\bbZ^2)}\leq Z$ has finite index, that $Z_0 \leq \ol{\phi(\bbZ^2)}$ is a further finite-index subgroup, that $\tau_0:\bbZ^2\times Z\to \Sone$ is a cocycle over $R_\phi$ and that $\bf{n}\in \bbZ^2\setminus \{\bs{0}\}$ is such that $\phi(\bf{n}) \in Z_0$ and $\tau_0(\bf{n},\cdot)$ is a $Z_0$-local nil-cocycle over $R_{\phi(\bf{n})}$. Then there are an extension $q:(\t{Z},\t{\phi})\to (Z,\phi)$ and a cocycle $\tau:\bbZ^2\times \t{Z}\to \Sone$ over $R_{\t{\phi}}$ such that $\tau_0(\bf{n},q(\cdot)) = \tau(\bf{n},\cdot)$ and $\tau$ is a $\ol{\t{\phi}(\bbZ^2)}$-local nil-cocycle.

Moreover, if $\phi$ has dense image then $\t{\phi}$ may also be chosen with dense image.
\end{prop}

\textbf{Remark}\quad In case $\phi$ has dense image we simply obtain a global nil-cocycle $\tau$ over $R_{\t{\phi}}$, since $\ol{\t{\phi}(\bbZ^2)} = \t{Z}$. \fin

\textbf{Proof}\quad It suffices to find
\begin{itemize}
\item an extension $q:(\t{Z},\t{\phi})\to (Z,\phi)$, ergodic if $R_\phi$ is ergodic,
\item a finite-index subgroup $\t{Z}_0\leq \ol{\t{\phi}(\bbZ^2)}$ such that $q(\t{Z}_0) = Z_0$ and $\t{\phi}(\bf{n}) \in \t{Z}_0$,
\item and a cocycle $\tau:\bbZ^2\times \t{Z}\to \Sone$ over $R_{\t{\phi}}$
\end{itemize}
such that the restriction $\tau|_{\L\times \t{Z}_0}$ is a nil-cocycle on $\t{Z}_0$ over $R_{\t{\phi}|_\L}$ for some finite-index $\L\leq \t{\phi}^{-1}(\t{Z}_0)$, and $\tau_0(\bf{n},q(\cdot)) = \tau(\bf{n},\cdot)$: given this, an application of Proposition~\ref{prop:loc-nil-to-nil-1} on each of the cosets of $\ol{\t{\phi}(\bbZ^2)}$ completes the proof.

First suppose that $\bf{n}' \in \phi^{-1}(Z_0)$ is linearly independent from $\bf{n}$, set $\L := \bbZ\bf{n} + \bbZ\bf{n}'$, and let $\O \subset \bbZ^2$ be a fundamental domain for $\G$.  By making a first extension of $(Z,\phi)$ by multiplying by the finite group rotation $\bbZ^2 \actson \bbZ^2/\L$ and replacing $Z_0$ with $Z_0\times \{\L\}$ if necessary, we may assume that $\L = \phi^{-1}(Z_0)$.

Now apply Lemma~\ref{lem:making-nilrtns-commute} to the $\L$-action $(Z_0,\phi|_\L)$ and cocycle $\tau_0(\bf{n},\cdot)|_{Z_0}$, which we know is lifted from some Lie group factor of $Z_0$ up to cohomology.  This gives some new rotation action $(Z',\phi')$ of $\L$ such that if $(\t{Z}_0,\t{\phi}_0):= (Z_0\times Z',(\phi|_\L,\phi'))$ and $q_0:\t{Z}_0 \to Z_0$ is the coordinate projection then there is a transformation $R_{\t{\phi}_0(\bf{n}')}\ltimes \s_2$ that commutes with $R_{\t{\phi}_0(\bf{n})}\ltimes \tau_0(\bf{n},q(\cdot))$ and such that $\s_2$ and $\tau_0(\bf{n},q(\cdot))$ together define a nil-cocycle $\tau':\L\times \t{Z}_0\to \Sone$.

Next, we may easily adjoin roots to $Z'$ in order to assume that $\phi':\L\to Z'$ is the restriction of some homomorphism $\bbZ^2\to Z'$, which we now also denote by $\phi'$.  Having done this, let $(\t{Z},\t{\phi}) := (Z\times Z',(\phi,\phi'))$ and $q:(\t{Z},\t{\phi})\to (Z,\phi)$ be the coordinate projection.  This locates $\t{Z}_0 \leq \t{Z}$ as a finite-index subgroup and realizes $\tau'$ as a nil-cocycle $\L\times \t{Z}_0\to \Sone$ over $R_{\t{\phi}|_\L}$ with $\tau'(\bf{n},\cdot) = \tau_0(\bf{n},q(\cdot))|_{\t{Z}_0}$.  Clearly we still have $\L = \t{\phi}^{-1}(\t{Z}_0)$.  In case $R_\phi$ was ergodic, we may also restrict to an $R_{\t{\phi}}$-ergodic component of this extension without losing any of these properties.

Finally we can extend $\tau'$ to a cocycle $\tau:\bbZ^2\times \t{Z}\to \Sone$ over $R_{\t{\phi}}$ in the only way permitted up to cohomology by the demands of the cocycle condition: that is, for $\bf{q} \in \bbZ^2$ and $z \in \t{Z}_0$ we set
\[\tau(\bf{q},z) := \tau'(\lfloor \bf{q}\rfloor,z)\cdot\tau_0(\{\bf{q}\},\phi(\lfloor\bf{q}\rfloor)q(z)),\]
and now for $\bf{k} \in \O$ we set
\[\tau(\bf{q},\t{\phi}(\bf{k})z) := \tau(\bf{q},z)\cdot \DDelta_{\t{\phi}(\bf{q})}\tau_0(\bf{k},q(z)).\]

This is easily checked to be a cocycle over $R_{\t{\phi}}$, and it manifestly agrees with $\tau'(\bf{q},z)$ when $\bf{q} \in \G$, since in that case $\lfloor \bf{q} + \bf{k}\rfloor = \bf{q}$. Since it satisfies the conditions of Proposition~\ref{prop:loc-nil-to-nil-1}, it is actually a nil-cocycle over $R_{\t{\phi}}$, so the proof is complete. \qed

\subsection{Direct integrals and inverse limits of nilsystems}

The central r\^ole played by pro-nilsystems in~\cite{HosKra05} will be reprised here, but only after we set up a suitable formalism for `direct integrals' of such systems.  Here we do this by building on the theory of extensions of systems by measurably-varying compact homogeneous space data in~\cite{Aus--ergdirint}.  A rather different approach to handling such systems has recently been used by Chu, Frantzikinakis and Host in~\cite{ChuFraHos09}, but their formalism seems to lend itself less readily to the kinds of detailed structural analysis we will need to perform later.

\begin{dfn}[Direct integrals of nilsystems and pro-nilsystems]
For a discrete Abelian group $\G$, a \textbf{direct integral of two-step $\G$-nilsystems} is a $\G$-system of the form
\[(\id_S \ltimes R_{\phi_\star})\ltimes \s:\G \actson ((S\ltimes Z_\star)\ltimes A_\star,\nu\ltimes m_{Z_\star\times A_\star})\]
for some invariant base space $(S,\nu)$, motionless compact Abelian Lie group data $Z_\star$ and $A_\star$, a measurable family of homomorphisms $\phi_s:\G\to Z_s$ and a cocycle-section $\s:\G\times Z_\star \to A_\star$, such that in addition the system $(Z_s\times A_s,m_{Z_s\times A_s},R_{\phi_s}\ltimes \s_s)$ is an ergodic two-step nilsystem for $\nu$-almost every $s \in S$.

Slightly abusively, a \textbf{two-step $\G$-pro-nilsystem} is an inverse limit of an inverse sequence of direct integrals of two-step $\G$-nilsystems.  We write $\sfZ_{\nil,2}^\G$ for the class of all such systems.
\end{dfn}

\textbf{Remark}\quad It is easy to see how this definition could be extended to higher-step nilsystems, but we do not make use of this in the present paper. \fin

\begin{lem}[Pro-nilsystems form an idempotent class]\label{lem:pronil-idemp}
The class $\sfZ_{\nil,2}^{\bbZ^2}$ of inverse limits of direct integrals of $\bbZ^2$-nilrotations forms an idempotent class (in the sense of Section 3 of~\cite{Aus--lindeppleasant1}).
\end{lem}

\textbf{Proof}\quad We must show closure of this class under joinings and inverse limits.  Since the class is defined by taking the completion under inverse limits of the class of direct integrals of $\bbZ^2$-nilsystems and any inverse limit of inverse limits can be written as a single inverse limit by a simple diagonal argument, the second of these properties is immediate.

To prove closure under joinings, observe first that if $\bfX_1,\bfX_2 \in \sfZ_{\nil,2}^{\bbZ^2}$ are generated by the sequences of factors $\pi^{(i)}_n:\bfX_i\to \bfX_{i,n}$ each of whose targets is a direct integral of nilsystems, then any joining of $\bfX_1$ and $\bfX_2$ is generated by the induced joinings of the pairs of factors $\bfX_{1,n}$ and $\bfX_{2,n}$, so it suffices to show that these are still direct integrals of nilsystems.  Therefore we may suppose that $\bfX_1$, $\bfX_2$ are themselves direct integrals of nilsystems, say with invariant base spaces $(S_i,\nu_i)$ for $i=1,2$.

If now $\l$ is a joining of these two systems on the product space $X_1\times X_2$, then it induces a $(\nu_1,\nu_2)$-coupling $\nu$ on $S := S_1\times S_2$, and each fibre measure $\l_s$ of its disintegration over $S$ is a joining of two ergodic $\bbZ^2$-nilsystems, say on the two nilmanifolds $G_{i,s_i}/\G_{i,s_i}$ for $i=1,2$.  As is standard (see, for instance, Leibman~\cite{Lei06}), this implies that $\l_s$ is a direct integral of the Haar measures on some family of sub-nilmanifolds of $G_{1,s_1}\times G_{2,s_2}/(\G_{1,s_1}\times \G_{2,s_2})$ invariant under the product nilrotations. Therefore by disintegrating each component $\l_s$ further, the action $T_1\times T_2$ on each of the fibres $\l_s$ may be expressed as a direct integral of some joinings that are themselves ergodic two-step $\bbZ^2$-nilsystems, and now the observation that the measurability of $\l_s$ with $s$ implies the measurability of the corresponding cocycle-sections completes the proof. \qed

\section{Characteristic factors for three directions in general position}\label{sec:lindep-1}

We henceforth assume the basic theory of Furstenberg self-joinings
and characteristic tuples of factors, referring where necessary to
the results of~\cite{Aus--lindeppleasant1} (particularly Theorem 1.1
and the results of Subsection 4.1 of that paper).

Now we will focus on the averages
\[S_N(f_1,f_2,f_3) = \frac{1}{N}\sum_{n=1}^N(f_1\circ T^{n\bf{p}_1})(f_2\circ T^{n\bf{p}_2})(f_3\circ T^{n\bf{p}_3})\]
for a $\bbZ^2$-action $T$ and three directions $\bf{p}_1$,
$\bf{p}_2$ and $\bf{p}_3$ in $\bbZ^2$ that are in general position
with $\bs{0}$. We will first show that any FIS$^+$ extension already has
characteristic factors with a structure we can describe quite
precisely, and will then turn to a finer analysis (and pass to some
further extensions) to obtain the more explicit picture of Theorem~\ref{thm:char-three-lines-in-2D}.  We will work with FIS$^+$
extensions rather than just FIS extensions for the sake of an
important application of fibre-normality in the proof of
Lemma~\ref{lem:2D-proj-full}.

\subsection{Overview and first
results}\label{subs:first-main-results}

The first tool at our disposal is the fact that FIS extensions are
pleasant for linearly independent tuples of directions (Proposition
4.5 of~\cite{Aus--lindeppleasant1}). This guarantees that after
ascending to an FIS$^+$ (and so certainly FIS) extension, our system
is at least pleasant and isotropized for any two of our $\bf{p}_i$.  To proceed further, we will need to understand the structure of the
Furstenberg self-joining
$\mu^{\rm{F}}_{T^{\bf{p}_1},T^{\bf{p}_2},T^{\bf{p}_3}}$ in much
greater detail.  Let us now agree to abbreviate this particular
Furstenberg self-joining to $\mu^{\rm{F}}$.

Our next steps are still quite routine. A standard re-arrangement
(see Section 4.1 of~\cite{Aus--lindeppleasant1}) gives
\begin{multline*}
\int_{X^3}f_1\otimes f_2\otimes f_3\,\d\mu^{\rm{F}}\\
= \int_Xf_i\cdot\Big(\lim_{N\to\infty}\frac{1}{N}\sum_{n=1}^N
(f_j\circ T^{n(\bf{p}_j - \bf{p}_i)})\cdot(f_k\circ T^{n(\bf{p}_k -
\bf{p}_i)})\Big)\,\d\mu
\end{multline*}
for any permutation $(i,j,k)$ of $(1,2,3)$.  It follows that
\[\int_{X^3}\bigotimes_{i=1}^3 f_i\,\d\mu^{\rm{F}} = \int_{X^3}\bigotimes_{i=1}^3\sfE_\mu(f_i\,|\,\b_i)\,\d\mu^{\rm{F}}\]
where for each $i=1,2,3$ and $\{j,k\} = \{1,2,3\}\setminus \{i\}$,
$\b_i:\bfX \to \bfV_i$ is the factor generated by all the double
nonconventional averages
\[\lim_{N\to\infty}\frac{1}{N}\sum_{n=1}^N
(f_j\circ T^{n(\bf{p}_j - \bf{p}_i)})\cdot(f_k\circ T^{n(\bf{p}_k -
\bf{p}_i)}).\] (Naturally, these re-arrangement games have analogs
for any linear nonconventional averages.)  It follows at
once that
\[\int_{X^3}\bigotimes_{i=1}^3 f_i\,\d\mu^{\rm{F}} = \int_{X^3}\bigotimes_{i=1}^3\sfE_\mu(f_i\,|\,\b'_i)\,\d\mu^{\rm{F}}\]
whenever $\b'_i \succsim \b_i$ for $i=1,2,3$.

Now, each of these double averages corresponds to a pair of linearly
independent directions (because $\bf{p}_j - \bf{p}_i$, $\bf{p}_k -
\bf{p}_i$ are linearly independent, by our assumption of general
position), and so falls within the scope of the Pleasant Extension
Theorem for linearly independent double averages (Theorem 1.1
of~\cite{Aus--lindeppleasant1}). This tells us that for any FIS$^+$
system the characteristic factors are simply composed of the
relevant isotropy factors, so that the above limit is equal to
\begin{multline*}
\lim_{N\to\infty}\frac{1}{N}\sum_{n=1}^N
(\sfE_\mu(f_j\,|\,\zeta_0^{T^{\bf{p}_j} =
T^{\bf{p}_i}}\vee\zeta_0^{T^{\bf{p}_j} = T^{\bf{p}_k}}) \circ
T^{n(\bf{p}_j - \bf{p}_i)})\\
\cdot(\sfE_\mu(f_k\,|\,\zeta_0^{T^{\bf{p}_k} =
T^{\bf{p}_i}}\vee\zeta_0^{T^{\bf{p}_j} = T^{\bf{p}_k}})\circ
T^{n(\bf{p}_k - \bf{p}_i)}).
\end{multline*}
Now, if $g_{ij}$ is bounded and
$\zeta_0^{T^{\bf{p}_i}=T^{\bf{p}_j}}$-measurable, $g_{jk}$ and
$h_{jk}$ are both bounded and
$\zeta_0^{T^{\bf{p}_j}=T^{\bf{p}_k}}$-measurable, and $h_{ki}$ is
bounded and $\zeta_0^{T^{\bf{p}_i}=T^{\bf{p}_k}}$-measurable, then
by re-arranging and applying the classical mean ergodic theorem we
find that
\begin{eqnarray*}
&&\lim_{N\to\infty}\frac{1}{N}\sum_{n=1}^N ((g_{ij}\cdot g_{jk})
\circ T^{n(\bf{p}_j - \bf{p}_i)}) \cdot((h_{jk}\cdot h_{ik})\circ
T^{n(\bf{p}_k - \bf{p}_i)})\\
&&= g_{ij}\cdot h_{ik}\cdot\lim_{N\to\infty}\frac{1}{N}\sum_{n=1}^N
(g_{jk}\cdot h_{jk}) \circ T^{n(\bf{p}_j - \bf{p}_i)}\\ &&=
g_{ij}\cdot h_{ik}\cdot\sfE_\mu(g_{jk}\cdot
h_{jk}\,|\,\zeta_0^{T^{\bf{p}_i} = T^{\bf{p}_j} = T^{\bf{p}_k}}),
\end{eqnarray*}
which is manifestly $(\zeta_0^{T^{\bf{p}_i} = T^{\bf{p}_j}}\vee
\zeta_0^{T^{\bf{p}_i} = T^{\bf{p}_k}})$-measurable.  By linearity
and continuity, it follows that the same is true of the above double
nonconventional averages for any $f_j$ and $f_k$, so we deduce that
$\b_i \lesssim \zeta_0^{T^{\bf{p}_i} = T^{\bf{p}_j}}\vee
\zeta_0^{T^{\bf{p}_i} = T^{\bf{p}_k}}$.  On the other hand, by
making a free choice of $g_{ij}$ and $h_{ik}$ in the above
calculation the reverse containment is also clear, hence $\b_i
\simeq \zeta_0^{T^{\bf{p}_i} = T^{\bf{p}_j}}\vee
\zeta_0^{T^{\bf{p}_i} = T^{\bf{p}_k}}$, and in the future we can
simply take $\b_i$ to equal this joining of isotropy factors.

In summary we have proved the following.

\begin{lem}\label{lem:subchar-factors}
If $\bfX$ is FIS$^+$ then under $\mu^{\rm{F}}$ the three coordinate
projections $\pi_i:X^3 \to X$, $i=1,2,3$, are relatively independent
over their further factors $\b_i\circ\pi_i$, where $\b_i :=
\zeta_0^{T^{\bf{p}_i} = T^{\bf{p}_j}}\vee \zeta_0^{T^{\bf{p}_i} =
T^{\bf{p}_k}}$, and these $\b_i$ comprise the unique minimal triple
of factors with this property. \qed
\end{lem}

\begin{dfn}[Subcharacteristic factors]
We will henceforth refer to $\b_i$ as the $i^{\rm{th}}$
\textbf{subcharacteristic factor} corresponding to the triple of
directions $\bf{p}_1$, $\bf{p}_2$ and $\bf{p}_3$.
\end{dfn}

Now we recall the basic criterion for characteristicity in terms of
$\mu^{\rm{F}}$, given, for example, as Corollary 4.2
in~\cite{Aus--lindeppleasant1}. This tells us that a triple of
factors $\xi_1$, $\xi_2$, $\xi_3$ of $\bfX$ is characteristic if for
any $f_1,f_2,f_3 \in L^\infty(\mu)$ and $\vec{T}$-invariant $g\in
L^\infty(\mu^{\rm{F}})$ we have
\[\int_{X^3}\prod_{i=1}^3(f_i\circ\pi_i)\cdot g\,\d\mu^{\rm{F}} = \int_{X^3}\prod_{i=1}^3(\sfE_\mu(f_i\,|\,\xi_i)\circ\pi_i)\cdot g\,\d\mu^{\rm{F}}.\]
Clearly this assertion is stronger than the relative independence of
$\pi_i$ that characterizes the $\b_i$, so it requires that $\xi_i
\succsim \b_i$.  In addition, since any $g\in
L^\infty(\mu^{\rm{F}})$ can be $L^2$-approximated by finite sums of
tensor products of functions in $L^\infty(\mu)$, this property also
requires that any $\vec{T}$-invariant function on $X^3$ be almost
surely measurable with respect to $\xi_1\times\xi_2\times\xi_3$.  It
turns out that these two demands on $\xi_1$, $\xi_2$, $\xi_3$ are
also sufficient for characteristicity.

\begin{lem}\label{lem:from-subchar-to-char}
A triple of factors $\xi_1$, $\xi_2$, $\xi_3$ of an FIS$^+$
$\bbZ^2$-system is characteristic if and only if
\begin{itemize}
\item $\xi_i \succsim \b_i$ for $i=1,2,3$, and
\item any $\vec{T}$-invariant function on $X^3$ is
$\mu^{\rm{F}}$-almost surely $(\xi_1\times \xi_2\times
\xi_3)$-measurable.
\end{itemize}
\end{lem}

\textbf{Proof}\quad Let $f_1$, $f_2$, $f_3$ and $g$ be as above.
Then $g$ is $(\xi_1\times \xi_2\times \xi_3)$-measurable, so we may
approximate it in $L^2$ by a finite sum $\sum_p g_{1,p}\otimes
g_{2,p}\otimes g_{3,p}$ with each $g_{i,p}$ being bounded and
$\xi_i$-measurable. For these functions we have
\begin{multline*}
\int_{X^3}\prod_{i=1}^3(f_i\circ\pi_i)\cdot \Big(\sum_p
\prod_{i=1}^3(g_{i,p}\circ\pi_i)\Big)\,\d\mu^{\rm{F}} = \sum_p
\int_{X^3}\bigotimes_{i=1}^3 \sfE_\mu(f_i\cdot
g_{i,p}\,|\,\xi_i)\,\d\mu^{\rm{F}}\\ = \sum_p
\int_{X^3}\bigotimes_{i=1}^3 \sfE_\mu(f_i\,|\,\xi_i)\cdot
g_{i,p}\,\d\mu^{\rm{F}} =
\int_{X^3}\prod_{i=1}^3(\sfE_\mu(f_i\,|\,\xi_i)\circ\pi_i)\cdot
\Big(\sum_p \prod_{i=1}^3(g_{i,p}\circ\pi_i)\Big)\,\d\mu^{\rm{F}},
\end{multline*}
first because $\xi_i \succsim \b_i$ and then because each $g_{i,p}$
is $\xi_i$-measurable. By continuity this yields
\[\int_{X^3}\prod_{i=1}^3(f_i\circ\pi_i)\cdot
g\,\d\mu^{\rm{F}} =
\int_{X^3}\prod_{i=1}^3(\sfE_\mu(f_i\,|\,\xi_i)\circ\pi_i)\cdot
g\,\d\mu^{\rm{F}}\] as required. \qed

Lemmas~\ref{lem:subchar-factors} and~\ref{lem:from-subchar-to-char}
now put us into a position to apply the non-ergodic
Furstenberg-Zimmer Inverse Theorem~\ref{thm:rel-ind-joinings}, since
we need to control the $\vec{T}$-invariant factor of the joining
$\mu^{\rm{F}}$ of three copies of $\bfX$, and these three copies are
relatively independent over the subcharacteristic factors $\b_1$,
$\b_2$, $\b_3$. First, however, recall that that theory applies to a
joining of systems that is relatively independent over a collection
of \emph{relatively ergodic} factors of each system.  In the current
setting the coordinate projection $\pi_i:X^3 \to X$ intertwines
$\vec{T}$ with $T^{\bf{p}_i}$, and in general the factor
$\b_i:\bfX\to\bfV_i$ need not be relatively ergodic for the
transformation $T^{\bf{p}_i}$.  We therefore first extend each
$\b_i$ further to
\[\a_i := \zeta_0^{T^{\bf{p}_i}}\vee \b_i = \zeta_0^{T^{\bf{p}_i}}\vee \zeta_0^{T^{\bf{p}_i} = T^{\bf{p}_j}}\vee \zeta_0^{T^{\bf{p}_i} = T^{\bf{p}_k}},\]
and can now apply the Furstenberg-Zimmer Theory to the relatively
independent joining $\mu^{\rm{F}}$ of three copies of $\bfX$ over
the three factors $\a_i$, each of which is $T^{\bf{p}_i}$-relatively
ergodic.  Let us write $\bfW_i$ for some $\bbZ^2$-system that we
take for the target of $\a_i$, so $\bfW_i$ extends $\bfV_i$ through
$\b_i|_{\a_i}$.

We can easily check that we have lost no generality at this step, in
that any triple of characteristic factors satisfies $\xi_i \succsim
\a_i$. Indeed, for each $i=1,2,3$ and any $\a_i$-measurable function
$g_i \in L^\infty(\mu)$, the lifted function $g_i\circ\pi_i$ is
$\vec{T}$-invariant and so by Lemma~\ref{lem:from-subchar-to-char}
is necessarily measurable with respect to $\xi_1\times
\xi_2\times\xi_3$; this clearly requires that $\xi_i \succsim \a_i$.

\begin{dfn}[Proto-characteristic factors]
We will henceforth refer to $\a_i =
\zeta_0^{T^{\bf{p}_i}}\vee\b_i:\bfX\to\bfW_i$ as the $i^{\rm{th}}$
\textbf{proto-characteristic factor} corresponding to the triple of
directions $\bf{p}_1$, $\bf{p}_2$ and $\bf{p}_3$.
\end{dfn}

\textbf{Remark}\quad In case $\bf{p}_1$, $\bf{p}_2$, $\bf{p}_3$ are
three linearly independent directions in some $\bbZ^d$, $d \geq 3$,
the main results of~\cite{Aus--nonconv} tell us that for a suitable
extension (such as an FIS extension of our $\bbZ^d$-system) the triple
$\a_1$, $\a_2$, $\a_3$ is actually characteristic.  The above
discussion shows that these factors are at least obvious lower
bounds for the actual characteristic factors $\xi_1$, $\xi_2$,
$\xi_3$, and that the remaining gap between $\xi_i$ and $\a_i$
(after we ascend to as well-behaved an extension as we can build)
must be accounted for by some essential `interaction' between the
transformations $T^{\bf{p}_1}$, $T^{\bf{p}_2}$ and $T^{\bf{p}_3}$
that cannot be removed by extending further without disrupting the
linear dependence relations among $\bf{p}_1$, $\bf{p}_2$ and
$\bf{p}_3$. To be a little imprecise, it is this defect that is
accounted for by the extra ingredient of the pro-nilsystem that appears in Theorem~\ref{thm:char-three-lines-in-2D}.
\fin

Now applying the Furstenberg-Zimmer Inverse
Theorem~\ref{thm:rel-ind-joinings} to the invariant functions on
$(X^3,\mu^{\rm{F}})$ in view of the above-found relative
independence over $\a_1$, $\a_2$ and $\a_3$, and coupling its
conclusion with Lemma~\ref{lem:from-subchar-to-char}, we deduce that
the extension of factors
$\bfY_i\stackrel{\a_i|_{\xi_i}}{\longrightarrow} \bfW_i$ must be
$T^{\bf{p}_i}$-isometric:

\begin{lem}
The minimal characteristic factors $\xi_1$, $\xi_2$, $\xi_3$ of an
FIS$^+$ $\bbZ^2$-system $\bfX$ satisfy
\[\a_i \precsim \xi_i
\precsim \zeta_{1/\a_i}^{T^{\bf{p}_i}}\quad\quad\hbox{for
}i=1,2,3.\] \qed
\end{lem}

\textbf{Remark}\quad In general for groups $\L\leq \G$, an extension
of a $\G$-system that is relatively ergodic and isometric for $T^{\
\uhr\L}$ for some proper subgroup $\L \leq \G$ need not be isometric
for the rest of the $\G$-action. Indeed, it is this fundamental
difficulty that mandates the notion of `primitive extension',
allowing the juxtaposition of isometric behaviour in some directions
and relatively weak-mixing behaviour in others, in Furstenberg and
Katznelson's original work on the multidimensional Szemer\'edi
Theorem~\cite{FurKat78}.  For this reason, the above lemma by itself
tells us little about the behaviour of the transformations
$T^\bf{n}$ that are linearly independent from $T^{\bf{p}_i}$ on the
factors $\xi_i$. In fact we will find that after ascending to a
suitable extension, the extension $\a_i|_{\xi_i}$ must be isometric
--- and even Abelian
--- for the whole of the $\bbZ^2$-action, but we will need several
more steps before reaching this fact. \fin

It follows that in order to identify the $\vec{T}$-invariant factor
of $(X^3,\mu^{\rm{F}})$ as far as it extends above $\a_1\times
\a_2\times \a_3$, it suffices to consider the restriction
\[(\zeta_{1/\a_1}^{T^{\bf{p}_1}}\times
\zeta_{1/\a_2}^{T^{\bf{p}_2}}\times
\zeta_{1/\a_3}^{T^{\bf{p}_3}})_\#\mu^{\rm{F}}\] of the Furstenberg
self-joining to a joining of the factors
$\bfZ_1^{T^{\bf{p}_i}}(\bfX/\a_i)$ for $i=1,2,3$.

Let us temporarily introduce the abbreviations $\zeta_i :=
\zeta_{1/\a_i}^{T^{\bf{p}_i}}$ and $\bfZ_i :=
\bfZ_1^{T^{\bf{p}_i}}(\bfX/\a_i)$; and let us also write $\bfZ$ for
the joining of the $\bfZ_i$ obtained by restricting $\mu^{\rm{F}}$;
$\bfW$ for the joining of the $\bfW_i$ obtained by restricting it
further; and $\vec{\a}$, $\vec{\xi}$ and $\vec{\zeta}$ for the
factor maps $\a_1\times\a_2\times\a_3$,
$\xi_1\times\xi_2\times\xi_3$ and
$\zeta_1\times\zeta_2\times\zeta_3$ of $(X^3,\mu^{\rm{F}})$
respectively.

Now we make our first appeal to the fibre-normality contained in the
FIS$^+$ condition.  Since by assumption $\bfX$ is fibre-normal over
\[\a_i = \zeta^{\bfX}_{\sfZ_0^{\bf{p}_i}\vee \sfZ_0^{\bf{p}_i - \bf{p}_j}\vee \sfZ_0^{\bf{p}_i - \bf{p}_k}},\]
we can coordinatize
\begin{center}
$\phantom{i}$\xymatrix{ \bfZ_i^{\
\uhr\bf{p}_i}\ar[dr]_{\a_i|_{\zeta_i}}\ar@{<->}[rr]^-\cong &&
\bfW_i^{\ \uhr\bf{p}_i}\ltimes
(G_{i,\bullet},m_{G_{i,\bullet}},\s_i)\ar[dl]^{\rm{canonical}}\\
& \bfW_i^{\ \uhr\bf{p}_i} }
\end{center}
as extensions by compact group data $G_{i,\bullet}$ for some
cocycle-sections $\s_i:W_i\to G_{i,\bullet}$ over $T^{\bf{p}_i}$.

Of course, knowing only that the factors $\xi_i$ are intermediate
between $\a_i$ and $\zeta_i$, they might still require nontrivial
homogeneous space data in coordinatizations as extensions of $\a_i$.
The simpler structure of fibre-normal extensions will prove crucial
shortly (during the proof of
Proposition~\ref{prop:somewhat-better-Mackey} in the next
subsection), so now we turn our attention to these maximal isometric
extensions to gain further insight into the relatively
$\vec{T}$-invariant factor over $\a_1\times \a_2\times \a_3$.  We
will eventually deduce that the extensions $\a_i|_{\xi_i}$ must in
fact have their own fairly simple structure in an FIS$^+$ system.

The above coordinatizations of the extensions $\a_i|_{\zeta_i}$
combine to give a coordinatization of the action of $\vec{T}$ on the
extension
$\bfZ\stackrel{\vec{\a}|_{\vec{\zeta}}}{\longrightarrow}\bfW$ as
\begin{center}
$\phantom{i}$\xymatrix{
(Z,\vec{\zeta}_\#\mu^{\rm{F}},\vec{T}|_{\vec{\zeta}})\ar[dr]_{\vec{\a}|_{\vec{\zeta}}}\ar@{<->}[rr]^-\cong
&& (W,\vec{\a}_\#\mu^{\rm{F}},\vec{T}|_{\vec{\a}})\ltimes
(\vec{G}_\bullet,m_{\vec{G}_\bullet},\vec{\s})\ar[dl]^{\rm{canonical}}\\
& (W,\vec{\a}_\#\mu^{\rm{F}},\vec{T}|_{\vec{\a}}) }
\end{center}
by the compact group data $\vec{G}_\bullet :=
G_{1,\pi_1(\bullet)}\times G_{2,\pi_2(\bullet)}\times
G_{3,\pi_3(\bullet)}$ and the cocycle-section $\vec{\s} :=
(\s_1,\s_2,\s_3):W \to \vec{G}_\bullet$ over $\vec{T}|_{\vec{\a}}$.
Note that under this coordinatization the measure
$\vec{\zeta}_\#\mu^{\rm{F}}$, which we know is a
$(\vec{T}|_{\vec{\a}}\ltimes \vec{\s})$-invariant lift of
$\vec{\a}_\#\mu^{\rm{F}}$, must actually equal
$\vec{\a}_\#\mu^{\rm{F}}\ltimes m_{G_\bullet}$ since the three
coordinate projections on $Z$ are relatively independent over their
further factors $\a_1$, $\a_2$ and $\a_3$.

At this point the non-ergodic Mackey
Theorem~\ref{thm:homo-nonergMackey-1} (specialized to the case of a
fibre-normal extension) comes to bear, immediately giving the
following.

\begin{prop}\label{prop:initial-Mackey-description}
For an FIS$^+$ $\bbZ^2$-system $\bfX$ there are measurable compact
$\vec{T}|_{\vec{\a}}$-invariant subgroup data $M_\bullet \leq
\vec{G}_\bullet$ and a Borel section $b:W\to \vec{G}_\bullet$ such
the $\vec{T}|_{\vec{\zeta}}$-invariant factor of
$(Z,\vec{\zeta}_\#\mu^{\rm{F}})$ is coordinatized by the map
\begin{multline*}
((w_1,w_2,w_3),(g_1,g_2,g_3))\\ \mapsto
\big(\zeta_0^{\vec{T}|_{\vec{\a}}}(w_1,w_2,w_3),M_{(w_1,w_2,w_3)}\cdot
b(w_1,w_2,w_3)\cdot (g_1,g_2,g_3)\big)
\end{multline*}
from $Z$ to $Z_0^{\vec{T}|_{\vec{\a}}}\ltimes
M_\bullet\backslash\vec{G}_\bullet$, and if the probability kernel
$P:Z_0^{\vec{T}|_{\vec{\a}}}\stackrel{\rm{p}}{\to} W$ represents the
$\vec{T}|_{\vec{\a}}$-ergodic decomposition of
$\vec{\a}_\#\mu^{\rm{F}}$ then the probability kernel
\[P'\ :\ Z_0^{\vec{T}|_{\vec{\a}}}\ltimes M_\bullet\backslash\vec{G}_\bullet \stackrel{\rm{p}}{\to} Z\ :\ (s,M_s \vec{g}')\stackrel{\rm{p}}{\mapsto} P(s,\,\cdot\,)\ltimes m_{b(\bullet)^{-1}\cdot M_s\cdot \vec{g}'}\]
represents the $\vec{T}|_{\vec{\zeta}}$-ergodic decomposition of
$\vec{\zeta}_\#\mu^{\rm{F}}$.
\end{prop}

We will generally refer to $M_\bullet$ and $b$ as the
\textbf{joining Mackey group} and the \textbf{joining Mackey
section} respectively, and will refer to them together as the
\textbf{joining Mackey data}.

It is from this proposition that our finer analysis of $\xi_1$,
$\xi_2$ and $\xi_3$ will really commence. This gives us a picture of
the $\vec{T}$-invariant factor of $(X^3,\mu^{\rm{F}})$ over the
proto-characteristic factors $\a_1$, $\a_2$ and $\a_3$ in terms of
much more concrete data such as the joining Mackey group and
section, for whose analysis some much more delicate tools are
available.  After some further preliminary work in the next
subsection, we will begin this analysis of the Mackey data in
Subsection~\ref{subs:Mackey-full-twodims} by showing that in fact in
an FIS$^+$ system the joining Mackey group must be relatively
`large', in the sense that the relatively $\vec{T}$-invariant
subextension of
$\bfZ\stackrel{\vec{\a}|_{\vec{\zeta}}}{\longrightarrow}\bfW$ that
remains after quotienting by it is always describable in terms of
compact group data extensions of each individual $\a_i$ by
\emph{Abelian groups}, and with cocycles that must satisfy a certain
combined coboundary equation.  This will give a description of each
$\xi_i$ as an Abelian isometric extension of $\a_i$ for the
restriction of the transformation $T^{\bf{p}_i}$. From there we will
show that each of these extensions is actually Abelian isometric for
the whole $\bbZ^2$-action, and then give a much more careful
analysis of the consequences of the equation relating the different
cocycles until the particular structures of Theorem~\ref{thm:char-three-lines-in-2D} emerge.

\subsection{The joining of the proto-characteristic factors}

The following proposition will give some useful insight into the
structure of the join of the proto-characteristic factors $\a_i$
under $\mu^\rm{F}$.

\begin{prop}\label{prop:subchar-joint-dist}
Under $\mu^\rm{F}$ the factors
\[\zeta_0^{T^{\bf{p}_1} = T^{\bf{p}_2}}\circ\pi_1 \simeq \zeta_0^{T^{\bf{p}_1} = T^{\bf{p}_2}}\circ\pi_2,\]
\[\zeta_0^{T^{\bf{p}_1} = T^{\bf{p}_3}}\circ\pi_1 \simeq \zeta_0^{T^{\bf{p}_1} = T^{\bf{p}_3}}\circ\pi_3\]
and
\[\zeta_0^{T^{\bf{p}_2} = T^{\bf{p}_3}}\circ\pi_2 \simeq \zeta_0^{T^{\bf{p}_2} = T^{\bf{p}_3}}\circ\pi_3\]
are relatively independent over
\[\zeta_0^{T^{\bf{p}_1} = T^{\bf{p}_2} = T^{\bf{p}_3}}\circ\pi_1\simeq \zeta_0^{T^{\bf{p}_1} = T^{\bf{p}_2} = T^{\bf{p}_3}}\circ\pi_2\simeq \zeta_0^{T^{\bf{p}_1} = T^{\bf{p}_2} = T^{\bf{p}_3}}\circ\pi_3.\]
\end{prop}

\textbf{Proof}\quad Let $f_{ij}$ be $\zeta_0^{T^{\bf{p}_i} =
T^{\bf{p}_j}}$-measurable for each pair $\{i,j\}$.  Then
$f_{ij}\circ\pi_i = f_{ij}\circ\pi_j$ $\mu^\rm{F}$-almost surely,
and using this freedom and the observation that for just two
linearly independent directions $\bf{n}_1,\bf{n}_2 \in \bbZ^2$ we
have simply $\mu^\rm{F}_{T^{\bf{n}_1},T^{\bf{n}_2}} =
\mu\otimes_{\zeta_0^{T^{\bf{n}_1}=T^{\bf{n}_2}}}\mu$, we can
evaluate
\begin{eqnarray*}
&& \int_{X^3} (f_{12}\circ\pi_1)\cdot (f_{13}\circ\pi_3)\cdot
(f_{23}\circ\pi_2)\,\d\mu^\rm{F}\\
&&=
\int_{X^3}(f_{12}\circ\pi_1)\cdot (f_{13}\circ\pi_3)\cdot (f_{23}\circ\pi_3)\,\d\mu^{\rm{F}}\\
&&= \int_{X^2}f_{12}\otimes (f_{13}\cdot f_{23})\,\d\mu_{T^{\bf{p}_1},T^{\bf{p}_3}}^{\rm{F}}\\
&&=\int_X
\sfE_\mu(f_{12}\,|\,\zeta_0^{T^{\bf{p}_1}=T^{\bf{p}_3}})\cdot
(f_{13}\cdot f_{23})\,\d\mu.
\end{eqnarray*}
On the other hand we have that $\zeta_0^{T^{\bf{p}_1} =
T^{\bf{p}_2}}$ and $\zeta_0^{T^{\bf{p}_1} = T^{\bf{p}_3}}$ are
relatively independent under $\mu$ over their meet
$\zeta_0^{T^{\bf{p}_1} = T^{\bf{p}_2} = T^{\bf{p}_3}}$ (see, for
instance, Lemma 7.3 in~\cite{Aus--ergdirint}) and hence that
\[\sfE_\mu(f_{12}\,|\,\zeta_0^{T^{\bf{p}_1}=T^{\bf{p}_3}}) =
\sfE_\mu(f_{12}\,|\,\zeta_0^{T^{\bf{p}_1}=T^{\bf{p}_2}=T^{\bf{p}_3}}),\]
and so the last line above simplifies to
\[\int_X
\sfE_\mu(f_{12}\,|\,\zeta_0^{T^{\bf{p}_1}=T^{\bf{p}_2} =
T^{\bf{p}_3}})\cdot (f_{13}\cdot f_{23})\,\d\mu,\] and reversing our
steps we find that this is also equal to
\[\int_{X^3} (\sfE_\mu(f_{12}\,|\,\zeta_0^{T^{\bf{p}_1}=T^{\bf{p}_2} =
T^{\bf{p}_3}})\circ\pi_1)\cdot (f_{13}\circ\pi_3)\cdot
(f_{23}\circ\pi_2)\,\d\mu^\rm{F}.\]

Arguing similarly for the pairs $13$ and $23$ we obtain
\begin{multline*}
\int_{X^3} (f_{12}\circ\pi_1)\cdot (f_{13}\circ\pi_3)\cdot
(f_{23}\circ\pi_2)\,\d\mu^\rm{F} \\ = \int_{X^3}
(\sfE_\mu(f_{12}\,|\,\zeta_0^{T^{\bf{p}_1}=T^{\bf{p}_2} =
T^{\bf{p}_3}})\circ\pi_1)\cdot
(\sfE_\mu(f_{13}\,|\,\zeta_0^{T^{\bf{p}_1}=T^{\bf{p}_2} =
T^{\bf{p}_3}})\circ\pi_3)\\ \cdot
(\sfE_\mu(f_{23}\,|\,\zeta_0^{T^{\bf{p}_1}=T^{\bf{p}_2} =
T^{\bf{p}_3}})\circ\pi_2)\,\d\mu^\rm{F},
\end{multline*}
as required. \qed

\subsection{The joining Mackey group has full two-dimensional
projections}\label{subs:Mackey-full-twodims}

In this subsection and the next our main goal is to prove that any
system has an FIS$^+$ extension for which the characteristic factors
can be coordinatized so that the Mackey group data of the
Furstenberg self-joining must be particularly simple.

Let us first introduce another useful piece of notation.

\begin{dfn}[Motionless data]\label{dfn:motionless}
If $(X,\mu,T)$ is a $\bbZ^2$-system and $x\mapsto G_x$ is a
measurable assignment of compact Abelian groups (from some fixed
fibre repository, as in Definition 3.1 of~\cite{Aus--ergdirint}),
then we will say that this assignment is \textbf{motionless} if it
is invariant under the whole $\bbZ^2$-action.  This situation will
always and exclusively be denoted by the use of the notation
$G_\star$ in place of $G_\bullet$, in which case we will often omit
to mention the motionlessness by name.
\end{dfn}

We make such efforts to work with non-ergodic data in order to avoid assuming that the $\bbZ^2$-actions we handle are ergodic overall, which would introduce the difficulty of repeatedly ensuring that our various extensions retain this ergodicity.  The reader will lose nothing by
thinking of group data of the form $G_\star$ as `effectively
constant' (since all the constructions we perform with such data
will be manifestly measurable).  The quality of motionlessness will
contrast, however, with group data over a $\bbZ^2$-system that is
invariant only for certain subactions, which will occur repeatedly
in the following.

\begin{prop}\label{prop:better-Mackey}
Any $\bbZ^2$-system $\bfX_0$ admits an FIS$^+$ extension
$\pi:\bfX\to\bfX_0$ in which the factors $\xi_i:\bfX\to \bfY_i$,
$i=1,2,3$, of the minimal characteristic triple can be coordinatized
over the proto-characteristic factors as
\begin{center}
$\phantom{i}$\xymatrix{
\bfY_i\ar[dr]_{\a_i|_{\xi_i}}\ar@{<->}[rr]^-\cong & &
\bfW_i\ltimes (A_\star,m_{A_\star},\s_i)\ar[dl]^{\rm{canonical}}\\
& \bfW_i}
\end{center}
for some compact Abelian group data $A_\star$ and cocycle-sections
$\s_i:\bbZ^2\times W_i \to A_\star$ over $T|_{\a_i}$ in such a way
that the resulting joining Mackey group data is
\[M_\star = \{(a_1,a_2,a_3) \in A_\star^3:\ a_1\cdot a_2\cdot a_3 = 1_{A_\star}\}\]
(noting that $\zeta_0^T\circ\pi_1 \simeq \zeta_0^T\circ\pi_2 \simeq
\zeta_0^T\circ\pi_3$ and so for $A_\star$ we have $A_{w_1} = A_{w_2}
= A_{w_3}$ $\mu^{\rm{F}}$-almost surely) and the joining Mackey
section may be expressed as some $b:W_1\times W_2\times W_3 \to
A_\star$ that satisfies
\[\s_1(\bf{p}_1,w_1)\cdot \s_2(\bf{p}_2,w_2)\cdot \s_3(\bf{p}_3,w_3) = \DDelta_{T|_{\a_1}^{\bf{p}_1}\times T|_{\a_2}^{\bf{p}_2}\times T|_{\a_3}^{\bf{p}_3}}b(w_1,w_2,w_3)\]
at $\vec{\a}_\#\mu^{\rm{F}}$-almost every $(w_1,w_2,w_3)$.

We will refer to the above group $M_\star$ as the \textbf{zero-sum}
subgroup of $A_\star^3$ (this is slightly abusive, since we write
the Abelian operation of $A_\star$ multiplicatively in this
subsection, but it should cause no confusion).
\end{prop}

\textbf{Remarks}\quad\textbf{1.}\quad Recall from Section 4.1
of~\cite{Aus--lindeppleasant1} that $\mu^{\rm{F}}$ is
$(T^{\bf{p}_1}\times T^{\bf{p}_2}\times T^{\bf{p}_3})$-invariant, so
that the appearance of a coboundary over $T|_{\a_1}^{\bf{p}_1}\times
T|_{\a_2}^{\bf{p}_2}\times T|_{\a_3}^{\bf{p}_3}$ above should cause
no concern.

\quad\textbf{2.}\quad At this stage our results are still geared
towards understanding the three factors $\xi_i:\bfX\to\bfY_i$,
$i=1,2,3$ separately. They do not immediately tell us anything about
the joint distribution of these factors under $\mu$.  This question
can be rather subtle, but we will learn a little more later in Propositions~\ref{prop:joining-Kron-to-isotropies} and~\ref{prop:jointdists}. \fin

Proposition~\ref{prop:better-Mackey} asserts that the whole
$\bbZ^2$-action $T|_{\xi_i}$ can be coordinatized as an extension of
$T|_{\a_i}$ in terms of an Abelian cocycle, rather than just the
$(\bbZ\bf{p}_i)$-subaction as discussed previously. We will prove it via a weaker result which gives a similar coordinatization of the extension
$\bfY_i\to\bfW_i$, but allows some additional `twisting' in the
joining Mackey group and does not yet give isometricity for the
actions of the whole of $\bbZ^2$.

\begin{prop}\label{prop:somewhat-better-Mackey}
If $\bfX$ is FIS$^+$ then the factors $\xi_i:\bfX\to\bfY_i$ have
$(\bbZ\bf{p}_i)$-subactions that can be coordinatized over $\a_i$ as
\begin{center}
$\phantom{i}$\xymatrix{
\bfY^{\ \uhr\bf{p}_i}_i\ar[dr]_{\a_i|_{\xi_i}}\ar@{<->}[rr]^-\cong & & \bfW^{\ \uhr\bf{p}_i}_i\ltimes (A_\star,m_{A_\star},\s_i)\ar[dl]^{\rm{canonical}}\\
&\bfW^{\ \uhr\bf{p}_i}_i }
\end{center}
for some compact Abelian group data $A_\star$ and cocycle-sections
$\s_i:W_i\to A_\star$ over $T|^{\bf{p}_i}_{\a_i}$ in such a way that
there are measurable families of isomorphisms
$\Theta_{i,\bullet}:W_1\times W_2\times W_3 \to \rm{Aut}\,A_\star$
such that the joining Mackey group data is
\[M_{\vec{w}} = \{(a_1,a_2,a_3)\in A_s^3:\ \Theta_{1,\vec{w}}(a_1)\cdot\Theta_{2,\vec{w}}(a_2)\cdot\Theta_{3,\vec{w}}(a_3) = 1_{A_s}\}\]
at $\vec{\a}_\#\mu^{\rm{F}}$-almost every $\vec{w} =
(w_1,w_2,w_3)\in W$, where $s := \zeta_0^T(w_1)$.
\end{prop}

\textbf{Remark}\quad Proposition~\ref{prop:somewhat-better-Mackey}
deduces some properties of the joining Mackey group merely from the
FIS$^+$ property.  By contrast, we will find that `straightening
out' the families of automorphisms $\Theta_{i,\bullet}$ to obtain
Proposition~\ref{prop:better-Mackey} will generally require a
further extension even if the original system was already FIS$^+$,
hence the form in which Proposition~\ref{prop:better-Mackey} is
phrased. \fin

This subsection will be dedicated to the proof of
Proposition~\ref{prop:somewhat-better-Mackey}, and we will then
deduce Proposition~\ref{prop:better-Mackey} from it in the next
subsection.

The technical result that really underlies
Proposition~\ref{prop:somewhat-better-Mackey} is the following. Part
of its interest is that its proof will use satedness in a new way,
not seen in the simpler arguments of~\cite{Aus--lindeppleasant1}.

\begin{lem}\label{lem:2D-proj-full}
If $\bfX = (X,\mu,T)$ is FIS$^+$, then under any coordinatizations
of the extensions
\begin{center}
$\phantom{i}$\xymatrix{
(\bfZ_{1/\a_i}^{T^{\bf{p}_i}})^{\ \uhr\bf{p}_i}\ar[dr]_{\rm{restriction\ of\ }\a_i}\ar@{<->}[rr]^-\cong & & \bfW_i^{\ \uhr\bf{p}_i}\ltimes (G_{i,\bullet},m_{G_{i,\bullet}},\s_i)\ar[dl]^{\rm{canonical}}\\
&\bfW_i^{\ \uhr\bf{p}_i} }
\end{center}
the joining Mackey group data $M_{(w_1,w_2,w_3)}$ has full
two-dimensional projections onto $G_{i,w_i}\times G_{j,w_j}$ for
$1\leq i < j\leq 3$ for $\vec{\a}_\#\mu^\rm{F}$-almost every
$(w_1,w_2,w_3)$.
\end{lem}

\textbf{Proof}\quad By symmetry it suffices to treat the case of
\[M_{12,(w_1,w_2,w_3)} := \{(g_1,g_2):\ \exists g_3 \in G_{3,w_3}\ \hbox{s.t.}\ (g_1,g_2,g_3)\in M_{(w_1,w_2,w_3)}\}.\]

Let us abbreviate $\bfZ_{1/\a_i}^{T^{\bf{p}_i}} =: \bfZ_i$ for
$i=1,2$, and now let $\vec{\bfZ}$ be the factor of the Furstenberg
self-joining $\bfX^\rm{F}$ generated by the factor maps
\[\bfX^\rm{F}\stackrel{\pi_1}{\longrightarrow}\bfX\to \bfZ_1,\]
\[\bfX^\rm{F}\stackrel{\pi_2}{\longrightarrow}\bfX\to \bfZ_2\]
and \[\bfX^\rm{F}
\stackrel{\pi_3}{\longrightarrow}\bfX\stackrel{\a_3}{\longrightarrow}
\bfW_3\] (so we do not keep the whole of $\bfZ_3$ in the third
factor). As a factor of $\bfX^\rm{F}$ this extends
$\vec{\a}:\bfX^\rm{F}\to\vec{\bfW}$, and the above coordinatizations
of $\bfZ_i^{\ \uhr\bf{p}_i}\to\bfW_i^{\ \uhr\bf{p}_i}$ for $i=1,2$
combine to coordinatize the action of the restriction of
$T^{\bf{p}_1}\times T^{\bf{p}_2}\times T^{\bf{p}_3}$ on
$\vec{\bfZ}\to\vec{\bfW}$ as an extension by the product group data
$G_{1,\pi_1(\bullet)}\times G_{2,\pi_2(\bullet)}$ with the above
product cocycle and with Mackey data $M_{12,(w_1,w_2,w_3)}$.

Let $\sfC := \sfZ_0^{\bf{p}_1}\vee \sfZ_0^{\bf{p}_1 - \bf{p}_2}\vee
\sfZ_0^{\bf{p}_1 - \bf{p}_3}$ and $\sfD := \sfZ_0^{\bf{p}_1}\vee
\sfZ_0^{\bf{p}_1 - \bf{p}_2}$. We will construct an extension of
$\t{\bfX}\to \bfX$ to which we can apply the assumption of
satedness.  In fact, letting $\L:=\bbZ\bf{p}_1 + \bbZ\bf{p}_2$ (a
full-rank sublattice of $\bbZ^2$), we will first use $\bfX^{\rm{F}}$
to construct an extension of the subaction system $\bfX^{\ \uhr\L}$,
then extend this further to recover an action of the whole of
$\bbZ^2$, and then argue that the maximal $\sfC$-factor of this
further extension forces us to the desired conclusion.

To extend $\bfX^{\ \uhr\L}$ let $\bfX'$ be the $\L$-system
constructed on the Furstenberg self-joining $(X^3,\mu^{\rm{F}})$ by
lifting $T^{\bf{p}_1}$ to $\t{T}^{\bf{p}_1}:= T^{\bf{p}_1}\times
T^{\bf{p}_2}\times T^{\bf{p}_3}$ and $T^{\bf{p}_2}$ to
$\t{T}^{\bf{p}_2} := (T^{\bf{p}_2})^{\times 3}$. Then
$\t{T}^{\bf{p}_1}$ and $\t{T}^{\bf{p}_2}$ both act as $T^{\bf{p}_2}$
on the second coordinate in $X^3$, so
\[\pi_2 \precsim \zeta_{0}^{\t{T}^{\bf{p}_1} =
\t{T}^{\bf{p}_2}};\] and also, we clearly have
\[(\zeta_0^{T^{\bf{p}_1}}\circ\pi_1)\vee (\zeta_0^{T^{\bf{p}_2}}\circ\pi_2)\vee (\zeta_0^{T^{\bf{p}_3}}\circ\pi_3)\precsim \zeta_0^{\t{T}^{\bf{p}_1}}.\]

On the other hand, under $\mu^\rm{F}$ we have
\[\zeta_0^{T^{\bf{p}_1}=T^{\bf{p}_3}}\circ\pi_1 \simeq
\zeta_0^{T^{\bf{p}_1}=T^{\bf{p}_3}}\circ\pi_3\] and
\[\zeta_0^{T^{\bf{p}_2}=T^{\bf{p}_3}}\circ\pi_2 \simeq
\zeta_0^{T^{\bf{p}_2}=T^{\bf{p}_3}}\circ\pi_3,\] so overall these
relations give
\begin{eqnarray*}
\a_3\circ\pi_3 &\simeq&
(\zeta_0^{T^{\bf{p}_1}=T^{\bf{p}_3}}\circ\pi_3)\vee(\zeta_0^{T^{\bf{p}_2}=T^{\bf{p}_3}}\circ\pi_3)\vee
(\zeta_0^{T^{\bf{p}_3}}\circ\pi_3)\\ &\simeq&
(\zeta_0^{T^{\bf{p}_1}=T^{\bf{p}_3}}\circ\pi_1)\vee(\zeta_0^{T^{\bf{p}_2}=T^{\bf{p}_3}}\circ\pi_2)\vee
(\zeta_0^{T^{\bf{p}_3}}\circ\pi_3)\\ &\precsim&
(\zeta_0^{T^{\bf{p}_1}=T^{\bf{p}_3}}\circ\pi_1)\vee
\zeta_0^{\t{T}^{\bf{p}_1}=\t{T}^{\bf{p}_2}}\vee
\zeta_0^{\t{T}^{\bf{p}_1}}
\end{eqnarray*}
and so also
\[(\a_1\circ\pi_1)\vee(\a_2\circ\pi_2)\vee(\a_3\circ\pi_3)\precsim
(\a_1\circ\pi_1)\vee
\zeta_0^{\t{T}^{\bf{p}_1}=\t{T}^{\bf{p}_2}}\vee\zeta_0^{\t{T}^{\bf{p}_1}}
\simeq (\zeta_{\sfC}^\bfX\circ\pi_1) \vee \zeta_{\sfD^{\
\uhr\L}}^{\bfX'}.\]

Now let $\pi:\t{\bfX}\to\bfX'\to\bfX$ be any further extension that
recovers an action of the whole of $\bbZ^2$ (this can always be
done: see, for instance, Subsection 3.2
in~\cite{Aus--lindeppleasant1}), so we must still have
\[(\a_1\circ\pi_1)\vee(\a_2\circ\pi_2)\vee(\a_3\circ\pi_3)\precsim
(\zeta_{\sfC}^\bfX\circ\pi) \vee \zeta_{\sfD}^{\t{\bfX}}.\]

Finally, the projection $M_{12,\bullet}$ is the Mackey group data
for the group data extension
\[\t{\bfW}^{\ \uhr\bf{p}_1}\ltimes (G_{1,\pi_1(\bullet)}\times G_{2,\pi_2(\bullet)},m_{G_{1,\pi_1(\bullet)}\times G_{2,\pi_2(\bullet)}},(\s_{1,\pi_1(\bullet)},\s_{2,\pi_2(\bullet)})).\]
The above construction locates this group data extension as a factor
of $\t{\bfX}$ that is contained within the joining of
$\pi:\t{\bfX}\to\bfX$ and $(\zeta_{\sfC}^\bfX\circ\pi) \vee
\zeta_{\sfD}^{\t{\bfX}} \precsim \zeta_{\sfC}^{\t{\bfX}}$. By
$\sfC$-satedness these two factors of $\t{\bfX}$ must be relatively
independent over
\[\zeta_\sfC^\bfX\circ \pi = \a_1\circ \pi.\]
In terms of the above coordinatizations, this implies that if $\t{x}$ is drawn from the probability distribution $\t{\mu}$, then its image $\zeta_\sfC^{\t{\bfX}}(\t{x})$ exactly determines the points
\begin{itemize}
\item $\big(\zeta_0^{T|_{\a_1}^{\bf{p}_1}\times T|_{\a_2}^{\bf{p}_2}\times T|_{\a_3}^{\bf{p}_3}}(w_1,w_2,w_3),M_{12,(w_1,w_2,w_3)}(g_1,g_2)\big)\ \ $ (because this is given by the restriction of
$\zeta_0^{\t{T}^{\bf{p}_1}} \precsim \zeta_\sfC^{\t{\bfX}}$ to
$\vec{\bfZ}$),
\item $(w_1,w_2,w_3)\ \ $ (because we have seen that $(\a_1\circ\pi_1)\vee(\a_2\circ\pi_2)\vee(\a_3\circ\pi_3)\precsim
\zeta_{\sfC}^{\t{\bfX}}$)
\item and $(w_2,g_2)\ \ $ (because $\pi_2
\lesssim \zeta_0^{\t{T}^{\bf{p}_1} = \t{T}^{\bf{p}_2}}$),
\end{itemize}
but all this information is conditionally independent from $(w_1,g_1)$ given $w_1 = \zeta_\sfC^\bfX\circ\pi(\t{x})$. This is possible only if
$M_{12,(w_1,w_2,w_3)} = G_{1,w_1}\times G_{2,w_2}$ almost surely, as
required. \qed

\textbf{Remark}\quad It is worth noting that although the
contradiction we obtain above is with isotropy-satedness, we have
used the full FIS$^+$ assumption because we have worked throughout
with an extension by group data.  In fact the above argument runs
into difficulties if we try to work with general homogeneous space
extensions, say by $G_{i,\bullet}/H_{i,\bullet}$, because in that
setting we cannot rule out that the group $M_{12,(w_1,w_2,w_3)}$ is
not the whole of $G_{1,w_1}\times G_{2,w_2}$ but is nevertheless
large enough that \[M_{12,(w_1,w_2,w_3)}(H_{1,w_1}\times H_{2,w_2})
= G_{1,w_1}\times G_{2,w_2}\] almost surely (which latter conclusion
is too weak for the next step of our argument below). This makes an
interesting contrast with the study of characteristic factors (even
without the freedom to pass to extensions) for just two commuting
transformations given in~\cite{Aus--ergdirint}.  There the relevant
joining Mackey group could be shown always to have full
one-dimensional projections, essentially because in that case the
joining of the proto-characteristic factors underneath this Mackey
group data is so simple that the one-dimensional projections of the
joining Mackey group data can easily be related to Mackey group data
for the isometric extensions in the original system (without
constructing an extension). It seems that matters become genuinely
more complicated for three-fold or higher Furstenberg self-joinings,
and some extra procedure such as the passage to fibre-normal
extensions is needed. \fin

Moving onwards, we will now make use of the following
group-theoretic lemma from Furstenberg and Weiss~\cite{FurWei96}.

\begin{lem}[Lemma 9.1
in~\cite{FurWei96}]\label{lem:group-th-from-FurWei} If $G_1$, $G_2$
and $G_3$ are compact metrizable groups and $M \leq G_1\times
G_2\times G_3$ has full two-dimensional projections then there are a
compact metrizable Abelian group $A$ and continuous epimorphisms
$\Psi_i:G_i\longrightarrow A$ (so that, in particular,
$[G_i,G_i]\leq \ker\Psi_i$) such that
\[M = \{(g_1,g_2,g_3):\ \Psi_1(g_1)\cdot\Psi_2(g_2)\cdot\Psi_3(g_3) = 1_A\}.\]
\qed
\end{lem}

In order to use this lemma, we need just a little more information
on the structure of the slices of $M_\bullet$, which we now acquire
in a few more short steps.

\begin{lem}\label{lem:invariances-leave-first-coord}
For an FIS$^+$ system we have
\[\zeta_0^{\vec{T}|_{\vec{\a}}}\wedge
\zeta_0^{(T^{\bf{p}_1})^{\times 3}|_{\vec{\a}}} \simeq
\zeta_0^{T^{\bf{p}_1}}\circ\pi_1|_{\vec{\a}}:\] that is, any
measurable subset of $W$ that is both $(T|_{\a_1}^{\bf{p}_1}\times
T|_{\a_2}^{\bf{p}_2}\times T|_{\a_3}^{\bf{p}_3})$-invariant and
$(T|_{\a_1}^{\bf{p}_1}\times T|_{\a_2}^{\bf{p}_1}\times
T|_{\a_3}^{\bf{p}_1})$-invariant is equal up to an
$\vec{\a}_\#\mu^\rm{F}$-negligible set to a
$T|_{\a_1}^{\bf{p}_1}$-invariant subset of $W_1$ lifted through the
first coordinate projection $W \to W_1$.
\end{lem}

\textbf{Proof}\quad The relation
$\zeta_0^{\vec{T}|_{\vec{\a}}}\wedge \zeta_0^{(T^{\bf{p}_1})^{\times
3}|_{\vec{\a}}} \succsim
\zeta_0^{T^{\bf{p}_1}}\circ\pi_1|_{\vec{\a}}$ is clear, so we focus
on its reverse.

Recall that for an FIS$^+$ system we have $\a_i = \b_i\vee
\zeta_0^{T^{\bf{p}_i}}$ with $\b_i = \zeta_0^{T^{\bf{p}_i}=
T^{\bf{p}_j}}\vee \zeta_0^{T^{\bf{p}_i} = T^{\bf{p}_k}}$.  Therefore
\[\vec{\a} \simeq (\zeta_0^{T^{\bf{p}_1}}\times
\zeta_0^{T^{\bf{p}_2}}\times \zeta_0^{T^{\bf{p}_3}})\vee (\b_1\times
\b_2\times \b_3).\]  The first of these factors is already invariant
under the restriction of $\vec{T}$ and so we have
\[\zeta_0^{\vec{T}|_{\vec{\a}}}\circ\vec{\a} \simeq (\zeta_0^{T^{\bf{p}_1}}\times
\zeta_0^{T^{\bf{p}_2}}\times \zeta_0^{T^{\bf{p}_3}})\vee
\big(\zeta_0^{\vec{T}|_{\vec{\a}}}\wedge (\b_1\times \b_2\times
\b_3)\big)\] (since the invariant factor of a joining in which the
first coordinate factor has trivial action is simply generated by
the first coordinate factor and the invariant sets of the second
coordinate factor). Let us next identify the second factor in the
join on the right-hand side of this equation.

Since $\zeta_0^{T^{\bf{p}_i} = T^{\bf{p}_j}}\circ\pi_i \simeq
\zeta_0^{T^{\bf{p}_i} = T^{\bf{p}_j}}\circ\pi_j$ under
$\mu^{\rm{F}}$, the factor $\b_1\times \b_2\times \b_3$ is actually
$\mu^{\rm{F}}$-almost surely determined by the first two coordinates
in $X^3$, and so it will suffice to identify
$\zeta_0^{T|_{\a_1}^{\bf{p}_1}\times T|_{\a_2}^{\bf{p}_2}}\wedge
(\b_1\times \b_2)$. Now, an easy calculation shows that that the
two-dimensional Furstenberg self-joining
$\mu^{\rm{F}}_{T^{\bf{p}_1},T^{\bf{p}_2}}$ is just the relatively
independent product $\mu\otimes_{\zeta_0^{T^{\bf{p}_1} =
T^{\bf{p}_2}}}\mu$ on $X^2$; and in view of the FIS property and the
consequent pleasantness of our system for all linearly independent
double ergodic averages (see Proposition 4.5
of~\cite{Aus--lindeppleasant1} and
Lemma~\ref{lem:from-subchar-to-char}), we have further that under
$\mu^{\rm{F}}_{T^{\bf{p}_1},T^{\bf{p}_2}}$ all $(T^{\bf{p}_1}\times
T^{\bf{p}_2})$-invariant subsets are measurable up to negligible
sets with respect to the factor $\zeta_0^{T^{\bf{p}_1}}\times
\zeta_0^{T^{\bf{p}_2}}$.

This therefore also applies to any $\vec{T}|_{\vec{\a}}$-invariant
measurable subset of $V_1\times V_2\times V_3$, and so the second
factor in the above join can actually be subsumed into the first to
give
\[\zeta_0^{\vec{T}|_{\vec{\a}}}\circ\vec{\a} \simeq (\zeta_0^{T^{\bf{p}_1}}\times
\zeta_0^{T^{\bf{p}_2}}\times \zeta_0^{T^{\bf{p}_3}}).\]

Finally, we observe similarly that the first coordinate factor of
$(\zeta_0^{T^{\bf{p}_1}}\times \zeta_0^{T^{\bf{p}_2}}\times
\zeta_0^{T^{\bf{p}_3}})$ is already invariant for the restriction of
$(T^{\bf{p}_1})^{\times 3}$, and so to find all sets that are
invariant for this transformation and measurable with respect to
this factor it suffices to consider the second and third
coordinates. Once again we have that the two-dimensional projection
$(\pi_1\times \pi_2)_\#\mu^{\rm{F}} =
\mu^{\rm{F}}_{T^{\bf{p}_2},T^{\bf{p}_3}}$ must simply equal
$\mu\otimes_{\zeta_0^{T^{\bf{p}_2} = T^{\bf{p}_3}}}\mu$, and the FIS
property implies that up to
$\mu^{\rm{F}}_{T^{\bf{p}_2},T^{\bf{p}_3}}$-negligible sets the only
$(T^{\bf{p}_1})^{\times 2}$-invariant sets in this space are
accounted for by the factor $\zeta_0^{T^{\bf{p}_1} = T^{\bf{p}_2} =
\id}\times \zeta_0^{T^{\bf{p}_1} = T^{\bf{p}_3} = \id}$.  Since
under $\mu^{\rm{F}}$ this product is clearly determined by the
$T^{\bf{p}_1}$-invariant factor of the first coordinate, the proof
is complete. \qed

\begin{lem}\label{lem:slices-normal}
If $G_1$ and $G_2$ are compact groups and $M\leq G_1\times G_2$ has
full one-dimensional projections (in the sense that for any $g_1 \in
G_1$ there exists $g_2 \in G_2$ such that $(g_1,g_2)\in M$, and
vice-versa), then the one-dimensional slice of $M$
\[L_1 := \{g_1 \in G_1:\ (g_1,1_{G_2})\in M\}\]
is a closed normal subgroup of $G_1$, and similarly for $L_2 \unlhd
G_2$.
\end{lem}

\textbf{Proof}\quad This is routine except for the conclusion of
normality.  By symmetry it suffices to treat the case $i=1$. Let
$r_1 \in G_1$.  Since $M$ has full one-dimensional projections we
can find $r_2 \in G_2$ such that $(r_1,r_2)\in M$. It is now easy to
check that
\begin{eqnarray*}r_1L_1 &=& \{g\in G_1:\ (r_1^{-1}g,e) \in M\}\\
&=& \{g\in G_1:\ (r_1,r_2)(r_1^{-1}g,e) \in M\}\\
&=&\{g\in G_1:\ (g,r_2) \in M\}\\
&=& \{g\in G_1:\ (gr_1^{-1},e)(r_1,r_2) \in M\}\\
&=& \{g\in G_1:\ (gr_1^{-1},e) \in M\} = L_1r_1.
\end{eqnarray*}
Since $r_1$ was arbitrary, $L_1$ is normal, as required. \qed

\begin{lem}[Deconstructing a relation between two group
correspondences]\label{lem:slices-respected} Suppose that $G_1$,
$G_2$ and $G_3$ are compact groups and that $M_1,M_2\leq G_1\times
G_2\times G_3$ are two subgroups that both have full two-dimensional
projections, and let their one-dimensional slices be
\[L_{i,1} := \{g \in G_1:\ (g,1_{G_2},1_{G_3})\in M_i\}\quad\quad\hbox{for}\ i=1,2\]
and similarly $L_{i,2}$, $L_{i,3}$.  Suppose further that
$\Phi_i:G_i\stackrel{\cong}{\longrightarrow} G_i$ and $h_i,k_i \in
G_i$ for $i=1,2,3$ satisfy
\[(h_1,h_2,h_3)\cdot(\Phi_1\times \Phi_2\times \Phi_3)(M_1)\cdot(k_1,k_2,k_3) = M_2.\]
Then $\Phi_i(L_{1,i}) = L_{2,i}$ for $i=1,2,3$.
\end{lem}

\textbf{Proof}\quad This follows fairly automatically upon checking
the above equation for different particular members of the relevant
group.  Clearly we may assume $i=1$ by symmetry.

Suppose that $g \in L_{1,1}$. Then the given equation tells us that
\[(h_1\cdot \Phi_1(g)\cdot k_1,h_2\cdot k_2,h_3\cdot k_3) = (m_1,m_2,m_3)\]
for some $(m_1,m_2,m_3) \in M_2$, and here, in particular, we have
that $m_2 = h_2\cdot k_2$ and $m_3$ do not depend on $g$.  Since the
above must certainly hold if $g = 1_{G_1}$, applying it also for any
other $g$ and differencing gives
\[(h_1\cdot \Phi_1(g)\cdot k_1)\cdot(h_1\cdot \Phi_1(1_{G_1})\cdot k_1)^{-1} = h_1\cdot\Phi_1(g)\cdot h_1^{-1} \in L_{2,1},\]
so $\Phi_1(L_{1,1}) \subseteq h_1^{-1}\cdot L_{2,1}\cdot h_1$.  An
exactly symmetric argument gives the reverse inclusion, so in fact
$\Phi_1(L_{1,1})$ is a conjugate of $L_{2,1}$.  However, since $M_1$
and $M_2$ have full coordinate projections onto $G_1$ and onto
$G_2\times G_3$, by Lemma~\ref{lem:slices-normal} it follows that in
fact $\Phi_1(L_{1,1}) = L_{2,1}$, as required. \qed

\begin{lem}\label{lem:slices-invar}
The one-dimensional slices of $M_\bullet$,
\[L_{1,(w_1,w_2,w_3)} := \{g_1 \in G_{1,w_1}:\ (g_1,1_{G_{2,w_2}},1_{G_{3,w_3}})\in M_{(w_1,w_2,w_3)}\}\]
and similarly $L_{2,(w_1,w_2,w_3)}$ and $L_{3,(w_1,w_2,w_3)}$, are
virtually functions of $w_1$ (respectively $w_2$, $w_3$) alone.
Also, under the above coordinatizations, for each $i = 1,2,3$ the
map
\[W_i\ltimes G_{i,\bullet}\to W_i\ltimes
(G_{i,\bullet}/L_{i,\bullet}): (w_i,g)\mapsto (w_i,gL_{i,w_i})\]
defines a factor of $\zeta_{1/\a_i}^{T^{\bf{p}_i}}$ for the whole
$\bbZ^2$-action $T$ (that is, it is respected by the restrictions of
every $T^\bf{n}$, not just of $T^{\bf{p}_i}$).
\end{lem}

\textbf{Proof}\quad Clearly by symmetry it suffices to treat the
case of $L_{1,(w_1,w_2,w_3)}$.  The crucial fact here is the
presence of the additional transformations of the factor
$\zeta_0^{\vec{T}}$ given by $(T^{\bf{n}})^{\times 3}$, $\bf{n} \in
\bbZ^2$. We can describe the restriction of $T^{\bf{n}}$ to the
tower of factors $\zeta_{1/\a_j}^{T^{\bf{p}_j}}\succsim \a_j$ for $j
= 1,2,3$ using the Relative Automorphism Structure
Theorem~\ref{thm:RAST}: under the above coordinatization we obtain
\[T|_{\zeta_{1/\a_j}^{T^{\bf{p}_j}}}^\bf{n} \cong T|_{\a_j}^\bf{n}\ltimes
(L_{\rho_{\bf{n},j}(\bullet)}\circ\Phi_{\bf{n},j,\bullet})\] for
some $\rho_{\bf{n},j}:W_j\to G_{j,\bullet}$ and
$T|_{\a_j}^{\bf{p}_j}$-invariant $\Phi_{\bf{n},j,\bullet}:W_j\to
\rm{Isom}(G_{j,\bullet},G_{j,T|_{\a_j}^\bf{n}(\bullet)})$.

Let us now phrase the condition that $S := T^{\bf{n} -
\bf{p}_1}\times T^{\bf{n} - \bf{p}_2}\times T^{\bf{n} - \bf{p}_3}$
respects $\zeta_0^{\vec{T}}$ in terms of these expressions and the
Mackey data. This requires that
$\zeta_0^{\vec{T}}|_{\vec{\zeta}}(S|_{\vec{\zeta}}(\vec{y}))$ depend
only on $\zeta_0^{\vec{T}}|_{\vec{\zeta}}(\vec{y})$, or equivalently
that $S$ almost surely carry the fibres of $\zeta_0^{\vec{T}}$ onto
themselves. In terms of the above Mackey description of
$\zeta_0^{\vec{T}}$ given by
Proposition~\ref{prop:initial-Mackey-description} this asserts that
for Haar-almost every $(g_1',g_2',g_3') \in G_{1,w_1}\times
G_{2,w_2}\times G_{3,w_3}$ there is some
\[(g_1'',g_2'',g_3'') \in G_{1,T|_{\a_1}^{\bf{n} -
\bf{p}_1}(w_1)}\times G_{2,T|_{\a_2}^{\bf{n} - \bf{p}_2}(w_2)}\times
G_{3,T|_{\a_3}^{\bf{n} - \bf{p}_3}(w_3)}\] such that
\begin{multline*}
\Big(\prod_{i=1}^3(L_{\rho_{\bf{n}-\bf{p}_i,i}(w_i)}\circ\Phi_{\bf{n}-\bf{p}_i,i,w_i})\Big)\big(b(w_1,w_2,w_3)^{-1}\cdot
M_{(w_1,w_2,w_3)}\cdot (g_1',g_2',g_3')\big)\\
= b(S|_{\vec{\a}}(w_1,w_2,w_3))^{-1}\cdot
M_{S|_{\vec{\a}}(w_1,w_2,w_3)}\cdot (g_1'',g_2'',g_3''),
\end{multline*}
or, re-arranging, that
\begin{eqnarray*}
&&b(S|_{\vec{\a}}(w_1,w_2,w_3))\cdot (\rho_{\bf{n}-\bf{p}_1,1}(w_1),\rho_{\bf{n}-\bf{p}_2,2}(w_2),\rho_{\bf{n}-\bf{p}_3,3}(w_3))\\
&&\quad\quad\cdot
\Big(\prod_{i=1}^3\Phi_{\bf{n}-\bf{p}_i,i,w_i}\Big)(b(w_1,w_2,w_3)^{-1})\cdot
\Big(\prod_{i=1}^3\Phi_{\bf{n}-\bf{p}_i,i,w_i}\Big)(M_{(w_1,w_2,w_3)})\\
&&\quad\quad\cdot (\Phi_{\bf{n}-\bf{p}_1,1,w_1}(g_1')(g_1'')^{-1},\Phi_{\bf{n}-\bf{p}_2,2,w_2}(g_2')(g_2'')^{-1},\Phi_{\bf{n}-\bf{p}_3,3,w_3}(g_3')(g_3'')^{-1})\\
&&\quad\quad\quad\quad\quad\quad\quad\quad\quad\quad\quad\quad\quad\quad\quad\quad\quad\quad\quad\quad\quad\quad\quad\quad=
M_{S|_{\vec{\a}}(w_1,w_2,w_3)}.
\end{eqnarray*}

We will now deduce the two desired conclusions from treating the
first coordinate projection in this equation using
Lemma~\ref{lem:slices-respected} for different values of $\bf{n}$.
By that lemma the above implies that
\[\Phi_{\bf{n}-\bf{p}_1,1,w_1}(L_{1,(w_1,w_2,w_3)}) = L_{1,S|_{\vec{\a}}(w_1,w_2,w_3)}.\]

If we first specialize this equation to $\bf{n} := \bf{p}_1$, then
of course we simply have $\Phi_{\bf{n}-\bf{p}_1,1,w_1} =
\id_{G_{1,w_1}}$, so the above equation tells us that the subgroup
$L_{1,(w_1,w_2,w_3)} \leq G_{1,w_1}$ is invariant under
$(\id_{W_1}\times T^{\bf{p}_1 - \bf{p}_2}\times T^{\bf{p}_1 -
\bf{p}_3})$.  Since we already know that it is $\vec{T}$-invariant
(since this holds for $M_\bullet$),
Lemma~\ref{lem:invariances-leave-first-coord} tells us that
$L_{1,(w_1,w_2,w_3)}$ virtually depends only on $w_1$, as required.

On the other hand, for any $\bf{m}$ we can set $\bf{n} := \bf{m} +
\bf{p}_1$ and find that the above equation expresses precisely the
condition that follows from the Relative Automorphism Structure
Theorem~\ref{thm:RAST} for
$(T^\bf{m})|_{\zeta^{T^{\bf{p}_1}}_{1/\a_1}}$ to respect the factor
corresponding to fibrewise quotienting by $L_{1,(w_1,w_2,w_3)}$
(which we have just seen virtually depends only on $w_1$).  This
completes the proof. \qed

\textbf{Proof of Proposition~\ref{prop:somewhat-better-Mackey}}\quad
If $\bfX$ is FIS$^+$ and we coordinatize
\begin{center}
$\phantom{i}$\xymatrix{
(\bfZ_{1/\a_i}^{T^{\bf{p}_i}})^{\,\uhr\bf{p}_i}\ar[dr]_{\a_i|_{\zeta_{1/\a_i}^{T^{\bf{p}_i}}}}\ar@{<->}[rr]^-\cong & & \bfW_i^{\,\,\uhr\bf{p}_i}\ltimes (G_{i,\bullet},m_{G_{i,\bullet}},\s_i)\ar[dl]^{\rm{canonical}}\\
&\bfW_i^{\,\,\uhr\bf{p}_i} }
\end{center}
then Lemma~\ref{lem:2D-proj-full} tells us that the associated
joining Mackey group data $M_\bullet$ has full two-dimensional
projections, and hence by Lemma~\ref{lem:group-th-from-FurWei} each
$M_{(w_1,w_2,w_3)}$ takes the form
\begin{multline*}
\{(g_1,g_2,g_3) \in \vec{G}_{(w_1,w_2,w_3)}:\\
\Psi_{1,(w_1,w_2,w_3)}(g_1)\cdot \Psi_{2,(w_1,w_2,w_3)}(g_2)\cdot
\Psi_{3,(w_1,w_2,w_3)}(g_3) = 1_{A_{(w_1,w_2,w_3)}}\}
\end{multline*}
for some compact metrizable Abelian group data $A_{(w_1,w_2,w_3)}$
and continuous epimorphisms
$\Psi_{i,(w_1,w_2,w_3)}:G_{i,w_i}\longrightarrow A_{(w_1,w_2,w_3)}$.
Moreover, by taking $A_{(w_1,w_2,w_3)}$ to be itself the quotient
$\vec{G}_{(w_1,w_2,w_3)}/M_{(w_1,w_2,w_3)}$, it is clear that we may
take $A_{(w_1,w_2,w_3)}$ and $\Psi_{i,(w_1,w_2,w_3)}$ to depend
measurably on $M_{(w_1,w_2,w_3)}$, and hence to vary measurably with
$(w_1,w_2,w_3)$.

Now Lemma~\ref{lem:slices-invar} gives that the one-dimensional
slices,
\[L_{1,(w_1,w_2,w_3)} := \{g_1 \in G_{1,w_1}:\ (g_1,1_{G_{2,w_2}},1_{G_{3,w_3}})\in M_{(w_1,w_2,w_3)}\}\]
and similarly $L_{2,(w_1,w_2,w_3)}$ and $L_{3,(w_1,w_2,w_3)}$, are
normal, are virtually functions of $w_1$ (respectively, $w_2$ and
$w_3$) and that the factors of the restriction of $T^{\bf{p}_i}$
given by fibrewise quotienting by these measurably-varying normal
subgroups are actually factor maps for the whole $\bbZ^2$-action
$T$.  Writing $A_{i,w_i}$ for the resulting quotient fibre group
$G_{i,w_i}/L_{i,w_i}$ and observing from
Lemma~\ref{lem:group-th-from-FurWei} that these are Abelian, these
intermediate systems are in fact the minimal characteristic factors
$\bfY_i$ and can be located according to another commutative diagram
\begin{center}
$\phantom{i}$\xymatrix{
(\bfZ_{1/\a_i}^{T^{\bf{p}_i}})^{\ \uhr\bf{p}_i}\ar[d]\ar@{<->}[rr]^-\cong & & \bfW_i^{\ \uhr\bf{p}_i}\ltimes (G_{i,\bullet},m_{G_{i,\bullet}},\s_i)\ar[d]^{\rm{canonical}}\\
\bfY_i^{\ \uhr\bf{p}_i}\ar[dr]\ar@{<->}[rr]^-\cong & & \bfW_i^{\ \uhr\bf{p}_i}\ltimes (A_{i,\bullet},m_{A_{i,\bullet}},\s'_i := \s_i\cdot L_{i,\bullet})\ar[dl]^{\rm{canonical}}\\
&\bfW_i^{\ \uhr\bf{p}_i}. }
\end{center}

Finally, observe from the definition of $L_{i,\bullet}$ that the
epimorphisms $\Psi_{i,\bullet}$ must factorize to give continuous
isomorphisms $\Theta_{i,(w_1,w_2,w_3)}:A_{i,w_i}\longrightarrow
A_{(w_1,w_2,w_3)}$ almost everywhere, and so it follows that
\[A_{1,w_1}\cong A_{2,w_2} \cong A_{3,w_3}\]
for almost all $(w_1,w_2,w_3)$ by some measurably-varying continuous
isomorphisms. On the other hand, $A_{i,\bullet}$ is also
$T^{\bf{p}_i}$-invariant, and so since the factors
$\zeta_0^{T^{\bf{p}_1}}\circ\pi_1$ and
$\zeta_0^{T^{\bf{p}_2}}\circ\pi_2$ of $\bfX^\rm{F}$ are relatively
independent over $\zeta_0^{T^{\bf{p}_1} = T^{\bf{p}_2} =
\id}\circ\pi_1\simeq \zeta_0^{T^{\bf{p}_1} = T^{\bf{p}_2} =
\id}\circ\pi_2$, it follows that we can adjust $A_{1,w_1}$ by a
measurably-varying family of continuous isomorphisms so that (up to
a negligible set) it depends only on $\zeta_0^{T^{\bf{p}_1} =
T^{\bf{p}_2} = \id}|_{\a_1}(w_1)$, and similarly for $A_{2,w_2}$ and
$A_{3,w_3}$.

To finish the proof we need only show that even this can be reduced
to a dependence only on $\zeta_0^T$.  This now follows because the
extension $\zeta_0^{T^{\bf{p}_1} = T^{\bf{p}_2} = \id} \gtrsim
\zeta_0^T$ is effectively a relatively ergodic extension of actions
of the finite group $\bbZ^2/(\bbZ\bf{p}_1 + \bbZ\bf{p}_2)$ with the
base action trivial, and so rather trivial application of the
non-ergodic Furstenberg-Zimmer Theory shows that each fibre of this
extension is a finite set, and that the transformations
$T^{\bf{n}}|_{\zeta_0^{T^{\bf{p}_1} = T^{\bf{p}_2} = \id}}$ simply
permute transitively the finitely many points of each fibre. Lifting
this picture, we see that the fibrewise actions of the
transformations $T^{\bf{n}}|_{\xi_i}$ must implicitly give
isomorphisms between each of the (finitely many) groups appearing as
$A_{1,w_1}$ for $w_1$ in a given fibre over $\zeta_0^T$, and so all
these compact Abelian groups coming from the same fibre are
isomorphic and these isomorphisms may be chosen measurably (since
there are only finitely many of them in question).  Therefore one
further re-coordinatization leads to $A_{1,w_1} = A_{2,w_2} =
A_{3,w_2} = A_s$ for some motionless data $A_\star$ and $s =
\zeta_0^T|_{\a_1}(w_1) = \zeta_0^T|_{\a_2}(w_2) =
\zeta_0^T|_{\a_3}(w_3)$, completing the proof of
Proposition~\ref{prop:somewhat-better-Mackey}. \qed

\textbf{Remark}\quad The main result we are working towards,
Theorem~\ref{thm:char-three-lines-in-2D}, itself tells us that for
an arbitrary system the characteristic factor $\bfY_i$ can
eventually be expressed as a subjoining from $\sfZ_0^{\bf{p}_i}$,
$\sfZ_0^{\bf{p}_i - \bf{p}_j}$, $\sfZ_0^{\bf{p}_i - \bf{p}_k}$ and
$\sfZ_{\Ab,2}^{\bf{p}_i,\bf{p}_i - \bf{p}_j,\bf{p}_i - \bf{p}_k}$ (two-step Abelian distal systems), and any joining
from these classes can easily be shown to have a further extension
that is simply an Abelian isometric extension of a
$(\sfZ_0^{\bf{p}_i}\vee \sfZ_0^{\bf{p}_i - \bf{p}_j}\vee
\sfZ_0^{\bf{p}_i - \bf{p}_k})$-system.  Intuitively, this suggests
that it should be possible to prove Abelianness of the
coordinatizing fibres of $\bfY_i\to \bfW_i$ after making only the
FIS assumption. Indeed, that implication could fail only if a system
could be found for which the coordinatizing fibres of $\bfY_i\to
\bfW_i$ are nontrivial homogeneous spaces
$G_{i,\bullet}/H_{i,\bullet}$, but such that to produce a further
nontrivial joining with a $(\sfZ_0^{\bf{p}_i}\vee \sfZ_0^{\bf{p}_i -
\bf{p}_j}\vee \sfZ_0^{\bf{p}_i - \bf{p}_k})$-system really requires
that we also involve a system from class
$\sfZ_{\Ab,2}^{\bf{p}_i,\bf{p}_i - \bf{p}_j,\bf{p}_i - \bf{p}_k}$, for
which the fibres over the Kronecker factor (which is always another
$(\sfZ_0^{\bf{p}_i},\sfZ_0^{\bf{p}_i - \bf{p}_j},\sfZ_0^{\bf{p}_i -
\bf{p}_k})$-subjoining) are Abelian. Presumably this would require
in turn that the Abelian fibres of the latter correspond to closed
Abelian subgroups $A_{i,\bullet} \leq G_{i,\bullet}$ with the
property that $A_{i,\bullet}H_{i,\bullet} = G_{i,\bullet}$ --- it is
this that would prevent the existence of a nontrivial joining to a
$(\sfZ_0^{\bf{p}_i}\vee \sfZ_0^{\bf{p}_i - \bf{p}_j}\vee
\sfZ_0^{\bf{p}_i - \bf{p}_k})$-system without also involving a
$\sfZ_{\Ab,2}^{\bf{p}_i,\bf{p}_i - \bf{p}_j,\bf{p}_i -
\bf{p}_k}$-system, because the whole extension $\bfY_i\to \bfW_i$
would still be relatively independent from the newly-adjoined system
even if this latter only failed to capture the subgroups
$A_{i,\bullet}$. This possibility seems remote, but I have not been
able to rule it out, and it seems to be rather easier to prove first
the abstract existence of FIS$^+$ extensions as in
Subsection~\ref{subs:fibre-normal} and then enjoy the simplification
of working with groups in places of homogeneous spaces above. \fin

\subsection{A zero-sum form for the joining Mackey group}

If we could take the isomorphisms $\Theta_{i,(w_1,w_2,w_3)}$
obtained in Proposition~\ref{prop:somewhat-better-Mackey} to depend
only on $w_i$, then we could simply use them to make one last
recoordinatization of the extensions $\bfY_i \to \bfW_i$ to complete
the proof of Proposition~\ref{prop:better-Mackey}.  I have not been
able to prove that this is possible in general, and here we will go around this problem by passing to a
further extension.

We begin this step with a few quite general lemmas.

\begin{lem}[Virtual isometricity implies isometricity]\label{lem:virtually-isom-implies-isom}
Suppose that $\L \leq \bbZ^d$ is a finite-index subgroup and that
$\pi:\bfX\to\bfY$ is an extension of $\bbZ^d$-systems such that the
extension of subactions $\pi:\bfX^{\ \uhr\L}\to \bfY^{\ \uhr\L}$ is
relatively ergodic and Abelian isometric with coordinatization
\begin{center}
$\phantom{i}$\xymatrix{ \bfX^{\
\uhr\L}\ar[dr]_{\pi}\ar@{<->}[rr]^-\cong & & \bfY^{\ \uhr\L}\ltimes
(A_\bullet,m_{A_\bullet},\s)\ar[dl]^{\rm{canonical}}\\
& \bfY^{\ \uhr\L}. }
\end{center}
Then there is some recoordinatization by an $S^{\ \uhr\L}$-invariant
measurable family of fibrewise automorphisms so that the whole
extension of $\bbZ^d$-actions can be coordinatized as isometric with
this compact Abelian group data.
\end{lem}

\textbf{Proof}\quad This is an easy consequence of the Relative
Automorphism Structure Theorem~\ref{thm:RAST}.  Applying that
theorem to the action $T$ regarded as itself an automorphic
$\bbZ^d$-action on the extension of the $\L$-subactions, we see that
the coordinatization of $T^{\ \uhr\L}$ as $S^{\ \uhr\L}\ltimes \s$
implies a coordinatization of $T^{\bf{n}}$ for each $\bf{n} \in
\bbZ^d$ as $S^{\bf{n}}\ltimes (L_{\rho_{\bf{n}}(\bullet)}\circ
\Phi_{\bf{n},\bullet})$ for some sections $\rho_\bf{n}:Y\to
A_\bullet$ and some measurable families of fibre-isomorphisms
$\Phi_{\bf{n},\bullet}:A_\bullet\to A_{S^{\bf{n}}(\bullet)}$. In
addition, each family of isomorphisms $\Phi_{\bf{n},\bullet}$ is
$S^{\ \uhr\L}$-invariant.

Of course, we must have $\rho_\bf{n} = \s(\bf{n},\,\cdot\,)$ and
$\Phi_{\bf{n},\bullet}\equiv \id_{A_\bullet}$ whenever $\bf{n}\in
\L$. Now consider the further factors
\[\bfY \stackrel{\zeta_0^{S^{\ \uhr\L}}}{\longrightarrow} \bfZ_0^{S^{\ \uhr\L}}\stackrel{\zeta_0^S|_{\zeta_0^{S^{\ \uhr\L}}}}{\longrightarrow}\bfZ_0^S.\]
Since the extension $\pi$ is relatively ergodic for the
$\L$-subactions, the compositions with $\pi$ of these two isotropy
factors coincide with those of the larger system $\bfX$. Also, the
restriction $S|_{\zeta_0^{S^{\ \uhr\L}}}$ can be identified with an
action of the finite quotient group $\bbZ^d/\L$ that is relatively
ergodic for the further factor map $\zeta_0^S|_{\zeta_0^{S^{\
\uhr\L}}}$, and so $S$ is actually transitive within almost all of
the fibres of $\zeta_0^S|_{\zeta_0^{S^{\ \uhr\L}}}$, which are
therefore identified as homogeneous spaces of this finite quotient
group. It follows that for almost every $s \in Z_0^S$, for almost
all pairs of points $y_1,y_2 \in (\zeta_0^S)^{-1}\{s\}$ there is
some $\bf{n}$ for which $\Phi_{\bf{n},y_1}$ carries $A_{y_1}$ (which
actually depends only on $\zeta_0^{S^{\ \uhr\L}}(y_1)$)
isomorphically onto $A_{y_2}$.  Therefore letting $\eta:Z_0^S\to
Z_0^{S^{\ \uhr\L}}$ be a measurable selector (Theorem~\ref{thm:meas-select}) for
$\zeta_0^S|_{\zeta_0^{S^{\ \uhr\L}}}$ we can make a simple
automorphism recoordinatization within each fibre of
$\zeta_0^S|_{\zeta_0^{S^{\ \uhr\L}}}$ to replace $A_\bullet$ with
$A_{\eta(\zeta_0^S(\bullet))}$, and hence assume that the group data
$A_\bullet$ is actually $S$-invariant and that
$\Phi_{\bf{n},\bullet}$ forms an $\rm{Aut}(A_\bullet)$-valued
cocycle-section.

Let us now write $R := S|_{\zeta_0^{S^{\ \uhr\L}}}$ and $\zeta :=
\zeta_0^S|_{\zeta_0^{S^{\ \uhr\L}}}$ for brevity and regard
$A_\bullet$ and each $\Phi_{\bf{n},\bullet}$ as a function defined
on $Z_0^{S^{\ \uhr\L}}$ rather than $Y$ (as we may by the
invariances established above). In this notation the trivial
requirement that $T$ commute with $T^{\ \uhr\L}$ shows that in fact
we must have $\Phi_{\bf{n},z} = \id_{A_z}$ whenever $R^{\bf{n}}(z) =
z$. We complete the proof by showing that there is a measurable
family $z\mapsto \Theta_z:Z_0^{S^{\ \uhr\L}}
\mapsto\rm{Aut}(A_\bullet)$ such that
\[\Theta_{R^{\bf{n}}(z)}\circ \Phi_{\bf{n},z}\circ \Theta_z^{-1} = \id_{A_z}\]
for $(\zeta_0^{S^{\ \uhr\L}})_\#\nu$-almost every $z\in Z_0^{S^{\
\uhr\L}}$ for every $\bf{n} \in \bbZ^d$. Let $z\mapsto \bf{m}(z)
\in\bbZ^d$ be a measurable selection such that $R^{\bf{m}(z)}(z) =
\eta(\zeta(z))$ (again, this is clearly possible from the
transitivity of $R$ on the fibres of $\zeta$), and now set
\[\Theta_z := \Phi_{\bf{m}(z),z}.\] We can compute from the fact that
$\Phi_{\cdot,\bullet}$ is an $\rm{Aut}(A_\bullet)$-valued
cocycle-section that
\begin{eqnarray*}
\Theta_{R^\bf{n}(z)}\circ \Phi_{\bf{n},z}\circ \Theta_z^{-1} &=&
\Phi_{\bf{m}(R^{\bf{n}}(z)),R^{\bf{n}}(z)}\circ
\Phi_{\bf{n}-\bf{m}(z),R^{\bf{m}(z)}(z)}\circ
\Phi_{\bf{m}(z),z}\circ (\Phi_{\bf{m}(z),z})^{-1}\\
&=&
\Phi_{\bf{m}(R^{\bf{n}}(z)),R^{\bf{n}-\bf{m}(z)}(R^{\bf{m}(z)}z)}\circ
\Phi_{\bf{n}-\bf{m}(z),R^{\bf{m}(z)}(z)}\\
&=& \Phi_{\bf{m}(R^{\bf{n}}(z))+\bf{n} - \bf{m}(z),R^{\bf{m}(z)}(z)}
= \id_{A_z},
\end{eqnarray*}
because
\begin{multline*}
R^{\bf{m}(R^{\bf{n}}(z))+\bf{n} - \bf{m}(z)}(R^{\bf{m}(z)}z) =
R^{\bf{m}(R^{\bf{n}}(z))+\bf{n}}(z)\\ =
R^{\bf{m}(R^{\bf{n}}(z))}(R^{\bf{n}}(z)) = \eta(\zeta(z)) =
R^{\bf{m}(z)}(z)
\end{multline*}
so the last cocycle appearing above must be
trivial. \qed

The r\^ole of the following lemma will be somewhat analogous to that
of Lemma~\ref{lem:slices-respected} in the previous subsection.

\begin{lem}\label{lem:equal-Mackeys-to-equal-isos}
Suppose that $A$ is a compact Abelian group and that
\[M_i = \{(a_1,a_2,a_3)\in A^3:\ \Theta_{i,1}(a_1)\cdot\Theta_{i,2}(a_2)\cdot\Theta_{i,3}(a_3) = 1_A\}\]
for $i=1,2$ are subgroups of $A^3$ with full two-dimensional
projections and trivial one-dimensional slices, and suppose also
that $\Phi_j:A\stackrel{\cong}{\longrightarrow} A$ for $j=1,2,3$ and
$(b_1,b_2,b_3)\in A^3$ are such that
\[(b_1,b_2,b_3)\cdot (\Phi_1\times \Phi_2\times \Phi_3)(M_1) = M_2.\]
Then
\[\Theta_{1,1}\circ\Phi^{-1}_1\circ\Theta_{2,1}^{-1} = \Theta_{1,2}\circ\Phi^{-1}_2\circ\Theta_{2,2}^{-1} = \Theta_{1,3}\circ\Phi^{-1}_3\circ\Theta_{2,3}^{-1}.\]
\end{lem}

\textbf{Proof}\quad First the condition that $(b_1,b_2,b_3)\cdot
(\Phi_1(1_A),\Phi_2(1_A),\Phi_3(1_A))\in M_2$ simplifies to
$(b_1,b_2,b_3) \in M_2$, and so we can multiply the given equation
by $(b_1,b_2,b_3)^{-1}$ to obtain simply
\[(\Phi_1\times \Phi_2\times \Phi_3)(M_1) = M_2.\]
We can now write this out more explicitly as
\begin{multline*}
\Theta_{1,1}(\Phi_1^{-1}(a_1))\cdot
\Theta_{1,2}(\Phi_2^{-1}(a_2))\cdot \Theta_{1,3}(\Phi_3^{-1}(a_3)) =
1_A\\ \Leftrightarrow\quad\quad \Theta_{2,1}(a_1)\cdot
\Theta_{2,2}(a_2)\cdot \Theta_{2,3}(a_3) = 1_A
\end{multline*}
for all $(a_1,a_2,a_3) \in A^3$.

Restricting first to the special case $a_1 = 1_A$ this now
re-arranges to give
\[\Theta_{1,3}(\Phi_3^{-1}(\Theta_{2,3}^{-1}(\Theta_{2,2}(a_2)))) = \Theta_{1,2}(\Phi_2^{-1}(a_2))\quad\quad\forall a_2\in A\]
and hence
\[\Theta_{1,3}\circ\Phi_3^{-1}\circ\Theta_{2,3}^{-1} = \Theta_{1,2}\circ\Phi_2^{-1}\circ\Theta_{2,2}^{-1},\]
and arguing similarly with $a_3 = 1_A$ shows that these are both
also equal to $\Theta_{1,1}\circ\Phi_1^{-1}\circ\Theta_{2,1}^{-1}$,
as required. \qed

A similar argument gives the forward implication of the following
lemma, while the reverse implication is an immediate check.

\begin{lem}\label{lem:when-good-subgroups-equal}
The groups $M_1$ and $M_2$ of the previous lemma are equal if and
only if
\[\Theta_{1,1}\circ\Theta_{2,1}^{-1} = \Theta_{1,2}\circ\Theta_{2,2}^{-1} = \Theta_{1,3}\circ\Theta_{2,3}^{-1}.\]
\qed
\end{lem}

Using the above results we can now show that, having once found the
extension $\bfY_i\to \bfW_i$ and the coordinatization of its
$(\bbZ\bf{p}_i)$-subaction promised by
Proposition~\ref{prop:somewhat-better-Mackey}, then after adjoining
a new $\sfW_i$-system if necessary we can render this extension
Abelian isometric for the whole $\bbZ^2$-action.

\begin{lem}[Making all transformations isometric]\label{lem:everybody-becomes-isometric}
Let $i\in \{1,2,3\}$ and $\sfW_i$ be the idempotent class
$\sfZ_0^{\bf{p}_i}\vee\sfZ_0^{\bf{p}_i -
\bf{p}_j}\vee\sfZ_0^{\bf{p}_i - \bf{p}_k}$. In the notation of
Proposition~\ref{prop:somewhat-better-Mackey}, any FIS$^+$
$\bbZ^2$-system $\bfX$ admits an FIS$^+$ extension
$\pi:\t{\bfX}\to\bfX$ such that we can coordinatize
\begin{center}
$\phantom{i}$\xymatrix{ (\zeta_{\sfW_i}^{\t{\bfX}}\vee
(\xi_i\circ\pi))(\tilde{\bfX})\ar[dr]_{\zeta_{\sfW_i}^{\t{\bfX}}|_{\zeta_{\sfW_i}^{\t{\bfX}}\vee
(\xi_i\circ\pi)}}\ar@{<->}[rr]^-\cong & & \t{\bfW_i}\ltimes
(A_\star,m_{A_\star},\s_i)\ar[dl]^{\rm{canonical}}\\
& \t{\bfW_i} }
\end{center}
for some compact Abelian group data $A_\star$ and cocycle sections
$\s_i:\bbZ^2\times \t{W_i}\to A_\star$.
\end{lem}

\textbf{Remarks}\quad It is very important to bear in mind that this
result gives an extension of $\bfX$ such that the extension
$\xi_i:\bfY_i\to \bfW_i$ may be lifted and then usefully
re-coordinatized for each $i$ \emph{separately}.  In general it
seems that the joint distribution of the systems $\bfY_1$, $\bfY_2$
and $\bfY_3$ as factors of the single system $\bfX$ can be extremely
complicated, and here we make no requirement that the
re-coordinatizations we obtain should enjoy any `consistency' in
terms of this joint distribution. \fin

\textbf{Proof}\quad By symmetry we may assume $i=1$. Let $\bfX'$ be
the extension of $\bfX^{\ \uhr(\bbZ\bf{p}_1 + \bbZ\bf{p}_2)}$ with
underlying space
\[(X\times W_2\times W_3,(\id_X\times \a_2\times \a_3)_\#\mu^\rm{F}),\]
with factor map $\pi'$ onto $X$ given by $\pi_1$ and with lifted
transformations \[(T')^{\bf{p}_1} := T^{\bf{p}_1}\times
T|_{\a_2}^{\bf{p}_2}\times T|_{\a_3}^{\bf{p}_3}\] and
\[(T')^{\bf{p}_2} := T^{\bf{p}_2}\times T|_{\a_2}^{\bf{p}_2}\times T|_{\a_3}^{\bf{p}_2}\]
(in what follows we could have exchanged the roles of $\bf{p}_2$ and
$\bf{p}_3$ in the above construction).  Now let $\pi:\t{\bfX}\to
\bfX$ be a further extension recovering an action of the whole of
$\bbZ^2$ (for example, an FP extension of $\bfX'$ as in Subsection
3.2 of~\cite{Aus--lindeppleasant1}):
\begin{center}
$\phantom{i}$\xymatrix{ \t{\bfX}^{\ \uhr(\bbZ\bf{p}_1 +
\bbZ\bf{p}_2)}\ar[rr]^{\pi}\ar[dr]& & \bfX^{\ \uhr(\bbZ\bf{p}_1 +
\bbZ\bf{p}_2)}\\ & \bfX'\ar[ur]_{\pi'}. }
\end{center}

Write $\t{\bfW}_1 := \sfW_1\t{\bfX}$ and let $A_\star$, $\s_i$ and
$\Theta_{i,\bullet}$ be as given by
Proposition~\ref{prop:somewhat-better-Mackey}. Now consider the new
extension
\[\zeta_{\sfW_1}^{\t{\bfX}}|_{\zeta_{\sfW_1}^{\t{\bfX}}\vee
(\xi_1\circ\pi)}:(\zeta_{\sfW_1}^{\t{\bfX}}\vee
(\xi_1\circ\pi))(\tilde{\bfX}) \to \t{\bfW_1}.\] This is obtained
from the original extension $\xi_1:\bfY_1\to \bfW_1$ by adjoining
the new $\sfW_1$-system $\t{\bfW}_1$ overall, and from the
$\sfW_1$-satedness of $\bfX$ we know that this adjoining is
relatively independent from $\bfY_1$ over $\bfW_1$, so that the
above enlarged extension has the same fibres $A_\star$ as the
original extension. From the way we have constructed $\t{\bfW}_1$ in
terms of the Furstenberg self-joining we see that we may insert the
extension of $(\bbZ\bf{p}_1 + \bbZ\bf{p}_2)$-systems
\begin{multline*}
\vec{\a}|_{\xi_1\times\a_1\times\a_2}:(\xi_1\times\id_{W_2}\times\id_{W_3})(\bfX')
= (Y_1\times W_2\times W_3,\mu',T')\\ \to
(W,\vec{\a}_\#\mu^\rm{F},T'|_{\vec{\a}})
\end{multline*}
into a commutative diagram of factors of $\t{\bfX}^{\
\uhr(\bbZ\bf{p}_1 + \bbZ\bf{p}_2)}$ as follows:
\begin{center}
$\phantom{i}$\xymatrix{ (\zeta_{\sfW_1}^{\t{\bfX}}\vee
(\xi_1\circ\pi))(\tilde{\bfX})^{\ \uhr(\bbZ\bf{p}_1 +
\bbZ\bf{p}_2)}\ar[d]\ar[rr]^-{\zeta_{\sfW_1}^{\t{\bfX}}|_{\zeta_{\sfW_1}^{\t{\bfX}}\vee
(\xi_1\circ\pi)}} & & \t{\bfW_1}^{\ \uhr(\bbZ\bf{p}_1 +
\bbZ\bf{p}_2)}\ar[d]\\
(\xi_1\times \id_{W_2}\times
\id_{W_3})(\bfX')\ar[d]\ar[rr]^-{\vec{\a}|_{\xi_1\times \a_2\times
\a_3}} & &
(W,\vec{\a}_\#\mu^\rm{F},T'|_{\vec{\a}})\ar[d]\\
\bfY_1^{\ \uhr(\bbZ\bf{p}_1 + \bbZ\bf{p}_2)}\ar[rr]_{\xi_1} & &
\bfW_1^{\ \uhr(\bbZ\bf{p}_1 + \bbZ\bf{p}_2)}. }
\end{center}

Appealing again to $\sfW_1$-satedness, each of the horizontal
extensions in this diagram inherits a coordinatization in terms of
$A_\star$ and $\s_1$ from the coordinatization of the bottom row. We
need to show that we can trivialize the isomorphism sections
associated with the restrictions of each $\t{T}^{\bf{n}}$ to the
extension of the top row. Letting $\L := \bbZ\bf{p}_1 + \bbZ\bf{p}_2
\leq \bbZ^2$ and noting that $\sfZ_0^\L \leq \sfW_1$, we deduce also
from $\sfW_1$-satedness that the above horiztonal extensions are all
still relatively ergodic for the $\L$-subactions, and now by
Lemma~\ref{lem:virtually-isom-implies-isom} it will suffice to
trivialize the isomorphism sections associated with $\t{T}^{\bf{n}}$
for $\bf{n} \in \L$.  This, in turn, may be done for the extension
of the middle row of the above diagram instead, since then lifting
the $(T')^{\bf{p}_1}$-invariant measurable family of fibrewise
automorphisms that we use to the top row completes the proof.

On the middle-row extension $\vec{\a}|_{\xi_1\times \a_2\times
\a_3}$ we can re-coordinatize the fibre-copies of $A_\star$ by the
fibrewise automorphisms
$\Theta_{2,\bullet}^{-1}\circ\Theta_{1,\bullet}$ (recalling that
this is a function of $(w_1,w_2,w_3)\in W$). We will show that this
trivializes the relevant isomorphism sections using the existence of
the additional commuting transformations $(T^\bf{n})^{\times 3}$ on
$(X^3,\mu^\rm{F})$. The Relative Automorphism Structure
Theorem~\ref{thm:RAST} tells us that for each $\bf{n} \in \L$ and
$i=1,2,3$ we can coordinatize
\[T^\bf{n}|_{\xi_i} \cong T^{\bf{n}}|_{\a_i}\ltimes (L_{\rho_{\bf{n},i}(\bullet)}\circ\Psi_{\bf{n},i,\bullet}),\]
and in case $i=1$ this coordinatization can be lifted to give
\[\tilde{T}^\bf{n}|_{\zeta_{\sfW_1}^{\tilde{\bfX}}\vee \xi_1} \cong \tilde{T}^{\bf{n}}|_{\zeta_{\sfW_1}^{\tilde{\bfX}}}\ltimes (L_{\t{\rho}_{\bf{n},1}(\bullet)}\circ\t{\Psi}_{\bf{n},1,\bullet})\]
with $\t{\rho}_{\bf{n},1} = \rho_{\bf{n},1}\circ
(\a_1\circ\pi)|_{\zeta_{\sfW_1}^{\t{\bfX}}}$ and similarly for
$\t{\Phi}_{\bf{n},1,\bullet}$. Now the condition that
$(T^\bf{n})^{\times 3}$ respect $\zeta_0^{\vec{T}}$ as a factor map
gives that for $\vec{\a}_\#\mu^\rm{F}$-almost every
$(w_1,w_2,w_3)\in W$, writing $s:= \zeta_0^{T}|_{\a_1}(w_1)$ as
before, for every $(a_1',a_2',a_3') \in A_s^3$ there is
$(a_1'',a_2'',a_3'') \in A_s^3$ such that
\begin{multline*}
(a_1',a_2',a_3')\cdot
(\rho_{\bf{n},1}(w_1),\rho_{\bf{n},2}(w_2),\rho_{\bf{n},3}(w_3))\cdot
(\Psi_{\bf{n},1,w_1}\times \Psi_{\bf{n},2,w_2}\times
\Psi_{\bf{n},3,w_3})(M_{(w_1,w_2,w_3)})\\
= (a_1'',a_2'',a_3'')\cdot
M_{(T|_{\a_1}^{\bf{n}}(w_1),T|_{\a_2}^{\bf{n}}(w_2),T|_{\a_3}^{\bf{n}}(w_3))}.
\end{multline*}

Applying Lemma~\ref{lem:equal-Mackeys-to-equal-isos} when $\bf{n} =
\bf{p}_2$ (and recalling that $\Psi_{\bf{p}_2,2,\bullet}\equiv
\id_{A_\star}$) now gives
\begin{multline*}
\Theta_{1,(w_1,w_2,w_3)}\circ\Psi_{\bf{p}_2,1,w_1}^{-1}\circ
\Theta_{1,(T|_{\a_1}^{\bf{p}_2}(w_1),T|_{\a_2}^{\bf{p}_2}(w_2),T|_{\a_3}^{\bf{p}_2}(w_3))}^{-1}\\
=\Theta_{2,(w_1,w_2,w_3)}\circ
\Theta_{2,(T|_{\a_1}^{\bf{p}_2}(w_1),T|_{\a_2}^{\bf{p}_2}(w_2),T|_{\a_3}^{\bf{p}_2}(w_3))}^{-1}
\end{multline*}
and hence
\begin{multline*}
\Theta_{2,(T|_{\a_1}^{\bf{p}_2}(w_1),T|_{\a_2}^{\bf{p}_2}(w_2),T|_{\a_3}^{\bf{p}_2}(w_3))}^{-1}\circ\Theta_{1,(T|_{\a_1}^{\bf{p}_2}(w_1),T|_{\a_2}^{\bf{p}_2}(w_2),T|_{\a_3}^{\bf{p}_2}(w_3))}\\
\circ\Psi_{\bf{p}_2,1,w_1}\circ
(\Theta_{2,(w_1,w_2,w_3)}^{-1}\circ\Theta_{1,(w_1,w_2,w_3)})^{-1} =
\id_{A_\star}
\end{multline*}
$\a_\#\mu^\rm{F}$-almost surely.

This implies that upon re-coordinatizing the fibre copies of
$A_\star$ by $\Theta_{2,\bullet}^{-1}\circ\Theta_{1,\bullet}$ the
family of isomorphisms $\Psi_{\bf{p}_2,1,\bullet}$ trivializes.
Since the re-coordinatizing fibrewise isomorphisms are invariant for
the restriction of $\t{T}^{\bf{p}_1}$, under the new
coordinatization that results the $(\bbZ\bf{p}_1)$-subaction is also
still coordinatized simply by an $A_\star$-valued cocycle-section,
and so we have obtained isometricity for the whole $\L$-subaction,
as desired. \qed

\begin{cor}\label{cor:somewhat-straightened-Mackeygroup}
Let $\sfW_i$ be the idempotent class
$\sfZ_0^{\bf{p}_i}\vee\sfZ_0^{\bf{p}_i -
\bf{p}_j}\vee\sfZ_0^{\bf{p}_i - \bf{p}_k}$. In the notation of
Proposition~\ref{prop:somewhat-better-Mackey}, any FIS$^+$
$\bbZ^2$-system $\bfX$ admits an FIS$^+$ extension
$\pi:\t{\bfX}\to\bfX$ such that we can coordinatize
\begin{center}
$\phantom{i}$\xymatrix{ (\zeta_{\sfW_i}^{\t{\bfX}}\vee
(\xi_i\circ\pi))(\tilde{\bfX})\ar[dr]_{\zeta_{\sfW_i}^{\t{\bfX}}|_{\zeta_{\sfW_i}^{\t{\bfX}}\vee
(\xi_i\circ\pi)}}\ar@{<->}[rr]^-\cong & & \t{\bfW_i}\ltimes
(A_\star,m_{A_\star},\s_i)\ar[dl]^{\rm{canonical}}\\
& \t{\bfW}_i }
\end{center}
for some compact Abelian group data $A_\star$ and cocycle sections
$\s_i:\bbZ^2\times \t{W_i}\to A_\star$ so that the resulting Mackey
group data for the joining of the above extensions under
$\t{\mu}^\rm{F}$ is
\[M_\bullet \cong \{(a_1,a_2,a_3)\in A_\star^3:\ a_1\cdot a_2\cdot a_3 = 1_{A_\star}\}.\]
$\vec{\a}_\#\t{\mu}^\rm{F}$-almost everywhere.
\end{cor}

\textbf{Remark}\quad Note that this is not yet a result describing
the overall joining Mackey group for the new system $\t{\bfX}$, but
only the Mackey group data for the joining $\t{\mu}^\rm{F}$
restricted to the subextensions of $\t{\bfY}_i\to \t{\bfW}_i$
obtained from $(\zeta_{\sfW_i}^{\t{\bfX}}\vee
(\xi_i\circ\pi))(\tilde{\bfX})$. \fin

\textbf{Proof}\quad For each $i=1,2,3$ let $\pi_{(i)}:\bfX_{(i)}\to
\bfX$ be an extension as given by
Lemma~\ref{lem:everybody-becomes-isometric}, and now let $\bfX'\to
\bfX$ be the relatively independent product of the extensions
\begin{center}
$\phantom{i}$\xymatrix{ \bfX_{(1)}\ar[dr]_{\pi_{(1)}} &
\bfX_{(2)}\ar[d]^{\pi_{(2)}} & \bfX_{(3)}\ar[dl]^{\pi_{(3)}}\\ &
\bfX }
\end{center}
and let $\t{\bfX}\to \bfX'$ any FIS$^+$ extension of this to give
the overall extension $\pi:\t{\bfX}\to \bfX$ by composition.

It is clear that the isometricity of the whole $\bbZ^2$-actions
obtained in Lemma~\ref{lem:everybody-becomes-isometric} persists
under passing to a further extension such as $\t{\bfX}$, since by
$\sfW_i$-satedness we may simply lift the group and cocycle data
describing the extensions
$\zeta_{\sfW_i}^{\bfX_{(i)}}|_{\zeta_{\sfW_i}^{\bfX_{(i)}}\vee
(\xi_i\circ\pi_{(i)})}$ further to give a coordinatization of
$\zeta_{\sfW_i}^{\t{\bfX}}|_{\zeta_{\sfW_i}^{\t{\bfX}}\vee
(\xi_i\circ\pi)}$.

We will deduce that $\t{\bfX}$ admits the desired simple form for
$M_\bullet$ by using again the presence of the automorphisms
$\t{T}^\bf{n}$ of the Furstenberg self-joining $\t{\mu}^\rm{F}$. As
a result of the simple coordinatization of each
$\t{T}|_{\zeta^{\t{\bfX}}_{\sfW_i}\vee (\xi_i\circ\pi_i)}$ as
$\t{T}|_{\zeta^{\t{\bfX}}_{\sfW_i}} \ltimes \s_i$ obtained from
Lemma~\ref{lem:everybody-becomes-isometric}, the condition that each
$(\t{T}^\bf{n})^{\times 3}$ respect $\zeta_0^{\vec{\t{T}}}$ now
becomes that for $\vec{\a}_\#\t{\mu}^\rm{F}$-almost every
$(\t{w}_1,\t{w}_2,\t{w}_3)\in W$, writing $s:=
\zeta_0^{\t{T}}|_{\t{\a}_1}(\t{w}_1)$ as before, for every
$(a_1',a_2',a_3') \in A_s^3$ there is $(a_1'',a_2'',a_3'') \in
A_s^3$ such that
\begin{multline*}
(a_1',a_2',a_3')\cdot
(\s_1(\bf{n},\t{w}_1),\s_2(\bf{n},\t{w}_2),\s_3(\bf{n},\t{w}_3))\cdot M_{(\t{w}_1,\t{w}_2,\t{w}_3)}\\
= (a_1'',a_2'',a_3'')\cdot
M_{(\t{T}|_{\t{\a}_1}^{\bf{n}}(\t{w}_1),\t{T}|_{\t{\a}_2}^{\bf{n}}(\t{w}_2),\t{T}|_{\t{\a}_3}^{\bf{n}}(\t{w}_3))}.
\end{multline*}

Since $M_\bullet$ still has trivial slices and full two-dimensional
projections (indeed, it cannot be larger for $\t{\mu}^\rm{F}$ than
the joining Mackey group data for $\mu^\rm{F}$, and if it no longer
had full two-dimensional projections then we could derive a
contradiction with satedness just in Lemma~\ref{lem:2D-proj-full}),
we may invoke its representation in the form
\[M_\bullet = \{(a_1,a_2,a_3)\in A_\star^3:\ \Theta_{1,\bullet}(a_1)\cdot\Theta_{1,\bullet}(a_2)\cdot\Theta_{3,\bullet}(a_3) = 1_{A_\star}\},\]
and now apply Lemma~\ref{lem:equal-Mackeys-to-equal-isos} to deduce
that
\begin{multline*}
\Theta_{1,(\tilde{T}|_{\tilde{\a}_1}^{\bf{n}}(\t{w}_1),\tilde{T}|_{\tilde{\a}_2}^{\bf{n}}(\t{w}_2),\tilde{T}|_{\tilde{\a}_3}^{\bf{n}}(\t{w}_3))}\circ(\Theta_{1,(\t{w}_1,\t{w}_2,\t{w}_3)})^{-1}\\
=
\Theta_{2,(\tilde{T}|_{\tilde{\a}_1}^{\bf{n}}(\t{w}_1),\tilde{T}|_{\tilde{\a}_2}^{\bf{n}}(\t{w}_2),\tilde{T}|_{\tilde{\a}_3}^{\bf{n}}(\t{w}_3))}\circ(\Theta_{2,(\t{w}_1,\t{w}_2,\t{w}_3)})^{-1}\\
=
\Theta_{3,(\tilde{T}|_{\tilde{\a}_1}^{\bf{n}}(\t{w}_1),\tilde{T}|_{\tilde{\a}_2}^{\bf{n}}(\t{w}_2),\tilde{T}|_{\tilde{\a}_3}^{\bf{n}}(\t{w}_3))}\circ(\Theta_{3,(\t{w}_1,\t{w}_2,\t{w}_3)})^{-1}
\end{multline*}
$\vec{\tilde{\a}}_\#\mu^\rm{F}$-almost everywhere.  From this
Lemma~\ref{lem:when-good-subgroups-equal} gives
\[M_{(\tilde{T}|_{\tilde{\a}_1}^{\bf{n}}(\t{w}_1),\tilde{T}|_{\tilde{\a}_2}^{\bf{n}}(\t{w}_2),\tilde{T}|_{\tilde{\a}_3}^{\bf{n}}(\t{w}_3))} = M_{(\t{w}_1,\t{w}_2,\t{w}_3)}\]
$\vec{\tilde{\a}}_\#\mu^\rm{F}$-almost everywhere for every $\bf{n}
\in \bbZ^2$. Since $M_\bullet$ is already $\vec{\t{T}}$-invariant,
Lemma~\ref{lem:invariances-leave-first-coord} now gives that it is
virtually measurable with respect
$\zeta_0^{\t{T}^{\bf{p}_1}}\circ\pi_1$ and hence in fact with
respect to $\zeta_0^{\t{T}}\circ \pi_1 \simeq \zeta_0^{\t{T}}\circ
\pi_2\simeq \zeta_0^{\t{T}}\circ \pi_3$.  Therefore, in particular,
we can actually choose $\Theta_{i,\bullet}$ depending only on
$\zeta_0^{\t{T}}|_{\t{\a}_i}(\t{w}_i)$ to represent this Mackey
data. One further fibrewise recoordinatization by the
$\t{T}$-invariant automorphisms $\Theta_{i,\bullet}$ of $A_\star$,
which by $\t{T}$-invariance does not disrupt the coordinatization of
our extensions by $A_\star$-valued cocycle-sections, now clearly
straightens out the joining Mackey group completely to give the
desired zero-sum form
\[\{(a_1,a_2,a_3)\in A_\star^3:\ a_1\cdot a_2\cdot a_3 = 1_{A_\star}\}\]
everywhere. \qed

\textbf{Remark}\quad Notice that in the above proof, when we form a
$\sfW_i$-adjoining of a $\sfW_i$-sated system $\bfX$ this preserves
that instance of satedness, but will typically disrupt
$\sfW_j$-satedness for any other $j$.  After three different
extensions for $i=1,2,3$ we cannot be sure that our new larger
system retains any satedness (or, similarly, any fibre-normality),
hence our need to form another FIS$^+$ extension to recover these
valuable properties that we assumed initially. \fin

\textbf{Proof of Proposition~\ref{prop:better-Mackey}}\quad Let
$\bfX_{(0)} := \bfX_0$ and let $\psi^{(1)}_{(0)}:\bfX_{(1)}\to
\bfX_{(0)}$ be an FIS$^+$ extension. We will extend this to an
inverse sequence of FIS$^+$ systems $(\bfX_{(m)})_{m\geq 0}$,
$(\psi^{(m)}_{(k)})_{m\geq k\geq 0}$ and then show that the inverse
limit has the desired property.

Given $m\geq 1$ and $(\bfX_{(k)})_{m\geq k\geq 0}$,
$(\psi^{(k)}_{(\ell)})_{m\geq k\geq \ell\geq 0}$ we construct
$\psi^{(m+1)}_{(m)}:\bfX_{(m+1)}\to\bfX_{(m)}$ as follows. Since
$\bfX_{(m)}$ is FIS$^+$, by
Proposition~\ref{prop:somewhat-better-Mackey} we can choose
coordinatizations
\begin{center}
$\phantom{i}$\xymatrix{
\bfY_{(m),i}^{\ \uhr\bf{p}_i}\ar[dr]_{\a_{(m),i}|_{\xi_{(m),i}}}\ar@{<->}[rr]^-\cong & & \bfW_{(m),i}^{\ \uhr\bf{p}_i}\ltimes (A_{(m),\star},m_{A_{(m),\star}},\s_{(m),i})\ar[dl]^{\rm{canonical}}\\
&\bfW_{(m),i}^{\ \uhr\bf{p}_i} }
\end{center}
of the minimal characteristic factors $\xi_{(m),i}$, with associated
joining Mackey group
\begin{multline*}
M_{(m),(w_1,w_2,w_3)}\\ = \{(a_1,a_2,a_3)\in A_{(m),\star}^3:\
\Theta_{(m),1,\vec{w}}(a_1)\cdot\Theta_{(m),2,\vec{w}}(a_2)\cdot\Theta_{(m),3,\vec{w}}(a_3)
= 1_{A_{(m),\star}}\}.
\end{multline*}
Now let $\psi^{(m+1)}_{(m)}:\bfX_{(m+1)}\to \bfX_{(m)}$ be the
FIS$^+$ extension of $\bfX_{(m)}$ given by
Corollary~\ref{cor:somewhat-straightened-Mackeygroup}.

Having formed this inverse sequence, let $\bfX_{(\infty)}$,
$(\psi_{(m)})_{m\geq 0}$ be its inverse limit.  We will show this
has the desired properties.

We know that the minimal characteristic factors of $\bfX_{(\infty)}$
satisfy $\xi_{(\infty),i} \succsim \a_{(\infty),i}$. On the other
hand a simple check (see Lemma 4.4 in~\cite{Aus--lindeppleasant1})
shows that
\[\xi_{(\infty),i}\simeq \bigvee_{m\geq 1}\xi_{(m),i}\circ\psi_{(m)},\]
so by sandwiching we also have
\[\xi_{(\infty),i}\simeq \bigvee_{m\geq 1}(\a_{(\infty),i}\vee(\xi_{(m),i}\circ\psi_{(m)})).\]

Thus each $\xi_{(\infty),i}$ is generated by all the intermediate
factors
\[\xi_{(\infty),i}\succsim
(\a_{(\infty),i}\vee(\xi_{(m),i}\circ\psi_{(m)})) \succsim
\a_{(\infty),i}.\] Moreover,
Corollary~\ref{cor:somewhat-straightened-Mackeygroup} gives us a
coordinatization of the restriction of the whole $\bbZ^2$-action
$T_{(\infty)}$ to each
$(\a_{(\infty),i}\vee(\xi_{(m),i}\circ\psi_{(m)})) \succsim
\a_{(\infty),i}$ as an Abelian isometric extension, and so in fact
the restriction of $T_{(\infty)}$ is Abelian isometric for the whole
extension $\xi_{(\infty),i}\succsim \a_{(\infty),i}$.

Next, since each individual system $\bfX_{(m)}$ is FIS$^+$, we must
have that $\a_{(\infty),i}$ and $\psi_{(m)}$ are relatively
independent over $\a_{(m),i}\circ\psi_{(m)}$.  Therefore the
property that the Abelian extension $\xi_{(m),i}\succsim \a_{(m),i}$
can be `untwisted' when we lift to
$\a_{(m+1),i}\vee(\xi_{(m),i}\circ\psi^{(m+1)}_{(m)}) \succsim
\a_{(m+1),i}$ to have a coordinatization enjoying the simple
zero-sum form for its Mackey group data given by
Corollary~\ref{cor:somewhat-straightened-Mackeygroup} lifts to the
extensions $\a_{(\infty),i}\vee(\xi_{(m),i}\circ\psi_{(m)})\succsim
\a_{(\infty),i}$.

In terms of these data, the Relative Factor Structure
Theorem~\ref{thm:RFST} now gives us an explicit description of the
extension $\xi_{(\infty),i}\succsim \a_{(\infty),i}$ inside the
inverse limit: it tells us that for each $m \geq k\geq 0$ there is a
$T_{(\infty)}|_{\a_{(\infty),i}}$-invariant family of continuous
epimorphisms
$\Phi^{(m)}_{(k),i,\bullet}:A_{(m),\star}\longrightarrow
A_{(k),\star}$ on $W_{(\infty),i}$ such that the canonical factor
map from $\a_{(\infty),i}\vee(\xi_{(m),i}\circ\psi_{(m)})$ onto
$\a_{(\infty),i}\vee(\xi_{(k),i}\circ\psi_{(k)})$ is coordinatized
as
\[\phi^{(m)}_{(k)} = \id_{W_{(\infty),i}}\ltimes (L_{\rho_{(k),i}(\bullet)}\circ \Phi^{(m)}_{(k),i,\bullet}).\]
Combining these data now gives a coordinatization of
$\xi_{(\infty),i}\succsim \a_{(\infty),i}$ as
\begin{center}
$\phantom{i}$\xymatrix{
\bfY_{(\infty),i}\ar[dr]_{\a_{(\infty),i}|_{\xi_{(\infty),i}}}\ar@{<->}[rr]^-\cong
& & \bfW_{(\infty),i}\ltimes
(A_{(\infty),i,\star},m_{A_{(\infty),i,\star}},\s_{(\infty),i})\ar[dl]^{\rm{canonical}}\\
& \bfW_{(\infty),i}. }
\end{center}
with fibres the inverse limit groups
\[A_{(\infty),i,\star} := \lim_{m\leftarrow}\,(A_{(m),\phi_{(m)}(\star)},\Phi^{(m+1)}_{(m),i,\phi_{(m+1)}(\bullet)}),\]
which are still compact Abelian, and are invariant for the whole
action $T_{(\infty)}$ because this is so of the groups
$A_{(m),\phi_{(m)}(\star)}$ and the epimorphisms
$\Phi^{(m+1)}_{(m),i,\phi_{(m+1)}(\bullet)}$.  The cocycle
$\s_{(\infty),i}$ is given by the simultaneous lift to
$A_{(\infty),i,\star}$ of the sequence of cocycles
$(\s_{(m),i})_{m\geq 1}$ (which exists by the construction of the
inverse limit groups).  Let
$\Phi_{(m),i,\bullet}:A_{(\infty),i,\star}\to A_{(m),\star}$ be the
canonical continuous epimorphisms associated to this inverse limit
group.

Finally, letting $M_{(\infty),\bullet}$ be the joining Mackey group
of these resulting coordinatizations of $\xi_{(\infty),i}\succsim
\a_{(\infty),i}$ we see that this must be the intersection of the
lifted Mackey groups $(\Phi_{(m),1,\bullet}\times
\Phi_{(m),2,\bullet}\times
\Phi_{(m),3,\bullet})^{-1}(M_{(m),\bullet})$, and so it still has
trivial one-dimensional slices and full two-dimensional projections,
implying that
\[A_{(\infty),1,\star} = A_{(\infty),2,\star} = A_{(\infty),3,\star}\]
(so we may drop the superfluous subscript), and in fact it is now
clear that $M_{(\infty),\bullet}$ has the simple zero-sum form.

Since $\bfX_{(\infty)}$ is still FIS$^+$ by
Proposition~\ref{prop:fibre-normals-stable}, this completes the
proof of Proposition~\ref{prop:better-Mackey} save for exhibiting
the cocycle equation
\begin{multline*}
\s_{(\infty),1}(\bf{p}_1,w_1)\cdot
\s_{(\infty),2}(\bf{p}_2,w_2)\cdot \s_{(\infty),3}(\bf{p}_3,w_3)\\ =
\DDelta_{T_{(\infty)}|_{\a_{(\infty),1}}^{\bf{p}_1}\times
T_{(\infty)}|_{\a_{(\infty),2}}^{\bf{p}_2} \times
T_{(\infty)}|_{\a_{(\infty),3}}^{\bf{p}_3}}b(w_1,w_2,w_3)
\end{multline*}
for some $b:W_{(\infty),1}\times W_{(\infty),2}\times
W_{(\infty),3}\to A_{(\infty),\star}$. Given the zero-sum form of $M_{(\infty),\bullet}$ this is now immediate from the introductory
discussion of Subsection~\ref{subs:first-main-results}. \qed

\subsection{First cocycle factorization}\label{subs:factorizing-cocycles}

Following the work of the preceding two sections we will now
consider an FIS$^+$ system $\bfX$ that satisfies in addition the
conclusions of Proposition~\ref{prop:better-Mackey}, and will next
begin to put the cocycles $\s_i$ into a more convenient form.

Our first step is to cut down the individual dependence of the
cocycle $\s_i(\bf{p}_i,\,\cdot\,)$ for $T|_{\a_i}^{\bf{p}_i}$ from
the proto-characteristic factor $\a_i$ to the subcharacteristic
factor $\b_i$ (we will not obtain any similar simplification for
$\s_i(\bf{n},\,\cdot\,)$ for any $\bf{n} \not\in\bbZ\bf{p}_i$, since
the coboundary equation obtained in
Proposition~\ref{prop:better-Mackey} does not give any immediate
information for these other $\bf{n}$). This relies on a fairly
simple measurable selection argument, but depends crucially on the
relative invariance of the restriction of $T^{\bf{p}_i}$ to
$\b_i|_{\a_i}:\bfW_i\to \bfV_i$. After this we will show how the
resulting cocycle $\s_i$ can be factorized as a product of even
simpler cocycles.

\begin{prop}\label{prop:no-extraneous-dependences}
Every system $\bfX_0$ has an extension $\pi:\bfX\to\bfX_0$ that is
FIS$^+$ and for which
\begin{center}
$\phantom{i}$\xymatrix{
\bfY_i\ar[dr]_{\a_i|_{\xi_i}}\ar@{<->}[rr]^-\cong & &
\bfW_i\ltimes (A_\star,m_{A_\star},\s_i)\ar[dl]^{\rm{canonical}}\\
& \bfW_i }
\end{center}
for some compact Abelian group data $A_\star$ and some cocycles
$\s_i$ such that the associated coordinatization by group data of
the subextension $\vec{\a}|_{\vec{\xi}}:\bfY\to \bfW$ inside the
Furstenberg self-joining has Mackey group data
\[M_\bullet = \{(a_1,a_2,a_3)\in A_\star^3:\ a_1\cdot a_2\cdot a_3 = 1_{A_\star}\}\quad\quad\vec{\a}_\#\mu^\rm{F}\hbox{-a.s.},\]
and also such that $\s_i(\bf{p}_i,\,\cdot\,)$ is measurable with
respect to $\beta_i|_{\a_i}$.
\end{prop}

\textbf{Proof}\quad Proposition~\ref{prop:better-Mackey} already
gives an FIS$^+$ extension $\bfX$ satisfying all of the desired
conditions except for the restricted dependence of
$\s_i(\bf{p}_i,\,\cdot\,)$. Proposition~\ref{prop:better-Mackey}
also gives the joint coboundary equation
\begin{multline*}
\s_1(\bf{p}_1,w_1)\cdot \s_2(\bf{p}_2,w_2)\cdot
\s_3(\bf{p}_3,w_3) = \DDelta_{\vec{T}|_{\vec{\xi}}}b(w_1,w_2,w_3)\\
\hbox{for }\vec{\a}_\#\mu^{\rm{F}}\hbox{-a.e. }(w_1,w_2,w_3)
\end{multline*}
for the corresponding Mackey section $b:W_1\times W_2\times W_3 \to
A$.

Consider the factor
\[\vec{\b}|_{\vec{\a}}:W_1\times W_2\times W_3 \to V_1\times V_2\times V_3.\]
We know from the discussion of
Subsection~\ref{subs:first-main-results} that the coordinate
projections $\pi_1$, $\pi_2$, $\pi_3$ on $W_1\times W_2\times W_3$
are relatively independent over their further factors
$\b_1\circ\pi_1$, $\b_2\circ\pi_2$, $\b_3\circ\pi_3$ under
$\mu^{\rm{F}}$, and so, choosing $T$-equivariant probability kernels
$P_i:V_i\stackrel{\rm{p}}{\longrightarrow} W_i$ representing the
disintegrations of $(\a_i)_\#\mu$ over $\b_i|_{\a_i}$, we can
express \[\vec{\a}_\#\mu^{\rm{F}} = \int_{V_1\times V_2\times
V_3}P_1(v_1,\,\cdot\,)\otimes P_2(v_2,\,\cdot\,)\otimes
P_3(v_3,\,\cdot\,)\,\vec{\b}_\#\mu^{\rm{F}}(\d(v_1,v_2,v_3)).\]

In conjunction with the above cocycle equation, we conclude from
this that:
\begin{quote}
for $\vec{\b}_\#\mu^{\rm{F}}$-a.e. $(v_1,v_2,v_3)$,\\ it holds that
for
$(P_2(v_2,\,\cdot\,)\otimes P_3(v_3,\,\cdot\,))$-a.e. $(w_2,w_3)$,\\
it holds that for $P_1(v_1,\,\cdot\,)$-a.e. $w_1$ we have
\[\s_1(\bf{p}_1,w_1)\cdot \s_2(\bf{p}_2,w_2)\cdot
\s_3(\bf{p}_3,w_3) = \DDelta_{\vec{T}|_{\vec{\xi}}}b(w_1,w_2,w_3).\]
\end{quote}
In addition, the above condition on $(w_2,w_3)$ is easily seen to be
measurable, and the extension $\vec{\b}|_{\vec{\a}}:W_1\times W_2\times W_3\to V_1\times V_2\times V_3$ is relatively $\vec{T}$-invariant (simply from the definition of the $\a_i$), and therefore by Proposition~\ref{prop:invar-meas-select} we can
choose a $\vec{T}$-equivariant measurable selector $\eta =
(\eta_2,\eta_3):V_1\times V_2\times V_3 \to W_2\times W_3$ such that
\begin{quote}
for $\vec{\b}_\#\mu^{\rm{F}}$-a.e. $(v_1,v_2,v_3)$,\\ it holds that
for $P_1(v_1,\,\cdot\,)$-a.e. $w_1$ we have
\[\s_1(\bf{p}_1,w_1)\cdot \s_2(\bf{p}_2,\eta_2(\vec{v}))\cdot
\s_3(\bf{p}_3,\eta_3(\vec{v})) =
\DDelta_{\vec{T}|_{\vec{\xi}}}b(w_1,\eta_2(\vec{v}),\eta_3(\vec{v})).\]
\end{quote}

Now let $\pi':\bfX'\to \bfV_1^{\ \uhr(\bbZ\bf{p}_1 +
\bbZ\bf{p}_2)}$ be the extension given by extending $\b_1(\bfX)$ to
a system on $V_1\times V_2\times V_3$ through the first coordinate
projection, lifting $(\b_1)_\#\mu$ to $\vec{\b}_\#\mu^\rm{F}$,
$T|_{\b_1}^{\bf{p}_1}$ to $\vec{T}|_{\vec{\b}}$ and
$T|_{\b_1}^{\bf{p}_2}$ to $(T^{\bf{p}_2})^{\times 3}|_{\vec{\b}}$. Let $\pi'':\bfX''\to \bfV_1$ be an extension over $\pi'$ that
recovers the action of the whole group $\bbZ^2$.  Finally, let
\[\tilde{\bfX} := \bfX\otimes_{\{\b_1= \pi''\}}\bfX''\]
regarded as an extension of $\bfX$ through the first coordinate
projection.

Under the measure $\mu^\rm{F}$ we have $\id_{V_1\times V_2\times
V_3} \simeq (\b_1\circ\pi_1) \vee (\zeta_0^{\bf{p}_2 -
\bf{p}_3}\circ \pi_2) \precsim (\b_1\circ\pi_1) \vee \pi_2$, and
$\b_1(\bfX)$ is a $(\sfZ_0^{\bf{p}_1 - \bf{p}_2}\vee
\sfZ_0^{\bf{p}_1 - \bf{p}_3})$-system and $\pi_2$ is manifestly
invariant under $(T')^{\bf{p}_1 - \bf{p}_2} =
\vec{T}((T^{\bf{p}_2})^{\times 3})^{-1}$.  Therefore the map
$\bfX''\to\bfX'$ (although it will typically not be a factor map for
the whole $\bbZ^2$-action) is nevertheless contained in
$\zeta_{\sfZ_0^{\bf{p}_1 - \bf{p}_2}\vee\sfZ_0^{\bf{p}_1 -
\bf{p}_3}}^{\bfX''}$, and so now we can simply reinterpret the above
cocycle equation in $\tilde{\bfX}$ as asserting that
$\s_1(\bf{p}_1,w_1)$ is cohomologous to a cocycle (given by
$\s_2(\bf{p}_2,\eta_2(\vec{v}))^{-1}\cdot
\s_3(\bf{p}_3,\eta_3(\vec{v}))^{-1}$) that is measurable with
respect to $\tilde{\b}_1$.

Clearly we can perform similar extensions to the end of enlarging
$\b_2$ and $\b_3$, and now alternately combining this kind of
extension and extensions obtained by re-implementing
Proposition~\ref{prop:better-Mackey}, the resulting inverse sequence
has an inverse limit that still enjoys all of the properties
guaranteed by Proposition~\ref{prop:better-Mackey} (by just the same
reasoning as for that proposition itself) and also enjoys the
restricted dependence of the newly-obtained cocycles
$\s_1(\bf{p}_i,\,\cdot\,)$. \qed

Our next trick will be to decompose the cocycles
$\s_i(\bf{p}_i,\,\cdot\,)$ obtained above into products of
simpler factorizing cocycles.  This is the first in a sequence of such factorizations that will eventually lead to Theorem~\ref{thm:char-three-lines-in-2D}.

\begin{prop}\label{prop:cocycle-prod-decomp}
For the extended system obtained above the cocycles $\s_i$ admit factorizations
\[\s_i(\bf{p}_i,\,\cdot\,) =
(\DDelta_{T|_{\a_i}^{\bf{p}_i}}b_i)\cdot\rho_{i,j}\cdot
\rho_{i,k}^{-1}\cdot \tau_i\] in which
\begin{itemize}
\item $\rho_{i,j}$ is $T^{\bf{p}_i - \bf{p}_j}$-invariant,
\item $\rho_{i,j} = \rho_{j,i}^{-1}$,
\item $\tau_i$ is measurable with respect to $(\zeta_1^T\wedge
\zeta_0^{T^{\bf{p}_i - \bf{p}_j}})\vee
(\zeta_1^T\wedge\zeta_0^{T^{\bf{p}_i - \bf{p}_k}})$,
\item the cocycles $\tau_i$ satisfy
\[(\tau_1\circ\pi_1)\cdot (\tau_2\circ\pi_2)\cdot (\tau_3\circ\pi_3) \in \B^1(\vec{T}|_{\bigvee_{ij}(\zeta_1^T\wedge\zeta_0^{T^{\bf{p}_i} = T^{\bf{p}_j}})\circ\pi_i};A_\star).\]
\end{itemize}
\end{prop}

Our approach here will be to first deduce 
something about the necessary structure of the transfer function
$b$, and then infer from this the desired structure for $\s_i$.

The subcharacteristic factor $\b_i$ is given by
$\zeta_0^{T^{\bf{p}_i} =
T^{\bf{p}_j}}\vee\zeta_0^{T^{\bf{p}_i}=T^{\bf{p}_k}}$, and the two
isotropy factors contributing to this join are relatively
independent under $\mu$ over $\zeta_0^{T^{\bf{p}_1} = T^{\bf{p}_2} =
T^{\bf{p}_3}}$, and so we can sensibly write points of $V_i$ as $v_i
= (v_{ij},v_{ik})$, where the two coordinates are independent random
variables after conditioning on $\zeta_0^{T^{\bf{p}_1} =
T^{\bf{p}_2} = T^{\bf{p}_3}}|_{\zeta_0^{T^{\bf{p}_i} =
T^{\bf{p}_j}}}(v_{ij}) = \zeta_0^{T^{\bf{p}_1} = T^{\bf{p}_2} =
T^{\bf{p}_3}}|_{\zeta_0^{T^{\bf{p}_i} = T^{\bf{p}_k}}}(v_{ik})$.

\begin{lem}
If $\s_1$, $\s_2$ and $\s_3$ are as output by
Proposition~\ref{prop:no-extraneous-dependences} and $b:W_1\times
W_2\times W_3 \to A_\star$ is a choice of joining Mackey section, so
\[\s_1(\bf{p}_1,\b_1|_{\a_1}(w_1))\cdot\s_2(\bf{p}_2,\b_2|_{\a_2}(w_2))\cdot\s_3(\bf{p}_3,\b_3|_{\a_3}(w_3)) = \DDelta_{\vec{T}|_{\vec{\a}}}b(w_1,w_2,w_3),\]
then there is a (possibly different) choice of $b$ satisfying this
equation such that
\begin{itemize}
\item[(1)] $b$ is measurable with respect to $\vec{\b}|_{\vec{\a}}$,
\item[(2)] and, writing our cocycles as functions of $v_{12}$, $v_{13}$ and
$v_{23}$, we have that $b$ takes the form
\[b(v_{12},v_{13},v_{23}) = b_1(v_{12},v_{13})\cdot b_2(v_{12},v_{23})\cdot b_3(v_{13},v_{23})\cdot c(z_{12},z_{13},z_{23})\]
where $z_{ij} := \zeta_1^T|_{\zeta_0^{T^{\bf{p}_i} =
T^{\bf{p}_j}}}(v_{ij})$, so in particular $c$ depends only on the
join under $\mu^\rm{F}$ of the group rotation factors
$\zeta_1^T\wedge \zeta_0^{T^{\bf{p}_i} = T^{\bf{p}_j}}$.
\end{itemize}
\end{lem}

\textbf{Proof}\quad\textbf{(1)}\quad The extension $\vec{\b}|_{\vec{\a}}:\bfW\to \bfV$ is
relatively $\vec{T}|_{\vec{\a}}$-invariant, and so given the
$\vec{\a}_\#\mu^\rm{F}$-almost sure equation
\[\s_1(\bf{p}_1,\b_1|_{\a_1}(w_1))\cdot\s_2(\bf{p}_2,\b_2|_{\a_2}(w_2))\cdot\s_3(\bf{p}_3,\b_3|_{\a_3}(w_3)) = \DDelta_{\vec{T}|_{\vec{\a}}}b(w_1,w_2,w_3),\]
by Proposition~\ref{prop:invar-meas-select} we can choose a
$\vec{T}$-equivariant measurable selector $\eta:V_1\times V_2\times
V_3\to W_1\times W_2\times W_3$ such that
\[\s_1(\bf{p}_1,v_1)\cdot\s_2(\bf{p}_2,v_2)\cdot\s_3(\bf{p}_3,v_3) = \DDelta_{\vec{T}|_{\vec{\a}}}b(\eta(v_1,v_2,v_3))\]
for $\vec{\b}_\#\mu^\rm{F}$-almost every $(v_1,v_2,v_3)\in V_1\times
V_2\times V_3$.  Now simply replacing $b$ by
$b\circ\eta\circ\vec{\b}|_{\vec{\a}}$ proves the first conclusion.

\quad\textbf{(2)}\quad We will show that the second conclusion
already holds for the transfer function $b$ output by part (1)
above.  This will rely on the following trick (first shown to me by
Bernard Host). First we agree to write our cocycles as functions on
$V_{12}\times V_{13}\times V_{23}$ instead of $V_1\times V_2\times
V_3$, since by Proposition~\ref{prop:subchar-joint-dist} this is
equivalent up to negligible sets.  Now for $l=1,2$ and
$ij\in\binom{\{1,2,3\}}{2}$ let $V^l_{ij}$ be a copy of $V_{ij}$ and
form the relatively independent product
\[(\tilde{V},\tilde{\l}) := \Big(\prod_{ij\in \binom{\{1,2,3\}}{2}} V^1_{ij}\times \prod_{ij\in \binom{\{1,2,3\}}{2}} V^2_{ij},(\vec{\b}_\#\mu^\rm{F})\otimes_{\zeta_0^{T^{\bf{p}_1} = T^{\bf{p}_2} = T^{\bf{p}_3}}\circ\pi_1}(\vec{\b}_\#\mu^\rm{F})\Big),\]
and note that by Proposition~\ref{prop:subchar-joint-dist} for this
space the natural projection factor maps onto the spaces $V_{i,j}$
are all relatively independent over a single factor map onto
$Z_0^{T^{\bf{p}_1} = T^{\bf{p}_2} = T^{\bf{p}_3}}$ (since
\[\zeta_0^{T^{\bf{p}_1} = T^{\bf{p}_2} = T^{\bf{p}_3}}\circ\pi_1\simeq \zeta_0^{T^{\bf{p}_1} = T^{\bf{p}_2} = T^{\bf{p}_3}}\circ\pi_2\simeq \zeta_0^{T^{\bf{p}_1} = T^{\bf{p}_2} = T^{\bf{p}_3}}\circ\pi_3\quad\hbox{).}\]

We now consider the given combined cocycle equation on each
subproduct of the form $V^{l_{12}}_{12}\times V^{l_{13}}_{13}\times V^{l_{23}}_{23}$ for $l_{12}$, $l_{13}$, $l_{23}\in \{1,2\}$.
Multiplying these equations with alternating sign gives
\begin{eqnarray*}
&&\DDelta_{\vec{T}|_{\vec{\b}}\times
\vec{T}|_{\vec{\b}}}\Big(\prod_{(l_{12},l_{13},l_{23})}b(v_{12}^{l_{12}},v_{13}^{l_{13}},v_{23}^{l_{23}})^{(-1)^{l_{12} + l_{13} + l_{23}}}\Big)\\
&&=
\prod_{(l_{12},l_{13},l_{23})}\big(\s_1(\bf{p}_1,v_{12}^{l_{12}},v_{13}^{l_{13}})
\cdot\s_2(\bf{p}_2,v_{12}^{l_{12}},v_{23}^{l_{23}}) \cdot
\s_3(\bf{p}_3,v^{l_{13}}_{13},v^{l_{23}}_{23})\big)^{(-1)^{l_{12} +
l_{13} + l_{23}}}\\ &&= 1_{A_\star},
\end{eqnarray*}
since the terms on the right-hand side here cancel completely.

It follows that the function
\[\prod_{(l_{12},l_{13},l_{23})}b(v_{12}^{l_{12}},v_{13}^{l_{13}},v_{23}^{l_{23}})^{(-1)^{l_{12} + l_{13} + l_{23}}}\]
on $(\tilde{V},\tilde{\l})$ is
$(\vec{T}|_{\vec{\b}}\times\vec{T}|_{\vec{\b}})$-invariant.
Moreover, $(\tilde{V},\tilde{\l})$ is a relatively independent
product of two copies of $(V,\vec{\b}_\#\mu^\rm{F})$ over a copy of
$Z_0^{T^{\bf{p}_1} = T^{\bf{p}_2} = T^{\bf{p}_3}}$, on which
$T^{\bf{p}_1}$, $T^{\bf{p}_2}$ and $T^{\bf{p}_3}$ all act by the
same rational rotation, since $\bf{p}_1 - \bf{p}_2$ and $\bf{p}_1 -
\bf{p}_2$ together generate a finite-index sublattice of $\bbZ^2$,
and over which each fibre copy of $(V,\vec{\b}_\#\mu^\rm{F})$
carries an action of $\vec{T}$ that is relatively ergodic over the
common copy of $\bfZ_0^{T^{\bf{p}_1} = T^{\bf{p}_2} = T^{\bf{p}_3}}$
up to another rational rotation factor (by
Lemma~\ref{lem:invariances-leave-first-coord} and since each pair
$\bf{p}_i$, $\bf{p}_i - \bf{p}_j$ also generate a finite-index
sublattice of $\bbZ^2$). It follows that the above product function
must actually be measurable with respect to the join of all the
relevant copies of the group rotation factors $\zeta_1^T\wedge
\zeta_0^{T^{\bf{p}_i}= T^{\bf{p}_j}}$, and so we may write it as
$c^\circ(z_{12}^1,z_{12}^2,\ldots,z_{23}^2)^{-1}$ in the obvious
notation.  Now we can simply re-arrange the definition of this
function to obtain
\[b(v^1_{12},v^1_{13},v^1_{23})
= \Big(\prod_{(l_{12},l_{13},l_{23})\neq
(1,1,1)}b(v^{l_{12}}_{12},v^{l_{13}}_{13},v^{l_{23}}_{23})^{(-1)^{l_{12} + l_{13} + l_{23}}}\Big)\cdot
c^\circ(z_{12}^1,\ldots,z_{23}^2).\]

Finally we choose a measurable selector $\eta:V^1_{12}\times
V^1_{13}\times V^1_{23}\to\tilde{V}$ so that the above equation is
satisfied at
$(v^1_{12},v^1_{13},v^1_{23},\eta(v^1_{12},v^1_{13},v^1_{23}))$ for
$\vec{\b}_\#\mu^\rm{F}$-a.e. $(v^1_{12},v^1_{13},v^1_{23})$.
Composed with this measurable selector, the function
\[b(v^1_{12},v^1_{13},v^2_{23})\]
virtually becomes a function of $v^1_{12}$ and $v^1_{13}$ alone, and
similarly for all other contributions to the product on the
right-hand side above except the last. Hence by suitably grouping
these together the above equation is now itself in the form
\[b(v_{12},v_{13},v_{23}) = b_1(v_{12},v_{13})\cdot b_2(v_{12},v_{23})\cdot b_3(v_{13},v_{23}) \cdot c(z_{12},z_{13},z_{23})\]
for suitable measurable functions $b_1$, $b_2$, $b_3$ and $c$, as
required. \qed

\begin{cor}
If $\s_1$, $\s_2$ and $\s_3$ are as output by
Proposition~\ref{prop:no-extraneous-dependences} then there are
sections $b_i:V_{ij}\times V_{ik}\to A_\star$ such that the
cohomologous cocycles $\s'_i := \s_i \cdot
\DDelta_{T|_{\a_i}}(b_i\circ\b_i|_{\a_i})$ are such that each
$\s'_i(\bf{p}_i,\,\cdot\,)$ is $\b_i|_{\a_i}$-measurable and these
satisfy
\[\s'_1(\bf{p}_1,v_1) \cdot\s'_2(\bf{p}_2,v_2) \cdot\s'_2(\bf{p}_3,v_3) = \DDelta_{\vec{T}|_{\vec{\b}}}c(v_1,v_2,v_3)\]
for some section $c:V_1\times V_2\times V_3\to A_\star$ that depends
only on the join of the group rotation factors $(\zeta_1^T\wedge
\zeta_0^{T^{\bf{p}_i}=T^{\bf{p}_j}})\circ\pi_i$.
\end{cor}

\textbf{Proof}\quad Let
\[b(v_{12},v_{13},v_{23}) = b_1(v_{12},v_{13}) \cdot b_2(v_{12},v_{23}) \cdot b_3(v_{13},v_{23}) \cdot c(z_{12},z_{13},z_{23})\]
be the factorization of $b$ obtained in the preceding lemma, and now
let $\s'_i := \s_i \cdot \DDelta_{T|_{\a_i}}(b_i\circ\b_i|_{\a_i})$
for these $b_i$.  In these terms the combined coboundary equation
simply re-arranges to give precisely
\[\s'_1(\bf{p}_1,v_1) \cdot \s'_2(\bf{p}_2,v_2) \cdot \s'_2(\bf{p}_3,v_3) = \DDelta_{\vec{T}|_{\bigvee_{ij}(\zeta_1^T\wedge\zeta_0^{T^{\bf{p}_i} = T^{\bf{p}_j}})\circ\pi_i}}c(z_{12},z_{13},z_{23}),\]
which is the required equation upon lifting $c$ to be a function on
$V_1\times V_2\times V_3$. \qed

\textbf{Proof of Proposition~\ref{prop:cocycle-prod-decomp}}\quad
Considering the equation
\[\s'_1(\bf{p}_1,v_1) \cdot \s'_2(\bf{p}_2,v_2) \cdot \s'_2(\bf{p}_3,v_3) = \DDelta_{\vec{T}|_{\vec{\b}}}c(v_1,v_2,v_3)\]
obtained from the preceding corollary, and recalling again the
relative independence of $v_{12}$, $v_{13}$ and $v_{23}$ under
$\vec{\b}_\#\mu^\rm{F}$ promised by
Proposition~\ref{prop:subchar-joint-dist}, we see that we can make a
measurable selection $\eta:V_{12}\times V_{13}\to V_{23}$ that
actually depends only on $\zeta_0^{T^{\bf{p}_1} = T^{\bf{p}_2} =
T^{\bf{p}_3}}(v_{12}) = \zeta_0^{T^{\bf{p}_1} = T^{\bf{p}_2} =
T^{\bf{p}_3}}(v_{13})$ such that
\[\s'_1(\bf{p}_1,v_{12},v_{13}) \cdot \s'_2(\bf{p}_2,v_{12},\eta(v_{12})) \cdot \s'_2(\bf{p}_3,v_{13},\eta(v_{13})) = (\DDelta_{\vec{T}|_{\vec{\b}}}c)(v_{12},v_{13},\eta(v_{12}))\]
almost surely, and so subtracting the second and third left-hand
terms from both sides gives an explicit equation for
$\s'_1(\bf{p}_1,\,\cdot\,)$ as a cocycle of the form
$\rho^\circ_{12}\cdot\rho^\circ_{13}\cdot\tau^\circ_1$ with
$\rho^\circ_{ij}$ a function only of $v_{ij}$ and $\tau^\circ_1$
measurable with respect to the join of its permitted group rotation
factors (although we must be careful: $\tau_1^\circ:(v_{12},v_{13})
\mapsto (\DDelta_{\vec{T}|_{\vec{\b}}}c)(v_{12},v_{13},\eta(v_{12}))$
is \emph{not} usually a coboundary, in spite of appearances, since
in this case $\eta$ is not a selector for a relatively
$\vec{T}$-invariant extension and so cannot necessarily be made
$\vec{T}$-equivariant).

The same is true of $\s'_2$ and $\s'_3$ by symmetry, and so we can
now substitute the resulting form for each
$\s'_i(\bf{p}_i,\,\cdot\,)$ once again into the combined cocycle
equation to obtain
\begin{multline*}
\tau^\circ_1(v_{12},v_{13})\cdot\tau^\circ_2(v_{12},v_{23})\cdot\tau^\circ_3(v_{13},v_{23})\\
\cdot ((\rho^\circ_{12}\cdot \rho^\circ_{21})(v_{12}))\cdot
((\rho^\circ_{13}\cdot \rho^\circ_{31})(v_{13}))\cdot
((\rho^\circ_{23}\cdot \rho^\circ_{32})(v_{23}))\\ =
\DDelta_{\vec{T}|_{\vec{\b}}}c(v_{12},v_{13},v_{23}).
\end{multline*}
Since $v_{12}$, $v_{13}$ and $v_{23}$ are certainly relatively
independent under $\mu$ over their factor-map images
$\zeta_1^T|_{\zeta_0^{T^{\bf{p}_1}=T^{\bf{p}_2}}}(v_{12})$,
$\zeta_1^T|_{\zeta_0^{T^{\bf{p}_1}=T^{\bf{p}_3}}}(v_{13})$ and
$\zeta_1^T|_{\zeta_0^{T^{\bf{p}_2}=T^{\bf{p}_3}}}(v_{23})$, it
follows that each $(\rho^\circ_{ij}\cdot \rho^\circ_{ji})(v_{ij})$
is virtually a function only of $z_{ij} =
\zeta_1^T|_{\zeta_0^{T^{\bf{p}_i}=T^{\bf{p}_j}}}(v_{ij})$. We now
define
\[\rho_{12} := \rho_{21}^{-1} = \rho^\circ_{12},\quad\quad \rho_{31} = \rho_{13}^{-1} := \rho^\circ_{31},\quad\quad\hbox{and }\rho_{23} = \rho_{32}^{-1} := \rho^\circ_{23}\]
and
\[\tau_1 := \tau^\circ_1\cdot (\rho^\circ_{13}\cdot \rho^\circ_{31}),\quad\quad\tau_2 := \tau^\circ_2\cdot (\rho^\circ_{12}\cdot \rho^\circ_{21}) \quad\quad\hbox{and }\tau_3 := \tau^\circ_3\cdot (\rho^\circ_{23}\cdot \rho^\circ_{32}),\]
so that $\tau_i\cdot \rho_{ij}\cdot\rho_{ik} =
\tau^\circ_i\cdot\rho^\circ_{ij}\cdot\rho^\circ_{ik}$ for each $i$,
to obtain an equivalent factorization of each
$\s'_i(\bf{p}_i,\cdot\,)$ in terms of which the combined cocycle
equation now simplifies to
\[(\tau_1\circ\pi_1)\cdot(\tau_2\circ\pi_2)\cdot(\tau_3\circ\pi_3) = \DDelta_{\vec{T}|_{\vec{\b}}}c\]
with all of these function now actually depending only on the join
of the relevant group rotation factors, as required. \qed

\subsection{Reduction to another proposition on factorizing cocycles}

The final proof of Theorem~\ref{thm:char-three-lines-in-2D} will follow from an enhancement of the cocycle factorization of Proposition~\ref{prop:cocycle-prod-decomp}.  In the present subsection we introduce this enhancement and show how it leads to the full theorem.

\textbf{Notation}\quad Extending Definition~\ref{dfn:motionless}, we will henceforth write
$(Z_\star,R_{\phi_\star})$ to denote a $\bbZ^2$-system whose
underlying space is the direct integral of some measurably-varying
family of compact Abelian groups $Z_\star$, indexed by some other
standard Borel probability space $(S,\nu)$ on which the action is
trivial, with the overall action a fibrewise rotation defined by a
measurable selection for each $s$ of a dense homomorphism
$\phi_s:\bbZ^2\to Z_s$: that is, $R_{\phi_\star}$ is
given by
\[R_{\phi_\star}^{\bf{n}}(s,z) := (s,z\cdot \phi_s(\bf{n}))\quad\quad\hbox{for\ }s\in S,\ z\in Z_s\ \hbox{and}\ \bf{n} \in \bbZ^2.\]
Although we sometimes omit to mention it, the measure on this system is the integral of the Haar measures $m_{Z_\star}$. We will refer to such a system as a \textbf{direct integral of
ergodic group rotations} and to $(S,\nu)$ as its \textbf{invariant
base space}. Sometimes we omit the base space $(S,\nu)$ from mention
completely, since once again the forthcoming arguments will all
effectively be made fibrewise, just taking care that all
newly-constructed objects can still be selected measurably. In
particular, we will often write just $Z_\star$ in place of $S\ltimes Z_\star$. \fin

\begin{prop}\label{prop:cocyc-factorize}
Let $\bfX$ be a system as output by
Proposition~\ref{prop:cocycle-prod-decomp} with the factorization of $\s_i(\bf{p}_i,\cdot)$ given there and let $\zeta_1^T:\bfX\to (Z_\star,m_{Z_\star},R_{\phi_\star})$ be a coordinatization of its Kronecker factor.  In addition let
\[\tau_i^{(m)}:= \tau_i\cdot (\tau_i\circ R_{\phi_\star(\bf{p}_i)})\cdot\cdots\cdot (\tau_i\circ R_{\phi_\star((m-1)\bf{p}_i)})\]
for any integer $m\geq 1$. For any motionless selection of characters $\chi_\star \in \widehat{A_\star}$ there are a fibrewise extension of ergodic group rotations $q_\star:(\t{Z}_\star,R_{\t{\phi}_\star})\to (Z_\star,R_{\phi_\star})$, a motionless selection of integers $m_\star \geq 1$ and a motionless family of finite-index subgroups $Z_{0,\star}\leq \t{Z}_\star$ such that $\t{\phi}_\star(m_\star\bf{p}_i)\in Z_{0,\star}$ and
\[\chi_\star\circ\tau^{(m_\star)}_i\circ q_\star = \tau_{i,j}\cdot\tau_{i,k}\cdot\tau_{\rm{nil},i}\cdot \DDelta_{\t{\phi}_\star(\bf{p}_i)}\b\]
for some $\b \in \C(\t{Z}_\star)$, where
\begin{itemize}
\item $\tau_{i,j}$ is $R_{\t{\phi}_\star(\bf{p}_i - \bf{p}_j)}$-invariant and $\tau_{i,k}$ is $R_{\t{\phi}_\star(\bf{p}_i - \bf{p}_k)}$-invariant, and
\item the two-step Abelian distal transformation $R_{\t{\phi}_s(m_s\bf{p}_i)}\ltimes \tau_{\nil,i}\actson (z_0Z_{0,s})\times \Sone$ is a two-step nilrotation for every coset $z_0Z_{0,s}\leq \t{Z}_s$ for $\nu$-almost every $s \in S$.
\end{itemize}
\end{prop}

In this section we deduce Theorem~\ref{thm:char-three-lines-in-2D} from Proposition~\ref{prop:cocyc-factorize} using three smaller lemmas.

\begin{lem}\label{lem:make-nilrtn-global}
Suppose that $\bf{n} \in \bbZ^2\setminus \{\bs{0}\}$, that $\phi_\star:\bbZ^2\to Z_\star$ is a motionless measurable family of dense homomorphisms, that $Z_{0,\star} \leq Z_\star$ has finite index almost surely and that $\s:Z_\star\to\Sone$ is a Borel map that restricts to a $Z_{0,\star}$-local nil-cocycle over $R_{\phi_\star(\bf{n})}$ almost surely.  Then there are a fibrewise extension of rotations $q:(\t{Z}_\star,\t{\phi}_\star)\to (Z_\star,\phi_\star)$, some tuple $\s_1,\s_2,\ldots,\s_{D_\star}:\bbZ^2\times \t{Z}_\star\to \Sone$ of cocycles over $R_{\t{\phi}_\star}$ such that each $\s_i$ restricts to a global nil-cocycle on $\t{Z}_s$ almost surely, and a $\t{\phi}_\star(\bf{n})$-invariant map $d:Z_\star\to \{1,2,\ldots,D_\star\}$ such that $\s(q(\t{z})) = \s_{d(\t{z})}(\bf{n},\t{z})$ for $m_{\t{Z}_s}$-almost every $\t{z}$ for $\nu$-almost every $s$.
\end{lem}

\textbf{Remark}\quad It is very important that for each fibre $Z_s$ we may need to introduce several different global nil-cocycles $\s_i$ to choose from in our representation of $\s$, since a priori we have no information at all that relates the behaviour of $\s$ on different cosets of $Z_{0,s}\cdot\ol{\phi_s(\bbZ\bf{n})}$, which may still be smaller than the whole of $Z_s$. \fin

\textbf{Proof}\quad First note that $\G_s := \phi_s^{-1}(Z_{0,s})$ is a subgroup of $\bbZ^2$ that varies measurably in $s$ (with the obvious discrete measurable structure on the countable set of subgroups of $\bbZ^2$).  By assumption it always contains $\bf{n}$.

Let $\O_s \subset \bbZ^2$ be a fundamental domain for $\G_s$, chosen to
contain $\bs{0}$, and let $\lfloor\cdot\rfloor_s +
\{\cdot\}_s$ be the resulting decomposition of $\bbZ^2$ into integer- and fractional-parts modulo $\G_s$.

Extending the fibre groups $Z_\star$ if necessary, we may apply Lemma~\ref{lem:making-nilrtns-commute} for some direction $\bf{m}_s \in \G_s$ linearly independent from $\bf{n}$ and then apply Proposition~\ref{prop:loc-nil-to-nil-2} to the finite-index inclusion $\bbZ\bf{m}_s + \bbZ\bf{n} \leq \G_s$ in order to assume that in place of our initially-posited $\s$ we actually have $\s(\bf{n},\cdot)$ for some map $\s:\G_\star\times Z_\star\to\Sone$ that restricts to a $Z_{0,s}$-local $\G_s$-nil-cocycle over $R_{\phi_s}$ for almost every $s$.  (It is easy to see that the selections in Lemma~\ref{lem:making-nilrtns-commute} and Proposition~\ref{prop:loc-nil-to-nil-2} can be made measurably in $s$.)

Now we form another extension $Z_\star\times (\bbZ^2/\G_\star)$ with the
measurable family of homomorphisms
\[\t{\phi}_\star:\bf{p}\mapsto (\phi_\star(\bf{p}),\bf{p} + \G),\]
and over this we consider the Abelian extension with toral fibres $(\Sone)^{\O_\star}$ and
with cocycle $\vec{\s}:\bbZ^2 \times \t{Z}_\star\to (\Sone)^{\O_\star}$ given by
\[\vec{\s}(\bf{p},(z,\bf{k} + \G)) := \Big(\s \big(\lfloor \{\omega + \bf{k}\}_s + \bf{p}\rfloor_s,\phi_s(-\{\omega + \bf{k}\}_s)\cdot z\big)\Big)_{\omega \in \O_s}.\]
By restricting to $\t{Z}_\star:= \ol{\t{\phi}_\star(\bbZ^2)}$, which must still cover the whole of $Z_\star$ through the restriction $q_\star$ of the coordinate projection $Z_\star\times (\bbZ^2/\G_\star)\to Z_\star$ because $\phi_\star$ has dense image almost surely, we may assume that $(\t{Z}_s,R_{\t{\phi}_s})$ is ergodic for almost every $s$, and re-interpret $\vec{\s}$ as a cocycle on this space.

For any fixed $\omega \in \O_s$ the associated coordinate of this cocycle is
\[\s_\omega(\bf{p},(z,\bf{k} + \G)) := \s \big(\lfloor \{\omega + \bf{k}\}_s + \bf{p}\rfloor_s,\phi_s(-\{\omega + \bf{k}\}_s)\cdot z\big).\]

We can now check the following:
\begin{itemize}
\item $\vec{\s}$ is a cocycle over $R_{\t{\phi}_\star}$: it suffices to check this for each $\omega \in \O_\star$ separately, to which end the cocycle equation for $\s$ gives
\begin{eqnarray*}
&&\s_\omega(\bf{p},\t{\phi}_s(\bf{q})\cdot(z,\bf{k} + \G))\cdot\s_\omega(\bf{q},(z,\bf{k} + \G))\\
&&= \s \big(\lfloor \{\omega + \bf{k} + \bf{q}\}_s + \bf{p}\rfloor_s,\\
&&\quad\quad\quad\quad\quad\quad \phi(\bf{q} + \{\omega + \bf{k}\}_s - \{\omega + \bf{k} + \bf{q}\}_s)\cdot \phi_s(-\{\omega + \bf{k}\}_s)\cdot z\big)\\
&&\quad\quad\quad\quad\cdot \s \big(\lfloor \{\omega + \bf{k}\}_s + \bf{q}\rfloor_s,\phi_s(-\{\omega + \bf{k}\}_s)\cdot z\big)\\
&&= \s \big(\lfloor \{\omega + \bf{k} + \bf{q}\}_s + \bf{p}\rfloor_s,\\
&&\quad\quad\quad\quad\quad\quad \phi(\{\omega + \bf{k}\}_s + \bf{q} - \{\{\omega + \bf{k}\}_s + \bf{q}\}_s)\cdot \phi_s(-\{\omega + \bf{k}\}_s)\cdot z\big)\\
&&\quad\quad\quad\quad\cdot \s \big(\lfloor \{\omega + \bf{k}\}_s + \bf{q}\rfloor_s,\phi_s(-\{\omega + \bf{k}\}_s)\cdot z\big)\\
&&= \s \big(\lfloor \{\omega + \bf{k} + \bf{q}\}_s + \bf{p}\rfloor_s, \phi(\lfloor\{\omega + \bf{k}\}_s + \bf{q}\rfloor_s)\cdot \phi_s(-\{\omega + \bf{k}\}_s)\cdot z\big)\\
&&\quad\quad\quad\quad\cdot \s\big(\lfloor \{\omega + \bf{k}\}_s + \bf{q}\rfloor_s,\phi_s(-\{\omega + \bf{k}\}_s)\cdot z\big)\\
&&= \s \big(\lfloor \{\omega + \bf{k}\}_s + \bf{q}\rfloor_s + \lfloor \{\omega + \bf{k} + \bf{q}\}_s + \bf{p}\rfloor_s,\phi_s(-\{\omega + \bf{k}\}_s)\cdot z\big)\\
&&= \s_\omega(\bf{q} + \bf{p},(z,\bf{k} + \G))
\end{eqnarray*}
for any $\bf{q},\bf{p} \in \bbZ^2$, as required;
\item $\vec{\s}$ is a nil-cocycle: once again, it suffices to check this coordinatewise, but if $\bf{p} \in \G_s$ then
\[\s\big(\lfloor \{\omega + \bf{k}\}_s + \bf{p}\rfloor_s,\phi_s(-\{\omega + \bf{k}\}_s)\cdot z\big) = \s\big(\bf{p},\phi_s(-\{\omega + \bf{k}\}_s)\cdot z\big),\]
and so our assumptions give that this is a $q^{-1}(Z_{0,s})$-local $\G_s$-nil-cocycle, and therefore since we have extended it to a cocycle for the whole of $\bbZ^2$ this extension must in fact be a global nil-cocycle by Proposition~\ref{prop:loc-nil-to-nil-1}.
\end{itemize}

Finally, in view of the identity
\[\{\{-\bf{p}\} + \bf{p}\} = \{\{-\bf{p}\} - (\lfloor-\bf{p}\rfloor + \{-\bf{p}\})\} = \{-\lfloor - \bf{p}\rfloor\}\equiv \bs{0},\]
we have
\begin{multline*}
\s(\bf{n},(z,\bf{k} + \G)) = \s(\lfloor\{\{-\bf{k}\} + \bf{k}\} + \bf{n}\rfloor,\phi(-\{\{-\bf{k}\} + \bf{k}\})\cdot z)\\ = \s_{\{-\bf{k}\}}(\bf{n},(z,\bf{k} + \G)).
\end{multline*}
Hence if we set $D_\star := |\O_\star|$ and so regard the coordinates $\s_\omega$ as indexed by $1,2,\ldots,D_\star$, then composing the map $(z,\bf{k} + \G)\mapsto \{-\bf{k}\}$ with this enumeration of $\O_\star$ gives a function $d$ with the desired properties.
\qed

\begin{lem}\label{lem:make-subjoining}
Suppose that $\bf{n}_1,\bf{n}_2,\bf{n}_3\in \bbZ^2$ are such that any two are linearly independent, that $\bfY$ is a $\bbZ^2$-system that is a joining of $\bfY_2 \in \sfZ_0^{\bf{n}_2}$, $\bfY_3 \in \sfZ_0^{\bf{n}_3}$ and $\bfZ = (Z_\star,m_{Z_\star},R_{\phi_\star}) \in \sfZ_1^{\bbZ^2}$, that $C_\star \leq \Sone$ is a motionless selection of closed subgroups, and that $\bfX = \bfY\ltimes (C_\star,m_{C_\star},\s)$ is an extension with a cocycle-section $\s:\bbZ^2\times Y\to C_\star$ that admits a factorization
\[\s(\bf{n}_1,\cdot) = \DDelta_{\phi_\star(\bf{n}_1)}b\cdot\s_2\cdot \s_3\cdot \s_\nil\]
in which
\begin{itemize}
\item $b$ takes values in $\Sone$,
\item $\s_j$ takes values in $\Sone$ and is lifted from $\bfY_j$ for $j=2,3$, and
\item $\s_\nil$ takes values in $\Sone$ and is such that the two-step Abelian transformation $R_{\phi_s(\bf{n}_1)}\ltimes \s_\nil \actson (z_0Z_{0,s})\times \Sone$ is a two-step nilrotation for every coset $z_0Z_{0,s}\leq Z_s$ for almost every $s$.
\end{itemize}
Then $\bfX$ is a subjoining of $\sfZ_0^{\bf{n}_1}$, $\sfZ_0^{\bf{n}_2}$, $\sfZ_0^{\bf{n}_3}$ and $\sfZ_{\nil,2}^{\bbZ^2}$.
\end{lem}

\textbf{Proof}\quad Letting $\bfY_j = (Y_j,\nu_j,S_j)$ for $j=2,3$, we may coordinatize $\bfY$ by some invariant $(\nu_2,\nu_3,m_{Z_\star})$-coupling $\l$ on $Y_2\times Y_3\times Z_\star$.  Also let $\pi:\bfX\to\bfY$ be the canonical factor.

Applying Lemma~\ref{lem:make-nilrtn-global} to the Borel map $\s_\nil$ we can find some relatively ergodic extension $q:(\t{Z}_\star,R_{\t{\phi}_\star})\to (Z_\star,R_{\phi_\star})$, some cocycle-section \[\vec{\s} = (\t{\s}_1,\t{\s}_2,\ldots,\t{\s}_{D_\star}):\bbZ^2\times \t{Z}_\star\to (\Sone)^{D_\star}\]
over $R_{\t{\phi}_\star}$ that restricts to a global nil-cocycle almost surely, and some $\t{\phi}(\bf{n}_1)$-invariant selection $d:Z_\star\to \{1,2,\ldots,D_\star\}$ such that
\[\s_\nil(q(\t{z})) = \s_{d(\t{z})}(\bf{n}_1,\t{z})\quad\quad\hbox{almost surely}.\]
Defining $\t{\bfZ} := (\t{Z}_\star,m_{\t{Z}_\star},R_{\t{\phi}_\star})$ and
\[\bfX_\nil:= \t{\bfZ}\ltimes ((\Sone)^{D_\star},m_{(\Sone)^{D_\star}},\vec{\s}),\]
it follows that this is a direct integral of two-step $\bbZ^2$-nilsystems.  Let $\pi_\nil:\bfX_\nil\to \t{\bfZ}$ be the canonical factor.

On the other hand, for $j=2,3$ we can extend $\s_j:Y_j\to\Sone$ to a cocycle $\s_j':(\bbZ\bf{n}_1 + \bbZ\bf{n}_j)\times Y_j\to\Sone$ over $S_j^{\ \uhr\ \bbZ\bf{n}_1 + \bbZ\bf{n}_j}$ simply by setting $\s'_j(\bf{n}_j,\cdot)\equiv 1$, where the cocycle equation for $\s_j'$ follows from the assumption that $\s_j$ depends only on a factor which has trivial $(\bbZ\bf{n}_j)$-subaction. This defines an extension
\[\bfX^\circ_j := \bfY_j^{\ \uhr\ \bbZ\bf{n}_1 + \bbZ\bf{n}_j}\ltimes (\Sone,m_{\Sone},\s'_j)\]
of $(\bbZ\bf{n}_1 + \bbZ\bf{n}_j)$-subactions in which $\bf{n}_j$ still acts trivially.  We may therefore interpret this as an extension of $((\bbZ\bf{n}_1 + \bbZ\bf{n}_j)/\bbZ\bf{n}_j)$-systems, and now constructing a further FP extension (Definition 3.17 in~\cite{Aus--lindeppleasant1}) for the inclusion of groups
\[\bbZ^2/\bbZ\bf{n}_j \geq (\bbZ\bf{n}_1 + \bbZ\bf{n}_j)/\bbZ\bf{n}_j\]
we obtain $\pi_j:\bfX_j\to\bfY_j$ such that $\bfX_j$ still has trivial $(\bbZ\bf{n}_j)$-subaction, and where $\bfX^\circ_j\to \bfY_j^{\ \uhr\ \bbZ\bf{n}_1 + \bbZ\bf{n}_j}$ appears as an intermediate extension of the $(\bbZ\bf{n}_1 + \bbZ\bf{n}_j)$-subactions.

Now, the system $\bfY$ is a joining of the targets of the factor maps $\pi$, $\pi_2$, $\pi_3$ and $q\circ \pi_\nil:\bfX_\nil\to\bfZ$, and so we can define $\t{\bfX}$ to be the joining of $\bfX_2$, $\bfX_3$, $\bfX_\nil$ and $\bfX$ that extends $\bfY$ and under which these four factor maps are relatively independent.  The above descriptions of these individual factors give a coordinatization of $\t{\bfX}$ on some space extending
\[(Y_2\times \Sone)\times (Y_3\times \Sone)\times (\t{Z}_\star\ltimes (\Sone)^{D_\star})\ltimes C_\star\]
(where the further extension of this space needed to describe the whole of $\t{\bfX}$ results from the FP extension $\bfX_j$ constructed over $\bfX_j^\circ$).  In particular, the explicit product space above carries an action of $\bbZ\bf{n}_1$ that is a factor of $\t{\bfX}^{\ \uhr\bbZ\bf{n}_1}$, and which is explicitly coordinatized as
\[\big(S_2^{\bf{n}_1}\times S_3^{\bf{n}_1}\times R_{\t{\phi}_\star(\bf{n}_1)}\big)\ltimes \big(\s_2(\bf{n}_1,y_2),\ \s_3(\bf{n}_1,y_3),\ \vec{\s}(\bf{n}_1,\t{z}),\ \s(\bf{n}_1,(y_2,y_3,q(\t{z})))\big).\]
By our assumptions and the output of Lemma~\ref{lem:make-nilrtn-global} this cocycle for the action of $\bf{n}_1$ satisfies
\[\s(\bf{n}_1,(y_2,y_3,q(\t{z}))) = \DDelta_{S^{\bf{n}_1}}b(y_2,y_3,q(\t{z}))\cdot \s_2(\bf{n}_1,y_2)\cdot \s_3(\bf{n}_1,y_3)\cdot \vec{\s}_{d(\t{z})}(\bf{n}_1,\t{z})\]
almost surely.  By the conjugate-minimality of the Mackey group data $M_\bullet$ of this subaction (see Theorem~\ref{thm:homo-nonergMackey-1}), it must satisfy
\[M_\bullet\leq \big\{(s_2,s_3,t_1,t_2,\ldots,t_{D_\star},z)\in (\Sone)^{D_\star + 2}\times C_\star:\ z = s_2 s_3 t_{d(\bullet)}\big\},\]
bearing in mind that this right-hand group data varies measurably and is invariant under $S_2^{\bf{n}_1}\times S_3^{\bf{n}_1}\times R_{\t{\phi}_\star(\bf{n}_1)}$ (although not necessarily under the whole $\bbZ^2$-action $S_2\times S_3\times R_{\t{\phi}_\star}$, because the selection map $d$ is known to be invariant only under the $(\bbZ\bf{n}_1)$-subaction).

In particular, we see that the coordinates $s_2$, $s_3$, $t_1$, \ldots, $t_{D_\star}$ together with the coset $(s_2,\ldots,z)\cdot M_\bullet$ together determine the value of the coordinate $z$.  This implies that under the joining $\t{\bfX}$ the coordinate projections onto $\bfX_2\in \sfZ_0^{\bf{n}_2}$, $\bfX_3 \in \sfZ_0^{\bf{n}_3}$, $\bfX_\nil \in \sfZ_{\nil,2}^{\bbZ^2}$ together with the factor $\sfZ_0^{\bf{n}_1}\t{\bfX}$ determine $\bfX$, and so this explicitly witnesses $\bfX$ as a $(\sfZ_0^{\bf{n}_1},\sfZ_0^{\bf{n}_2},\sfZ_0^{\bf{n}_3},\sfZ_{\nil,2}^{\bbZ^2})$-subjoining, as required. \qed

\begin{lem}\label{lem:simplifying-part-invar}
If $\bf{n}_1,\bf{n}_2 \in \bbZ^2$ are linearly independent and $m\geq 1$ then any system in the class $\sfZ_0^{\bbZ\bf{n}_1}\vee \sfZ_0^{m\bbZ\bf{n}_2}$ is a factor of a system in the class $\sfZ_0^{\bbZ\bf{n}_1}\vee \sfZ_0^{\bbZ\bf{n}_2}$.
\end{lem}

\textbf{Proof}\quad Since we already have $\sfZ_0^{\bbZ\bf{n}_1} \subseteq \sfZ_0^{\bbZ\bf{n}_1}\vee \sfZ_0^{\bbZ\bf{n}_2}$ it suffices to show that any $(X,\mu,T) \in \sfZ_0^{m\bbZ\bf{n}_2}$ has an extension in the class $\sfZ_0^{\bbZ\bf{n}_1}\vee \sfZ_0^{\bbZ\bf{n}_2}$.

To show this, let us first treat the case in which $\bf{n}_1$, $\bf{n}_2$ comprise a basis of $\bbZ^2$. Observe that the extension $\zeta:\bfX \to \sfZ_0^{\bf{n}_2}\bfX$ must be a direct integral of group rotations for the subaction of $\bbZ\bf{n}_2$ over the $\bf{n}_2$-invariant system $\sfZ_0^{\bf{n}_2}\bfX$, where almost all of the fibre groups are quotients of $\bbZ/m\bbZ$ (hence finite).  Let $A_\bullet$ be these finite measurable $T|_\zeta$-invariant group data, and let $\bfY := (Y,\nu,S)$ be a coordinatization of $\sfZ_0^{\bf{n}_2}\bfX$.

Since the whole action $T$ must respect this factor and commute with $T^{\bf{n}_2}$, it follows from Theorem~\ref{thm:RAST} that $T$ can be coordinatized over the factor $\zeta$ as $S\ltimes \s$ for $\s$ a cocycle section $\bbZ^2\times Y\to A_\bullet$.

Now let $\t{\bfX} := \bfY\ltimes (A_\bullet^2,m_{A_\bullet^2},(\s_1,\s_2))$ with
\[\s_1(p\bf{n}_1 + q\bf{n}_2,\cdot) := \s(q\bf{n}_2,\cdot),\quad \s_2(p\bf{n}_1 + q\bf{n}_2,\cdot) := \s(p\bf{n}_1,\cdot).\]
From the vanishing $S^{\bf{n}_2} = \id_Y$ we can deduce firstly that $\s_1$, $\s_2$ satisfy the equations of a cocycle over $S$ and secondly that
\begin{multline*}
\s_1(p\bf{n}_1 + q\bf{n}_2,\cdot) \cdot \s_2(p\bf{n}_1 + q\bf{n}_2,\cdot) = \s(q\bf{n}_2,\cdot) \cdot \s(p\bf{n}_1,\cdot)\\
= \s(q\bf{n}_2,\cdot) \cdot \s(p\bf{n}_1,S^{q\bf{n}_2}(\cdot)) = \s(p\bf{n}_1 + q\bf{n}_2,\cdot).
\end{multline*}

Now the map $(y,a,a')\mapsto (y,a\cdot a')$ defines a factor map
$\t{\bfX}\to \bfX$, and on the other hand the two coordinate
projections $(y,a,a')\mapsto (y,a)$ and $\mapsto (y,a')$ yield
factors that are two different extensions of $\bfY$, the first being
$\bf{n}_1$-invariant and the second $\bf{n}_2$-invariant, and so
since also $\bfY \in \sfZ_0^{\bbZ\bf{n}_2}$ it follows that the
joining $\t{\bfX}$ of these two systems lies in
$\sfZ_0^{\bf{n}_1}\vee \sfZ_0^{\bf{n}_2}$, as required.

Finally, in case $\bf{n}_1, \bf{n}_2$ do not span the whole of $\bbZ^2$, we first extend the $(\bbZ\bf{n}_1 + \bbZ\bf{n}_2)$-subaction to obtain $\t{\bfX}$ as above, and then form a further extension $\bfX'$ of $\t{\bfX}$ that recovers the action of the whole of $\bbZ^2$ (for example, an FP extension, as in Definition 3.17 of~\cite{Aus--lindeppleasant1}), for which we then still have that $\bfX$ is a factor of $(\sfZ_0^{\bf{n}_1}\vee \sfZ_0^{\bf{n}_2})\bfX'$. \qed

\textbf{Proof of Theorem~\ref{thm:char-three-lines-in-2D} from Proposition~\ref{prop:cocyc-factorize}}\quad This is most easily phrased using satedness and an argument by contradiction.  In view of Lemma~\ref{lem:pronil-idemp} and the existence of multiply-sated extensions (Theorem 3.11 in~\cite{Aus--lindeppleasant1}) we may assume that $\bfX_0$ itself is sated relative to all joins of isotropy factors and two-step pro-nilsystem factors.  We will show by contradiction that under this assumption $\bfX_0$ must itself admit the characteristic factors described in Theorem~\ref{thm:char-three-lines-in-2D}.

Thus, suppose that Theorem~\ref{thm:char-three-lines-in-2D} fails for some triple of directions $\bf{p}_1$, $\bf{p}_2$, $\bf{p}_3$. We know that $\bfX_0$ does admit some minimal characteristic factors $\xi_{i,0}$, $i=1,2,3$, for these directions, so our supposition implies that
\begin{eqnarray}\label{eq:noncont}
\xi_{i,0}\not\precsim \zeta_0^{T_0^{\bf{p}_i}}\vee \zeta_0^{T_0^{\bf{p}_i}=T_0^{\bf{p}_j}}\vee\zeta_0^{T_0^{\bf{p}_i}=T_0^{\bf{p}_k}}\vee \zeta^{T_0}_{\nil,2}
\end{eqnarray}
for at least one $i \in \{1,2,3\}$.

Let $\bfX\to\bfX_0$ be a further extension as given by Proposition~\ref{prop:cocycle-prod-decomp}.  Since the minimal characteristic factors $\xi_i$ of $\bfX$ must certainly contain those of $\bfX_0$, that theorem now implies that 
\begin{itemize}
\item $\xi_{i,0}\circ\pi$ is contained in some Abelian isometric extension of
\[\zeta_0^{T^{\bf{p}_i}}\vee \zeta_0^{T^{\bf{p}_i}=T^{\bf{p}_j}}\vee\zeta_0^{T^{\bf{p}_i}=T^{\bf{p}_k}}\vee \zeta^{T}_1:\bfX\to \bfW = (W,\theta,R),\]
where $\bfW$ is a suitable choice of target system for this factor map, and
\item the extension is given by a cocycle $\s_i:\bbZ^2\times W\to A_\star$ that admits a factorization
\[\s_i(\bf{p}_i,\,\cdot\,) = (\DDelta_{T|_{\a_i}^{\bf{p}_i}}b_i)\cdot\rho_{i,j}\cdot
\rho_{i,k}^{-1}\cdot \tau_i\]
as in Proposition~\ref{prop:cocyc-factorize}.
\end{itemize}
Let us write $\bfY\to \bfW$ for this Abelian extension that coordinatizes $\xi_i$ over the above factor map.

Now let $U$ be a compact metrizable Abelian fibre repository for the groups $A_\star$, let $\chi_n$, $n\geq 1$ be an enumeration of $\hat{U}$, and for each $n$ let $\chi_{n,\star} := \chi_n|_{A_\star}$ be the resulting measurable family of characters in $\hat{A_\star}$.  Since the restrictions $\chi_n|_{A_s}$ span the whole of $L^2(m_{A_s})$ for each $s$ (because $L^2(m_{A_s}) = \ol{\rm{span}}\,\hat{A_s}$ and any member of $\hat{A_s}$ is the restriction of some character on $U$), it follows that the family of fibrewise factors
\begin{multline*}
\id_W\ltimes \chi_{n,\star}:(W \ltimes A_\star, \theta \ltimes m_{A_\star},R\ltimes \s_\star)\\ \to (W \ltimes \chi_{n,\star}(A_\star), \theta \ltimes m_{\chi_{n,\star}(A_\star)},R\ltimes \chi_{n,\star}(\s_\star))
\end{multline*}
generates the whole of $\bfY \cong (W \ltimes A_\star,\theta \ltimes m_{A_\star}, R \ltimes \s_\star)$.  This defines a family of factors $\k_n:\bfY\to \bfY_n$ generating $\bfY_n$, and hence we obtain $\xi_i \simeq \bigvee_{n\geq 1}(\k_n\circ \xi_i)$.

We complete the proof by showing that each $\bfY_n$ is a factor of a member of $\sfZ_0^{\bf{p}_i}\vee\sfZ_0^{\bf{p}_i - \bf{p}_j}\vee\sfZ_0^{\bf{p}_i - \bf{p}_k}\vee \sfZ_{\nil,2}^{\bbZ^2}$: since this is an idempotent class, it then follows that $\bfY$ is also a factor of a member of this class, and hence that there is a further extension $\t{\pi}:\t{\bfX}\to \bfX_0$ such that $\xi_{i,0}\circ\t{\pi}$ is contained in
\[(\sfZ_0^{\bf{p}_i}\vee \sfZ_0^{\bf{p}_i - \bf{p}_j}\vee\sfZ_0^{\bf{p}_i - \bf{p}_k}\vee \sfZ_{\nil,2}^{\bbZ^2})\t{\bfX}.\]
In view of~(\ref{eq:noncont}) this gives the desired contradiction with the $(\sfZ_0^{\bf{p}_i}\vee \sfZ_0^{\bf{p}_i - \bf{p}_j}\vee\sfZ_0^{\bf{p}_i - \bf{p}_k}\vee \sfZ_{\nil,2}^{\bbZ^2})$-satedness of $\bfX_0$.

To find our extension of $\bfY_n$ we now call on the further factorization given by  Proposition~\ref{prop:cocyc-factorize}. After adjoining the enlarged Kronecker system $(\t{Z}_\star,m_{\t{Z}_\star},R_{\t{\phi}_\star})\to (Z_\star,m_{Z_\star},R_{\phi_\star})$ relatively independently to $\bfX$ if necessary, that proposition gives a motionless selection of integers $m_\star\geq 1$ and a new factorization of $\chi_{n,\star}\circ\tau_i^{(m_\star)}$. From the ingredients of that new factorization, we can combine $\beta$ with $\chi_{n,\star}\circ b_i$ and $\tau_{i,j}$ with $\chi_{n,\star}\circ\rho_{i,j}^{(m_\star)}$ to give
\[\chi_{n,\star}\circ\s_i(m_\star\bf{p}_i,\,\cdot\,) = (\DDelta_{T|_{\a_i}^{\bf{p}_i}}b'_i)\cdot\rho'_{i,j}\cdot
(\rho'_{i,k})^{-1}\cdot \tau'_i,\]
where now
\begin{itemize}
\item $\rho'_{i,j}$ is $T^{\bf{p}_i - \bf{p}_j}$-invariant and $\rho'_{i,k}$ is $T^{\bf{p}_i - \bf{p}_k}$-invariant, and
\item $\tau'_i$ is lifted from $Z_\star$, takes values in $\Sone$ and is such that the two-step Abelian transformation $R_{\phi_s(m_s\bf{p}_i)}\ltimes \tau'_i\actson (z_0Z_{0,s})\times \Sone$ is a two-step nilrotation for every coset $z_0Z_{0,s}\leq Z_s$ for almost every $s$.
\end{itemize}

To use this factorization of $\chi_{n,\star}\circ \s_i(m_\star\bf{p}_i,\cdot)$ we must decompose $\bfY_n$ a little further.  For each $m\geq 1$ let $S_m := \{s \in S:\ m_s = m\}$,
\[C_{s,m} := \left\{\begin{array}{ll}\chi_{n,s}(A_s)&\quad\hbox{if}\ s\in S_m\\ (0)&\quad\hbox{else,}\end{array}\right.\]
and let $r_{s,m}:\chi_s(A_s) \longrightarrow C_{s,m}$ be the fibrewise quotient map which for each $s$ equals either $\id_{\chi_s(A_s)}$ or the zero map accordingly.  Since $\bigcup_{m=1}^\infty S_m$ is a disjoint union of measurable sets that is $\nu$-conegligible in $S$, the fibrewise quotient factors $\id_W\ltimes r_{\star,m}:W\ltimes \chi_{n,\star}(A_\star)\to W\ltimes C_{\star,m}$ together generated $\bfY_n$.  Letting $\bfY_{n,m}$ be the targets of these factors, it will therefore suffices to show that each of these is individually a factor of a member of $\sfZ_0^{\bf{p}_i}\vee\sfZ_0^{\bf{p}_i - \bf{p}_j}\vee\sfZ_0^{\bf{p}_i - \bf{p}_k}\vee \sfZ_{\nil,2}^{\bbZ^2}$.

However, now an appeal to Lemmas~\ref{lem:make-nilrtn-global} and~\ref{lem:make-subjoining} with $C_\star:= C_{\star,m}$ gives that each $\bfY_{n,m}$ is a factor of a member of
$\sfZ_0^{m\bf{p}_i}\vee \sfZ_0^{\bf{p}_i - \bf{p}_j}\vee\sfZ_0^{\bf{p}_i - \bf{p}_k}\vee \sfZ_{\nil,2}^{\bbZ^2}$, and Lemma~\ref{lem:simplifying-part-invar} shows that in fact this class is simply equal to $\sfZ_0^{\bf{p}_i}\vee \sfZ_0^{\bf{p}_i - \bf{p}_j}\vee\sfZ_0^{\bf{p}_i - \bf{p}_k}\vee \sfZ_{\nil,2}^{\bbZ^2}$, so this completes the proof. \qed

\subsection{Reduction of the cocycle factorization to the ergodic case}\label{subs:reduce-to-erg}

\begin{prop}\label{prop:ergodic-case-enough}
If the conclusion of Proposition~\ref{prop:cocyc-factorize} holds for almost every individual fibre group, then it holds in general.
\end{prop}

\textbf{Proof}\quad Let $(S,\nu)$ be the standard Borel probability space indexing the family $Z_\star$.  Our assumption is that a suitable extension $q_s:(\t{Z}_s,R_{\t{\phi}_s})\to (Z_s,R_{\phi_s})$, integer $m_s$ and factorization
\[\chi_s\circ\tau_i^{(m_s)}\circ q_s = \tau_{s,i,j}\cdot\tau_{s,i,k}\cdot\tau_{s,\rm{nil},i}\cdot \DDelta_{\t{\phi}_s(\bf{p}_i)}\b_s\]
exist for $\nu$-almost all $s$ separately, and we must show that these data can be chosen so that the family $(\t{Z}_s,\t{\phi}_s)$ and integers $m_s$ are measurable in $s$ and each of $\tau_{s,i,j}$, $\tau_{s,i,k}$, $\tau_{s,\nil,i}$ and $\b_s$ is $\nu$-almost surely the restriction to the fibre $\{s\}\times \t{Z}_s$ of some measurable map $S\ltimes \t{Z}_\star\to\Sone$ or $\bbZ^2\times (S\ltimes \t{Z}_\star) \to\Sone$.

In essence this follows from an appeal to the Measurable Selection Theorem (in the form of Theorem 2.2 in~\cite{Aus--ergdirint}, for example), but we must be quite careful in how we handle the measurability issues resulting from the additional variability in the domain $Z_s$.  The key idea is to identify the data
\[q_s:(\t{Z}_s,\t{\phi}_s)\to (Z_s,\phi_s),\quad \b_s:\t{Z}_s\to\Sone,\quad\tau_{s,i,j},\tau_{s,i,k},\tau_{s,\nil,i}:\bbZ^2\times \t{Z}_s\to \Sone\]
with a sequence of approximating data that involve only compact Abelian Lie groups.  This will carry the advantage that compact Abelian Lie groups have only countably many closed subgroups, which will clarify some of those measurability issues.

By symmetry, let us now assume that $(i,j,k) = (1,2,3)$, and to lighten notation in the remainder of this subsection let $\s_s := \chi_s\circ \tau_i\circ q_s$ and write $\tau_{s,\ell}$ in place of $\tau_{s,1,\ell}$ for $\ell =2,3$ and $\tau_{s,\nil}$ in place of $\tau_{s,\nil,1}$.

\quad\textbf{Step 1}\quad Fix some compact metrizable Abelian fibre repository $U$ for the family $Z_\star$.  By embedding $U$ into $(\Sone)^\bbN$ using an enumeration of $\hat{U}$, we may simply assume that $U = (\Sone)^\bbN$. Since any compact metrizable Abelian extension of $Z_s$ may be written as a closed subgroup of $(\Sone)^\bbN\times Z_s$ with extension epimorphism given by the second coordinate projection (this time by enumerating $\hat{\t{Z}_s}\cap Z_s^\perp$, for example), it follows from our assumptions that for each $s$ separately we may realize $\t{Z}_s$ as a closed subgroup of $\t{U} := (\Sone)^\bbN\times (\Sone)^\bbN$ which projects onto $Z_s \leq (\Sone)^\bbN$ under the second coordinate projection.

Let $Q_{1,N}:(\Sone)^\bbN\to (\Sone)^N$ be the projection onto the first $N$ coordinates and let $Q_N := (Q_{1,N},Q_{1,N}):\t{U}\to (\Sone)^N\times (\Sone)^N$, so that this is an inverse sequence of quotients generating the whole of $\t{U}$.  Given these we may also define finite-level connecting epimorphisms $Q^N_K = (Q^N_{1,K},Q^N_{1,K}):(\Sone)^N\times (\Sone)^N\to (\Sone)^K\times (\Sone)^K$ when $N\geq K$.

\quad\textbf{Step 2}\quad In order to make use of these approximating Lie groups, let $\L_N$ be the set of all nonuples
\[(m,Z,\t{Z},\t{\phi},\rho,\b,\tau_2,\tau_3,\tau_\nil)\]
such that
\begin{itemize}
\item $m$ is a positive integer;
\item $Z\leq (\Sone)^N$, $\t{Z}\leq (\Sone)^N\times (\Sone)^N$ and $\t{Z}$ projects onto $Z$ under the second coordinate projection;
\item $\t{\phi}:\bbZ^2\to \t{Z}$ is a homomorphism;
\item $\rho:Z\to \Sone$, $\b:\t{Z}\to \Sone$ and $\tau_2,\tau_3,\tau_\nil:\bbZ^2\times \t{Z}\to\Sone$ are Haar-a.e. equivalence classes of Borel maps.
\end{itemize}
Since $(\Sone)^N\times (\Sone)^N$ has only countably many closed subgroups, we may partition $\L_N$ into countably many subsets according to the subgroup $\t{Z}$ appearing in the sequence, and now for each possible $\t{Z}$ the corresponding subset of $\L_N$ may simply be identified with
\[\bbN\times \Hom(\bbZ^2,\t{Z})\times \C(Z)\times \C(\t{Z})\times \C(\bbZ^2\times \t{Z})^3.\]
Regarding $\C(Z)$ and its cousins with their usual Polish topologies, we consider this product endowed with the product topology and its Borel $\s$-algebra and now piece these $\s$-algebras together to obtain a Borel structure on the whole of $\L_N$.  This defines a standard Borel structure on $\L_N$ because there are only countably many pieces.

Let $\O := \prod_{N\geq 1}\L_N$ and $\O_N := \prod_{K\leq N}\L_K$, so there is a natural projection map $\O\to\O_N$ for each $N$.  Consider $\O$ and each $\O_N$ endowed with its product Borel structure.

\quad\textbf{Step 3}\quad Inside $\O$ we now define the subset $\O_{\rm{conv}}\subseteq \O$ to comprise those sequences
\[\big((m_N,Z_N,\t{Z}_N,\t{\phi}_N,\rho_N,\b_N,\tau_{N,2},\tau_{N,3},\tau_{N,\nil})\big)_{N\geq 1}\]
such that the following hold:
\begin{itemize}
\item $m_N$ does not depend on $N$;
\item we have $Q^N_K(\t{Z}_N) = \t{Z}_K$ for each $N\geq K$, so that we may define $\t{Z}_\infty := \bigcap_{N\geq 1}Q_N^{-1}(\t{Z}_N)\leq \t{U}$ and $Z_\infty := \bigcap_{N\geq 1}Q_{1,N}^{-1}(Z_N)\leq U$, and observe that $\t{Z}_\infty$ projects onto $Z_\infty$ under the second coordinate projection and that $Q_N(\t{Z}_\infty) = \t{Z}_N$ and $Q_{1,N}(Z_\infty) = Z_N$ for all $N$;
\item similarly, we have $Q^N_K\circ\t{\phi}_N = \t{\phi}_K$ whenever $N\geq K$, so that we may unambiguously define $\t{\phi}_\infty\in\Hom(\bbZ^2,\t{Z}_\infty)$ by letting $\t{\phi}_\infty(\bf{n})$ be the unique element of $\bigcap_{N\geq 1}Q_N^{-1}\{\t{\phi}_N(\bf{n})\}$, and also define $\phi_\infty\in\Hom(\bbZ^2,Z_\infty)$ to be the composition of $\t{\phi}_\infty$ with the projection $\t{Z}_\infty\to Z_\infty$;
\item each of the function sequences converges in probability, in the sense that for any $\eps > 0$ there is some $K\geq 1$ such that
\[m_{\t{Z}_N}\{|\b_N - \b_K\circ Q^N_K| > \eps\} < \eps\quad\quad\forall N\geq K,\]
where the additive difference $\b_N - \b_K\circ Q^N_K$ is understood as a difference of two complex numbers, and similarly for the $\rho_N$, $\tau_{N,i}$ and $\tau_{N,\nil}$;
\item the functions $\bbZ^2\times \t{Z}_N\to\Sone$ are asymptotically cocycles, in the sense that for any $\bf{n},\bf{m} \in \bbZ^2$ and $\eps > 0$ there is some $K\geq 1$ such that
\[m_{\t{Z}_N}\{|\DDelta_{\t{\phi}_N(\bf{n})}\tau_{N,i}(\bf{m},\cdot) - \DDelta_{\t{\phi}_N(\bf{m})}\tau_{N,i}(\bf{n},\cdot)| > \eps\} < \eps\quad\quad\forall N\geq K,\]
and similarly for $\tau_{N,\nil}$;
\item moreover, the functions $\tau_{N,\nil}$ actually stabilize at some finite level $K$, in the sense that $\tau_{N,\nil} = \tau_{K,\nil}\circ Q^N_K$ Haar-a.e. for all $N\geq K$, and the stable value $\tau_{K,\nil}$ is a cocycle $\bbZ^2\times \t{Z}_K\to\Sone$ over $R_{\t{\phi}_N}$.
\end{itemize}
Clearly $\O_{\rm{conv}}$ is a Borel subset of $\O$.  The third condition implies that the lifts $\b_N\circ Q_N:\t{Z}_\infty\to\Sone$ converge in probability on $\t{Z}_\infty$, and similarly for the other maps, to some Borel maps $\rho_\infty:Z_\infty\to\Sone$, $\b_\infty:\t{Z}_\infty\to\Sone$ and $\tau_{\infty,2},\tau_{\infty,3}:\bbZ^2\times \t{Z}_\infty\to\Sone$, where these latter are cocycles over $R_{\t{\phi}_\infty}$.  In addition, simply by lifting the stable value $\tau_{K,\nil}:\bbZ^2\times \t{Z}_K\to\Sone$ we obtain a cocycle $\tau_{\infty,\nil}:\bbZ^2\times \t{Z}_\infty\to\Sone$ over $R_{\t{\phi}_\infty}$.  We will refer to these groups, homomorphisms and maps as the \textbf{limit data} of the sequence $\big((m_N,Z_N,\t{Z}_N,\t{\phi}_N,\rho_N,\b_N,\tau_{N,2},\tau_{N,3},\tau_{N,\nil})\big)_{N\geq 1}$.

Conversely, given a subgroup $\t{Z}_\infty \leq (\Sone)^\bbN\times (\Sone)^\bbN$, a homomorphism $\t{\phi}_\infty:\bbZ^2\to \t{Z}_\infty$ and maps $\rho_\infty$, $\b_\infty$, $\tau_{\infty,2}$, $\tau_{\infty,3}$ and $\tau_{\infty,\nil}$ satisfying all the conditions listed above, then since any Borel map on a compact group may be approximated in probability by a map lifted from a Lie quotient group it follows that there is some sequence \[\big((m_N,Z_N,\t{Z}_N,\t{\phi}_N,\rho_N,\b_N,\tau_{N,2},\tau_{N,3},\tau_{N,\nil})\big)_{N\geq 1} \in \O_{\rm{conv}}\] giving rise to them as its limit data.

\quad\textbf{Step 4}\quad We now make use of the Borel sets of canonical nil-cocycles $\A(\cdot)$ obtained in Lemma~\ref{lem:canon-nil-cocyc-measble}. We consider the further subset $\O_{\rm{final}} \subseteq \O_{\rm{conv}}$ comprising those convergent sequences $\big((m_N,Z_N,\t{Z}_N,\t{\phi}_N,\rho_N,\b_N,\tau_{N,2},\tau_{N,3},\tau_{N,\nil})\big)_{N\geq 1}$ such that the limiting cocycles satisfy
\[\tau_{\infty,2}(\bf{p}_1 - \bf{p}_2,\cdot) = \tau_{\infty,3}(\bf{p}_1 - \bf{p}_3,\cdot) \equiv 1\]
and that if $\tau_{K,\nil}$ has stabilized at level $K$ then $(\t{\phi}_K,\tau_{K,\nil}) \in \A(\t{Z}_K)$, the class of pairs introduced in Lemma~\ref{lem:canon-nil-cocyc-measble}. It follows from Lemma~\ref{lem:canon-nil-cocyc-measble} that $\O_{\rm{final}}$ is a Borel subset of $\O_{\rm{conv}}$.

\quad\textbf{Step 5}\quad The importance of these preliminaries is that they provide a standard Borel space $\O_{\rm{final}}$ from which we can make measurable selections.  Let $E \subseteq S\times \O_{\rm{final}}$ comprise those pairs
\[\big(s,\big((m_N,Z_N,\t{Z}_N,\t{\phi}_N,\rho_N,\b_N,\tau_{N,2},\tau_{N,3},\tau_{N,\nil})\big)_{N\geq 1}\big)\]
such that
\begin{itemize}
\item (smaller subgroups behave as they should) for each $N$ we have $Z_N = Q_{1,N}(Z_s)$ and $\t{\phi}_N$ projects to $Q_{1,N}\circ \phi_s$ under composition with $\t{Z}_N\to Z_N$;
\item (cocycles have desired limits) we have $\rho_\infty = \s_s(m_N\bf{p}_1,\cdot)$.
\end{itemize}

Once again, $E$ is easily seen to be Borel in $S\times \O_{\rm{final}}$, and our initial assumption pertaining to $\nu$-a.e. individual $s$ promises that the fibre $E\cap (\{s\}\times\O_{\rm{final}})$ is $\nu$-almost surely nonempty.  Hence an appeal to the Measurable Selection Theorem now gives a $\nu$-conegligible subset $S_0\subseteq S$ and a Borel selection
\[s\mapsto \big(s,\big((m_s,Z_{s,N},\t{Z}_{s,N},\t{\phi}_{s,N},\rho_{s,N},\b_{s,N},\tau_{s,N,2},\tau_{s,N,3},\tau_{s,N,\nil})\big)_{N\geq 1}\big)\]
defined for $s \in S_0$.

Finally we let $\t{Z}_\star$ be the measurable family $s\mapsto \bigcap_{N\geq 1}Q_N^{-1}(\t{Z}_{s,N})$, $q_\star:(\t{Z}_\star,\t{\phi}_\star)\to (Z_\star,\phi_\star)$ the coordinate projection, and observe that the sequences
\[\b_{\star,N}\circ Q_N,\quad \tau_{\star,N,i}\circ Q_N\quad\hbox{and}\quad \tau_{\star,N,\nil}\circ Q_N,\]
regarded as maps $S\ltimes \t{Z}_\star\to \Sone$ and $\bbZ^2\times (S\ltimes \t{Z}_\star)\to \Sone$, also converge in probability for $\nu\ltimes m_{\t{Z}_\star}$ and their limits define the map $\b_\star$ and cocycles $\tau_{\star,i}$ and $\tau_{\star,\nil}$ required for Proposition~\ref{prop:cocyc-factorize}, together with the measurable selections of integers $m_\star$.  This completes the proof. \qed

\subsection{Another consequence of satedness}

We now make a slight detour to introduce a property of certain direct integrals of Kronecker systems that we will need later and that seems to merit its own subsection, and show how it can be deduced from
the FIS property.

\begin{dfn}[DIO system]
A direct integral of $\bbZ^d$-group rotations
$(U_\star,m_{U_\star},\psi_\star)$ with invariant base space
$(S,\nu)$ has the \textbf{disjointness of independent orbits
property}, or is \textbf{DIO}, if for subgroups $\G_1,\G_2 \leq
\bbZ^d$ we have
\[\G_1\cap \G_2 =
\{\bs{0}\}\quad\quad\Rightarrow\quad\quad
\overline{\phi_s(\G_1)}\cap \overline{\phi_s(\G_2)} =
\{1_{Z_s}\}\quad \hbox{for}\ \nu\hbox{-a.e.}\ s.\]
\end{dfn}

\begin{prop}\label{prop:FISimpliesDIO}
If $\bfX$ is an FIS $\bbZ^d$-system then the factor $\bfZ_1^T =
\sfZ_1^{\bbZ^d}\bfX$, the maximal factor of $\bfX$ that can be
coordinatized as a direct integral of group rotations, is such that
for any subgroups $\G_1,\G_2 \leq \bbZ^d$ with trivial intersection
we have
\[\zeta_1^T \leq (\zeta_1^T\wedge \zeta_0^{T^{\ \uhr\G_1}})\vee(\zeta_1^T\wedge \zeta_0^{T^{\ \uhr\G_2}}),\]
and $\bfZ_1^T$ is DIO.
\end{prop}

\textbf{Proof}\quad Let $\pi:\bfX \to
(Z_\star,m_{Z_\star},R_{\phi_\star})$ be a coordinatization of
$\zeta_1^T:\bfX\to\bfZ_1^T$, say with invariant base space
$(S,\nu)$. Fix $\G_1, \G_2 \leq \bbZ^d$ with trivial intersection
and let $\G:= \G_1 + \G_2$. First note that if $\G$ has infinite
index in $\bbZ^d$ then we can choose another subgroup $\L \leq
\bbZ^d$ that is a complement to the radical
\[\rm{rad}\,\G:= \{\bf{n} \in \bbZ^d:\ k\bf{n} \in \G\ \hbox{for some}\ k\in\bbZ\setminus\{0\}\},\]
so that now $\G_1 \cap (\G_2 + \L) = \{\bs{0}\}$ and $\G_1 + \G_2 +
\L$ has finite index in $\bbZ^d$; and so simply by replacing $\G_2$
with $\G_2 + \L$ if necessary it suffices to treat the case in which
$\G$ has finite index in $\bbZ^d$.

The remainder of the proof breaks into two steps.

\quad\textbf{Step 1}\quad We first observe that any direct integral
of $\bbZ^d$-group rotations $(U_\star,m_{U_\star},\psi_\star)$
(which we may assume has ergodic fibres) is a
$(\sfZ_0^{\G_1},\sfZ_0^{\G_2})$-subjoining.

Let us first see this when $\G_1 + \G_2 = \bbZ^d$, so that we may
express $\bbZ^d = \G_1 \oplus \G_2$ and let $\rm{proj}_i:\bbZ^d\to
\G_i$ be the resulting coordinate projections. In this case the
construction is very simple: for by the ergodicity of the fibres we
have $\overline{\psi_s(\G_1)} + \overline{\psi_s(\G_2)} = U_s$
almost surely, and now we can define the extension of direct
integrals of group rotations
\[(U_{1,\star},m_{U_{1,\star}},\psi_{1,\star})\to (U_\star,m_{U_\star},\psi_\star)\]
with the same invariant base space $(S,\nu)$ by setting \[U_{1,s} :=
\overline{\psi_s(\G_1)} \times \overline{\psi_s(\G_2)}\] and
$q_{1,s}:U_{1,s}\onto U_{0,s}:(u,v)\mapsto uv$ and defining the
extended homomorphism by
\[\psi_{1,s}(\bf{n}) = (\psi_s(\rm{proj}_1(\bf{n})),\psi_s(\rm{proj}_2(\bf{n})))\]
(all of these specifications being manifestly still measurable in
$s$).  The extended system is now clearly a joining of the systems
\[(\overline{\psi_\star(\G_i)},m_{\overline{\psi_\star(\G_i)}},\psi_\star\circ\rm{proj}_i)\]
for $i=1,2$, each of which has trivial $\G_{3-i}$-subaction.

If $\G_1 + \G_2$ is a proper subgroup of $\bbZ^d$ then we must work
a little harder. We can treat this case abstractly by first
constructing a suitable extension for the subaction $R_{\psi}^{\
\uhr\G_1 + \G_2}$ using the argument above, and then constructing a further extension to recover an action of the whole of $\bbZ^d$, such as an FP extension as in Subsection 3.2 of~\cite{Aus--lindeppleasant1}, which is easily seen to retain the desired disjointness of orbit closures and to give another direct integral of group rotations. However, for clarity let us describe a suitable construction a little more explicitly in the present setting.

Let $K_{i,s}:=\overline{\psi_s(\G_i)}$ and $K_s :=
\overline{\psi_s(\G)} = K_{1,s} \times K_{2,s}$, let $\O \subseteq
\bbZ^d$ be a fundamental domain for the finite-index subgroup $\G$
and let $\{\cdot\}:\bbZ^d\to \O$, $\lfr\cdot\rfr:\bbZ^d\to \G$ be
the corresponding fractional- and integer-part maps. Let us also decompose $\lfr\cdot\rfr$ further as $\lfr\cdot\rfr_1 + \lfr\cdot\rfr_2$ with $\lfr\cdot\rfr_i:\bbZ^d\to
\G_i$ (clearly having chosen $\O$ there is a unique such
decomposition). Finally let $w_{s,\omega} := \psi_s(\omega) \in U_s$
for $\omega \in \O$.

Now consider the map $q_{1,s}:U_{1,s} := K_{1,s}\times K_{2,s}\times
(\bbZ^d/\G)\to U_s$ given by
\[(u,v,\bf{m} + \G) \mapsto u\cdot v\cdot w_{s,\{\bf{m}\}}.\]
This is easily seen to be onto, because the original homomorphism
$\psi_s$ was dense.  On $S\ltimes U_{1,\star}$ we define the
$\bbZ^d$-action $R_1$ by
\begin{multline*}
R_1^\bf{n}:(s,u,v,\bf{m} + \G)\\ \mapsto (s,\psi_s(\lfr\bf{n} +
\bf{m}\rfr_1 - \lfr\bf{m}\rfr_1)u,\psi_s(\lfr\bf{n}+\bf{m}\rfr_2 -
\lfr\bf{m}\rfr_2)v,\bf{m} + \bf{n} + \G)
\end{multline*}
(it is easily checked that the right-hand side here depends only on
the class $\bf{m} + \G$, so this is a well-defined action), and now
we see that
\begin{eqnarray*}
&&q_{1,s}(R_1^\bf{n}(s,u,v,\bf{m} + \G))\\ &&= \psi_s(\lfr\bf{n} +
\bf{m}\rfr_1 - \lfr\bf{m}\rfr_1)\cdot u\cdot
\psi_s(\lfr\bf{n}+\bf{m}\rfr_2 - \lfr\bf{m}\rfr_2)\cdot v\cdot
w_{s,\{\bf{m} + \bf{n}\}}\\
&&= \psi_s(\lfr\bf{n} + \bf{m}\rfr - \lfr\bf{m}\rfr)\cdot u\cdot
v\cdot \psi_s(\{\bf{n} + \bf{m}\} - \{\bf{m}\})\cdot w_{s,\{\bf{m}\}}\\
&&= \psi_s(\bf{n})\cdot(u\cdot v\cdot w_{s,\{\bf{m}\}}).
\end{eqnarray*}
Thus $q_{1,\star}:(U_{1,\star},m_{U_{1,\star}},R_1)\to
(U_\star,m_{U_\star},\psi_\star)$ defines an extension of $\bbZ^d$-systems. Since the subaction of the finite-index
subgroup group $\G_1 + \G_2$ simply acts by rotations inside each of
the $[\G_1 +\G_2 : \bbZ^d]$-many fibres of $K_{1,\star}\times
K_{2,\star}$ in $U_{1,\star}$, this subaction is actually a direct
integral of direct sums of group rotation actions and hence the
overall action is also a direct integral of group rotations. Finally
we observe that the fibrewise restriction of $R_1$ to the canonical
factor with fibres $K_{1,\star}\times (\bbZ^d/\G)$ has trivial
$\G_2$-subaction and its fibrewise restriction to the canonical
factor with fibres $K_{2,\star}\times (\bbZ^d/\G)$ has trivial
$\G_1$-subaction, so this extended system is a member of
$\sfZ_0^{\G_1}\vee \sfZ_0^{\G_2}$.

\quad\textbf{Step 2}\quad Since we assume that $\bfX$ is
$(\sfZ_0^{\G_1}\vee \sfZ_0^{\G_2})$-sated, Step 1 now implies that
$\pi \precsim \zeta_0^{T^{\ \uhr\G_1}}\vee \zeta_0^{T^{\
\uhr\G_2}}$.  On the other hand $\zeta_0^{T^{\ \uhr\G_1}}$ and
$\zeta_0^{T^{\ \uhr\G_2}}$ are relatively independent over
$\zeta_0^{T^{\ \uhr(\G_1 + \G_2)}}$, and since $\G_1 + \G_2$ has
finite index in $\bbZ^2$ this in turn is simply an extension of
$\zeta_0^T$ by finite group rotations that factorize through the
quotient map $\bbZ^d\to \bbZ^d/\G$. By the non-ergodic
Furstenberg-Zimmer Theorem~\ref{thm:rel-ind-joinings} it follows
that $\pi$ and $\zeta_0^{T^{\ \uhr\G_1}}\vee \zeta_0^{T^{\
\uhr\G_2}}$ are relatively independent under $\mu$ over
\[(\zeta^T_{1/\zeta_0^{T^{\ \uhr\G}}}\wedge \zeta_0^{T^{\ \uhr\G_1}})\vee (\zeta^T_{1/\zeta_0^{T^{\ \uhr\G}}}\wedge\zeta_0^{T^{\ \uhr\G_2}}),\]
so the above containment implies that $\pi$ is actually contained in
this join.

However, again since $\G$ has finite index in $\bbZ^d$ and any
compact extension of a \emph{finite} group rotation system is still
compact, we must in fact have $\zeta^T_{1/\zeta_0^{T^{\ \uhr\G}}} =
\zeta_1^T$, and so we have deduced the first desired conclusion that
$\pi$ is contained in the join of its further $\G_1$- and
$\G_2$-invariant factors. Since these are coordinatized by
the fibrewise quotient maps
\[S\ltimes Z_\star \to S\ltimes (Z_\star/\overline{\phi_\star(\G_i)}):(s,z)\mapsto (s,z\overline{\phi_s(\G_i)})\quad\quad\hbox{for}\ i=1,2,\]
in order for these to generate the whole of $Z_s$ above $\nu$-almost
every $s$ it must hold that the cosets $z\overline{\phi_s(\G_1)}$
and $z\overline{\phi_s(\G_2)}$ together uniquely determine $z \in
Z_s$ for almost every $s$, or equivalently that
\[\overline{\phi_s(\G_1)}\cap \overline{\phi_s(\G_2)} = \{1_{Z_s}\}\quad\quad\hbox{for}\ \nu\hbox{-almost every}\ s,\]
as required. \qed

\textbf{Example}\quad Although the DIO property will shortly prove very useful, ascending to a DIO extension can make
a very simple initially-given group rotation system $(U_\star,m_{U_\star},\phi_\star)$
into a very much more complicated extension
$(\t{U}_\star,m_{\t{U}_\star},\t{\phi}_\star)$.  For example,
letting $w \in \Sone$ be an irrational rotation and
$\phi:\bbZ^2\longrightarrow (\Sone)^2 =: U_0$ be the homomorphism
$(m,n)\mapsto (w^m,w^n)$, we can build a DIO extension by choosing a sequence
$((m_{i1},m_{i2}),(n_{i1},n_{i2}))_{i\geq 1}$ of linearly
independent pairs of members of $\bbZ^2$, and then constructing the
inverse sequence of systems
$\big((U_{(i)},m_{U_{(i)}},\phi_{(i)})\big)_{i\geq 0}$,
$(q^{(i)}_{(j)})_{i\geq j\geq 0}$ recursively so that given
$U_{(i)}$ the map $q^{(i+1)}_{(i)}$ sends $(s,t)$ to
$(s^{m_{i1}}t^{m_{i2}},s^{n_{i1}}t^{n_{i2}})$.  It is easy to see
that this construction gives rise to a sequence of surjective
endomorphisms of $U_{(i)}\cong (\Sone)^2$ (in which
$q^{(i+1)}_{(i)}$ has covering number
$\Big|\det\Big(\begin{array}{cc}m_{i1}& m_{12}\\
n_{i1}& n_{i2}\end{array}\Big)\Big|$, in particular), but that the
resulting inverse limit group is an extremely complicated beast
indeed.  A more detailed discussion of such inverse limit
constructions can be found in Rudolph's paper~\cite{Rud93}. \fin

The importance of Proposition~\ref{prop:FISimpliesDIO} for our study
of $\bbZ^2$-systems is that it substantially simplifies our picture
of the joinings of direct integrals of group rotations
\[(\zeta_1^T\wedge \zeta_0^{T^{\bf{p}_i - \bf{p}_j}})\vee (\zeta_1^T\wedge \zeta_0^{T^{\bf{p}_i - \bf{p}_k}}),\]
and also their overall joining under $\mu^{\rm{F}}$, that underly
the new maps $\tau_i$ that appear in the factorization of
Proposition~\ref{prop:cocycle-prod-decomp}.

Indeed, we have just seen that for an FIS system each of the above
factors simply equals $\zeta_1^T$. Letting $\zeta_1^T:\bfX\to
(Z_\star,m_{\star},R_{\phi_\star})$ be a coordinatization as above with
invariant base space $(S,\nu)$ (so $(S,\nu)$ can
be identified with $\bfZ_0^T$), it follows from the definition of $\mu^\rm{F}$ that its restriction to the factor $(\zeta_1^T)^{\times 3}$ is the joining-limit as $N\to\infty$ of the measures
\begin{multline*}
\frac{1}{N}\sum_{n=1}^N m_{(\phi(n\bf{p}_1),\phi(n\bf{p}_2),\phi(n\bf{p}_3))\cdot\{(z,z,z):\ z \in Z\}}\\ = \frac{1}{N}\sum_{n=1}^N m_{(1,\phi(n(\bf{p}_2 - \bf{p}_1)),\phi(n(\bf{p}_3-\bf{p}_1)))\cdot\{(z,z,z):\ z \in Z\}}\\ \to m_{\ol{\{(1,\phi(n(\bf{p}_2 - \bf{p}_1)),\phi(n(\bf{p}_3 - \bf{p}_1))):\ n\in\bbZ\}}\cdot \{(z,z,z):\ z\in Z\}}.
\end{multline*}

Clearly if
\[(1,u_2,u_3) \in \ol{\{(1,\phi(n(\bf{p}_2 - \bf{p}_1)),\phi(n(\bf{p}_3 - \bf{p}_1))):\ n\in\bbZ\}}\leq Z^3\]
then $u_i \in \ol{\phi(\bbZ(\bf{p}_i - \bf{p}_1))}$ for $i=2,3$ and $u_2u_3^{-1}\in \ol{\phi(\bbZ(\bf{p}_3 - \bf{p}_2))}$.  On the other hand, given any $(1,u_2,u_3)$ satisfying these constraints, if we choose $n_i \in \bbZ$ so that $\phi(n_i(\bf{p}_2 - \bf{p}_1))\to u_2$ as $i\to\infty$ and pass to a subsequence so that $\phi(n_i(\bf{p}_3 - \bf{p}_1))$ also converges, say to $v$, then we see that $vu_3^{-1} = (vu_2^{-1})(u_2u_3^{-1})$ must lie in $\ol{\phi(\bbZ(\bf{p}_3 - \bf{p}_1)}\cap \ol{\phi(\bbZ(\bf{p}_3 - \bf{p}_2))}$.  By the DIO property this is $\{1\}$, so $v = u_3$ and $(1,u_2,u_3)$ is in our subgroup.

Hence given the DIO property the restriction of $\mu^{\rm{F}}$ is simply $m_{\vec{Z}}$ for
\[\vec{Z} = \{(z_1,z_2,z_3)\in Z^3:\ z_i\ol{\phi(\bbZ(\bf{p}_i - \bf{p}_j))} = z_j\ol{\phi(\bbZ(\bf{p}_i - \bf{p}_j))}\ \forall i\neq j\}.\]

We can now deduce another important consequence of the DIO property.

\begin{cor}\label{cor:simple-rewriting-group-rotations}
If $(Z,m_Z,\phi)$ is a $\bbZ^2$-group rotation having the DIO property and $\bf{n}_1,\bf{n}_2 \in \bbZ^2$ are linearly independent then the extension of group rotations
\[q:(Z,m_Z,R_{\phi(\bf{n}_1)}) \to (Z/\ol{\phi(\bbZ\bf{n}_2)},m_{Z/\ol{\phi(\bbZ\bf{n}_2)}},R_{\phi(\bf{n}_1)\ol{\phi(\bbZ\bf{n}_2)}})\]
is relatively invariant.
\end{cor}

\textbf{Proof}\quad Since $\bf{n}_1$ and $\bf{n}_2$ are linearly independent the DIO property gives \[\ol{\phi(\bbZ\bf{n}_1)}\cap \ker q = \ol{\phi(\bbZ\bf{n}_1)}\cap \ol{\phi(\bbZ\bf{n}_2)} = \{1\},\] and hence $q$ restricts to an isomorphism
\[q|_{\ol{\phi(\bbZ\bf{n}_1)}}:\ol{\phi(\bbZ\bf{n}_1)}\to \big(\ol{\phi(\bbZ\bf{n}_1)}\cdot \ol{\phi(\bbZ\bf{n}_2)}\big)/\ol{\phi(\bbZ\bf{n}_2)}.\]
Since the individual ergodic components of $R_{\phi(\bf{n}_1)}$ in $Z$ are simply the cosets of $\ol{\phi(\bbZ\bf{n}_1)}$, it follows that $q$ maps each of these isomorphically to a corresponding ergodic fibre of $R_{\phi(\bf{n}_1)\ol{\phi(\bbZ\bf{n}_2)}}$, as required. \qed

\subsection{Completion of the cocycle factorization in the ergodic case}\label{subs:factorizing-in-erg-case}

Let us finally pick up the thread that we
set down at the end of the
Subsection~\ref{subs:factorizing-cocycles}.

After ascending to an extended system $\bfX$ as given by
Proposition~\ref{prop:no-extraneous-dependences} and adopting the
factorization of Proposition~\ref{prop:cocycle-prod-decomp} for a
given triple $\bf{p}_1$, $\bf{p}_2$, $\bf{p}_3$, the resulting
cocycles $\tau_i$ satisfy a combined coboundary equation with
transfer function depending only on the joining under $\mu^\rm{F}$
of the direct integrals of group rotations
$(\zeta_1^T\wedge\zeta_0^{T^{\bf{p}_i} = T^{\bf{p}_j}})\circ\pi_i$. In the notation introduced at the end of the preceding subsection this equation reads
\begin{eqnarray}\label{eq:combcoc}
\tau_1(z_1)\cdot
\tau_2(z_2)\cdot \tau_3(z_3) =
\DDelta_{(w_1,w_2,w_3)}c(z_1,z_2,z_3)\quad\quad m_{\vec{Z}_\star}\hbox{-a.e.}\
(z_1,z_2,z_3)
\end{eqnarray}
for some Borel $c:\vec{Z}_\star\to A_\star$, where $w_i := \phi_\star(\bf{p}_i)$.

The above gives us an equation relating the restrictions of the cocycles $\tau_i$ to $Z_s$ for almost every $s$.  On the other hand, Proposition~\ref{prop:FISimpliesDIO} promises that the ergodic group rotation action $(Z_s,m_{Z_s},\phi_s)$ has the DIO property for almost every $s$; and by Proposition~\ref{prop:ergodic-case-enough} the conclusion of Proposition~\ref{prop:cocyc-factorize} will hold if we merely prove it for almost every $s$ individually, ignoring issues of measurability in $s$.  Therefore we may now assume that our overall system is ergodic, and that there is simply a single DIO group rotation $(Z,m_Z,R_\phi)$ in play.  This will both lighten the notation in the arguments to come, which will involve moving quickly among various subgroups of $Z$, and will save us the trouble of re-proving `measurably-varying' versions of a host of standard results from Moore's cohomology theory of locally compact groups (see Appendix~\ref{app:cohom}).

We will therefore now drop all mention of the invariant base space $(S,\nu)$, and will omit the subscript $_\star$ indicating motionless dependence from data such as the fibre group $A$.

Also, our proof of Proposition~\ref{prop:cocyc-factorize} will make use of equation~(\ref{eq:combcoc}) only after composing with some character $\chi \in \hat{A}$, after which the same equation is obtained for the resulting $\Sone$-valued maps. Hence it will suffice from this point on to consider $\Sone$-valued maps, and so to lighten the notation further we will henceforth treat each $\tau_i$ and $c$ as themselves $\Sone$-valued.

Our first step is to obtain solutions to an analog of the Conze-Lesigne equations, but which make only a weaker demand in that they are `directional'.  This step has close parallels with the extraction of the Conze-Lesigne equations in the study of characteristic factors for $\bbZ$-actions and related problems: see, in particular, Meiri~\cite{Mei90}, Rudolph~\cite{Rud93} and Furstenberg and Weiss~\cite{FurWei96}.

\begin{lem}\label{lem:lift-u}
If $u \in \ol{\phi(\bbZ(\bf{p}_1 - \bf{p}_3))}$ then there is some $u_2\in Z$ such that $(u,u_2,1) \in \vec{Z}$.
\end{lem}

\textbf{Proof}\quad Clearly 
\[\ol{\phi(\bbZ(\bf{p}_1 - \bf{p}_3))} \leq \ol{\phi(\bbZ(\bf{p}_1 - \bf{p}_2))} \cdot \ol{\phi(\bbZ(\bf{p}_2 - \bf{p}_3))},\]
and so $u$ may be expressed accordingly as $u_1u_2$.  This gives $u_2\ol{\phi(\bbZ(\bf{p}_2 - \bf{p}_3))} = \ol{\phi(\bbZ(\bf{p}_2 - \bf{p}_3))}$ and $u\ol{\phi(\bbZ(\bf{p}_1 - \bf{p}_2))} = u_2\ol{\phi(\bbZ(\bf{p}_1 - \bf{p}_2))}$, and so $(u,u_2,1) \in \vec{Z}$ as required. \qed

\begin{lem}\label{lem:processing-cocyc-over-Kron}
In the notation explained above, for each
ordering $\{i,j,k\} = \{1,2,3\}$ the Borel map $\tau_i$ has the property that for every $u \in \ol{\phi(\bbZ(\bf{p}_i - \bf{p}_j))}$ there are Borel maps $b_u:Z\to\Sone$ and $c_u:Z/\ol{\phi(\bbZ(\bf{p}_i - \bf{p}_k))}\to\Sone$ such that
\[\DDelta_u\tau_i(z) = \DDelta_{\phi(\bf{p}_i)}b_u(z)\cdot c_u(z\cdot \ol{\phi(\bbZ(\bf{p}_i - \bf{p}_k))})\]
Haar almost surely.  In which case we write that $b_u$ \textbf{solves equation} E$(i,j,u)$ and term $c_u$ its \textbf{one-dimensional auxiliary}.
\end{lem}

\textbf{Remark}\quad Observe that if $\ol{\phi(\bbZ(\bf{p}_i - \bf{p}_j))} = \ol{\phi(\bbZ(\bf{p}_i - \bf{p}_k))} = Z$ then the above conclusions simply promise solutions to the classical Conze-Lesigne equations.  However, the DIO property prevents these subgroups of $Z$ from being dense in all cases except $Z = \{0\}$, so for us this is only a first step on route to Proposition~\ref{prop:cocyc-factorize}. \fin

\textbf{Proof}\quad By symmetry it suffices to show that when $i=1$ the equation E$(1,3,u)$ admits a solution for every $u \in \ol{\phi(\bbZ(\bf{p}_1 - \bf{p}_3))}$.  Recall that we write $w_i := \phi(\bf{p}_i)$, and let us also set $w_{ij} := w_i\ol{\phi(\bbZ(\bf{p}_i - \bf{p}_j))} = w_j\ol{\phi(\bbZ(\bf{p}_i - \bf{p}_j))}$ when $i \neq j$.

By Lemma~\ref{lem:lift-u} we have $(u,u_2,1) \in \vec{Z}$ for some $u_2 \in Z$.  We may therefore consider equation~(\ref{eq:combcoc}) shifted by this element of $\vec{Z}$, and now dividing the shifting equation by the original gives
\begin{eqnarray}\label{eq:combcoc2}
\DDelta_u\tau_1(z_1)\cdot \DDelta_{u_2} \tau_2(z_2) = \DDelta_{(w_1,w_2,w_3)}\t{b}_u(z_1,z_2,z_3)
\end{eqnarray}
$m_{\vec{Z}}$-almost surely, where $\t{b}_u(\vec{z}) := c((u,u_2,1)\cdot \vec{z})\cdot\ol{c(\vec{z})}$.

However, by
Corollary~\ref{cor:simple-rewriting-group-rotations} we can
re-coordinatize $(Z,m_Z,R_{w_i})$ as a relatively invariant extension of $(Z/\ol{\phi(\bbZ(\bf{p}_1 - \bf{p}_2))},m_{Z/\ol{\phi(\bbZ(\bf{p}_1 - \bf{p}_2))}},R_{w_{12}})$ for $i=1,2$, say as
\[(Z,m_Z,R_{w_i}) \cong (S_i,\nu_i,\id)\otimes (Z/\ol{\phi(\bbZ(\bf{p}_1 - \bf{p}_2))},m_{Z/\ol{\phi(\bbZ(\bf{p}_1 - \bf{p}_2))}},R_{w_{12}})\]
for some auxiliary standard Borel spaces
$(S_i,\nu_i)$. In these new coordinatizations equation~(\ref{eq:combcoc2}) reads
\[(\DDelta_u\tau_1)(s_1,z_{12})\cdot(\DDelta_{u_2}\tau_2)(s_2,z_{12})
= \DDelta_{\id\times \id\times R_{w_{12}}} \t{b}_u(s_1,s_2,z_{12})\] $(\nu_1\otimes \nu_2\otimes
m_{Z/\ol{\phi(\bbZ(\bf{p}_1 - \bf{p}_2))}})$-almost surely (where we have been a little casual in identifying the differenced function $\DDelta_u\tau_1$ as a function on $S_1\times (Z/\ol{\phi(\bbZ(\bf{p}_1 - \bf{p}_2))})$). This we can re-arrange to give
\[(\DDelta_u\tau_1)(s_1,z_{12}) = \DDelta_{\id\times\id\times R_{w_{12}}} \t{b}_u(s_1,s_2,z_{12})\cdot \overline{(\DDelta_{u_2}\tau_2)(s_2,z_{12})},\]
so picking some $s_2 = s_2^\circ$ for which this holds for almost
every $(s_1,z_{12})$ we deduce that
\[\DDelta_u\tau_1(s_1,z_{12}) = \DDelta_{\id\times R_{w_{12}}}b_0(s_1,z_{12})\cdot c_0(z_{12})\]
with $b_0(s_1,z_{12}):= \t{b}_u(s_1,s_2^\circ,z_{12})$ and
$c_0(z_{12}) :=
(\DDelta_{u_2}\tau_2)(s_2^\circ,z_{12})$. Recalling our identification $(Z,R_{w_1})\cong (S_1\times Z/\ol{\phi(\bbZ(\bf{p}_1 - \bf{p}_2))},\id\times R_{w_{12}})$, we recognize this as equation E$(1,3,u)$, so $b_0$ is a solution. \qed

\textbf{Remark}\quad Although very simple, the above analysis of equation~(\ref{eq:combcoc}) was possible only in light of the DIO property and its consequence Corollary~\ref{cor:simple-rewriting-group-rotations}, which in turn hold only because of the very strong FIS assumption.  I suspect that without the DIO property there may be instances of the combined cocycle equation for which the above conclusion fails. \fin

The equations E$(i,j,u)$ solved by the preceding lemma are already suggestively close to the Conze-Lesigne equations of Proposition~\ref{prop:CL}.  That propositions makes it clear that any nil-cocycle admits solutions to every E$(i,j,u)$, but owing to the restriction $u \in \ol{\phi(\bbZ(\bf{p}_i - \bf{p}_k))}$ in Lemma~\ref{lem:processing-cocyc-over-Kron} there are other examples of cocycles $\tau_i$ that admit solutions to these equations.

In particular, this is trivially so if $\tau_i$ is invariant in either of the directions $\phi(\bf{p}_i - \bf{p}_j)$ or $\phi(\bf{p}_i - \bf{p}_k)$, since if $\tau_i$ is $\phi(\bf{p}_i - \bf{p}_j)$-invariant then $\DDelta_u\tau_i \equiv 1$ for any $u \in \ol{\phi(\bf{p}_i - \bf{p}_j)}$, and if it is $\phi(\bf{p}_i - \bf{p}_k)$-invariant then we may simply let $b_u\equiv 1$ and $c_u := \DDelta_u\tau_i$.

This sheds some light on the point of Proposition~\ref{prop:cocyc-factorize}: we will prove that any cocycle $\tau_i$ admitting solutions to all of the equations E$(i,j,u)$ must factorize into examples of these different kinds (partially invariant cocycles and nil-cocycles, or more precisely local nil-cocycles).  That will give the factorization needed for Proposition~\ref{prop:cocyc-factorize}, and will take up the remainder of this section.

In pursuit of this goal we will first prove a factorization result for the functions $b_u$, $b_v$ that solve the equations E$(i,j,u)$ and E$(i,k,v)$, and then use that to factorize $\tau_i$ itself.

For brevity let us now set $\bf{n}_1 := \bf{p}_1$ and $\bf{n}_j:= \bf{p}_1 - \bf{p}_j$ for $j=2,3$, and let $K_i := \ol{\phi(\bbZ\bf{n}_i)}\leq Z$ for $i=1,2,3$ and $Z_{ij}:= K_i\cdot K_j\leq Z$ for $i\neq j$. Since $\phi$ has dense image and $\bbZ\bf{n}_i + \bbZ\bf{n}_j$ has finite index in $\bbZ^2$, each $Z_{ij}$ has finite-index in $Z$.

Let $u\mapsto b_u$ be a measurable selection of solutions to the equations E$(1,2,u)$ or E$(1,3,u)$; note that this is unambiguous since the DIO property gives that $K_2\cap K_3 = \{0\}$, so for $u\neq 0$ at most one of the equations can apply. We now extend the definition of the measurable selection
$b_\bullet$ from $K_2\cup K_3$ to $Z_{23} := K_2\cdot K_3$ by
setting
\[b_{uv} := b_v\cdot(b_u\circ R_v)\]
when $u \in K_2$ and $v \in K_3$.

Our analysis will rely on the cohomological results of Appendix~\ref{app:cohom}.  In order to bring these to bear we must first perform various further manipulations on the Borel selection $z\mapsto b_z$.

Using $b$ we define the map $\k:Z_{23}\times Z_{23}\to \C(Z)$ by
\[\k(z,z'):= (b_z\circ R_{z'})\cdot \ol{b_{zz'}}\cdot b_{z'}.\]
This is a Borel $2$-cocycle in the sense of Moore's cohomology theory for locally compact groups (see Appendix A), where we endow the Polish Abelian group $\C(Z)$ with the obvious rotation action of $Z_{23}$ restricted from that of $Z$. Note the important triviality that here and henceforth we write our cocycles with the order of the arguments reversed, thus:
\[(\k(z,z')\circ R_{z''}) \cdot\ol{\k(z,z'z'')}\cdot\k(zz',z'') \cdot\ol{\k(z',z'')} = 1,\]
since this convention seems slightly more natural in our present setting.  Of course this theory behaves exactly as does the conventional version, since they are isomorphic under writing the arguments in reverse order.

We will first show that $\k$ takes values in the rotation-invariant subgroup
\[\W = \C(Z)^{K_1}\cdot \C(Z)^{K_2}\cdot \C(Z)^{K_3} \leq \C(Z),\]
where $\C(Z)^K$ denotes the subgroup of maps in $\C(Z)$ that are invariant under translation by the subgroup $K$. This will need the following toy result on reducing the dependences of coboundary
equations (see, for example, Moore and
Schmidt~\cite{MooSch80}).

\begin{lem}\label{lem:baby-cocyc-factorization}
Suppose that $q:(\t{Z},R_{\t{w}})\to (Z,R_w)$ is an
extension of ergodic group rotations and that $\s:Z\to \Sone$ is a
measurable function such that $\s\circ q$ is a coboundary over
$R_{\t{w}}$. Then there is some $\theta \in \Sone$
such that $\theta\cdot\s$ is a coboundary over $R_w$. \qed
\end{lem}

\begin{lem}\label{lem:cocycle-basics}
The following hold:
\begin{enumerate}
\item If $b$ and $b'$ both satisfy E$(1,2,u)$ then
$b \cdot \ol{b'} \in \W$.
\item If $u\in K_2$, $v\in K_3$ and $b$ and $b'$ satisfy respectively E$(1,2,u)$ and E$(1,3,v)$ then
\[[R_v\ltimes b',R_u\ltimes b] = \id_Z\ltimes (\DDelta_vb \cdot\ol{\DDelta_ub'})\\
\in \id_Z\ltimes\C(Z)^{K_1}.\]
\item The cocycle $\k$ almost surely takes values in $\W$.
\end{enumerate}
\end{lem}

\textbf{Proof}\quad\textbf{1.}\quad Letting $c$ and $c'$ be the
respective one-dimensional auxiliaries of $b$ and $b'$ and dividing
the resulting instances of equation E$(1,2,u)$ gives
\[\DDelta_{\phi(\bf{n}_1)}(b \cdot\ol{b'})(z) = (c' \cdot\ol{c})(z K_3).\]
Now we can apply Lemma~\ref{lem:baby-cocyc-factorization} to each of the finitely many ergodic components of $R_{\phi(\bf{n}_1)K_3}\actson Z/K_3$ (which are just the cosets of $K_1K_3/K_3$ in $Z/K_3$), choosing for each of them some ergodic component of $R_{\phi(\bf{n}_1)}$ that covers it, and so deduce that $(c'\cdot\ol{c})(zK_3) = \theta(z)\cdot\DDelta_{\phi(\bf{n}_1)}g(z)$ for some $\theta \in
\C(Z)^{K_1K_3}$ and $g\in \C(Z)^{K_3}$.

Substituting this expression now gives
\[\DDelta_{\phi(\bf{n}_1)}(b \cdot\ol{b'} \cdot\ol{g}) = \theta,\]
so that on each of the finitely many cosets of $K_1K_3$ the map $b \cdot\ol{b'}\cdot\ol{g}$ agrees with the restriction of the product of a character and a $K_1$-invariant function.

If $\chi\in \C(Z)$ is such that its restriction to each coset of $K_1K_3$ agrees with the restriction of some character, then since $K_1\cap K_2 = \{0\}$ we may factorize $\chi$ within each coset into a product of $K_1$- and $K_2$-invariant characters.  These factorizations then combine to give one factorization of $\chi$ as a member of $\C(Z)^{K_1}\cdot\C(Z)^{K_2}$.  Overall, this shows that $b\cdot \ol{b'}$ is the product of $g \in \C(Z)^{K_3}$ with members of $\C(Z)^{K_1}$ and $\C(Z)^{K_2}$, as required.

\quad\textbf{2.}\quad This follows similarly: if $c,c'$ are the
respective one-dimensional auxiliaries of $b,b'$, then the
assumed equations give
\[\DDelta_u\tau_1(z) = \DDelta_{\phi(\bf{n}_1)}b(z) \cdot c(zK_3)\]
and
\[\DDelta_v\tau_1(z) = \DDelta_{\phi(\bf{n}_1)}b'(z) \cdot c'(zK_2),\]
and so if we now difference the first of these equations by $v$ and the second by $u$ and then divide the two new equations that result, we are left with
\[\DDelta_{\phi(\bf{n}_1)}(\DDelta_v b \cdot\ol{\DDelta_u b'}) = 1.\]
Hence $\DDelta_v b \cdot\ol{\DDelta_u b'} \in \C(Z)^{K_1}$, and a direct
calculation gives that
\[[R_v\ltimes b',R_u\ltimes b] = \id_Z\ltimes (\DDelta_v b \cdot \ol{\DDelta_u b'}).\]

\quad\textbf{3.}\quad If $u_i \in K_2$ and $v_i \in
K_3$ for $i=1,2$, then by definition
\begin{eqnarray*}
\k(u_1v_1,u_2v_2) &=& ((b_{v_1} \cdot b_{u_1}\circ R_{v_1})\circ R_{u_2 v_2})\cdot \ol{b_{v_1v_2} \cdot (b_{u_1u_2}\circ
R_{v_1v_2})}\\
&& \cdot (b_{v_2}\cdot (b_{u_2}\circ R_{v_2}))\\
&=& (b_{v_1}\circ R_{u_2v_2})\cdot\ol{b_{v_1v_2}}\cdot b_{v_2}\\
&& \cdot (b_{u_1}\circ R_{v_1u_2v_2}) \cdot \ol{b_{u_1u_2}\circ R_{v_1v_2}} \cdot (b_{u_2}\circ R_{v_2})\\
&=& (b_{v_1}\circ R_{v_2})\cdot\ol{b_{v_1v_2}}\cdot b_{v_2}\\
&& \cdot (b_{v_1}\circ R_{u_2v_2})\cdot\ol{b_{v_1}\circ R_{v_2}}\\
&& \cdot ((b_{u_1}\circ R_{u_2}) \cdot \ol{b_{u_1u_2}} \cdot b_{u_2})\circ R_{v_1v_2}\\
&& \cdot b_{u_2}\circ R_{v_2} \cdot\ol{b_{u_2}\circ R_{v_1v_2}}\\
&=& (b_{v_1}\circ R_{v_2})\cdot\ol{b_{v_1v_2}}\cdot b_{v_2}\\
&& \cdot ((b_{u_1}\circ R_{u_2}) \cdot \ol{b_{u_1u_2}} \cdot b_{u_2})\circ R_{v_1v_2}\\
&& \cdot (\DDelta_{u_2}b_{v_1}\cdot \ol{\DDelta_{v_1}b_{u_2}})\circ R_{v_2}.
\end{eqnarray*}
This expression now contains two kinds of factor, which we can show
must almost surely lie in $\W$ using two separate arguments:
\begin{itemize}
\item Simply by multiplying equation E$(1,3,v_2)$
\[\DDelta_{v_2}\tau_1(z) = \DDelta_{\phi(\bf{n}_1)}b_{v_2}(z)\cdot c_{v_2}(zK_2)\]
and the shifted equation
E$(1,3,v_1)\circ R_{v_2}$
\[\DDelta_{v_1}\tau_1(zv_2) = \DDelta_{\phi(\bf{n}_1)}(b_{v_1}\circ R_{v_2})(z) \cdot c_{v_1}(zv_2K_3)\]
we see that $(b_{v_1}\circ R_{v_2}) \cdot b_{v_2}$ is a solution of equation E$(1,3,v_1v_2)$ for almost
every $v_1$ and $v_2$, and hence by part 1 above that
\[(b_{v_1}\circ R_{v_2})\cdot \ol{b_{v_1v_2}}\cdot b_{v_2} \in \W,\]
and similarly
\[((b_{u_1}\circ R_{u_2})\cdot \ol{b_{u_1u_2}}\cdot b_{u_2})\circ R_{v_1v_2} \in \W\]
for almost every $u_1$ and $u_2$.
\item The expression $\DDelta_{u_2}b_{v_1} \cdot\ol{\DDelta_{v_1}b_{u_2}}$ has
been shown to lie almost surely in $\W$ in part
2 above. \qed
\end{itemize}

This target module $\W$ for $\k$ is huge, and so it seems this cocycle may be extremely complicated.  But at least by another measurable selection we can factorize it as
\[\k = \k_1\cdot \k_2\cdot \k_3\]
where each $\k_i$ is a $2$-cochain taking values in $\C(Z)^{K_i}$.  The point will be that we can modify this to obtain a factorization of $\k$ into pieces that individually behave well.

To do this we will first pass to one higher degree of cohomology.  Applying the coboundary operator to the above gives
\[1 = d\k_1\cdot d\k_2\cdot d\k_3.\]

In order to discuss the individual factors $d\k_i$ we need another piece of notation.

\begin{dfn}[Locally affine functions]
If $Z$ is a compact Abelian group and $Z_0 \leq Z$ a finite-index subgroup, then a function $f \in \C(Z)$ is \textbf{$Z_0$-locally affine} if its restriction to every coset of $Z_0$ agrees with the restriction of some affine function (that is, a constant multiple of a character).  We write $\E(Z;Z_0) \leq \C(Z)$ for the subgroup of $Z_0$-locally affine functions.
\end{dfn}

\begin{lem}\label{lem:hard-to-factorize}
There is a finite-index subgroup $Z_0\leq Z$ such that if $\g_i \in \C(Z)^{K_i}$ for $i=1,2,3$ and $\g_1\cdot \g_2\cdot \g_3 \equiv 1$ then in fact $\g_i \in \E_i(Z;Z_0)^{K_i}$ for each $i$.
\end{lem}

\textbf{Proof}\quad Choose $Z_0 := K_1K_2\cap K_1K_3\cap K_2K_3$, so this has index at most
\[[Z:K_1K_2][Z:K_1K_3][Z:K_2K_3].\]
Now let $u\in K_1$ and take a difference to obtain
$\DDelta_u\g_2\cdot \DDelta_u\g_3 \equiv 1$.  Clearly $\DDelta_u\g_i \in \C(Z)^{K_i}$, so this equation implies that each of these two factors is actually invariant under the whole of $Z_{23} = K_2K_3$, and hence certainly on $Z_0$.  Therefore $\g_i|_{zZ_0}$ must be an affine function on $zZ_0$ for each coset $zZ_0\leq Z$ for $i=2,3$: that is, $\g_i \in \E(Z;Z_0)^{K_i}$.  Differencing in a different direction we can treat $\g_1$ similarly. \qed

The lemma implies that in fact the $3$-cocycle $d\k_i:Z_{23}\times Z_{23}\times Z_{23}\to \C(Z)$ takes values in the smaller group of maps $\E_i(Z;Z_0)^{K_i}$ for each $i$.  Observe also that if $\g \in \E(Z;Z_0)$ is $K_i$-invariant then it must actually reside in $\E(Z;Z_0K_i)^{K_i}$.  Let us now write $d\k_i|_{Z_0\times Z_0\times Z_0}$ for the map
\[Z_0\times Z_0\times Z_0\to \E(Z;Z_0)^{K_i}\]
obtained by restricting the domain on which $d\k_i$ is defined (\emph{not} restricting the individual functions in its target module). It follows that we may identify $d\k_i|_{Z_0\times Z_0\times Z_0}$ with a cocycle taking values in a direct sum of copies of $\E(Z_0K_i)^{K_i}$, one for each coset of $Z_0K_i$ in $Z$.

Now, since $Z_0 \leq K_iK_j \cong K_i \times K_j$, clearly $(Z_0\cap K_i)\leq K_i$ has finite index and similarly for $j$, and so replacing $Z_0$ by $(Z_0\cap K_i)(Z_0\cap K_j)$ if necessary we may apply the virtual vanishing result of Lemma~\ref{lem:deg-3-vanishing} to deduce that each of the above-mentioned components of $d\k_i|_{Z_0\times Z_0\times Z_0}$ is actually an $\E(Z_0K_i)^{K_i}$-valued coboundary.  Repeating this for each $i$ and now letting $Z_0$ be the intersection of the various finite-index subgroups obtained in the process, we can put these components back together to obtain
$d\k_i|_{Z_0\times Z_0\times Z_0} = d\a_i$ for some $2$-cochains $\a_i:Z_0\times Z_0\to \E(Z;Z_0K_i)^{K_i}$.

It follows that
\[\k|_{Z_0\times Z_0} = (\k_1\cdot \k_2\cdot \k_3)|_{Z_0\times Z_0} = \k'_1\cdot \k_2'\cdot \k'_3\cdot \a\]
where
\begin{itemize}
\item $\k'_i = \k_i|_{Z_0\times Z_0}\cdot \ol{\a_i}$ is a $2$-cocycle with values in \[\C(Z)^{K_i}\cdot \E(Z;Z_0K_i)^{K_i} = \C(Z)^{K_i}\] for $i=1,2,3$,
\item and $\a := \a_1\a_2\a_3$ is a $2$-cocycle with values in
\[\E(Z;Z_0K_1)^{K_1}\cdot \E(Z;Z_0K_2)^{K_2}\cdot \E(Z;Z_0K_3)^{K_3} \leq \E(Z;Z_0).\]
\end{itemize}

The above factorization of $\k|_{Z_0\times Z_0}$ can now be unwoven into a useful factorization of $b_z$ (at least for $z \in Z_0$).

\begin{prop}\label{prop:factorizing-bs}
After possibly shrinking $Z_0$ further, there are Borel maps $b_i:Z_0\to \C(Z)^{K_i}$ for $i=1,2,3$, a $Z_0$-local nil-selector $b_\nil:Z_0\to \C(Z)$ (see Definition~\ref{dfn:local-nil-select}) and a Borel map $\b \in \C(Z)$ such that \[b_z = b_{1,z}b_{2,z}b_{3,z}b_{\nil,z}\DDelta_z\b\quad\quad\hbox{for Haar-a.e. }z \in Z_0.\]
\end{prop}

\textbf{Proof}\quad This will require an analysis of each $\k'_i$ and of $\a$; we break these into separate steps.

\quad\textbf{Step 1}\quad We can identify $\k'_i$ with a direct sum of $[Z:Z_0K_i]$-many $2$-cocycles taking values in the $Z_0$-module $\C(Z_0K_i)^{K_i}\cong \C(Z_0)^{K_i\cap Z_0}$.

Shrinking $Z_0$ further if necessary, we may assume that $Z_0$ is the product subgroup $(Z_0\cap K_1)\cdot (Z_0\cap K_2)$, and then by Lemma~\ref{lem:C(Z)K-vald-collapse} for each of the $\C(Z_0)^{Z_0\cap K_1}$-valued components $\l$ of $\k'_1$ we may write $\l = \l' \cdot da$ for some $a:Z_0\to \C(Z_0K_1)^{K_1}$ and $\l'\in \Z^2(Z_0,\Sone)$.  In addition, by Theorem~\ref{thm:cohomA} this $2$-cocycle $\l'$ is inflated up to cohomology from some finite-dimensional quotient $Z_0\to (\Sone)^D\times F$. Now the dimension shifting Proposition~\ref{prop:HofMoscalc} and the vanishing result Corollary~\ref{cor:vanishing-from-HofMos} show that $\rmH^2((\Sone)^D,\Sone) = (0)$ for any $D$, and hence that $\l'$ must trivialize upon restricting to some further finite-index subgroup of $Z_0$.  Arguing thus for each of the finitely many components of $\l$ and then reassigning the label $Z_0$ to the intersection of the finitely many finite-index subgroups so obtained, we may therefore assume that each of the $\C(Z_0K_1)^{K_1}$-valued components of $\k'_1$ is a coboundary, and now their primitives combine to show that $\k'_1$ itself is a coboundary.  Repeating this argument for $i=2,3$ shows that for some (perhaps much smaller) finite-index subgroup $Z_0$ we may write $\k'_i = db_i$ for some $b_i:Z_0\to\C(Z)^{K_i}$.

\quad\textbf{Step 2}\quad We next make a similar analysis of $\a:Z_0\times Z_0\to \E(Z;Z_0)$.  First observe that it, too, breaks into $[Z:Z_0]$-many components, each of which may be identified as a $2$-cocycle $Z_0\times Z_0\to \E(Z_0)$.  Let $\l$ now be one of these components of $\a$.

Applying Theorems~\ref{thm:cohomA} and~\ref{thm:cohomB} and the long exact sequence corresponding to the presentation
\[\Sone\into \E(Z_0)\onto \hat{Z_0},\]
we first deduce that $\l = (\l'\circ q^{\times 2})\cdot da$ for some $a:Z_0\to \E(Z_0)$ and the inflation through some finite-dimensional quotient $q:Z_0\onto Z_1 \cong (\Sone)^D\times F$ of a $2$-cocycle $\l' \in \Z^2(Z_1,\E(Z_1))$.  Shrinking $Z_0$ again if necessary we may assume that in fact $F = (0)$.

Now let $\ol{\l'}$ be the image of $\l'$ under the quotient map $\E(Z_1)\to \hat{Z_1} \stackrel{\cong}{\to} \bbZ^D$. Arguing coordinate-wise in $\bbZ^D$, the standard calculation of Corollary~\ref{cor:reps-for-Sone0} promises that each coordinate of $\ol{\l'}$ is cohomologous to a $2$-cocycle of the form
\[(u_1,u_2)\mapsto \lfloor \{\g(u_1)\} + \{\g(u_2)\}\rfloor\]
for some $\g \in \hat{(\Sone)^D}$.  Crucially, these are all \emph{symmetric} functions on $Z_1\times Z_1$, and so we can write $\ol{\l'} = \ol{\l''} + d\ol{\a'}$ for some symmetric $2$-cocycle $\ol{\l''}:Z_1\times Z_1\to \bbZ^D$ and some $\ol{\a'}:Z_1\to \bbZ^D$.  Choosing a measurable lift of $\ol{\a'}$ that takes values in $\E(Z_1)$, inflating it through $q$ and combining it with $a$, we may assume that in fact $\ol{\l''} = \ol{\l'}$ is itself a symmetric function on $Z_1\times Z_1$.

\quad\textbf{Step 3}\quad We will now make two uses of the standard identification of $\rmH^2(\cdot,\cdot)$ with equivalence classes of group extension (see Moore~\cite{Moo64(gr-cohomI-II)}). First we deduce that there is some extension
\[\bbZ^D \into A \onto Z_1\]
that gives rise to the $2$-cocycle $\ol{\l'}$, and since $\ol{\l'}$ is symmetric it follows that $A$ is still Abelian: that is, $A$ is a finite-dimensional locally compact Abelian group with a covering map $A\onto Z_1$.

On the other hand, the original $\E(Z_1)$-valued $2$-cocycle $\l'$ corresponds to a larger extension
\[\E(Z_1)\into G\onto Z_1,\] into which $A$ now fits as an intermediate quotient $G\onto A\onto Z_1$, where the kernel of the first of these two quotient maps is $\E(Z_1)/\hat{Z_1} \cong \Sone$.  Therefore $G$ is an extension of $\Sone$ by $A$ with the trivial action: in particular, $G$ is a two-step nilpotent group, and by standard classification results (it is clearly connected and without small subgroups) it must be a Lie group.

Thus we have shown that each component $\l$ of $\a$ is actually inflated from the cocycle describing some two-step nilpotent Lie group extension of $\E(Z_1)$ by $Z_1$, and so it equals the coboundary of some nil-selector $Z_1\to \C(Z_1)$.  Given such nil-selectors for all the components of $\a$, lifting them back up to $Z_0$ and combining them we obtain a $Z_0$-local nil-selector $b_\nil:Z_0\to \C(Z)$ such that $\a = db_\nil$.

\quad\textbf{Step 4}\quad It remains to put the above information together.  We have obtained a finite-index $Z_0\leq Z$ and $b_i$ for $i=1,2,3$ and $b_\nil$ such that
\[\k|_{Z_0\times Z_0} = d(b|_{Z_0}) = d(b_1\cdot b_2\cdot b_3\cdot b_\nil),\]
and hence
\[b|_{Z_0}\cdot \ol{b_1\cdot b_2\cdot b_3\cdot b_\nil}:Z_0\to \C(Z)\]
is a $1$-cocycle.  The triviality of $\rmH^1(Z_0,\C(Z))$ now gives some $\b \in \C(Z)$ such that this $1$-cocycle equals $\DDelta_\bullet\b$.
Re-writing this equality gives the desired conclusion. \qed

Since the conclusion of Proposition~\ref{prop:cocyc-factorize} is clearly invariant under modifying $\tau_1$ by a coboundary, we may at this stage divide it by $\DDelta_{\phi(\bf{n}_1)}\b$ and simultaneously divide each $b_z$ by $\DDelta_z\b$, and so henceforth assume that $\b \equiv 1$.

\begin{lem}\label{lem:locconst-to-eigen}
If $\phi:\bbZ^2\to Z$ is a dense homomorphism, $Z_0 \leq
Z$ has finite index, $\theta \in \C(Z)^{Z_0}$ and $\bf{n} \in
\bbZ^2$, then there is an extension $q:(\t{Z},\t{\phi})\to (Z,\phi)$
of ergodic rotations such that $\theta\circ q$ is a coboundary over
$R_{\t{\phi}(\bf{n})}$.
\end{lem}

\textbf{Proof}\quad Let $m\geq 1$ be minimal such that $\phi(m\bf{n})\in Z_0$ and let
\[\theta'(z) := \theta(z)\theta(z\phi(\bf{n}))\cdots\theta(z\phi((m-1)\bf{n})).\]
It will suffice to find an extension such that $\theta'\circ q$ is a
coboundary over $R_{\t{\phi}(m\bf{n})}$, since if $\theta'\circ q =
\DDelta_{\t{\phi}(m\bf{n})}g$ then $(\theta\circ
q)\cdot\ol{\DDelta_{\t{\phi}(\bf{n})}g}$ is a cocycle whose $m$-fold
composition vanishes, and so this may be shown to be a coboundary by hand using the fact that the cosets $\phi(\bf{n})Z_0$, \ldots, $\phi(m\bf{n})Z_0$ are distinct.

Now $R_{\phi(m\bf{n})}$ preserves each coset of $Z_0$ within $Z$.  Therefore we may simply form $Z' := Z\times (\Sone)^{[Z:Z_0]}$, let $q$ be the first coordinate projection and lift $\phi$ to some $\t{\phi} = (\phi,\phi')$ so that the finitely many values taken by $\theta$ are all eigenvalues of the rotation $R_{\phi'(m\bf{n})}$ on $(\Sone)^{[Z:Z_0]}$.  Having done this, patching together the corresponding eigenfunctions exhibits $\theta\circ q$ as the desired coboundary. The proof is completed by restricting from $Z'$ to the closed subgroup $\ol{\t{\phi}(\bbZ^2)}$, which still covers $Z$ because $\phi$ was assumed ergodic. \qed

\textbf{Proof of Proposition~\ref{prop:cocyc-factorize}}\quad We have already reduced to the ergodic case, and clearly it suffices to assume that $(i,j,k) = (1,2,3)$. Let $Z_0\leq Z$ be given by Proposition~\ref{prop:factorizing-bs}, and let $m\geq 1$ be minimal such that $\phi(m\bf{n}_1) \in Z_0$.

\quad\textbf{Step 1}\quad That proposition and the equations solved in Lemma~\ref{lem:processing-cocyc-over-Kron} give for almost every $z \in Z_0$ that
\[\DDelta_z\tau_1 = \DDelta_{\phi(\bf{n}_1)}b_z\cdot c_{2,z}\cdot c_{3,z}\]
for some $c_{2,z}\in \C(Z)^{K_2}$, $c_{3,z} \in \C(Z)^{K_3}$, where
\[b_z = b_{1,z}b_{2,z}b_{3,z}b_{\nil,z}.\]
Substituting this latter factorization and observing that $\DDelta_{\phi(\bf{n}_1)}b_{1,z}\equiv 1$ (because $b_{1,z}$ is $K_1$-invariant), we obtain
\[\DDelta_z\tau_1 = \DDelta_{\phi(\bf{n}_1)}b_{\nil,z}\cdot\big(\DDelta_{\phi(\bf{n}_1)}b_{2,z}\cdot c_{2,z}\big)\cdot \big(\DDelta_{\phi(\bf{n}_1)}b_{3,z}\cdot c_{3,z}\big).\]

If we consider this equation together with its shifts by $R_{\phi(\bf{n}_1)}$, $R_{\phi(2\bf{n}_1)}$, \ldots, $R_{\phi((m-1)\bf{n}_1)}$ and multiply them, we obtain
\[\DDelta_z\tau^{(m)}_1 = \DDelta_{\phi(m\bf{n}_1)}b_{\nil,z}\cdot\big(\DDelta_{\phi(m\bf{n}_1)}b_{2,z}\cdot c^{(m)}_{2,z}\big)\cdot \big(\DDelta_{\phi(m\bf{n}_1)}b_{3,z}\cdot c^{(m)}_{3,z}\big)\]
for the function
\[\tau^{(m)}_1 := \tau_1\cdot (\tau_1\circ R_{\phi(\bf{n}_1)})\cdot\cdots \cdot (\tau_1\circ R_{\phi((m-1)\bf{n}_1)})\]
and similarly-defined $c^{(m)}_{i,z}$ for $i=2,3$, which clearly still lie in their respective $\C(Z)^{K_i}$.

Now, since $b_\nil$ is a $Z_0$-local nil-selector, letting $\s_\nil\in \C(Z)$ be its value over $\phi(m\bf{n}_1)$ we also obtain
\[\DDelta_z\s_\nil = \DDelta_{\phi(m\bf{n}_1)}b_{\nil,z}\cdot \theta\]
for some $\theta\in \C(Z)^{Z_0}$ (this is an application of Proposition~\ref{prop:CL} on each coset of $Z_0$ separately).

Applying Lemma~\ref{lem:locconst-to-eigen} and replacing the finite-index containment of groups $Z_0 \leq Z$ by their resulting extensions if necessary, we may assume that $\theta$ is in fact a coboundary over $R_{\phi(m\bf{n}_1)}$, and now adjusting $\tau_1$ by this coboundary we may assume further that $\theta = 0$ to obtain
\[\DDelta_z((\tau_1^{(m)}\circ q)\cdot\ol{\s_\nil}) = c'_{2,z}\cdot c'_{3,z}\]
where
\[c'_{i,z}:= \DDelta_{\phi(m\bf{n}_1)}b_{i,z}\cdot c^{(m)}_{i,z}.\]

\quad\textbf{Step 2}\quad The above equation implies that the map $c':z\mapsto c_{2,z}'c_{3,z}'$ is a $1$-cocycle $Z_0\to \C(Z)^{K_2}\cdot \C(Z)^{K_3}$.  Any function in $\C(Z)^{K_2}\cdot \C(Z)^{K_3}$ can be factorized as a product of members of $\C(Z)^{K_i}$, $i=2,3$, and this factorization is unique up to a member of $\C(Z)^{K_2K_3}$.  Therefore the two components $z\mapsto c'_{i,z}$ must be $1$-cocycles individually up to an error which is captured by a $2$-cocycle $Z_0\times Z_0\to\C(Z)^{K_2K_3}$.

Since $Z_0\leq K_2K_3$, the $Z_0$-module $\C(Z)^{K_2K_3}$ decomposes into a direct sum of copies of $\Sone$ with trivial $Z_0$-action. However, any class in $\rmH^2(Z_0,\Sone)$ trivializes on restricting to some further finite-index subgroup (using again Theorem~\ref{thm:cohomA} and the vanishing $\rmH^2((\Sone)^D,\Sone) = (0)$), and so by shrinking $Z_0$ further we may assume that $dc_2' = \ol{dc_3'}$ is a $\C(Z)^{K_2K_3}$-valued $2$-coboundary.  Adjusting $c_2'$ and $c_3'$ by its primitive, it follows that we may in fact assume that each $c'_i$ is individually a $1$-cocycle.

\quad\textbf{Step 3}\quad We next solve this cocycle condition.  Noting that $c_i'$ takes values in $\C(Z)^{K_i}\subseteq \C(Z)$, by the triviality of $\rmH^1(Z_0,\C(Z))$ we can express $c_{i,z}' = \DDelta_z\b_i'$ for some $\b_i' \in \C(Z)$.  Hence we need only understand the condition that $\DDelta_z\b'_i$ be $K_i$-invariant for almost all $z\in Z_0$.  Since the crossed homomorphism $z\mapsto \DDelta_z\b'_i$ is automatically continuous, this invariance in fact holds for strictly all $z\in Z_0$.

Considering this condition first for $z \in Z_0\cap K_i$, we deduce that $\DDelta_z\b'_i$ is constant on each coset of $Z_0\cap K_i$. This requires that $\b'_i$ restrict to an affine map on each coset of $Z_0\cap K_i$ (which can always be extended to an affine map on $Z$).

Next, if $z,z'\in Z$ are such that $z(Z_0\cap K_i) \neq z'(Z_0\cap K_i)$ but $zK_i = z'K_i$, and if $\b'_i|_{z(Z_0\cap K_i)} = \theta\g|_{z(Z_0\cap K_i)}$ and $\b'_i|_{z'(Z_0\cap K_i)} = \theta'\g'|_{z'(Z_0\cap K_i)}$ for some $\theta,\theta'\in \Sone$ and $\g,\g'\in \hat{Z}$, then $\DDelta_w\b'_i$ for $w \in Z_0\cap K_i$ takes the constant values $\g(w)$ on $z(Z_0\cap K_i)$ and $\g'(w)$ on $z'(Z_0\cap K_i)$, and so since these two constants must agree it follows that $\g$ and $\g'$ must restrict to the same member of $\hat{Z_0\cap K_i}$ (and so in fact by adjusting the choice of constants $\theta,\theta'$ we may assume $\g = \g'$).

If instead $z' = zw$ for some $w \in Z_0$ and $\g$, $\g'$ are as above then the function $\DDelta_w\b'_i$ agrees with a constant multiple of $\g'\cdot \ol{\g}$ on $z(Z_0\cap K_i)$, and so since it is $(Z_0\cap K_i)$-invariant it follows again that $\g$ and $\g'$ must restrict to the same element of $\hat{Z_0\cap K_i}$.

Finally, suppose that $z' = zwk$ for some $w \in Z_0$ and $k\in K_i$ and that $\b'_i$ agrees with the restrictions of the affine maps $\theta_1\g$, $\theta_2\g$, $\theta_3\g$ and $\theta_4\g$ on the four cosets $z(Z_0\cap K_i)$, $zw(Z_0\cap K_i)$, $zk(Z_0\cap K_i)$ and $zwk(Z_0\cap K_i)$ respectively for some $\theta_1,\theta_2,\theta_3,\theta_4 \in \Sone$ and $\g\in \hat{Z}$.  Then $\DDelta_w \b'_i$ takes the value $\theta_2\ol{\theta_1}$ on $z(Z_0\cap K_i)$ and the value $\theta_4\ol{\theta_3}$ on $zk(Z_0\cap K_i)$, so since these must be equal we have $\theta_2\ol{\theta_1} = \theta_4\ol{\theta_3}$.

It follows that on each coset $zZ_0K_i$, $\b'_i$ must take the form $\theta\cdot \g|_{zZ_0K_i}$ for some fixed $\g \in \hat{Z}$ and some $\theta:zZ_0K_i\to \Sone$ which factorizes into a product of a $Z_0$-invariant function and a $K_i$-invariant function. Therefore overall we find that $\b'_i$ can be factorized as $\chi_i\cdot \theta_i\cdot \b''_i$ with $\theta_i\in \C(Z)^{Z_0}$, $\chi_i$ a map which is affine on each coset of $Z_0K_i$ and $\b''_i\in\C(Z)^{K_i}$.

Re-arranging everything we have so far, we obtain
\[\DDelta_z((\tau_1^{(m)}\circ q)\cdot\ol{\b_2''\cdot \b_3''}) = \DDelta_z(\s_\nil\cdot\theta_2\cdot \theta_3\cdot\chi_2\cdot\chi_3)\]
for all $z \in Z_0$, and hence
\[(\tau_1^{(m)}\circ q)\cdot\ol{\b_2''\cdot \b_3''} = \l\cdot\s_\nil\cdot\theta_2\cdot \theta_3\cdot\chi_2\cdot\chi_3\]
for some $\l\in \C(Z)^{Z_0}$.  Here the right-hand size is a product of a $Z_0$-local nil-cocycle with an element of $\E(Z;Z_0)$, and since affine functions clearly satisfy the requirement of Proposition~\ref{prop:CL} such a product is still a $Z_0$-local nil-selector.  Letting $\tau_{1,i}:= \b_i''$ for $i=2,3$ and $\tau_\nil$ be this local nil-cocycle, this completes the proof.  \qed

\section{Convergence for some quadratic averages}\label{sec:convproof}

We will now use Theorem~\ref{thm:char-three-lines-in-2D} to prove Theorem~\ref{thm:polyconv}.

\subsection{Joint distributions of one-dimensional isotropy factors}

The proof of convergence will require some basic
results on the possible distributions of collections of
one-dimensional isotropy factors of a $\bbZ^2$-system, which will again make use of the DIO property.

The following two propositions contain the extra control of joinings of one-dimensional isotropy factor that we need.

\begin{prop}\label{prop:joining-Kron-to-isotropies}
Suppose that $\bf{n}_1$, $\bf{n}_2$, $\bf{n}_3 \in \bbZ^2 \setminus
\{\bs{0}\}$ are three directions no two of which are parallel, that
$\bfX_1 = (X_1,\mu_1,T_1)\in \sfZ_0^{\bf{n}_1}$, $\bfX_2 =
(X_2,\mu_2,T_2)\in \sfZ_0^{\bf{n}_2}$, $\bfX_3 = (X_3,\mu_3,T_3)\in
\sfZ_0^{\bf{n}_3}$ and that $\bfZ = (Z,\nu,S)$ is a group rotation
$\bbZ^2$-system. Suppose further that $\bfX = (X,\mu,T)$ is a 
joining of these four systems through the factor maps
$\xi_i:\bfX\to\bfX_i$, $i=1,2,3$ and $\a:\bfX\to\bfZ$. Then
$(\xi_1,\xi_2,\xi_3,\a)$ are relatively independent under $\mu$ over
their further factors
$(\zeta_1^{T_1}\circ\xi_1,\zeta_1^{T_2}\circ\xi_2,\zeta_1^{T_3}\circ\xi_3,\a)$.
\end{prop}

\textbf{Remark}\quad In this proposition, the subscripts on `$T_i$', $i=1,2,3$, label different \emph{whole actions}: as usual, the individual transformations are indicated by a superscript, as in $T_i^\bf{v}$ for $\bf{v} \in \bbZ^2$. \fin

\textbf{Proof}\quad  We will prove that under $\bfX$ the factors
$\xi_1$, $\xi_2$, $\xi_3$ and $\a$ are relatively independent over
$\zeta_1^{T_1}\circ\xi_1$, $\xi_2$, $\xi_3$ and $\a$; repeating this
argument to handle $\xi_2$ and $\xi_3$ then gives the full result.

Letting $\bfY = (\xi_3\vee\a)(\bfX)$ be the factor of $\bfX$
generated by $\xi_3$ (which is $T^{\bf{n}_3}$-invariant) and $\a$
(which is isometric for $T$, hence certainly for $T^{\bf{n}_3}$), we
see that this is isometric for the $\bbZ\bf{n}_3$-subaction.  This implies
that its joining to any other system is relatively independent over
the maximal factor of that other system that is isometric for the $\bbZ\bf{n}_3$-subaction.

On the other hand, $\xi_1$ and $\xi_2$ must be relatively
independent over $\xi_1\wedge \xi_2$ under $\mu$ (simply by
averaging with respect to $\bf{n}_2$), and the subactions
generated by $\bf{n}_1$ and by $\bf{n}_2$ are both trivial on this
meet, so $\xi_1\wedge\xi_2 \precsim
\zeta_0^{T^{\bf{n}_1},T^{\bf{n}_2}}$.  Therefore the target system of $\xi_1\wedge \xi_2$ is a
direct integral of finite group rotations factoring through the
quotient $\bbZ^2/(\bbZ\bf{n}_1 + \bbZ\bf{n}_2)$.

Since $\xi_1\vee \xi_2$ must be joined to $\xi_3\vee\a$ relatively
independently over the maximal $T^{\bf{n}_3}$-isometric factor of
$\xi_1\vee \xi_2$, it follows from the Furstenberg-Zimmer Structure
Theorem~\ref{thm:rel-ind-joinings}
that $\xi_1\vee\xi_2$ is in particular joined to $\xi_3\vee\a$
relatively independently over the join of maximal isometric
subextensions
\[(\zeta_{1/(\xi_1\wedge\xi_2)|_{\xi_1}}^{T_1^{\bf{n}_3}}\circ\xi_1)\vee(\zeta_{1/(\xi_1\wedge\xi_2)|_{\xi_2}}^{T_2^{\bf{n}_3}}\circ\xi_2).\]
Since $\xi_1\wedge\xi_2$ has target a direct integral of
\emph{periodic} rotations, the maximal $T_i^{\bf{n}_3}$-isometric
subextension of $\xi_i\to (\xi_1\wedge\xi_2)|_{\xi_i}$ is simply the
maximal factor of $\xi_i$ that is coordinatizable as a direct
integral of group rotations for each $i=1,2$: that is, it is
$\zeta_1^{T_i}\circ\xi_i$. Hence we have shown that under $\mu$ the
factors $\xi_1\vee\xi_2$ and $\xi_3\vee\a$ are relatively
independent over
$(\zeta_1^{T_1}\circ\xi_1)\vee(\zeta_1^{T_2}\circ\xi_2)$ and
$\xi_3\vee\a$.  Thus whenever $f_i \in L^\infty(\mu_i)$ for
$i=1,2,3$ and $g\in L^\infty(\nu)$ we have
\begin{eqnarray*}
&&\int_X
(f_1\circ\xi_1)\cdot(f_2\circ\xi_2)\cdot(f_3\circ\xi_3)\cdot(g\circ\a)\,\d\mu\\
&&= \int_X
\sfE_\mu\big((f_1\circ\xi_1)\cdot(f_2\circ\xi_2)\,\big|\,(\zeta_1^{T_1}\circ\xi_1)\vee(\zeta_1^{T_2}\circ\xi_2)\big)\cdot(f_3\circ\xi_3)\cdot(g\circ\a)\,\d\mu\\
&&= \int_X
(\sfE_\mu(f_1\,|\,\zeta_1^{T_1})\circ\xi_1)\cdot(\sfE_\mu(f_2\,|\,\zeta_1^{T_2})\circ\xi_2)\cdot(f_3\circ\xi_3)\cdot(g\circ\a)\,\d\mu\\
&&= \int_X
(\sfE_\mu(f_1\,|\,\zeta_1^{T_1})\circ\xi_1)\cdot(f_2\circ\xi_2)\cdot(f_3\circ\xi_3)\cdot(g\circ\a)\,\d\mu,
\end{eqnarray*}
where the second equality follows from the relative independence of
$\xi_1$ and $\xi_2$ over $\xi_1\wedge \xi_2$, which is contained in
$\zeta_1^{T_i}\circ\xi_i$ for both $i=1,2$. This completes the
proof. \qed

Our second characterization of joint distributions of isotropy factors will require the following result from Furstenberg and Weiss (Lemma 10.3 of that paper).

\begin{lem}\label{lem:factorizing-transfer}
If $\bfX_1$, $\bfX_2$ are ergodic $\bbZ$-systems and
$f_i:X_i\to\Sone$, $i=1,2$, are Borel maps for which there is some
Borel $g:X_1\times X_2\to\Sone$ with $f_1\otimes f_2 =
\DDelta_{T_1\times T_2}g$, $(\mu_1\otimes \mu_2)$-a.s., then in fact
there are constants $c_i \in \Sone$ and Borel maps $g_i:X_i\to\Sone$
such that $f_i = c_i\cdot \DDelta_{T_i} g_i$. \qed
\end{lem}

\begin{prop}\label{prop:jointdists}
Suppose that $\bf{n}_1,\bf{n}_2,\bf{n}_3,\bf{n}_4 \in
\bbZ^2\setminus \{\bs{0}\}$ are directions no two of which are
parallel, that $\bfX_i = (X_i,\mu_i,T_i)\in \sfZ_0^{\bf{n}_i}$ for
$i=1,2,3,4$ and that $\bfY = (Y,\nu,S)$ is a two-step Abelian
isometric $\bbZ^2$-system. Suppose further that $\bfX = (X,\mu,T)$
is a joining of these five systems through the factor maps
$\xi_i:\bfX\to \bfX_i$, $i=1,2,3,4$ and $\eta:\bfX\to\bfY$, with the
maximality properties that $\xi_i = \zeta_0^{T^{\bf{n}_i}}$ for
$i=1,2,3,4$ and $\eta \succsim \zeta_1^T$.

\quad\textbf{(1)}\quad Under these assumptions the factor maps
$\xi_1,\xi_2,\xi_3,\xi_4,\eta$ are relatively independent under
$\mu$ over their further factors
\[\zeta_{\Ab,2}^{T_1}\circ
\xi_1,\ \zeta_{\Ab,2}^{T_2}\circ \xi_2,\ \zeta_{\Ab,2}^{T_3}\circ
\xi_3,\ \zeta_{\Ab,2}^{T_4}\circ\xi_4,\ \eta.\]

\quad\textbf{(2)}\quad If in addition we know that $\bfY$ is a
two-step $\bbZ^2$-pro-nilsystem whose Kronecker factor has the DIO property, then the five factors above are
actually relatively independent over \[\zeta_{\nil,2}^{T_1}\circ \xi_1,\
\zeta_{\nil,2}^{T_2}\circ \xi_2,\ \zeta_{\nil,2}^{T_3}\circ \xi_3,\
\zeta_{\nil,2}^{T_4}\circ\xi_4,\ \eta.\]
\end{prop}

\textbf{Proof}\quad\textbf{(1)}\quad First set $\b_i :=
\zeta_2^{T_i}\circ\xi_i$ and $\a_i := \zeta_{\Ab,2}^{T_i}\circ
\xi_i$ for $i=1,2,3,4$, so each $\a_i\succsim
\zeta_1^{T_i}\circ\xi_i$ is the maximal Abelian subextension of
$\b_i\succsim \zeta_1^{T_i}\circ\xi_i$.

We need to prove that
\[\int_X f_1f_2f_3f_4g\,\d\mu = \int_X \sfE_\mu(f_1\,|\,\a_1)\sfE_\mu(f_2\,|\,\a_2)\sfE_\mu(f_3\,|\,\a_3)\sfE_\mu(f_4\,|\,\a_4)g\,\d\mu\] for any
$\xi_i$-measurable functions $f_i$ and $\eta$-measurable function
$g$.  In fact it will suffice to prove that
\[\int_X f_1f_2f_3f_4g\,\d\mu = \int_X f_1f_2f_3\sfE_\mu(f_4\,|\,\a_4)g\,\d\mu,\]
since then repeating the same argument for the other three isotropy
factors in turn completes the proof.

By Proposition~\ref{prop:joining-Kron-to-isotropies} the three factors
$\zeta_1^T\vee\xi_1$, $\zeta_1^T\vee \xi_2$ and $\zeta_1^T\vee\xi_3$
must be joined relatively independently over $\zeta_1^T$.  On the
other hand, the factor $\xi_4\vee\eta$ is an extension of
$\zeta_1^T$ that is certainly still an Abelian isometric extension
for the $(\bbZ\bf{n}_4)$-subaction, and so $\xi_1\vee
\xi_2\vee\xi_3\vee\zeta_1^T$ must be joined to it relatively
independently over
\[\zeta_2^{T^{\bf{n}_4}}\wedge\big(\xi_1\vee
\xi_2\vee\xi_3\vee\zeta_1^T\big).\]

However, now the Furstenberg-Zimmer Structure Theorem tells us that
this last factor must be contained in
\[(\zeta_2^{T^{\bf{n}_4}}\wedge\xi_1)\vee(\zeta_2^{T^{\bf{n}_4}}\wedge \xi_2)\vee(\zeta_2^{T^{\bf{n}_4}}\wedge\xi_3)\vee\zeta_1^T\]
(using that $\zeta_2^{T^{\bf{n}_4}}\wedge (\xi_i\vee \zeta_1^T) =
(\zeta_2^{T^{\bf{n}_4}}\wedge \xi_i)\vee \zeta_1^T$, because
$\zeta_1^T$ is already one-step distal). Here the factors
$\zeta_2^{T^{\bf{n}_4}}\wedge \xi_i$ are actually isometric
extensions of $\zeta_1^T\wedge \xi_i$ (not just of
$\zeta_1^{T^{\bf{n}_4}}\wedge\xi_i$), since in each case
isometricity for the $(\bbZ\bf{n}_4)$-subaction and
\emph{invariance} for the $(\bbZ\bf{n}_i)$-subaction together imply
isometricity for the whole $\bbZ^2$-system
$\zeta_1^{T^{\bf{n}_4}}\wedge\xi_i$, since $\bbZ\bf{n}_i +
\bbZ\bf{n}_4$ has finite index in $\bbZ^2$ by the non-parallel
assumption.

Overall this tells us that $\xi_4\vee\eta$ is relatively independent
from the factors $\xi_1$, $\xi_2$ and $\xi_3$ over their further
factors $\b_1$, $\b_2$ and $\b_3$; and now applying the same
argument with any of the other isotropy factors as the distinguished
factor in place of $\xi_4$, we deduce that this latter is relatively
independent from all our other factors over $\b_4$.

By reducing to the factor of $\bfX$ generated by the $\b_i$ and
$\eta$, we may therefore assume that each $\bfX_i$ is itself a
two-step distal system (since the join $\b_1\vee
\b_2\vee\b_3\vee\b_4\vee\eta$ is still two-step distal, and so its
maximal isotropy factor in each direction $\bf{n}_i$ is also
two-step distal and hence equal to $\b_i$).

To make the remaining reduction to have $\a_i$ in place of $\b_i$,
now let $\bfZ_1^T = (Z_\star,m_{Z_\star},R_{\phi_\star})$ be some
coordinatization of the Kronecker factor $\zeta_1^T$ as a direct
integral of ergodic $\bbZ^2$-group rotations, and let us pick
coordinatizations
\begin{center}
$\phantom{x}$\xymatrix{\bfX_i\ar[dr]_{\zeta_1^T|_{\xi_i}}\ar@{<->}[rr]^-\cong
&& \bfZ_1^{T_i}\ltimes
(G_{i,\bullet}/H_{i,\bullet},m_{G_{i,\bullet}/H_{i,\bullet}},\s_i)\ar[dl]^{\rm{canonical}}\\
& \bfZ_1^{T_i}  }
\end{center}
and
\begin{center}
$\phantom{x}$\xymatrix{\bfY\ar[dr]_{\zeta_1^T|_\eta}\ar@{<->}[rr]^-\cong
&& \bfZ_1^T\ltimes
(A_\bullet,m_{A_\bullet},\tau)\ar[dl]^{\rm{canonical}}\\
& \bfZ_1^T.  }
\end{center}
As usual this may be done so that the $\s_i$ and $\tau$ are
relatively ergodic.

Let us first complete the proof in case the systems $\bfX_i$, $i=1,2,3$ are all fibre-normal over their Kronecker factors (Definition~\ref{dfn:fib-nor}), so that we may take $H_{i,\bullet} = \{1_{G_{i,\bullet}}\}$.

Given this, any joining of the above relatively ergodic
group extensions of $\bfZ_1^T$ is described by some
$T|_{\zeta_1^T}$-invariant measurable Mackey group data
\[M_z \leq \prod_{i=1}^4 G_{i,z_i}\times A_z\]
and a section $b:Z\to \prod_{i=1}^4 G_{i,z_i}\times A_z$, where $z
\in Z_\star$ and we write $z_i = \zeta_0^{T^{\bf{n}_i}}|_{\zeta_1^T}(z) = z\ol{\phi_\star(\bbZ\bf{n}_i)}$.  To
complete the proof of part (1) under our fibre-normality assumption we will show that
\[M_z \geq \prod_{i=1}^4[G_{i,z_i},G_{i,z_i}]\times \{1_{A_z}\}\]
almost surely, since in this case we may quotient out each extension
$\bfX_i\to \bfZ_1^{T_i}$ fibrewise by the normal subgroups
$[G_{i,\bullet},G_{i,\bullet}] \leq G_{i,\bullet}$ to obtain that
our joining is relatively independent over some Abelian
subextensions, as required.

The point is that for any three-subset $\{i_1,i_2,i_3\}\subset
\{1,2,3,4\}$ the projection of $M_\bullet$ onto the product of
factor groups $G_{i_j,z_{i_j}}$, $j=1,2,3$ is just the Mackey group
data of the joining of $\xi_{i_1}$, $\xi_{i_2}$, $\xi_{i_3}$ and
$\zeta_1^T$ as factors of $\bfX$.  By
Proposition~\ref{prop:joining-Kron-to-isotropies} these are relatively
independent over $\zeta_1^T$, so this coordinate projection of the
Mackey group must be the whole of $\prod_{j=1}^3G_{i_j,z_{i_j}}$.
Hence $M_\bullet$ has full projections onto any three of
the $G_{i,z_i}$, and so for any $g_1,h_1 \in G_{1,z_1}$ (say) we can
find $g_2 \in G_{2,z_2}$, $h_3 \in G_{3,z_3}$ and $a,b \in A_z$ such
that
\begin{multline*}
(g_1,g_2,1,1,a),(h_1,1,h_3,1,b) \in M_z\\
\Rightarrow\quad\quad [(g_1,g_2,1,1,a),(h_1,1,h_3,1,b)] =
([g_1,h_1],1,1,1,1)\in M_z.
\end{multline*}

Arguing similarly for the other $G_{i,z_i}$, we deduce that
$M_\bullet$ contains the Cartesian product of commutator subgroups,
as required.

Finally, if the systems $\bfX_i$ are not fibre-normal, then regarding them as systems with acting group $\bbZ^2/\bbZ\bf{n}_i$ and applying Proposition~\ref{prop:fibre-normals-exist} gives extensions $\t{\bfX}_i\to \bfX_i$ that are fibre-normal, and we may now extend $\bfX$ to a joining $\t{\bfX}$ of these systems with $\bfY$ simply by joining these new extensions relatively independently over $\bfX$.  Having done this the above argument shows that the factor maps onto the $\t{\bfX}_i$ and $\bfY$ are relatively independent over their maximal two-step Abelian subextensions, and hence the joining of the original systems $\bfX_i$ must be relatively independent over some factors that are simultaneously contained in these two-step Abelian distal factors of the $\t{\bfX}_i$.  Since the class of two-step Abelian distal systems is closed under taking factors (for example, by a simple appeal to the Furstenberg-Zimmer Inverse Theorem and the Relative Factor Structure Theorem~\ref{thm:RFST}), this completes the proof.

\quad\textbf{(2)}\quad We can prove this more delicate assertion by
considering the Mackey group data $M_\bullet$ and cocycle-section
obtained above more carefully. By part (1) we may reduce to the case in which each $\xi_i = \a_i$ for each $i$: that is, each $\xi_i\succsim \zeta_1^{T_i}\circ \xi_i$ is itself an Abelian isometric extension. Retaining the notation from part (1), this means we may take the group data $G_{i,\bullet}$ to be Abelian for each $i$, and now as before the joint distribution of the five factors is given by the Mackey group data and section.  As recalled in Theorem~\ref{thm:homo-nonergMackey-1}, these are characterized by the minimality of $M_\bullet$ subject to the cocycle equation
\begin{multline}\label{eq:cobdry}
(\s_1(\bf{n},z_1),\s_2(\bf{n},z_2),\s_3(\bf{n},z_3),\s_4(\bf{n},z_4),\tau(\bf{n},z))\\ \cdot (b\circ R_{\phi_\star(\bf{n})}(z))\cdot b(z)^{-1} \in M_z,
\end{multline}
where as previously we write $T|_{\zeta_1^T} = R_{\phi_\star}$.

The group data $M_\bullet$ is invariant under $T|_{\zeta_1^T} =
R_{\phi_\star}$, and so the same is true of its one-dimensional slices
such as
\[M_{1,z} := M_z\cap (G_{1,z_1}\times \{1\}\times \{1\}\times \{1\}\times \{1\}).\]
Identifying $M_{i,z}$ with a subgroup of $G_{i,z_i}$ by ignoring the restricted coordinates, we now note that the $R_{\phi_\star}$-invariance of $M_\bullet$ implies that this one-dimensional slice may be regarded as depending only on $z_i$ (and, of course, still being $T_i|_{\zeta_1^{T_i}}$-invariant).  We may
therefore consider the subextension of each $\xi_i\succsim \zeta_1^{T_i}\circ \xi_i$ corresponding to the fibrewise quotient
maps onto the quotient groups $A_{i,z_i} := G_{i,z_i}/M_{i,z_i}$.  Since $M_\bullet$ contains the product of the one-dimensional slices, the factors $\xi_i$ and $\eta$ are all relatively independent over the joining of these subextensions, and for that smaller joining the corresponding Mackey group data has trivial one-dimensional slices.

Let us now adjust our notation so that $M_\bullet \leq \prod_{i=1}^4A_{i,z_i}\times A_z$
and $b_\bullet$ are the Mackey group data and section for
the joining of these smaller Abelian extensions, so that we may now assume $M_{i,z} = \{1\}$ almost surely
for each $i$. We complete the proof by showing that under this
further assumption, the coboundary equation~(\ref{eq:cobdry})
implies that each of the cocycles $\s_i$ has one-dimensional
projections that all satisfy the conditions of
Proposition~\ref{prop:CL} over almost every ergodic component of $\zeta_1^{T_i}$.

Thus, now suppose that $\chi_\star \in \hat{A_{1,\star}}$ is a
motionless selection of characters. We will prove that the cocycle
$\chi_\star(\s_{1,\star})$ admits solutions to all of the
Conze-Lesigne equations, and hence defines a two-step nilsystem
factor of $\zeta_0^{T_1}$.  The argument for the other $T_i$ is
similar.

The point is that since the one-dimensional slices of $M_\star$ are
almost surely trivial, there are motionless selections
$\chi_{i,\star}\in \hat{A_{i,\star}}$ for $i=2,3,4$ and
$\chi_{\nil,\star}\in\hat{A_{\nil,\star}}$ such that
\[\vec{\chi}_\star := \chi_{1,\star}\otimes \chi_{2,\star}\otimes \chi_{3,\star}\otimes \chi_{4,\star}\otimes \chi_{\nil,\star} \in M_\star^\perp\]
almost surely (because now the composition of homomorphisms
\begin{multline*}
A_{1,z}\cong A_{1,z}\times\{(1,1,1,1)\}\subset A_{1,z_1}\times
A_{2,z_2}\times A_{3,z_3}\times A_{4,z_4}\times A_z\\
\onto (A_{1,z_1}\times A_{2,z_2}\times A_{3,z_3}\times
A_{4,z_4}\times A_z)/M_z
\end{multline*}
is an injection, and any character on $A_{1,z}$ can therefore be
extended to a character on the right-hand quotient group). Applying this product character to the relation~(\ref{eq:cobdry})
and setting $\s_i':= \chi_{i,\star}(\s_i)$, $\s'_\nil :=
\chi_{\nil,\star}(\s_\nil)$ and $b' := \vec{\chi}_\star(b)$ we
obtain
\[\s'_1\s'_2\s'_3\s'_4\s'_\nil = \DDelta_{\phi(\bullet)}b'.\]

We need to find solutions to the Conze-Lesigne equations for
$\s'_{1,s}$ for almost every point $s$ in the invariant base space of $\bfX$.  Since $\s'_{1,s}$ is a cocycle for an action of $\bbZ^2/\bbZ\bf{n}_1$ (that is, $\s'_{1,s}(\bf{n}_1,\cdot) \equiv 1)$,
by Proposition~\ref{prop:loc-nil-to-nil-1} it will suffice to do
this for some finite-index subgroup $Z_0\leq Z_s$ and for the subaction of some $\bf{n}' \in \bbZ^2$ which is linearly independent from $\bf{n}_1$ and such that $\phi_s(\bf{n}') \in Z_0$.

By assumption $\bf{n}_1$ and $\bf{n}_2$ are non-parallel, and so we will do this with $\bf{n}' \in \bbZ\bf{n}_2\setminus \{\bs{0}\}$.  Evaluating the above coboundary equation in
the direction $\bf{n}_2$ and using that $\s'_{2,s}(\bf{n}_2,\cdot)
\equiv 1$ gives
\[\s'_{1,s}(\bf{n}_2,\cdot)\s'_{3,s}(\bf{n}_2,\cdot)\s'_{4,s}(\bf{n}_2,\cdot)\s'_{\nil,s}(\bf{n}_2,\cdot) = \DDelta_{\phi(\bf{n}_2)}b_s'.\]

Let $Z_0 := \ol{\phi_s(\bbZ\bf{n}_1 +
\bbZ\bf{n}_3})$ and consider $z = z'z'' \in Z_0$ with $z' \in
\ol{\phi_s(\bbZ\bf{n}_1)}$ and $z'' \in \ol{\phi_s(\bbZ\bf{n}_3)}$.
Differencing the above equation by $z''$ and using that
$\s'_{1,s}(\bf{n}_1,\cdot)$ and $\s'_{3,s}(\bf{n}_2,\cdot)$ are respectively $\ol{\phi_s(\bbZ\bf{n}_1)}$- and $\ol{\phi_s(\bbZ\bf{n}_3)}$-invariant gives
\[\DDelta_z\s'_{1,s}(\bf{n}_2,\cdot)\cdot\DDelta_{z''}\s'_{4,s}(\bf{n}_2,\cdot)\cdot\DDelta_{z''}\s'_{\nil,s}(\bf{n}_2,\cdot) = \DDelta_{\phi_s(\bf{n}_2)}\DDelta_{z''}b'_s.\]

Since $\s'_{\nil,s}$ is already a nil-cocycle, we have
\[\DDelta_{z''}\s'_{\nil,s}(\bf{n}_2,\cdot) = \DDelta_{\phi_s(\bf{n}_2)}b''\cdot c(\bf{n}_2)\]
for some $b'' \in \C(Z)$ and $c \in \Hom(\bbZ^2,\Sone)$. Therefore
re-arranging gives
\begin{eqnarray}\label{eq:final-prod-cobdry}
\DDelta_z\s'_{1,s}(\bf{n}_2,\cdot)\cdot(c(\bf{n}_2)\DDelta_{z''}\s'_{4,s}(\bf{n}_2,\cdot)) = \DDelta_{\phi_s(\bf{n}_2)}((\DDelta_{z''}b'_s)\cdot \ol{b''})).
\end{eqnarray}

Finally we note that the two factors on the left-hand side of this
equation depend only on $z_1 = z\ol{\phi_s(\bbZ\bf{n}_1)}$ and $z_4 = z\ol{\phi_s(\bbZ\bf{n}_4)}$ respectively.  Since the DIO property of $\bfZ$ promises that $\ol{\phi_s(\bbZ\bf{n}_1)}\cap \ol{\phi_s(\bbZ\bf{n}_4)} = \{1\}$, these two projections of $z$ are independent when $z$ is chosen from the further finite-index
subgroup $Z_1 := (\ol{\phi_s(\bbZ\bf{n}_1)}\cap Z_0)\cdot (\ol{(\phi_s(\bbZ\bf{n}_4)}\cap Z_0)$.  Choosing $m\geq 1$ such that $\phi(m\bf{n}_2) \in Z_1$ and raising equation~(\ref{eq:final-prod-cobdry}) to the power $m$, an
application of Lemma~\ref{lem:factorizing-transfer} now shows that
$\DDelta_z\s'_{1,s}(\bf{n}',\cdot)$ is itself a quasi-coboundary when $\bf{n}' := m\bf{n}_2$, as
required for Proposition~\ref{prop:loc-nil-to-nil-1}. \qed

\subsection{A reduction to simpler averages}

We can now introduce our pleasant extensions for the averages of Theorem~\ref{thm:polyconv}.

\begin{thm}\label{thm:char-poly}
Any $\bbZ^2$-system $\bfX_0$ admits an extension
$\pi:\bfX\to \bfX_0$ in which the factor
\[\xi_1 = \xi_2 := \zeta_\pro^{T^{\bf{e}_1}}\vee \zeta_0^{T^{\bf{e}_2}}\vee \zeta_{\nil,2}^T\]
is characteristic for the averages $S_N(\cdot,\cdot)$ appearing in Theorem~\ref{thm:polyconv}, in the sense that
\[S_N(f_1,f_2)\sim S_N(\sfE_\mu(f_1\,|\,\xi_1),\sfE_\mu(f_2\,|\,\xi_2))\]
in $L^2(\mu)$ as $N\to\infty$ for any $f_1,f_2 \in L^\infty(\mu)$.
\end{thm}

\begin{lem}\label{lem:vdC-for-poly}
If $\bfX$ is as output by Theorem~\ref{thm:char-three-lines-in-2D} and \[\frac{1}{N}\sum_{n=1}^N (f_1\circ T_1^{n^2})(f_2\circ
T_1^{n^2}T_2^n) \not\to 0\] in $L^2(\mu)$ as $N\to\infty$ then there
are some $\eps > 0$ and an increasing sequence of integers $1 \leq
h_1 < h_2 < \ldots$ such that
\[\big\|\sfE_\mu(f_1\,|\,\zeta_0^{T_1^{2h_i}}\vee\zeta_0^{T_1^{2h_i}T_2^{-1}}\vee\zeta_0^{T_2^{-1}}\vee\zeta_{\nil,2}^T)\big\|_2^2 \geq \eps\]
and
\[\big\|\sfE_\mu(f_2\,|\,\zeta_0^{T_1^{2h_i}}\vee\zeta_0^{T_1^{2h_i}T_2}\vee\zeta_0^{T_2}\vee\zeta_{\nil,2}^T)\big\|_2^2 \geq \eps\]
for each $i\geq 1$.
\end{lem}

\textbf{Proof}\quad Setting $u_n := (f_1\circ T_1^{n^2})(f_2\circ
T_1^{n^2}T_2^n) \in L^2(\mu)$, the version of the classical van der
Corput estimate for bounded Hilbert space sequences (see, for
instance, Section 1 of Furstenberg and Weiss~\cite{FurWei96}) shows
that
\[\frac{1}{N}\sum_{n=1}^N (f_1\circ T_1^{n^2})(f_2\circ T_1^{n^2}T_2^n) \not\to 0\]
in $L^2(\mu)$ as $N\to\infty$ only if
\begin{eqnarray*}
&&\frac{1}{H}\sum_{h=1}^H\frac{1}{N}\sum_{n=1}^N\langle
u_n,u_{n+h}\rangle\\
&&=\frac{1}{H}\sum_{h=1}^H\int_Xf_1\cdot\frac{1}{N}\sum_{n=1}^N
((f_1\circ T_1^{h^2})\circ T_1^{2hn}) (f_2\circ
T_2^n) ((f_2\circ T_1^{h^2}T_2^h)\circ T_1^{2hn}T_2^n)\,\d\mu\\
&&\not\to 0,
\end{eqnarray*}
and hence, by the Cauchy-Schwartz inequality, only if $f_1 \neq 0$
and for some $\eps
> 0$ there is an increasing sequence $1 \leq h_2 < h_2 < \ldots$
such that
\begin{eqnarray*}
&&\|f_1\|_2\Big\|\lim_{N\to\infty}\frac{1}{N}\sum_{n=1}^N(f_1\circ
T_1^{h_i^2}\circ T_1^{2h_in})(f_2\circ T_2^n)(f_2\circ
(T_1^{h_i^2}T_2^{h_i})\circ (T_1^{2h_i}T_2)^n)\Big\|_2\\
&&\geq
\Big|\int_Xf_1\cdot\Big(\lim_{N\to\infty}\frac{1}{N}\sum_{n=1}^N
(f_1\circ T_1^{h_i^2}\circ T_1^{2h_in}) (f_2\circ T_2^n) (f_2\circ
(T_1^{h_i^2}T_2^{h_i})\circ
(T_1^{2h_i}T_2)^n)\Big)\,\d\mu\Big|\\
&&\geq \|f_1\|_2\eps.
\end{eqnarray*}
It follows that each of $f_1$, $f_2$ should have conditional expectation of norm at least $\sqrt{\eps}$ onto the corresponding factor in any characteristic triple of factors for the above linear averages in $n$, so by Theorem~\ref{thm:char-three-lines-in-2D} this translates into the desired assertion. \qed

This tells us that if $S_N(f_1,f_2) \not\to 0$ then each
of $f_1$ and $f_2$ must enjoy a large conditional expectation onto
not just one factor of $\bfX$ with a special structure, but a whole
infinite sequence of these factors.  We will now use this to cut
down the characteristic factors we need for the averages $S_N$
further by considering the possible joint distributions of the members
of these infinite families of factors. Crucially, we can make use of the relative independence studied in the previous subsection through the following simple lemma.

\begin{lem}\label{lem:rel-ind-factor-seq} Suppose that $(X,\mu)$ is a
standard Borel probability space, $\pi_n:X\to Y_n$ is a sequence of
factor maps of $X$ and $\a_n:Y_n\to Z_n$ is a sequence of further
factor maps of $Y_n$ such that $(\pi_n,\pi_m)$ are relatively
independent over $(\a_n\circ\pi_n,\a_m\circ\pi_m)$ whenever $n \neq
m$ (note that we assume only pairwise relative independence).  If $f \in L^\infty(\mu)$ is such
that $\limsup_{n\to\infty}\|\sfE_\mu(f\,|\,\pi_n)\|_2 > 0$, then
also $\limsup_{n\to\infty}\|\sfE_\mu(f\,|\,\a_n)\|_2 > 0$.
\end{lem}

\textbf{Proof}\quad By thinning out our sequence if necessary, we
may assume that for some $\eta > 0$ we have
$\|\sfE_\mu(f\,|\,\pi_n)\|_2 \geq \eta$ for all $n$.  Suppose, for
the sake of contradiction, that $\sfE_\mu(f\,|\,\a_n) \to 0$ as
$n\to\infty$. Consider the sequence of Hilbert subspaces $L_n \leq
L^2(\mu)$ comprising those functions that are $\pi_n$-measurable and
the further subspaces $K_n \leq L_n$ comprising those that are
$\a_n$-measurable.  Then by assumption all the subspaces $L_n\ominus
K_n$ are mutually orthogonal, but $f$ has orthogonal projection of
norm at least $\eta/2$ onto all but finitely many of them, which is
clearly impossible. \qed

\textbf{Proof of Theorem~\ref{thm:char-poly}}\quad Let $\pi:\bfX\to\bfX_0$ be an extension as given by Theorem~\ref{thm:char-three-lines-in-2D}. If $f_1,f_2 \in
L^\infty(\mu)$ have $S_n(f_1,f_2)\not\to 0$, then
Lemma~\ref{lem:vdC-for-poly} promises that $f_1$ has a uniformly
large conditional expectation onto each of the factors
$\zeta_0^{T_1^{2h_i}}\vee\zeta_0^{T_1^{2h_i}T_2^{-1}}\vee\zeta_0^{T_2^{-1}}\vee\zeta_{\nil,2}^T$
for some infinite sequence $h_1 < h_2 < \ldots$.  Since no two of
the four vectors $(2h,0)$, $(h_i,1)$, $(h_j,1)$, $(0,1) \in \bbZ^2$ are
parallel when $h_i\neq h_j$ and $h := \rm{l.c.m.}(h_i,h_j)$, by
Proposition~\ref{prop:jointdists} the factors in this sequence are
pairwise relatively independent over the further factor
$\zeta_\pro^{T_1}\vee \zeta_0^{T_2}\vee \zeta_{\nil,2}^T$,
and so Lemma~\ref{lem:rel-ind-factor-seq} shows that in fact $f_1$
must have a nonzero conditional expectation onto this latter factor
also. The argument for $f_2$ is similar. \qed

\subsection{Completion of the convergence proof}

\textbf{Proof of Theorem~\ref{thm:polyconv}}\quad By Theorem~\ref{thm:char-poly} this will follow if we prove that $S_N(f_1,f_2)$ converges whenever $f_i$
is $\xi_i$-measurable. By approximation in $L^2(\mu)$ and
multilinearity, it actually suffices to consider the averages
$S_N(f_{11} f_{12} g_1,f_{21}f_{22}g_2)$ in which each $f_{j1}$ is
$T^\ell_1$-invariant for some large $\ell\geq 1$, each
$f_{j2}$ is $T_2$-invariant and each $g_j$ is
$\zeta^T_{\nil,2}$-measurable.

Next, writing
\begin{eqnarray*}
&&S_N(f_{11} f_{12}g_1,f_{21}f_{22}g_2) =
\frac{1}{N}\sum_{n=1}^N((f_{11}\cdot f_{12}\cdot g_1)\circ
T_1^{n^2})((f_{21}\cdot f_{22}\cdot g_2)\circ T_1^{n^2}T_2^n)\\
&&\sim \frac{1}{\ell}\sum_{k = 0}^{\ell-1}\frac{1}{(N/\ell)}\sum_{n
= 1}^{\lfloor N/\ell\rfloor} ((f_{11}\cdot f_{12}\cdot g_1)\circ
T_1^{(\ell n+k)^2})((f_{21}\cdot
f_{22}\cdot g_2)\circ T_1^{(\ell n+k)^2}T_2^{\ell n+ k})\\
&&= \frac{1}{\ell}\sum_{k = 0}^{\ell-1}(f_{11}\circ
T_1^{k^2})\Big(\frac{1}{(N/\ell)}\sum_{n = 1}^{\lfloor
N/\ell\rfloor} (g_1\circ T_1^{(\ell
n+k)^2})((f_{12}\cdot f_{21}\cdot f_{22}\cdot g_2)\circ T_1^{(\ell
n+k)^2}T_2^{\ell n+ k})\Big)
\end{eqnarray*}
(recalling that $\sim$ denotes asymptotic agreement in $L^2(\mu)$ as
$N\to\infty$), we see that it will suffice to prove convergence in
$L^2(\mu)$ for all averages along infinite arithmetic progressions
of the form
\[\frac{1}{(N/\ell)}\sum_{n = 1}^{\lfloor
N/\ell\rfloor} (g\circ T_1^{(\ell n)^2
+2k(\ell n)})(f\circ T_1^{(\ell n)^2 + 2k(\ell n)}T_2^{\ell
n})\] for all $k \in
\{0,1,\ldots,\ell-1\}$, where $g$ is $\zeta_{\nil,2}^T$-measurable. Let us now re-label $T_i^\ell$ as $T_i$ (and so effectively
restrict our attention to the subaction of $\ell\bbZ^2$), $\lfloor N/\ell\rfloor$ as $N$ and set $a := 2k$ so that the
above averages can be written as
\[\frac{1}{N}\sum_{n = 1}^{N} (g\circ T_1^{\ell n^2
+ an})(f\circ T_1^{\ell n^2 + an}T_2^n).\]

If we now simply re-run the standard application of the van der Corput estimate for these averages and consider the resulting non-vanishing integral under the Furstenberg self-joining, the assumption that $g$ is $\zeta^T_{\nil,2}$-measurable enables us to condition $f$ also onto some more restricted factor.  Specifically, the van der Corput estimate implies that if the above averages do not asymptotically vanish in $L^2(\mu)$ then also
\begin{eqnarray*}
&&\frac{1}{H}\sum_{h=1}^H\frac{1}{N}\sum_{n=1}^N\int_X (g\circ T_1^{\ell(n+h)^2 + a(n+h)})(\ol{g}\circ T_1^{\ell n^2 + an})\\ &&\quad\quad\quad\quad\quad\quad\quad\quad\quad\quad\cdot(f\circ T_1^{\ell(n+h)^2 + a(n+h)}T_2^{n+h})(\ol{f}\circ T_1^{\ell n^2 + an}T_2^n)\,\d\mu\\
&&= \frac{1}{H}\sum_{h=1}^H\frac{1}{N}\sum_{n=1}^N\int_X (g\circ T_1^{2\ell hn + \ell h^2 + ah}T_2^{-n})(\ol{g}\circ T_2^{-n})\\ &&\quad\quad\quad\quad\quad\quad\quad\quad\quad\quad\quad\quad\quad\quad\cdot(f\circ T_1^{2\ell hn +\ell h^2 + ah}T_2^h)\ol{f}\,\d\mu\\ && \not\to 0
\end{eqnarray*}
as $N\to\infty$ and then $H\to\infty$.  Hence there must be infinitely many $h$ for which the linear averages
\[\frac{1}{N}\sum_{n=1}^N(g\circ T_1^{2\ell hn + \ell h^2 + ah}T_2^{-n})(\ol{g}\circ T_2^{-n})(f\circ T_1^{2\ell hn +\ell h^2 + ah}T_2^h)\]
do not tend to zero in $L^2(\mu)$.  Another appeal to the van der Corput estimate therefore gives that
\[\int_{X^3}((g\circ T_1^{\ell h^2 + ah})\otimes \ol{g}\otimes (f\circ T_1^{\ell h^2 + ah}T_2^h))\cdot G\,\d\mu^{\rm{F}}_{T_1^{2\ell h}T_2^{-1},T_2^{-1},T_1^{2\ell h}} \neq 0\]
for some $(T_1^{2\ell h}T_2^{-1}\times T_2^{-1}\times T_1^{2\ell h})$-invariant function $G$.

This latter non-vanishing asserts that through the third coordinate projection $X^3 \to X$, the lifted function $f\circ \pi_3$ enjoys a nonzero correlation with a product of the function $(g\circ T_1^{\ell h^2 + ah}\circ \pi_1)\otimes (\ol{g}\circ \pi_2)$, which is measurable with respect to a two-step pro-nilsystem factor, and the function $G$ which is invariant under a lift of the transformation $T_1^{2\ell h}$.  Using the satedness of $\bfX$ again this implies that $f$ must itself correlate with the join of $\zeta^T_{\nil,2}$ and $\zeta_{0}^{T_1^{2\ell h}}$.

We may therefore break up $f$ again as $g'f'$ with $g'$ being $\zeta_{\nil,2}^T$-measurable and $f'$ being $T_1^{2\ell h}$-invariant. Re-inserting this into the averages of interest and increasing the previously-used value of $\ell$ accordingly, it follows that we need only prove convergence of
\[\frac{1}{N}\sum_{n=1}^N (g\circ T_1^{\ell n^2 + an})(g'\circ T_1^{\ell n^2 + an}T_2^n)(f'\circ T_2^n)\]
with $g$, $g'$ and $f'$ as above.

Finally, let $\mu = \int_S\mu_s\,\nu(\d s)$ be the $T$-ergodic decomposition of $\mu$.  We will show that the above averages form a Cauchy sequence in $L^2(\mu_s)$ for $\nu$-almost every $s$ separately.

By definition, for $\nu$-almost every $\mu_s$ the projection of $\mu_s$ onto the factor $\zeta^T_{\nil,2}$ is concentrated on an inverse limit of ergodic two-step nilsystems.  Hence for each such fixed $s$ the functions $g$ and $g'$ may be approximated in $L^2(\mu_s)$ by lifts of functions that are continuous on some finite-dimensional nilsystem appearing in this inverse sequence.  Letting these approximating functions be $h$ and $h'$, it now suffices to prove the convergence of
\[\frac{1}{N}\sum_{n=1}^N (h\circ T_1^{\ell n^2 + an})(h'\circ T_1^{\ell n^2 + an}T_2^n)(f'\circ T_2^n)\]
in $L^2(\mu_s)$.

However, having reduced to this problem it turns out that we can appeal to pointwise convergence using a recent theorem of Host and Kra.  In Theorem 2.22 of~\cite{HosKra07} they show that in our setting there is some $\mu_s$-conegligible subset $X_0 \subseteq X$ such that the sequence
\[\frac{1}{N}\sum_{n=1}^N b_n\cdot f'(T_2^nx)\]
converges for every nilsequence $(b_n)_{n\geq 1}$ and every $x\in X_0$.  Since sampling continuous functions along polynomial orbits on a nilmanifold still produces nilseqeuences (see, for instance, Leibman's papers~\cite{Lei05,Lei05(2)}), it follows that 
\[\frac{1}{N}\sum_{n=1}^N h(T_1^{\ell n^2 + an}x)\cdot h'( T_1^{\ell n^2 + an}T_2^nx)\cdot f'(T_2^nx)\]
converges for every $x \in X_0$, and hence that these averages do converge in $L^2(\mu_s)$.  This completes the proof. \qed

\textbf{Remark}\quad Before leaving the quadratic averages of Theorem~\ref{thm:polyconv}, we note that in their recent preprint~\cite{ChuFraHos09} Chu, Frantzikinakis and Host have proven (as a corollary of a strong convergence result for some different nonconventional averages) that for the question of \emph{weak} convergence in $L^2(\mu)$, our averages admit the even smaller characteristic pair of factors $\zeta_\pro^{T_1}\vee \zeta_\pro^{T_2}$, $\zeta_\pro^{T_2}$.  Although their approach does not give strong convergence, combined with the fact of that convergence proved here it follows that those smaller factors are in fact characteristic for strong convergence.  However, it seems hard to prove this using only repeated appeals to the van der Corput estimate and results on linear averages.  This may indicate a higher-dimensional instance of a phenomenon that Leibman has studied in some detail for polynomial nonconventional averages associated to a  single $\bbZ$-action~(\cite{Lei10}): the pro-nilsystem characteristic factors indicated by the Host-Kra Theory in that setting can sometimes be reduced further, but (so far) only by using more detailed results about nilsystems.  While Leibman obtains a more-or-less complete characterization of the extent of this phenomenon in one dimension, our exploration of its higher-dimensional generalization is only just beginning.

It is also worth noting that our use of the structure of direct integrals of nilsystems is in many ways similar to that of Chu, Frantzikinakis and Host, once the relevance of that structure has been established; the principal difference in our situation is that we must take a very different route to the result that these and various isotropy systems are sufficient ingredients to describe the characteristic factors completely. \fin

\section{Closing remarks}\label{sec:end}

The strategy of passing to a pleasant extension of a system in order
to enable a simplified description of its nonconventional averages
seems to be quite a powerful one, and I suspect that it will have
much further-reaching consequences in this area in the future.

In view of Theorem~\ref{thm:char-three-lines-in-2D} and the various earlier results that are known about $\bbZ$-actions or actions of linearly independent subgroups, it is natural to attempt a conjecture about characteristic factors for general polynomial nonconventional averages.  To be a little vague, this would assert that for any polynomial mappings $p_i:\bbZ\to\bbZ^d$, $i=1,2,\ldots,k$, if $\bfX$ is a $\bbZ^d$-systems that is sated relative to a sufficiently large list of joins of different idempotent classes of system then the averages
\[\frac{1}{N}\sum_{n=1}^N\prod_{i=1}^k (f_i\circ T^{p_i(n)})\]
for $f_i\in L^\infty(\mu)$ admit a characteristic tuple of factors $\xi_1$, $\xi_2$, \ldots, $\xi_k$ each of which is a join of systems for which some nontrivial subgroup of $\bbZ^d$ acts as a pro-nilsystem of some finite step.  Such a result as this would not settle the convergence of the above averages immediately, but it would surely constitute a major reduction of the problem. However, just what satedness assumption is needed, and how the lists of partially pro-nilsystem factors that appear in each $\xi_i$ could be determined in terms of $p_1$, $p_2$, \ldots, $p_k$, remain unclear from the few special cases that are known.

Although I feel that the pursuit of such a more general result is perhaps the most pressing issue suggested by our work above, it seems worth mentioning a few more specific questions that may be within easier reach.

Firstly, I suspect that some generalization of Theorem~\ref{thm:polyconv} to commuting actions
of $\bbZ^r$ should lie fairly close at hand, with the
principal new difficulty being that of reigning in the complexity of
the notation: probably one could prove that for any two commuting
actions $T_i:\bbZ^{r_i}\curvearrowright (X,\mu)$, $i=1,2$, any
quadratic form $Q:\bbZ^s\to\bbZ^{r_1}$ and any homomorphism
$L:\bbZ^s\to\bbZ^{r_2}$ the averages
\[\frac{1}{N^s}\sum_{\bf{n} \in \{1,2,\ldots,N\}^s}(f_1\circ T_1^{Q(\bf{n})})(f_2\circ T_1^{Q(\bf{n})}T_2^{L(\bf{n})})\]
converge in $L^2(\mu)$.

The next simplest case to consider might be
that of the averages
\[\frac{1}{N}\sum_{n=1}^N (f_1\circ T_1^{n^2})(f_2\circ T_2^{n^2})\]
for commuting transformations $T_1$ and $T_2$, but already the approach to these via characteristic factors seems to require some substantial improvement on Theorem~\ref{thm:char-three-lines-in-2D}.  Repeatedly applying the van der
Corput estimate to these until we reach some linear nonconventional
averages now throws up such averages
corresponding to seven different directions in $\bbZ^2$, for which
the available pleasant extensions will surely be much more complicated than those of 
Theorem~\ref{thm:char-three-lines-in-2D}.

The above averages are also the subject of the following
less ambitious question (put to me by Vitaly Bergelson), which may
be within closer range of the methods currently available:

\begin{ques}
Is it true that if $T_1^{-1}T_2\curvearrowright (X,\mu)$ and
$T_1\times T_2\curvearrowright (X^2,\mu^{\otimes 2})$ are both
totally ergodic then we have
\[\frac{1}{N}\sum_{n=1}^N (f_1\circ T_1^{n^2})(f_2\circ T_2^{n^2}) \to \int_X f_1\,\d\mu\cdot \int_X f_2\,\d\mu\quad\quad \hbox{in}\ L^2(\mu)\ \forall f_1,f_2 \in L^\infty(\mu)?\]
\end{ques}

Note that this is true if we assume instead that \emph{every}
direction in our $\bbZ^2$-action is totally ergodic, as follows from
the extension of Host and Kra's nilsystem machinery to
higher-dimensional actions under this assumption worked out by
Frantzikinakis and Kra in~\cite{FraKra05}.

\appendix

\section{Background on Moore cohomology}\label{app:cohom}

The classical cohomology of discrete groups (see, for instance,
Weibel~\cite{Wei94}) was extended to the category of locally compact
groups acting on Polish Abelian groups by Moore in a far-reaching
sequence of
papers~\cite{Moo64(gr-cohomI-II),Moo76(gr-cohomIII),Moo76(gr-cohomIV)},
and it is his version of the theory that we use in this paper.  We refer
the reader to those papers for a clear introduction to the subject,
discussion of the various issues that arise in the attempt to take
the topologies of the groups into account, and also a discussion
with further references of the relation in which this theory stands
to various other cohomology theories that have been developed for
locally compact groups.  (Let us also remark in passing that these
measurable cohomology groups have already appeared in
ergodic-theoretic works from time to time in the past; consider, for
example, the paper~\cite{Lem97} of Lema\'nczyk.)  

In this appendix we recall some important properties of Moore's measurable cohomology groups for locally compact groups, including some continuity properties of these groups under forming inverse limits of compact base groups that were recently established in~\cite{AusMoo--cohomcty}, and also give the details of a few purely cohomological calculations that were needed in Subsection~\ref{subs:factorizing-in-erg-case}.

Given a compact Abelian group $Z$ and a
Polish Abelian $Z$-module $A$ we write $\cal{Z}^r(Z,A)$ to denote
the Borel $r$-cocycles $Z^r\to A$ that appear in the inhomogeneous bar resolution, $\cal{B}^r(Z,A)$ to denote the
subgroup of coboundaries, and $\rmH^r(Z,A) :=
\cal{Z}^r(Z,A)/\cal{B}^r(Z,A)$ to denote the resulting Moore cohomology
group (which we will not topologize here).

\textbf{Notation}\quad In keeping with the rest of this paper, we will use multiplicative notation for cochains taking values in $(\Sone)^D$ for some $D$, but additive notation for cochains and other maps into discrete Abelian groups. Also, we now let $\lfloor\cdot\rfloor:\bbR\to\bbZ$ be the usual integer-part map and $\{\cdot\}:\Sone\to [0,1)$ be the lift such that $\{\rm{e}^{2\pi\rm{i}s}\}= s$ when $s \in [0,1)$. \fin

\begin{thm}[Theorem B from~\cite{AusMoo--cohomcty}]\label{thm:cohomA}
If $(G_m)_{m\geq 1}$, $(\pi^m_k)_{m\geq k\geq 1}$ is an inverse sequence of compact groups with inverse limit $G$, $(\pi_m)_{m\geq 1}$ then
\[\rmH^p(G,A) \cong \lim_{m\rightarrow}\rmH^p(G_m,A)\]
under the inverse limit of the inflation maps $\rm{inf}^p_{\pi_m}:\rmH^p(G_m,A)\to \rmH^p(G,A)$ whenever $A$ a discrete Abelian group or a finite-dimensional torus. \qed
\end{thm}

\begin{thm}[Theorem C from~\cite{AusMoo--cohomcty}]\label{thm:cohomB}
If $G$ is any compact group and
\[A_1 \subseteq A_2 \subseteq \ldots\]
is an increasing sequence of discrete $G$-modules with union $A = \bigcup_{m\geq 1}A_m$, also equipped with the discrete topology, then
\[\rmH^p(G,A) \cong \lim_{m\rightarrow}\rmH^p(G,A_m)\]
under the direct limit of the maps on cohomology induced by the inclusions $A_m \subseteq A$.
\end{thm}

Theorems A and B are valuable in conjunction with the explicit calculations that are available for cohomology over compact Abelian Lie groups.

\begin{prop}\label{prop:HofMoscalc}
The graded Moore cohomology ring $\rmH^\ast((\Sone)^D,\bbZ)$ is isomorphic to the symmetric algebra
\[\rm{Sym}^{2\ast}\hat{(\Sone)^D}\cong \rm{Sym}^{2\ast}(\bbZ^D),\]
graded so that every individual element of $\hat{(\Sone)^D}$ has degree \emph{two}.  This graded cohomology ring is isomorphic to $\rmH^{\ast - 1}((\Sone)^D,\Sone)$ under the switchback isomorphisms given by the long exact sequence of the presentation $\bbZ\into\bbR\onto \Sone$ together with the vanishing $\rmH^\ast((\Sone)^D,\bbR) = (0)$.
\end{prop}

\textbf{Proof}\quad This follows from the identification of $\rmH^\ast((\Sone)^D,\bbZ)$ with the \v{C}ech cohomology $\rmH^\ast_{\rm{sp}}((\Sone)^D,\bbZ)$ of the Milnor classifying space $B((\Sone)^D)$ of $(\Sone)^D$, as proved by Wigner in Theorem 4 of~\cite{Wig73}.  These latter cohomology rings have been computed very explicitly by Hofmann and Mostert in their monograph~\cite{HofMos73}: the result we need follows from Theorem 1.9 in their Chapter V. \qed

\begin{cor}\label{cor:reps-for-Sone}
Every cohomology class in $\rmH^3((\Sone)^D,\Sone) \cong \rmH^4((\Sone)^D,\bbZ)$ contains a representative of the form
\begin{eqnarray*}
\psi_{\g_1,\g_2,\ldots,\g_M,\chi_1,\chi_2,\ldots,\chi_M}&:&(\Sone)^D\times (\Sone)^D\times (\Sone)^D\to \Sone\\
&:&(z_1,z_2,z_3)\mapsto \prod_{m=1}^M \chi_m(z_3)^{\lfloor\{\g_m(z_1)\} + \{\g_m(z_2)\}\rfloor}
\end{eqnarray*}
for some $\g_1,\g_2,\ldots,\g_M$, $\chi_1,\chi_2,\ldots,\chi_M\in \hat{(\Sone)^D}$ (indeed, if $D\geq 1$ they each contain infinitely many such representatives).  The map
\[\psi_{\g_1,\g_2,\ldots,\g_M,\chi_1,\chi_2,\ldots,\chi_M}\mapsto \sum_{m=1}^M \g_m\odot \chi_m\]
descends at the level of cohomology to the isomorphism
\[\rmH^3((\Sone)^D,\Sone)\stackrel{\rm{switchback}}{\cong}\rmH^4((\Sone)^D,\bbZ)\stackrel{\scriptsize{\hbox{Prop.\ref{prop:HofMoscalc}}}}{\cong} \hat{(\Sone)^D}\odot \hat{(\Sone)^D} \cong \rm{Sym}^{2\ast}(\bbZ^D)\big|_{\rm{deg}=4}.\]
\end{cor}

\textbf{Proof}\quad These representatives are precisely those obtained from all possible lists of elements of $\rmH^1((\Sone)^D,\Sone)$ --- that is, of single characters --- by moving these to $\rmH^2((\Sone)^D,\bbZ)$ through the switchback, forming cup products in $\rmH^4((\Sone)^D,\bbZ)$ and moving back to $\rmH^3((\Sone)^D,\Sone)$.  This generates a complete list of representatives and gives rise to the asserted isomorphism in view of the identification with $\rm{Sym}^{2\ast}(\bbZ^D)$ given by Proposition~\ref{prop:HofMoscalc}. \qed

\begin{cor}\label{cor:vanishing-from-HofMos}
The groups $\rmH^p((\Sone)^D,\bbZ)$ vanish for all odd $p$ and are torsion-free for all $p\geq 1$, and consequently the groups $\rmH^p((\Sone)^D,F)$ also vanish for all odd $p$ when $F$ is a finite Abelian group.
\end{cor}

\textbf{Proof}\quad The first two assertions follow at once from the identification
\[\rmH^\ast((\Sone)^D,\bbZ)\cong \rm{Sym}^{2\ast}(\bbZ^D).\]
For the third, note first that using the Structure Theorem for finite Abelian groups and arguing coordinate-wise it suffices to treat $F = \bbZ/n\bbZ$. Now the odd-degree vanishing for $\bbZ$-valued cocycles together with the presentation $\bbZ\stackrel{\times n}{\into}\bbZ\onto \bbZ/n\bbZ$ give switchback maps $\rmH^p((\Sone)^D,\bbZ/n\bbZ)\into \rmH^{p+1}((\Sone)^D,\bbZ)$ that are monomorphisms for all odd $p$, but also their images must vanish under multiplying by $n$ and so must take values among the torsion elements of $\rmH^{p+1}((\Sone)^D,\bbZ)$. Since there are none of these it follows that $\rmH^p((\Sone)^D,\bbZ/n\bbZ) = (0)$. \qed

\begin{cor}\label{cor:reps-for-Sone0}
If $Z$ is a compact Abelian group then each class in $\rmH^2(Z,\bbZ)$ contains exactly one representative of the form
\[Z\times Z\to \bbZ:(z_1,z_2) \mapsto \lfloor\{\g(z_1)\} + \{\g(z_2)\}\rfloor\]
for some $\g \in \hat{Z}$.  In addition, if $Z$ is connected then any class in $\rmH^2(Z,\bbZ/n\bbZ)$ also contains a (possibly non-unique) representative of the above form, where the integer on the right-hand side is now to be understood modulo $n$.
\end{cor}

\textbf{Proof}\quad The above representatives are precisely those obtained by implementing the switchback isomorphism
\[\rmH^1(Z,\Sone)\cong \hat{Z}\to \rmH^2(Z,\bbZ)\]
that arises from the presentation
\[0 \to \bbZ\to \bbR \stackrel{\rm{e}^{2\pi\rm{i}\cdot}}{\to} \Sone\to 0,\]
then using the vanishing result $\rmH^\ast(Z,\bbR) = (0)$ and the particular choice of lifting $\Sone \to \bbR:z\mapsto \{z\}$.

For the second conclusion the previous corollary gives $\rmH^1(Z,\bbZ) = \rmH^3(Z,\bbZ) = (0)$ for connected $Z$, and so in this case the quotient map $\rmH^2(Z,\bbZ)\to \rmH^2(Z,\bbZ/n\bbZ)$ is an isomorphism in view of the collapsing of the long exact sequence. \qed

In the remainder of this appendix we include the two chief cohomological vanishing results that were needed in the main text.

\begin{lem}\label{lem:C(Z)K-vald-collapse}
Suppose that $K_1,K_2$ are compact Abelian groups, $Z := K_1\times K_2$ and $q_i:Z\to K_i$ are the coordinate projections. Identify $\C(K_2)$ with $\C(Z)^{K_1}$, and consider it as a $Z$-submodule of $\C(Z)$ with the rotation action.  Then any $2$-cocycle $\theta:Z\times Z\to \C(Z)^{K_1}$ is of the form $(\theta_1\circ q_1^{\times 2})\cdot d\a$ for some measurable $\a:Z\to \C(Z)^{K_1}$ and some $2$-cocycle $\theta_1:K_1\times K_1 \to \Sone$.
\end{lem}

\textbf{Proof}\quad Consider the map
\[\Psi:\Z^2(K_1,\Sone)\to \Z^2(Z,\C(Z)^{K_1})\]
defined by
\[\Psi(\theta_1)((z_1,z_2),(z'_1,z'_2)) := \iota(\theta_1(z_1,z'_1)),\]
where $z = (z_1,z_2)$ denotes the coordinates in $K_1\times K_2$ and $\iota:\Sone \to \C(Z)^{K_1}$ is the constant-functions embedding.  This map $\Psi$ sends cocycles to cocycles and coboundaries to coboundaries, so descends to a homomorphism of cohomology groups
\[\rmH^2(K_1,\Sone)\to \rmH^2(Z,\C(Z)^{K_1}).\]

The present lemma is asserting that this homomorphism is surjective.  However, it is actually an isomorphism: if one recognizes $\C(Z)^{K_1}$ as the result of inducing the module $\Sone$ with trivial $K_1$-action to the larger group $Z$, this is the classical Shapiro Isomorphism.  In the setting of measurable group cohomology, it is constructed as an abstract isomorphism by Moore in Theorem 6 of~\cite{Moo76(gr-cohomIII)}.  It only remains to check that this abstract isomorphism is realized by the above map $\Psi$ at the level of cochains.  This follows from the proof in~\cite{Moo76(gr-cohomIII)} by a routine diagram chase. \qed

\begin{lem}\label{lem:getting-rep-with-anning-chars}
Suppose that $Z = K_1\times K_2$ is a product of two compact Abelian groups and that $\psi:Z\times Z\times Z\to \rmH^3(Z,\Sone)$ is the $3$-cocycle
\[(z_1,z_2,z_3)\mapsto \prod_{m=1}^M\chi_m(z_3)^{\lfloor \{\g_m(z_1)\} + \{\g_m(z_2)\}\rfloor}\]
corresponding to some choice of $\g_1,\g_2,\ldots,\g_M$, $\chi_1,\chi_2,\ldots,\chi_M \in \hat{Z}$.  Then the class $[\psi]$ trivializes under the inclusion of the constant-valued maps $\Sone \into \C(Z)^{K_1}$ if and only if $\psi$ is cohomologous in $\rmH^3(Z,\Sone)$ to a cocycle expressible as above with $\chi_1,\chi_2,\ldots,\chi_M \in K_1^\perp$.
\end{lem}

\textbf{Proof}\quad Observing that the $3$-cocycle equation holds \emph{strictly} everywhere for the above function $\psi$, and using again Theorem 5 in~\cite{Moo76(gr-cohomIII)}, we may assume that there is some $\k:Z\times Z\to \C(Z)^{K_1}$ such that the equation $\psi = d\k$ holds among elements of $\C(Z)^{K_1}$ \emph{strictly} everywhere on $Z^3$: that is, that
\[\psi(z_1,z_2,z_3) = d\k(z_1,z_2,z_3)(wK_1)\]
for Haar-a.e. $w\in Z$ for strictly every $(z_1,z_2,z_3)\in Z^3$.  It follows, in particular, that for every $(z_1,z_2,z_3) \in K_1^3$ the above holds for Haar-a.e. $w$, and hence by Fubini's Theorem we may find some $w \in Z$ for which the above holds for a.e. $(z_1,z_2,z_3)\in K_1^3$: that is, we have successfully restricted the given coboundary equation to the possibly-negligible subset $K_1^3\leq Z^3$.  However, now for this fixed $w$, using the fact that $z_3K_1 = K_1$, the restricted equation indicates that $\psi|_{K_1\times K_1\times K_1}$ must be an $\Sone$-valued $3$-coboundary on $K_1$.  This restriction is given explicitly by
\[\prod_{m\in J}\big(\chi_m|_{K_1}(z_3)\big)^{\lfloor \{\g_m|_{K_1}(z_1)\} + \{\g_m|_{K_1}(z_2)\}\rfloor}\]
where
\[J := \{m\in \{1,2,\ldots,M\}:\ \chi_m\not\in K_1^\perp\}.\]
(we choose to keep any terms with $\chi_m \not\in K_1^\perp$ but $\g_m \in K_1^\perp$, even though they also vanish on $K_1^3$). Now for each $\g_m$ let $\t{\g}_m$ denote the unique element of $K_2^\perp$ for which $\t{\g}_m|_{K_1} = \g_m|_{K_1}$, and similarly associate $\t{\chi}_m$ to each $\chi_m$. Given these, let $\psi_1$ be the inflation of $\psi|_{K_1\times K_1\times K_1}$ back up to $Z^3$ through the coordinate projection map $Z\to K_1$, so
\[\psi_1(z_1,z_2,z_3) := \prod_{m\in J}\t{\chi}_m(z_3)^{\lfloor \{\t{\g}_m(z_1)\} + \{\t{\g}_m(z_2)\} \rfloor}.\]
If $\psi|_{K_1\times K_1\times K_1} = d\l$ for some $\l:K_1\times K_1\to \Sone$, then also $\psi_1 = d\l_1$ where $\l_1$ is the inflation of $\l$, and hence $\psi$ is cohomologous to $\psi \cdot \psi_1^{-1}$.  However,
\begin{eqnarray*}
(\psi \cdot \psi_1^{-1})(z_1,z_2,z_3) &=& \prod_{m \in \{1,2,\ldots,M\}\setminus J}\chi_m(z_3)^{\lfloor \{\g_m(z_1)\} + \{\g_m(z_2)\}\rfloor}\\
&&\quad \cdot \prod_{m\in J}(\chi_m \cdot \t{\chi}_m^{-1})(z_3)^{\lfloor \{\g_m(z_1)\} + \{\g_m(z_2)\}\rfloor}\\ &&\quad  \cdot \prod_{m\in J}\t{\chi}_m(z_3)^{\lfloor \{\g_m(z_1)\} + \{\g_m(z_2)\} \rfloor - \lfloor \{\t{\g}_m(z_1)\} + \{\t{\g}_m(z_2)\} \rfloor}.
\end{eqnarray*}
Since $\chi_m|_{K_1} = \t{\chi}_m|_{K_1}$, the first and second terms of this right-hand side are already in the desired form.  On the other hand, owing to the symmetrization implied by the identification with $\rm{Sym}^\ast$ in Corollary~\ref{cor:reps-for-Sone}, and since $(\g_m \cdot \t{\g}_m^{-1})\odot \t{\chi}_m = \t{\chi}_m\odot (\g_m \cdot \t{\g}_m^{-1})$, the last sum is cohomologous to
\[\prod_{m\in J}(\g_m(z_3) \cdot \t{\g}_m(z_3)^{-1})^{\lfloor \{\t{\chi}_m(z_1)\} + \{\t{\chi}_m(z_2)\} \rfloor},\]
which is also of the desired form because $\g_m|_{K_1} = \t{\g}_m|_{K_1}$. \qed

\begin{lem}\label{lem:deg-3-vanishing}
If $Z = K_1\times K_2$ is a product of compact Abelian groups and $\psi:Z\times Z\times Z\to \E(Z)^{K_1}$ is a $3$-cocycle whose class trivializes under the inclusion $\E(Z)^{K_1}\subseteq \C(Z)^{K_1}$, then there is a finite-index subgroup $Z_0\leq Z$ such that $\psi|_{Z_0\times Z_0\times Z_0}$ is a coboundary.
\end{lem}

\textbf{Proof}\quad Most of the work here will go into reducing to the case in which $Z$ is a Lie group.

\quad\textbf{Step 1}\quad We have a natural presentation of $Z$-modules
\[\Sone\into \E(Z)^{K_1}\onto \hat{Z}^{K_1}\cong K_1^\perp\]
according which $\Sone$ is identified with the constant-valued members of $\E(Z)^{K_1}$, and both $\Sone$ and $K_1^\perp$ have the trivial $Z$-action (even though the action on $\E(Z)^{K_1}$ is not trivial).  To this presentation corresponds the long exact sequence
\begin{multline*}
\ldots\to \rmH^{p-1}(Z,K_1^\perp)\\
\stackrel{\rm{switchback}}{\to} \rmH^p(Z,\Sone)\to \rmH^p(Z,\E(Z)^{K_1})\to \rmH^p(Z,K_1^\perp)\\
\stackrel{\rm{switchback}}{\to} \rmH^{p+1}(Z,\Sone)\to \ldots,
\end{multline*}
and so the class $[\psi]\in \rmH^3(Z,\E(Z)^{K_1})$ is uniquely identified by its image $[\ol{\psi}]\in \rmH^3(Z,K_1^\perp)$, where $\ol{\psi}(z_1,z_2,z_3) := \psi(z_1,z_2,z_3)\cdot \Sone$ takes values in
\[\E(Z)^{K_1}/\Sone\cong K_1^\perp,\]
together with some element of $\rmH^3(Z,\Sone)$ parameterizing the location of $[\psi]$ in the fibre over $[\ol{\psi}]$.

We may express each $K_i$ as an inverse limit of an increasing sequence of Lie quotient groups $q_{i,m}:K_i\onto K_{i,m}$, which combine to give Lie quotient groups $q_m := q_{1,m}\times q_{2,m}:Z\onto Z_m := K_{1,m}\times K_{2,m}$. Pontrjagin duality gives $K_1^\perp = \bigcup_{m\geq 1}K_{1,m}^\perp\circ q_m$, where $K^\perp_{1,m}$ is understood as a subgroup of $\hat{Z_m}$ and each $K_{1,m}^\perp\circ q_m$ may be identified as a $Z$-module with trivial action. Given this, Theorem~\ref{thm:cohomB} implies that
\[\rmH^p(Z,K_1^\perp) \cong \lim_{m\rightarrow}\rmH^p(Z,K_{1,m}^\perp\circ q_m)\]
under the inclusion maps.  Therefore there are some $m\geq 1$ and $\ol{\phi} \in \Z^3(Z,K_{1,m}^\perp\circ q_m)$ such that $[\ol{\psi}] = [\ol{\phi}]$, so letting $\phi:Z^3\to \E(Z_m)^{K_{1,m}}\circ q_m$ be a measurable lift of $\ol{\phi}$ it follows that $\psi \cdot \phi^{-1} = d\k \cdot \phi'$ for some $\k:Z^2\to \E(Z)^{K_1}$ and $\Sone$-valued $3$-cocycle $\phi'$. Adjusting $\psi$ by $d\k$ (which does not effect the desired conclusion) if necessary we may assume $\k = 0$, after which this equation tells us that $\psi$ takes values in $\E(Z_m)^{K_{1,m}}\circ q_m$ for some finite $m$.  By omitting a finite initial segment of the sequence $(Z_m)_{m\geq 1}$ and re-labelling we may also assume that $m=1$.

\quad\textbf{Step 2}\quad The increasing sequences of epimorphisms $q_{i,m}:K_i\onto K_{i,m}$ also define families of intermediate connecting epimorphisms $q^m_{i,k}:K_{i,m}\onto K_{i,k}$, and hence $q^m_k := q^m_{1,k}\times q^m_{2,k}:Z_m\onto Z_k$, whenever $m\geq k \geq 1$. Since the $Z$-action on the module $\E(Z_1)^{K_{1,1}}\circ\pi_1$ actually factorizes through the quotient $Z\onto Z_1$, an application of Theorem~\ref{thm:cohomA} in conjunction with the long exact sequence above gives
\[\rmH^3(Z,\E(Z_1)^{K_{1,1}}\circ q_1)\cong \lim_{m\rightarrow}\rmH^3(Z_m,\E(Z_1)^{K_{1,1}}\circ q^m_1),\]
now under the direct limit of the inflation maps, and so by adjusting $\psi$ by another $(\E(Z_1)^{K_{1,1}}\circ q_1)$-valued coboundary  we may assume that it is itself lifted from a cocycle $\psi_1 \in \Z^3(Z_m,\E(Z_1)^{K_{1,1}}\circ q^m_1)$ for some finite $m$.  Re-labelling we may again assume that $m = 1$.

\quad\textbf{Step 3}\quad Since each $K_{i,1}$ is finite-dimensional it is isomorphic to $(\Sone)^{D_i}\times F_i$ for some $D_i\geq 0$ and finite Abelian group $F_i$, and consequently $Z_1\cong (\Sone)^{D_1 + D_2}\times (F_1\times F_2)$.  Let $Z_{1,0}\cong (\Sone)^{D_1 + D_2}$ be the identity component of $Z_1$ and $Z_0 := q_1^{-1}(Z_{1,0})$, so that $[Z:Z_0] = |F_1||F_2| < \infty$.  We will show that $\psi|_{Z_0\times Z_0\times Z_0}$ is an $\E(Z)^{K_1}$-valued $3$-coboundary.

Considering the presentation
\[\Sone\into \E(Z_1)^{K_{1,1}}\onto K_{1,1}^\perp,\]
we first observe that $\ol{\psi_1} := \psi_1 \cdot \Sone$ is a $3$-cocycle with values in $K_{1,1}^\perp \cong \bbZ^{D_1}\times \hat{F_1}$, so by Corollary~\ref{cor:vanishing-from-HofMos} its class must trivialize on $Z_{1,0}$. Therefore there is some $\bar{\k}:Z_{1,0}^2\to \bbZ^{D_1}\times \hat{F_1}$ such that
\[\ol{\psi_1}|_{Z_{1,0}\times Z_{1,0}\times Z_{1,0}} = d\bar{\k},\]
and so letting $\k$ be a measurable $\E(Z)^{K_1}$-valued lift of $\bar{\k}$ it follows that $\psi_1|_{Z_{1,0}\times Z_{1,0}\times Z_{1,0}} - d\k$ takes values in the subgroup of constant-valued functions $\Sone \subseteq \E(Z)^{K_1}$.  Now our initial assumptions promise that the class of this $\Sone$-valued $3$-cocycle still trivializes under the inclusion $\Sone\subseteq \C(Z_{1,0})^{K_{1,1,0}}$, where $K_{1,1,0}$ is the identity component of $K_{1,1}$.

However, in view of this Lemma~\ref{lem:getting-rep-with-anning-chars} gives some $\a_0:Z_{1,0}^2 \to \Sone$ such that
\[\big((\psi_1|_{Z_{1,0}\times Z_{1,0}\times Z_{1,0}}) \cdot d\k^{-1} \cdot d\a_0^{-1}\big)(z_1,z_2,z_3) = \prod_{m=1}^M \chi_m(z_3)^{\lfloor\{\g_m(z_1)\} + \{\g_m(z_2)\}\rfloor}\]
for some $\chi_1$, $\chi_1$, \ldots, $\chi_M \in \widehat{Z_{1,0}}\cap K_{1,1,0}^\perp$, and lifting these characters to elements of $K_{1,1}^\perp$ this is now explicitly the coboundary of the $\E(Z_1)^{K_{1,1}}$-valued $2$-cochain
\[(u_1,u_2)\mapsto \prod_{m=1}^M \chi_m(\cdot)^{\lfloor\{\g_m(u_1)\} + \{\g_m(u_2)\}\rfloor}.\]
Re-arranging completes the proof. \qed

\section{Measurable selectors}

At several points in this paper we needed to appeal to some basic results on the existence of measurable selectors, often as a means of making rigorous a selection of representatives of one or another kind of data above the ergodic components of a non-ergodic system.

\begin{thm}\label{thm:meas-select}
Suppose that $(X,\S_X)$ and $(Y,\S_Y)$ are standard Borel spaces,
that $A \subseteq X$ is Borel and that $\pi:X\to Y$ is a Borel
surjection. Then the image $\pi(A)$ lies in the
$\nu^\rm{c}$-completion of $\S_Y$ for every Borel probability
measure $\nu$ on $(Y,\S_Y)$ with completion $\nu^\rm{c}$, and for
any such $\nu$ there is a map $f:B\to A$ with domain $B \in \S_Y$
such that $B \subseteq \pi(A)$, $\nu^\rm{c}(\pi(A)\setminus B) = 0$
and $\pi\circ f = \id_B$. \qed
\end{thm}

\textbf{Proof}\quad See, for example, 423O and its consequence
424X(h) in Fremlin~\cite{FreVol4}. \qed

The above result prompts the following two standard definitions.

\begin{dfn}[Universal measurability]
Given a standard Borel space $(X,\S_X)$, a measurable subset $A$ of $X$ is \textbf{universally measurable} if for any Borel probability $\mu$ on $X$ there is some $A' \in \S_X$ such that $\mu^\rm{c}(A\triangle A') = 0$; thus, the first part of the above conclusion is that $\pi(A)$ is universally measurable.
\end{dfn}

\begin{dfn}[Measurable selectors]
We refer to a map $f$ as given by the above theorem as a
\textbf{measurable selector} for the set $A$.
\end{dfn}

\textbf{Remark}\quad We should stress that this is only one of
several versions of the `measurable selector theorem', due variously
to von Neumann, Jankow, Lusin and others. Note in particular that in
some other versions a map $f$ is sought that select points of $A$
for \emph{strictly} all points of $\pi(A)$. In the above generality
we cannot guarantee that a strictly-everywhere selector $f$ is
Borel, but only that it is Souslin-analytic and hence universally
measurable (of course, from this the above version follows at once).
On the other hand, if the map $\pi|_A$ is countable-to-one, then a
version of the result due to Lusin does guarantee a
strictly-everywhere Borel selector $f$. This version has already
played a significant r\^ole in our corner of ergodic theory in the
manipulation of the Conze-Lesigne equations (see, for
example,~\cite{ConLes84,FurWei96,BerTaoZie09}), and so we should be
careful to distinguish it from the above.  A thorough account of all
these different results and their proofs can be found in Sections
423, 424 and 433 of Fremlin~\cite{FreVol4}. \fin

In the right circumstances it is possible to use
Theorem~\ref{thm:meas-select} to obtain a Borel selector that is
equivariant for a group of transformations, by making use of a
coordinatization of the invariant factor.  We prove this only for countable groups for simplicity, but the general case of l.c.s.c. groups and jointly measurable actions should follow with a little more care. It will be used in Subsection~\ref{subs:factorizing-cocycles}. (Note that an incorrect version of the following proposition and proof appeared as Proposition 2.4 in~\cite{Aus--ergdirint}.)

\begin{prop}\label{prop:invar-meas-select}
Suppose that $T:\G\curvearrowright (X,\mu)$ and $S:\G \curvearrowright (Y,\nu)$ are actions of a countable group on standard Borel probability spaces; that $\pi:(X,\mu,T)\to (Y,\nu,S)$ is a factor map; that $A\subseteq X$ lies in the $\mu$-completion of $\S_X$ and is conegligible and $T$-invariant; and also that $\pi$ is \textbf{relatively invariant}, meaning that $\S_X$ is generated by $\pi^{-1}(\S_Y)$ together with $\S_X^T$. Then there are an $S$-invariant set $B\in \S_Y$ such that $B\subseteq \pi(A)$ and $\nu^\rm{c}(\pi(A)\setminus B) = 0$ and a Borel map $f:B\to
A$ such that $f\circ S^\g = T^\g\circ f$ and $\pi\circ f = \id_B$.
\end{prop}

\textbf{Proof}\quad Let $A_0 \subseteq A$ be a conegligible Borel subset, and let $A_1 := \bigcap_{\g \in \G}T^\g(A_0)$, so this is still conegligible and Borel (using the countability of $\G$) but also $T$-invariant.  By replacing $A$ with $A_1$ if necessary, we may simply assume that $A$ itself is Borel.  However, having done this, we may replace $X$ with $A$ (since $(A,\mu|_{\S_A})$ is still a standard Borel probability space), and so assume that $A = X$.

Next choose a factor $\zeta_0^S:\bfY\to \bfZ_0^S$ that coordinatizes $\S_Y^S$, and then choose another factor $\zeta_0^T:\bfX\to \bfZ_0^T$ that coordinatizes $\S_X^T$ and such that there is a factorizing map $\phi:\bfZ_0^T\to \bfZ_0^S$ for which $\zeta_0^S\circ \pi(x) = \phi\circ \zeta_0^T(x)$ for all $x$ (for example, this can be guaranteed by replacing an initial choice of $\zeta_0^T$ with the map $(\zeta_0^T,\zeta_0^S\circ \pi):X\to Z_0^T\times Z_0^S$, which generates the same $\s$-subalgebra of $\S$ up to negligible sets).

Now consider the condition that $\S_X^T$ and $\pi^{-1}(\S_Y)$ together generate the whole of $\S_X$.  This amounts to the assertion that the map
\[(\pi,\zeta_0^T):X\to Y\times Z_0^T:x\mapsto (\pi(x),\zeta_0^T(x))\]
defines a measure-theoretic isomorphism of systems, and hence by restricting $X$ to a further full-measure $T$-invariant Borel subset we may actually identify it with an $(S\times \id)$-invariant Borel subset of $Y\times Z_0^T$ under this map $(\pi,\zeta_0^T)$ (by the standard pointwise description of measure-theoretic isomorphisms between standard Borel systems -- see Theorem 2.15 in Glasner~\cite{Gla03}).

Moreover, the condition that $\zeta_0^S\circ \pi = \phi\circ \zeta_0^T$ implies that
\[X \subseteq Y \times_{\{\zeta_0^S = \phi\}} Z_0^T := \{(y,z):\ \zeta_0^S(y) = \phi(z)\},\]
and so $\mu$ is $(S\times \id)$-invariant and is supported on this set.  We will next argue that $\mu$ actually equals the relatively independent product of $\nu$ and $(\zeta_0^T)_\#\mu$ on this set.  In case $\G$ is amenable this follows easily from the Norm Ergodic Theorem, but we can make a related argument even if it is not. Suppose that $f \in L^2(\nu)$ and $g \in L^2((\zeta_0^T)_\#\mu)$, and consider the integral
\[\int_{Y\times Z_0^T}f(y)g(z)\,\mu(\d y,\d z).\]
By the $(S\times \id)$-invariance of $\mu$, this is equal to
\[\int_{Y\times Z_0^T}f(S^\g y)g(z)\,\mu(\d y,\d z)\quad\quad \hbox{for any}\ \g \in \G,\]
and hence is also equal to $\int_{Y\times Z_0^T}f'(y)g(z)\,\mu(\d y,\d z)$ for any function $f'$ in the closed convex hull of the norm-bounded set $\{f\circ S^\g:\ \g \in \G\} \subseteq L^2(\nu)$.  This closed convex set has a unique element $f'$ of minimal norm (because the norm in $L^2(\mu)$ is uniformly convex), and now its uniqueness implies that this $f'$ is also $S$-invariant.  From this it follows that we must have $f' = \sfE_\nu(f\,|\,\zeta_0^S)$, for otherwise we could pick some $S$-invariant $h \in L^2(\nu)$ such that $\langle h,f'\rangle > \langle h,\sfE_\nu(f\,|\,\zeta_0^S)\rangle$, and this latter quantity must equal $\langle h,S^\g f\rangle $ for every $\g$ (because $f\circ S^\g - \sfE_\nu(f\,|\,\zeta_0^S)$ is orthogonal to all $S$-invariant functions), so that $h$ would define a hyperplane separating $f'$ from the closed convex set we have constructed.

Therefore
\[\int_{Y\times Z_0^T}f(y)g(z)\,\mu(\d y,\d z) = \int_{Y\times Z_0^T}\sfE_\nu(f\,|\,\zeta_0^S)(y)g(z)\,\mu(\d y,\d z),\]
and by taking linear combinations of product functions this implies that $\mu = \nu\otimes_{\{\zeta_0^S = \phi\}}(\zeta_0^T)_\#\mu$.

Let $P:Z_0^S\stackrel{\rm{p}}{\longrightarrow}Y$ be a probability kernel representing the disintegration of $\nu$ over $\zeta_0^S$, and $Q:Z_0^S\stackrel{\rm{p}}{\longrightarrow} Z_0^T$ a kernel representing the disintegration of $(\zeta_0^T)_\#\mu$ over $\phi$.  In terms of these kernels we have
\[\nu\otimes_{\{\zeta_0^S = \phi\}}(\zeta_0^T)_\#\mu = \int_{Z_0^S}P(w,\,\cdot\,)\otimes Q(w,\,\cdot\,)\ (\zeta_0^S)_\#\nu(\d w).\]
Since $(\nu\otimes_{\{\zeta_0^S = \phi\}}(\zeta_0^T)_\#\mu)(X) = 1$, it follows from Fubini's Theorem that for $(\zeta_0^T)_\#\mu$-a.e. $z \in Z_0^T$ we have
\[P(\phi(z),\{y \in Y:\ (z,y) \in X\}) = 1.\]
Let $C \subseteq Z_0^T$ be a conegligible Borel subset of $z$ for which this holds.  By the previous Theorem~\ref{thm:meas-select} we may find a conegligible Borel subset $D \subseteq \phi(C)$ that admits a measurable selector $g:D\to C$.

Finally, let $B_0 := (\zeta_0^S)^{-1}(D) \subseteq Y$, so this is $S$-invariant and conegligible, and consider the map
\[f:B_0\to Y\times Z_0^T:y \mapsto (y,g(\zeta_0^S(y))).\]
By the definition of the set $C$ and selector $g$, for each $w \in D$ we know that $P(w,\,\cdot\,)$-a.e. $y \in Y$ is such that $(y,g(w)) \in X$.  This implies that the further subset $B := \{y \in B_0:\ f(y) \in X\}$ has
\[\nu(B) = \int_D P(w,\{y:\ \zeta_0^S(y) = w\}\cap \{y:\ (y,g(w)) \in X\})\,(\zeta_0^S)_\#\nu(\d w) = 1,\]
and now restricting $f$ to $B$ gives a Borel measurable selector $B\to X$ that manifestly satisfies $f\circ S^\g = (S\times \id)^\g\circ f = T^\g\circ f$, as required. \qed

\begin{dfn}[Equivariant measurable selectors]
We refer to a map $f$ as given by the above proposition as a
\textbf{$T$-equivariant measurable selector} for the set $A$.
\end{dfn}

\parskip 0pt

\bibliographystyle{abbrv}
\bibliography{bibfile}

\vspace{10pt}

\small{\textsc{Courant Institute, New York University, New York, NY 10012, USA}}

\vspace{5pt}

\small{Email: \verb|tim@cims.nyu.edu|}

\vspace{5pt}

\small{URL: \verb|http://www.cims.nyu.edu/~tim|}

\end{document}